\documentclass[reqno]{amsart}

\usepackage[english]{babel}
\usepackage{amsmath}
\usepackage{amsthm}
\usepackage{amsfonts}
\usepackage{amssymb}
\usepackage{anysize}
\usepackage{hyperref}
\usepackage{appendix}
\usepackage{mathtools}
\usepackage{graphicx}
\usepackage{xcolor}
\usepackage{subcaption}
\usepackage{pgfplots}
\pgfplotsset{compat=1.12}

\newtheorem{defi}{Definition}[section]
\newtheorem{lem}[defi]{Lemma}
\newtheorem{theo}[defi]{Theorem}
\newtheorem{cor}[defi]{Corollary}
\newtheorem{pro}[defi]{Proposition}

\newtheorem{rem}[defi]{Remark}

\DeclareMathOperator{\divop}{div}

\DeclareMathOperator{\supp}{supp}
\DeclareMathOperator{\diam}{diam}
\DeclareMathOperator{\dist}{dist}

\DeclareMathOperator{\RR}{\mathbb{R}}

\DeclareMathOperator{\cut}{cut}
\DeclareMathOperator{\co}{co}
\DeclareMathOperator*{\esssup}{ess\,sup}
\DeclareMathOperator{\Skew}{Skew}

\allowdisplaybreaks 

\title[Filippov flows and mean-field limits in the kinetic singular Kuramoto model]{Filippov flows and mean-field limits in the kinetic singular Kuramoto model}
\author[David Poyato]{David Poyato}
\address{Departamento de Matem\'{a}tica Aplicada and Excellence Research Unit ``Modeling Nature'' (MNat), Universidad de Granada, Granada, 18071, Spain} 
\email{davidpoyato@ugr.es}

\begin{document}

\date{\today}



\begin{abstract}
The agent-based singular Kuramoto model was proposed in \cite{P-P-S} as a singular version of the Kuramoto model of coupled oscillators that is consistent with Hebb's rule of neuroscience. In such paper, the authors studied its well-posedness via the concept of Filippov solutions. Interestingly, they found some new emergent phenomena in the paradigm of Kuramoto model: clustering into subgroups and emergence of global phase synchronization taking place at finite time. 

This paper aims at introducing the associated kinetic singular Kuramoto model along $\mathbb{T}\times \mathbb{R}$, that is reminiscent of the classical Kuramoto--Sakaguchi equation. Our main goal is to propose a well-posedness theory of measure-valued solutions that remains valid after eventual phase collisions. The results will depend upon the specific regime of singularity: subcritical, critical and supercritical. The cornerstone is the existence of Filippov characteristic flows for interaction kernels with jump discontinuities. Our second goal is to study stability with respect to initial data, that in particular will provide quantitative estimates for the mean-field limit of the agent-based model towards the kinetic equation in quadratic Wasserstein-type distances. Finally, we will recover global phase synchronization at the macroscopic scale under appropriate generic assumptions on the initial data. The most singular regime will be tackled separately with an alternative method, namely, a singular hyperbolic limit on a regularized second order model with inertia. Note that we will work within the manifold $\mathbb{T}\times \mathbb{R}$ and will avoid resorting on Euclidean approximations. Also, no gradient-type structure will be needed, as opposed to some preceding literature. 
\end{abstract}

\maketitle

\tableofcontents


\section{Introduction}

The interest on understanding collective dynamics models has notably raised during the recent years. Not only it has attracted the attention of the most applied scientific community, but it has also become an important subject of study for the most theoretic research. From the point of view of the applied side, one of the main goals is to describe how each couple of individuals in the system interact. Notice that working with the whole set of biological variables of a complex system is a convoluted task. Then, one has to reduce the set of parameters just to a few of them, the most relevant ones, carrying most of the important information about the dynamics of the complex living system. A large family of this type of biological systems arise from \textit{soft active matter}, see \cite{Mar}. In general, active matter consists of self-propelled units or agents with the ability of converting stored or ambient energy into movement. The direction of self-propulsion is determined by agents' anisotropy itself, rather than being prescribed by some external field, as it is the case for lifeless systems. Very recently, there has been an intense development on this field, specially with regards to new non-Newtonian dynamics, such as the Cucker--Smale, Kuramoto, Vicsec models or the aggregation equation \cite{Bo-C-M,Bu-C-M,C-S,C-U,D-F-MA,D-F-MA-T,De-M,H-T,Ku-1,Ku-2,M-EK,M-C-O,T-B,T-B-L,V-C-B-C-S}. Those models have proven useful to describe emergence of patterns in groups of cells, animals or other living systems. Specifically, the interaction of active particles gives rise to highly correlated collective motion that usually ends up with an eventual emergent behavior of the whole populations as a unique entity. Such collective motion is named differently according to the specific biological context, e.g., flocking, synchronization, schooling, swarming, etc. In other words, patterns arise as a self-organized dynamics of the whole group coming from agent-based microscopic interactions. Elucidating how it takes place is the main problem here and, from the theoretical point of view, it has proven a rich source of relevant problems in mathematics. Indeed, new methods and important strong techniques have been obtained as a consequence of the positive feedback between both communities. 

Here, we shall focus on a specific agent-based model that was proposed in \cite{P-P-S} as a singularly weighted version of the original Kuramoto model \cite{Ku-1,Ku-2}. The original system by Kuramoto is a paradigmatic model describing collective synchronization of oscillators, which refers to the remarkable phenomenon of a large population of coupled oscillators that spontaneously synchronize to oscillate at a common frequency. Although Kuramoto initially proposed it for synchronization of chemical reactions, such many-body cooperative effect can be observed in many other examples in nature like the flashing of fireflies, chorusing of crickets, beating of cardiac cells or metabolic synchrony in yeast cell suspensions, see \cite{Acebron}. For some applications to neuronal synchronization and how the realistic human connectome maps that are available in the literature affect the emergence of synchronization, see \cite{V-M-M} and references therein. Such ideas exploit neuronal connections in the brain turn out to be organized in moduli structured in a hierarchical nested fashion across many scales, and it affects the neural dynamics \cite{R-S-T-B,Z-Z-Z-H-K-1,Z-Z-Z-H-K-2}

Let us comment on the above-mentioned agent-based singularly weighted version of the Kuramoto model. It stands for the following system of $N$ coupled oscillators \cite{P-P-S}
\begin{equation}\label{E-kuramoto-discrete}
\left\{
\begin{array}{l}
\displaystyle\dot{\theta}_i=\Omega_i+\frac{K}{N}\sum_{j=1}^Nh(\theta_j-\theta_i),\\
\displaystyle\theta_i(0)=\theta_{i,0}.
\end{array}
\right.
\end{equation}
We can think of such model to describe the evolution of the phases $\theta_i=\theta_i(t)$ of neuron signals located at specific areas in the brain. Hence, $\dot{\theta}_i=\dot{\theta}_i(t)$ represent the firing frequencies of neurons and $\Omega_i$ appears as an heterogeneity that play the role of a biased tendency of agents to move at their own frequency while being influenced by their neighbors via the coupling force $h$. The proposed periodic singular force reads
\begin{equation}\label{E-kuramoto-kernel}
h(\theta)=\frac{\sin\theta}{\vert \theta\vert_o^{2\alpha}},\ \theta\in \mathbb{R},
\end{equation}
where $|\theta|_o$ is the Riemannian distance of $e^{i\theta}$ towards $1$ along the unit circle, that is
\[
|\theta|_o :=\vert\bar\theta\vert\quad \text{for}\quad \bar\theta \equiv \theta \ \text{mod} ~ 2\pi, \quad \bar\theta\in (-\pi, \pi].
\]
Note that if $h(\theta)= \sin\theta$, then we recover the classical Kuramoto model. Regading the singular coupling \eqref{E-kuramoto-kernel}, it was proposed in \cite{P-P-S} via a rigorous singular limit in the fast learning regime of the Kuramoto model coupled with a learning rule for the coupling weights, that is governed by a Hebbian-type plasticity function. Such plasticity function obeys Hebb's rule \cite{H}, that represents the fact that synchronous activation of cells (firing of neurons) must lead to selectively pronounced increases in synaptic strength between those cells. This is the foundation of associative or Hebbian learning and provides an explanation about the adaptation and synchronization of neurons in the brain during a learning process.

The parameter $\alpha$ leads to three different regimes of the singularity $\alpha\in (0,\frac{1}{2})$, $\alpha=\frac{1}{2}$ and $\alpha\in (\frac{1}{2},1)$ that we will respectively call the \textit{subcritical}, \textit{critical} and \textit{supercritical} cases. Note that in the subcritical case the coupling function is H\"older-continuos, in the critical case it is bounded but discontinuous, while the supercritical regime corresponds to an unbounded singular kernel. Therefore, the Cauchy--Lipschitz theory cannot guarantee existence and uniqueness of global-in-time solutions, then requiring the concept of \textit{Filippov trajectories} \cite{Au-Ce,F}. This problem has been extensively studied in \cite{P-P-S}. In particular, uniqueness of Filippov solutions was proven to hold forwards-in-time only. Although it might appear a deficiency on the model at first glance, it actually suggests a new relevant dynamics of oscillators: \textit{finite-time sticking} and \textit{clustering into groups}. Namely, after some phases eventually agree in finite time, there is a chance that they keep stuck together for all times or they instantaneously disassociate. The rule governing such behavior has been proved to depend only upon certain specific conditions on the natural frequencies of the formed cluster \cite{P-P-S}. Such sticky behavior suggests that small clusters can emerge in finite time ending up with the eventual global synchronizations of all the oscillators in a unique big group. In other words, that can be regarded as \textit{finite-time phase synchronization} of \eqref{E-kuramoto-discrete} that was guaranteed in \cite{P-P-S} under certain initial assumptions. A similar sticky dynamic was found for related models like the singular Cucker--Smale model or the aggregation equation with mildly singular potentials, see \cite{C-D-F-L-S,C-J-L-V,M-P,Pes-1,Pes-2}.

Closely related to the aforementioned issues, one of the most classical problems in physics and mathematics is to achieve the perfect fit between such individual-based descriptions and the associated coarse-grained continuous versions, where systems are explained in an averaged way though a probabilistic distribution function. In kinetic theory, it is called \textit{mean-field limit} and, when the amount of individuals is large enough, it yields accurate macroscopic descriptions of the system in terms of a kinetic Vlasov-type equation. As it is already known, such mean-field methods for particle systems give rise to plenty of relevant models in statistical physics and fluid mechanics like Boltzmann, Vlasov--Poisson, Euler or Navier--Stokes equations, that have become the motor of important advances in mathematics. In particular, one can try to extend the methods to the new class of non-Newtonian interactions coming from the active soft matter community. When the coupling force between any two agents is Lipschitz, the mean-field limit can be recovered from standard methods \cite{McK,N}. However, interactions are usually discontinuous or even singular for most of the real systems. In such cases, the mean-field approximation is not obvious and new methods are required \cite{C-C-R,G-SR-T,H-L,H-J,H-M,J,J-W,J-W-2,J-W-3,Ka,M-M,M-M-W,N,S}.

This will be the subject of study in this paper. Specifically, we shall focus on the rigorous derivation of the Vlasov-type equation with non-smooth interactions that arises as mean-field limit of the above-mentioned singular Kuramoto model of $N$ coupled oscillators as $N\rightarrow\infty$, see \cite{P-P-S}:
\begin{equation}\label{VE}
\left\{
\begin{array}{l}
\displaystyle\frac{\partial f}{\partial t}+\frac{\partial}{\partial \theta}\left(\Omega f-K(h*\rho)f\right)=0,\\
\displaystyle f(0,\cdot,\cdot)=f_0.
\end{array}
\right.
\end{equation}
Here, $f=f(t,\theta,\Omega)$ represents the probability of finding an oscillator at time $t \geq 0$, with  frequency  $\theta\in (-\pi,\pi]$ and natural frequency $\Omega \in \mathbb{R}$, endowed with periodic boundary conditions for the variable $\theta$. Due to the periodicity of phases, such kinetic equation will be identified in a natural way with a nonlinear transport equation along the Riemannian manifold $\mathbb{T}\times \mathbb{R}$. This changes the natural phase-space from the a standard Euclidean space to a non-Euclidean ambient space. Specifically \eqref{VE} can be restated as follows
$$
\left\{
\begin{array}{l}
\displaystyle\frac{\partial f}{\partial t}+\divop_{(z,\Omega)}(\mathcal{V}[f]f)=0,\\
\displaystyle f(0,\cdot,\cdot)=f_0,
\end{array}
\right.
$$
for some tangent transport field $\mathcal{V}[f]$ to be defined later in Section \ref{S-macroscopic-model}. We remark here that we will not use any extra gradient-type structure of the system, as it was the case in previous literature for other models, see e.g. \cite{C-D-F-L-S,C-H-L-X-Y,C-H-M,F-K-M,H-K-M-P,K-M}. Contrarily, we will just work with the non-smooth tangent transport field $\mathcal{V}[f]$ with any extra structure. Such point of view requires extending the analysis of a non-smooth transport field given in \cite{A-C, A-G-S} from the Euclidean setting to our new ambient space. In addition, the dynamic of the agent-based system suggests that global phase synchronization (emergence of Dirac masses) will also take place in finite time at the macroscopic scale. Then, we must deal with weak solutions that are merely measure-valued. Our approach will be supported by the methods of Filippov characteristics and is inspired in \cite{C-J-L-V,Pou-R}. Specifically, we will try to give some sense to the Filippov flow of this non-smooth tangent transport field. To such end, we will study conditions on $\mathcal{V}[f]$ in the different regimes of $\alpha$ that allows constructing the associated characteristic flow globally-in-time. The uniqueness of the flow is again one-sided and will be guaranteed by the internal structure of the kernel $h$ at points of loss of Lipschitz-continuity. In particular, we shall show that although non-smooth, the tangent vector field $\mathcal{V}[f]$ is one-sided Lipschitz-continuous. Then, the associated characteristic system
\begin{equation*}
\left\{\begin{array}{l}
\displaystyle\frac{dX}{dt}(t;t_0,x_0)\in\mathcal{K}[\mathcal{V}[f_t]](X(t;t_0,x_0)),\\
\displaystyle \hspace{0.1cm} X(t_0;t_0,x_0)=x_0,
\end{array}\right.
\end{equation*}
enjoys global solutions in Filippov's sense that are unique forwards-in-time \cite{Au-Ce,F,P-P-S}. Here $\mathcal{K}[\mathcal{V}[f]]$ will denote the Filippov set-valued tangent field associated with $\mathcal{V}[f]$ to be defined later in Section \ref{S-critical-regime}. This will be the first result, that becomes the cornerstone to construct measure-valued solutions in the subcritical and critical cases. Our second result is the rigorous derivation of the mean field limit. Specifically, consider empirical measures supported on Filippov solutions to \eqref{E-kuramoto-discrete}-\eqref{E-kuramoto-kernel}, that is,
$$
\mu^N_t:=\frac{1}{N}\sum_{i=1}^N\delta_{z_i^N(t)}(z)\otimes \delta_{\Omega_i^N}(\Omega),
$$
where $z_i^N(t):=e^{i\theta_i^N(t)}$. Then under appropriate assumptions supported by the law of large numbers we shall show that $\mu^N\rightarrow f$ in $C([0,T],\mathcal{P}_1(\mathbb{T}\times \mathbb{R})-W_1)$ and the limiting $f$ solves \eqref{VE} in Filippov flow sense
$$
X_f(t;0,\cdot)_{\#}f_0=f_t,\ t\geq 0,
$$
where $X_f=X_f(t;0,z,\Omega)$ is the Filippov flow associated with $\mathcal{V}[f]$. Here, $W_1$ means the Rubinstein--Kantorovich distance on $\mathbb{T}\times \mathbb{R}$ and $\mathcal{P}_1$ stands for the probability measures with finite first order moment. This is a special consequence of our third result, namely, stability with respect to initial data in Wasserstein-type distances. Specifically, we will introduce some \textit{Dobrushin}-type estimate for any two measure-valued solutions $f^1=f^1_t(\theta,\Omega)$ and $f^2=f^2_t(\theta,\Omega)$ to \eqref{VE} with respect to two different metrics on the set of probability measures: the standard $2$-Wasserstein distance $W_2$ and a new version $W_{2,g}$, specially designed for this problem, that will be called fibered Wasserstein distance:
$$
\begin{array}{ll}
\displaystyle W_2(f^1_t,f^2_t)\leq e^{Ct}W_2(f^1_0,f^2_0), & t\geq 0,\\
\displaystyle W_{2,g}(f^1_t,f^2_t)\leq e^{C't}W_{2,g}(f^1_0,f^2_0), & t\geq 0.
\end{array}
$$
The latter one only works for solutions with the same distribution $g$ of natural frequencies. In particular, it allows recovering uniqueness results for any general initial data. Our forth result is the study of asymptotic behavior of solutions. Namely, we will show that under appropriate conditions, the macroscopic system \eqref{VE} enjoys finite-time global phase synchronization. 

Unfortunately, notice that there is no way to extend any (generalized) characteristic method to the most singular regime $\alpha\in (\frac{1}{2},1)$ because $\mathcal{V}[f]$ lacks of sense. To cover such supercritical case, we shall develop an alternative method that is valid for all $\alpha\in (0,1)$ (at least for identical oscillators $g=\delta_0$). We will augment the first order singular Kuramoto system into an auxiliary second order regularized system with inertia, frequency damping and diffusion, see \cite{C-H-L-X-Y,C-H-M,C-H-Y-1,C-H-Y-2}. Under an appropriate scaling depending on a parameter $\varepsilon\searrow 0$, the inertia and noise terms will vanish and singularity of the coupling function will be recovered. This defines a singular hyperbolic hydrodynamic limit of vanishing inertia type like in \cite{G-N-P-S,N-Pou-S,P-S}. For the application of similar methods in other models both with regular and singular interaction kernels, see \cite{E-F-R,E-F-S,F-S,F-S-T}. Via a compactness method, the sequence of augmented regularized distribution functions $P^\varepsilon=P^\varepsilon(t,\theta,\omega,\Omega)$ can be shown to have bounded zeroth and first order frequency moments that weakly converge to a weak measure-valued solution of the macroscopic singular system. The main point here is an accurate a priori control on the hierarchy of frequency moments of the second order regularized kinetic description, that in particular, includes time-equicontinuity for the sequence of reduced regularized distribution functions $f^\varepsilon=f^\varepsilon(t,\theta,\Omega)$ in some negative Sobolev space, i.e.,
$$
\Vert f^\varepsilon_{t_1}-f^\varepsilon_{t_2}\Vert_{W^{-1,1}(\mathbb{T}\times \mathbb{R})}\leq C\vert t_1-t_2\vert^{1/2}.
$$

Although we will just illustrate this techniques in the particular setting of the singular Kuramoto model, the author believes that most of the tools can be pushed towards other non-gradient type systems whose internal variables belong, in a natural way, to more general Lie groups. In addition, interactions do not need to be necessarily smooth, but at least enjoy similar one-sided Lipschitz properties or Osgood-type moduli of continuity.

The rest of the paper is structured as follows. Sections \ref{S-macroscopic-model}, \ref{S-weak-solutions-existence}, \ref{S-uniqueness-mean-field} and \ref{S-synchro} are devoted to the subcritical regime $\alpha\in (0,\frac{1}{2})$. Specifically, in Section \ref{S-macroscopic-model} we will introduce the model \eqref{VE} as a nonlinear transport equation along $\mathbb{T}\times \mathbb{R}$, some regularity properties of the tangent transport field $\mathcal{V}[f]$ will be derived and we will also revisit the main concepts of measure-valued solutions in the literature. In Section \ref{S-weak-solutions-existence} we shall prove the existence of global-in-time weak measure-valued solutions to \eqref{VE}. The above-mentioned Dobrushin-type estimates will be discussed in Section \ref{S-uniqueness-mean-field}, in particular their applications to obtain uniqueness results of weak measure-valued solutions and the mean-field limit from the particle system \eqref{E-kuramoto-discrete}-\eqref{E-kuramoto-kernel} towards \eqref{VE}. Finally, we will prove emergence of global phase synchronization in finite time for this subcritical regime under appropriate initial assumptions. Section \ref{S-critical-regime} will focus on the critical case $\alpha=\frac{1}{2}$. Specifically, we will adapt the above methods to obtain measure-valued solutions in Filippov's sense and we will provide analogue results regarding stability, uniqueness, mean field limit and emergence of synchronization. In Section \ref{S-supercritical-regime} we will explore the aforementioned singular hyperbolic limit of vanishing inertia type for the supercritical regime $\alpha\in (\frac{1}{2},1)$. The appendices are devoted to some technical observations that we will use throughout the paper. In Appendix \ref{appendix-relation-measure-spaces} we will recall some notation about measures along $\mathbb{T}$ or periodic measures. In Appendix \ref{appendix-Prokhorov} we shall recall some weak Banach-valued version of the Ascoli--Arzel\`a theorem that we be used in Section \ref{S-supercritical-regime}. Finally, Appendix \ref{appendix-differentiability-distance} recalls the differentiability properties of the squared Riemannian distance in a complete Riemannian manifold. In particular, we will introduce the concept of \textit{one-sided upper Dini differentiability} that we will sometimes use to deal with quadratic Wasserstein distances along manifolds.

Let us end this part by summarizing some notation to be used along the paper:
\begin{enumerate}
\item We will set $\xi\in C^\infty_c([0,\infty))$ to be any non-increasing cut-off function  verifying the following properties
\begin{enumerate}
\item $0\leq \xi\leq 1$, 
\item $\left.\xi\right\vert_{[0,1]}\equiv 1$ and $\left.\xi\right\vert_{[2,\infty]}\equiv 0$,
\item $\int_0^\infty\eta'(r)\,dr=-1.$
\end{enumerate}
Also, we will use the following scaled cut-off functions for every $\varepsilon>0$
\begin{equation}\label{E-scaled-cut-off}
\xi_\varepsilon(r):=\xi(\varepsilon^{-1}r),\ r\in [0,+\infty).
\end{equation}
\item For the product Riemannian manifold $\mathbb{T}\times \mathbb{R}$ we will denote the projections
\begin{equation}\label{E-projections}
\begin{array}{cccccccccc}
\pi_z&:\mathbb{T}\times \mathbb{R}&\longrightarrow&\mathbb{T}, & & \pi_\Omega&:\mathbb{T}\times\mathbb{R}&\longrightarrow&\mathbb{R},\\
 & (z,\Omega)&\longmapsto &z, & & & (z,\Omega)&\longmapsto& \Omega.
 \end{array}
\end{equation}
\item For any locally compact topological space $E$, $\mathcal{M}(E)$ will denote the space of all finite Radon measures on $E$. In addition, $\mathcal{M}^+(E)$ stands for the subset of non-negative measures and $\mathcal{P}(E)$ is the subset of probability measures.
\item For any complete Riemannian manifold $(M,\left<\cdot,\cdot\right>)$ and any $x\in M$, we denote the tangent space by $T_xM$ and $\left<\cdot,\cdot\right>_x:T_xM\rightarrow \mathbb{R}$ is the metric at $x$. Also, we denote the Riemannian distance between any couple of points $x,y\in M$ by
\begin{equation}\label{E-riemannian-distance}
d(x,y):=\inf\left\{\mathcal{L}[\gamma]:=\int_0^1\vert \gamma'(s)\vert\,ds:\,\gamma\in C^1([0,1],M),\ \gamma(0)=x, \gamma(1)=y\right\}.
\end{equation}
By the Hopf--Rinow theorem, such infimum is attained by some minimizing geodesic, see \cite{DC}.
\end{enumerate}

\section{Tangent transport field along $\mathbb{T}\times \mathbb{R}$ and measure-valued solutions}\label{S-macroscopic-model}
In this section we will introduce some tools and notation for the model that we will use along the paper with regards to the subcritical regime of the singularity $\alpha\in (0,\frac{1}{2})$. On the one hand, we will reformulate the macroscopic system as a genuine transport equation along a manifold and will focus on deriving properties of the tangent transport field. That will be the cornerstone to derive the well-posedness of global-in-time measure-valued solutions to the system with $\alpha\in (0,\frac{1}{2})$ in the forthcoming Sections \ref{S-weak-solutions-existence} and \ref{S-uniqueness-mean-field}. On the other hand, we will revisit different equivalent concepts of measure-valued solutions ranging from weak measure-valued solutions to superposition solutions and solutions in the sense of the characteristic flow. Although the regularity of the tangent transport field will vary (or even may not make sense at all) in the critical and supercritical regimes, we will later try to adapt this ideas to those more singular cases in Sections \ref{S-critical-regime} and \ref{S-supercritical-regime}.

\subsection{Formal derivation of the Vlasov equation}\label{SS-formal-mean-field-limit}

In this paper we will not address the propagation of chaos problem. For that issue, the reader may be interested in the following related literature \cite{H-J,H-M,J,J-W,J-W-2,Ka,M-M,M-M-W,S}. Nevertheless, let us recall a  formal derivation of \eqref{VE} that may arise form the propagation of chaos in the system. Specifically, since the natural frequencies in the discrete model \eqref{E-kuramoto-discrete} are constant parameters, we can equivalently state the system as follows  
$$
\left\{
\begin{array}{l}
\displaystyle\dot{\theta}_i=\Omega_i+\frac{K}{N}\sum_{j=1}^Nh(\theta_j-\theta_i),\\
\displaystyle \dot{\Omega}_i=0,\\
\vspace{-0.2cm}\\
\displaystyle\theta_i(0)=\theta_{i,0},\ \Omega_i(0)=\Omega_{i,0}\equiv \Omega_i.
\end{array}
\right.
$$
Then, not only $\theta_i$ is regarded as a mechanical variable but also $\Omega_i$ is. Consequently, the BBGKY hierarchy of Liouville equations for the joint laws $F^N_t=F^N_t(x_1,\ldots,x_N,\Omega_1,\ldots,\Omega_N)\in \mathcal{P}_{sym}(\mathbb{R}^N\times \mathbb{R}^N)$ reads
\begin{equation}\label{E-BBGKY}
\frac{\partial F^N}{\partial t}+\sum_{i=1}^N\frac{\partial}{\partial\theta_i}\left(\left(\Omega_i+\frac{K}{N}\sum_{j=1}^Nh(\theta_j-\theta_i)\right)F^N\right)=0.
\end{equation}
Let us consider the projection map onto the first $k\in \{1,\ldots,N\}$ variables, that is,
$$\begin{array}{cccc}
\pi^{k,N}:&\mathbb{R}^N\times \mathbb{R}^N&\longrightarrow&\mathbb{R}^k\times \mathbb{R}^k,\\
 & (\Theta^N,\Omega^N) & \longmapsto & (\Theta^{k,N},\Omega^{k,N}),
\end{array}
$$
where we denote $\Theta^{k,N}:=(\theta_1,\ldots,\theta_k)$ and $\Omega^{k,N}:=(\Omega_1,\ldots,\Omega_k)$, for any $\Theta^N=(\theta_1,\ldots,\theta_N)\in \mathbb{R}^N$ and $\Omega^N=(\Omega_1,\ldots,\Omega_N)\in \mathbb{R}^N$. Then, we can consider the marginal measures $F^{k,N}_t:=\pi^{k,N}_{\#}(F^N_t)\in \mathcal{P}_{sym}(\mathbb{R}^k\times \mathbb{R}^k)$. Thanks to the assumed symmetry in the system, integration in \eqref{E-BBGKY} yields
\begin{align}
\frac{\partial F^{k,N}}{\partial t}+\sum_{i=1}^k\frac{\partial}{\partial \theta_i}&\left(\left(\Omega_i+\frac{K}{N}\sum_{j=1}^kh(\theta_j-\theta_i)\right)F^{k,N}\right.	\nonumber\\
&\hspace{0.3cm}\left.+K\frac{N-k}{N}\int_{\mathbb{R}\times\mathbb{R}}h(\theta_{k+1}-\theta_i)\,d_{(\theta_{k+1},\Omega_{k+1})}F^{k+1,N}\right)=0.\label{E-BBGKY-marginal}
\end{align}
Observe that the hierarchy is not necessarily closed yet. Via a diagonal argument we can obtain weak limits of an appropriate subsequence (that we denote in the same manner)
$$F^{k,\infty}:=\mbox{weak}\,*-\lim_{N\rightarrow\infty} F^{k,N}.$$
For simplicity of the notation, let us denote $F:=F^{1,\infty}$. Let us assume that all the initial values are tensorized, that is $F^{k,\infty}_0=F^{\otimes^k}_0$. Then, \textit{propagation of chaos} means that such tensorization remains true for all times, i.e., 
$$F^{k,\infty}_t=F^{\otimes^k}_t,\ \mbox{ for all }\ t\geq 0.$$
Conditionally under such property, we can pass to the limit as $N\rightarrow\infty$ in \eqref{E-BBGKY-marginal} for $F^{1,N}$ and close an equation for $F$ as follows
$$
\frac{\partial F}{\partial t}+\frac{\partial}{\partial\theta}\left(\left(\Omega+K\int_{\mathbb{R}^2}h(\theta'-\theta)\,d_{(\theta',\Omega')}F\right)F\right)=0.
$$
This is the classical Vlasov equation of the system. See the aforementioned references where propagation of chaos is proved for some systems. Finally, since the phases $\theta\in \mathbb{R}$ only make sense modulo $2\pi$, we can consider the following map projecting each $\theta\in \mathbb{R}$ into its representative $\bar \theta\in (-\pi,\pi]$ modulo $2\pi$
$$
\begin{array}{cccc}
\bar\pi:&\mathbb{R}&\longrightarrow & (-\pi,\pi],\\
 & \theta & \longmapsto & \bar \theta.
\end{array}
$$
It clearly maps $F_t$ into $f_t\in \mathcal{P}((-\pi,\pi]\times \mathbb{R})$ through $\bar\pi_{\#}F_t=f_t$, and $f=f_t$ obviously fulfils \eqref{VE} in distributional sense. As previously stated, we are not interested in the propagation of chaos topic, but rather on the mean field limit via the empirical measure technique that will be described throughout the paper.

\subsection{Reformulating a nonlinear transport equation along a manifold}
In the above part, the kinetic singular Kuramoto model \eqref{VE} was formally derived. Solutions are regarded as periodic measures $f_t\in \mathcal{P}((-\pi,\pi]\times \mathbb{R})$ with respect to $\theta$. It is clear that we can equivalently regard them as measures $f_t\in \mathcal{P}(\mathbb{T}\times \mathbb{R})$ by virtue of the identification of the interval $(-\pi,\pi]$ with the torus $\mathbb{T}$, see Theorem \ref{theo-relation-measure-spaces-isometry} in Appendix \ref{appendix-relation-measure-spaces}. Here on, we will think of both spaces as the same space and will change notation from one to another without any notice for simplicity of arguments. Before going further in our results, we will introduce some basic properties and notation with regards to the Riemannian manifold $\mathbb{T}\times \mathbb{R}$. In the first result we comment on the Riemannian distance, tangent space and vector fields in such Riemannian manifold.

\begin{lem}\label{L-TxR-notation-1}
Consider the Riemannian manifold $\mathbb{T}\times \mathbb{R}$ endowed with the standard metric. Then,
\begin{enumerate}
\item The Riemannian distance in $\mathbb{T}\times \mathbb{R}$ between any couple $(z_1,\Omega_1),(z_2,\Omega_2)\in\mathbb{T}\times \mathbb{R}$ is defined as
$$
d((z_1,\Omega_1),(z_2,\Omega_2))=\left(d(z_1,z_2)^2+(\Omega_1-\Omega_2)^2\right)^{1/2}=\left(\vert \theta_1-\theta_2\vert_o^2+(\Omega_1-\Omega_2)^2\right)^{1/2},
$$
for any $\theta_1,\theta_2\in\mathbb{R}$ such that $z_1=e^{i\theta_1}$ and $z_2=e^{i\theta_2}$. Here, $\vert \cdot\vert_o$ means the orthodromic distance in the unit torus, that is $\vert \theta\vert_o=\vert \bar\theta\vert$, where $\bar \theta$ is the representative modulo $2\pi$ of $\theta$ in $(-\pi,\pi]$.
\item The tangent space at $(z,\Omega)\in \mathbb{T}\times \mathbb{R}$ reads
$$T_{(z,\Omega)}(\mathbb{T}\times \mathbb{R})=T_{z}\mathbb{T}\times T_\Omega\mathbb{R}=\{(p\,iz,q):\,p,q\in\mathbb{R}\}.$$
As a consequence, the space $\mathfrak{X}\mathcal{C}(\mathbb{T}\times \mathbb{R})$ of $\mathcal{C}$-regular tangent vectors along $\mathbb{T}\times \mathbb{R}$, for any given regularity class $\mathcal{C}$, consist of the fields $V$ with components
\begin{equation}\label{E-fields}
V_{(z,\Omega)}=(P(z,\Omega)\,iz,Q(z,\Omega)),
\end{equation}
for a couple of scalar functions $P,Q:\mathbb{T}\times \mathbb{R}\longrightarrow\mathbb{R}$ with regularity in the class $\mathcal{C}$.
\item Consider $V\in \mathfrak{X}C^1(\mathbb{T}\times \mathbb{R})$ given by \eqref{E-fields} for some $P,Q\in C^1(\mathbb{T}\times \mathbb{R})$. Then,
\begin{equation}\label{E-fields-div}
\divop V= \frac{\partial P}{\partial \theta}+\frac{\partial Q}{\partial\Omega}\equiv iz\frac{\partial P}{\partial z}+\frac{\partial Q}{\partial\Omega},
\end{equation}
where we have used the identifications in the Appendix \ref{appendix-relation-measure-spaces} and the notation of the complex derivatives in Definition \ref{def-complex-derivatives}. The same formula \eqref{E-fields-div} holds true for distributional derivatives.
\end{enumerate}
\end{lem}

In the second lemma, we introduce the geodesics and parallel transport in $\mathbb{T}\times \mathbb{R}$ that will be of interest in some upcoming results.

\begin{lem}\label{L-TxR-notation-2}
Let us set any point $x=(z,\Omega)\in\mathbb{T}\times \mathbb{R}$, where $z=e^{i\theta}$ for some $\theta\in \mathbb{R}$, and any tangent vector $v=(p\,iz,q)\in T_{(z,\Omega)}\mathbb{T}\times \mathbb{R}$.
\begin{enumerate}
\item Define the geodesic $\gamma_{z,p}$ of $\mathbb{T}$ issued at $z$ in the direction $p\,iz$, i.e.,
$$\gamma_{z,p}(s):=e^{i(\theta+ps)}.$$
Then, the geodesic $\widehat{\gamma}_{x,v}$ of $\mathbb{T}\times \mathbb{R}$ issued at $x$ in the direction $v$ reads
$$\widehat{\gamma}_{x,v}(s)=(\gamma_{z,p}(s),\Omega+qs)=(e^{i(\theta+ps)},\Omega+qs).$$
In particular, the Riemannian exponential map reads
$$\widehat{\exp}_x(v)=(\exp_z(p\,iz),\Omega+q)=(e^{i(\theta+p)},\Omega+q).$$
\item Let $\widehat{\gamma}_{x,v}$ be the associated geodesic, and set $s_1,s_2\in \mathbb{R}$. Then, the parallel transport from $\widehat{\gamma}_{x,v}(s_1)$ to $\widehat{\gamma}_{x,v}(s_2)$ along $\widehat{\gamma}_{x,v}$ is the linear isometry (see \cite{DC})
$$\begin{array}{cccc}
\displaystyle\mathcal{\tau}[\widehat{\gamma}_{x,v}]_{s_1}^{s_2}:&T_{\widehat{\gamma}_{x,v}(s_1)}(\mathbb{T}\times \mathbb{R})&\longrightarrow &T_{\widehat{\gamma}_{x,v}(s_2)}(\mathbb{T}\times \mathbb{R}),\\
\displaystyle & (p'i\gamma_{z,p}(s_1),q') & \longmapsto & (p'i \gamma_{z,p}(s_2),q').
\end{array}$$

\end{enumerate}
\end{lem}

Since the proofs are standard, we omit them here. 

\begin{defi}\label{D-transport-field}
Consider $\alpha\in(0,\frac{1}{2})$ and $K>0$. We  (formally) define the function $\mathcal{P}[\mu]$ and the tangent vector field $\mathcal{V}[\mu]$ along the manifold $\mathbb{T}\times \mathbb{R}$  by
\begin{align*}
\mathcal{P}[\mu](\theta,\Omega)&:=\Omega-K\int_{\mathbb{T}}\int_{\mathbb{R}}h(\theta-\theta')\,d_{(\theta',\Omega')}\mu,\\
\mathcal{V}[\mu](z,\Omega)&:=(\mathcal{P}[\mu](z,\Omega)\,iz,0),
\end{align*}
where $\mu\in \mathcal{M}(\mathbb{T}\times \mathbb{R})$ is any finite Radon measure.
\end{defi}

Using such notation one can easily check that the kinetic singular Kuramoto model can be restated as the following nonlinear transport equation in conservative form along $\mathbb{T}\times\mathbb{R}$:
\begin{equation}\label{E-kuramoto-transport-TxR}
\left\{
\begin{array}{l}
\displaystyle\frac{\partial f}{\partial t}+\divop_{(z,\Omega)}(\mathcal{V}[f]f)=0,\\
\displaystyle f(0,\cdot,\cdot)=f^0.
\end{array}
\right.
\end{equation}
Of course, now there is no need to impose explicit periodicity conditions because they are actually considered  implicitly in the geometry of the space $\mathbb{T}\times \mathbb{R}$.

\subsection{Properties of the transport field}\label{SS-transport-field-properties}

Here on, most of our effort will be devoted to derive some one-sided modulus of continuity of $\mathcal{V}[f]$. In particular, such (weak) regularity will entail that the characteristic flow associated with $\mathcal{V}[f]$ is well--defined forwards-in-time, for any weak measure-valued solution to \eqref{E-10}. This will be the cornerstone to show well-posedness of global-in-time weak measure valued solutions to \eqref{E-kuramoto-transport-TxR} in Sections \ref{S-weak-solutions-existence} and \ref{S-uniqueness-mean-field}. Before moving to such regularity issues of the transport field $\mathcal{V}[f]$, let us set the appropriate spaces of time-dependent measures.

\begin{defi}\label{D-spaces-measures}
Fix $T>0$, then we will consider
\begin{align*}
\mathcal{C}_\mathcal{M}(0,T)&:=\left\{\mu \in L^\infty(0,T;\mathcal{M}(\mathbb{T}\times \mathbb{R})):\,t\longmapsto \int_{\mathbb{T}\times \mathbb{R}}\varphi\,d\mu_t\in C([0,T]),\ \forall\,\varphi\in C_c(\mathbb{T}\times \mathbb{R})\right\},\\
\widetilde{C}_\mathcal{M}(0,T)&:=\left\{\mu \in L^\infty(0,T;\mathcal{M}(\mathbb{T}\times \mathbb{R})):\,t\longmapsto \int_{\mathbb{T}\times \mathbb{R}}\varphi\,d\mu_t\in C([0,T]),\ \forall\,\varphi\in C_b(\mathbb{T}\times \mathbb{R})\right\},\\
\mathcal{T}_\mathcal{M}(0,T)&:=\left\{\mu\in L^\infty(0,T;\mathcal{M}(\mathbb{T}\times \mathbb{R})):\,(\mu_t)_{t\in [0,T]}\mbox{ is uniformly tight}\right\},\\
\mathcal{AC}_\mathcal{M}(0,T)&:=\left\{\mu \in L^\infty(0,T;\mathcal{M}(\mathbb{T}\times \mathbb{R})):\,t\longmapsto \int_{\mathbb{T}\times \mathbb{R}}\varphi\,d\mu_t\in AC(0,T),\ \forall\,\varphi\in C_c^\infty(\mathbb{T}\times \mathbb{R})\right\}.
\end{align*}
For simplicity of the notation, we will sometimes remove the dependence on $T$ when it is clear.
\end{defi}

Some properties about the preceding spaces of measures are in order:

\begin{pro}
For the spaces in Definition \ref{D-spaces-measures}, the following properties hold true:
\begin{enumerate}
\item The spaces $\mathcal{C}_\mathcal{M}$ and $\widetilde{\mathcal{C}}_\mathcal{M}$ can be represented as follows
\begin{align*}
\mathcal{C}_\mathcal{M}&=C([0,T],\mathcal{M}(\mathbb{T}\times \mathbb{R})-\mbox{weak}\,*),\\
\widetilde{\mathcal{C}}_\mathcal{M}&=C([0,T],\mathcal{M}(\mathbb{T}\times \mathbb{R})-\mbox{narrow}).
\end{align*}
\item In the above definition of $\mathcal{C}_\mathcal{M}$, test functions $\varphi\in C_c(\mathbb{T}\times \mathbb{R})$ can be replaced by any intermediate regularity class being dense in $C_0(\mathbb{T}\times \mathbb{R})$, e.g.,
$$C_c^\infty(\mathbb{T}\times \mathbb{R}),\ C^k_c(\mathbb{T}\times \mathbb{R})\mbox{ and }W^{k,p}(\mathbb{T}\times \mathbb{R}),\ \mbox{ with }k\in\mathbb{N},\ p>2.$$
\item The space $\widetilde{\mathcal{C}}_\mathcal{M}$ can be represented as follows
$$\widetilde{\mathcal{C}}_\mathcal{M}=\mathcal{C}_\mathcal{M}\cap\mathcal{T}_\mathcal{M}.$$
\item The next embeddings take place
$$\mathcal{AC}_\mathcal{M}\subseteq \mathcal{C}_\mathcal{M}\ \mbox{ and }\ \mathcal{AC}_\mathcal{M}\cap \mathcal{T}_\mathcal{M}\subseteq \widetilde{\mathcal{C}}_\mathcal{M}.$$
\item The next embedding takes place for any $1\leq p,q\leq \infty$ 
$$L^\infty(0,T;\mathcal{M}(\mathbb{T}\times \mathbb{R}))\cap W^{1,q}(0,T;W^{-1,p'}(\mathbb{T}\times \mathbb{R})-\mbox{weak}\,*)\subseteq \mathcal{AC}_\mathcal{M}.$$
\end{enumerate}
\end{pro}

\begin{proof}
The first assertion is clear and the second one is a straightforward density argument of the set of smooth and compactly supported functions $C^\infty_c(\mathbb{T}\times \mathbb{R})$ in $C_0(\mathbb{T}\times \mathbb{R})$. The third item is nothing but Prokhorov's compactness theorem, whilst the forth is clear by definition. The last claim follows from Sobolev's embedding theorem in one dimension.
\end{proof}

Let us now recall the concept of weak measure-valued solution to \eqref{E-kuramoto-transport-TxR}.

\begin{defi}\label{D-weak-measure-solution}
We will say that $f\in \mathcal{C}_\mathcal{M}$ is a weak measure-valued solution to \eqref{E-kuramoto-transport-TxR} when
$$\int_0^T\int_{\mathbb{T\times \mathbb{R}}}\frac{\partial\varphi}{\partial t}\,d_{(z,\Omega)}f_t\,dt+\int_0^T\int_{\mathbb{T}\times \mathbb{R}}\left<\mathcal{V}[f_t],\nabla_{(z,\Omega)}\varphi\right>\,d_{(z,\Omega)}f_t\,dt=-\int_{\mathbb{T}\times \mathbb{R}}\varphi(0,z,\Omega)d_{(z,\Omega)}f^0,$$
for every $\varphi\in C^1_c([0,T)\times \mathbb{T}\times \mathbb{R})$.
\end{defi}

Notice that the nonlinear term in Definition \ref{D-weak-measure-solution} makes sense for any $f\in \mathcal{C}_\mathcal{M}$, no matter whether $f$ also belongs to $\widetilde{C}_\mathcal{M}$ or $\mathcal{AC}_\mathcal{M}$. However, as it is the case for many other models, solutions end up being more regular in time than simply $\mathcal{C}_\mathcal{M}$. In such case, we can restate the above weak formulation for $f \in \mathcal{AC}_\mathcal{M}$ as follows.

\begin{pro}\label{P-weak-measure-solution-charact}
Consider $\alpha\in(0,1)$, $K>0$ and fix $f\in \mathcal{AC}_\mathcal{M}$. Then, the following two statement are equivalent:
\begin{enumerate}
\item $f$ is a weak measure-valued solution to \eqref{E-kuramoto-transport-TxR} in the sense of Definition \ref{D-weak-measure-solution}.
\item The following identity holds
\begin{equation}\label{E-weak-measure-solution-charact}
\displaystyle \frac{d}{dt}\int_{\mathbb{T}\times \mathbb{R}}\phi\,d_{(z,\Omega)}f_t=\int_{\mathbb{T}\times \mathbb{R}}\left<\mathcal{V}[f_t],\nabla_{(z,\Omega)}\phi\right>d_{(z,\Omega)}f_t,
\end{equation}
for a.e. $t\in [0,T]$ and $f(0,\cdot,\cdot)=f^0$, for every $\phi\in C^\infty_c(\mathbb{T}\times \mathbb{R})$.
\end{enumerate}
\end{pro}

\begin{proof}
First, let us assume that $f$ solves \eqref{E-kuramoto-transport-TxR} in the sense of Definition \ref{D-weak-measure-solution}. Take any $\phi\in C^\infty_c(\mathbb{T}\times\mathbb{R})$, consider $\eta\in C^\infty_c(0,T)$ and define the test function
$$\varphi(t,z,\Omega):=\eta(t)\phi(z,\Omega).$$
Since $\varphi\in C^\infty_c((0,T)\times \mathbb{T}\times \mathbb{R})$, then Definition \ref{D-weak-measure-solution} entails
$$\int_0^T\eta'(t)\left(\int_{\mathbb{T}\times\mathbb{R}}\phi\,d_{(z,\Omega)}f_t\right)\,dt+\int_0^T\eta (t)\left(\int_{\mathbb{T}\times \mathbb{R}}\left<\mathcal{V}[f_t],\nabla_{(z,\Omega)}\phi\right>\,d_{(z,\Omega)}f_t\right)\,dt=0.$$
By the arbitrariness of $\eta$, we can identify the weak derivative of $t\mapsto\int_{\mathbb{T}\times \mathbb{R}}\phi\,d_{(z,\Omega)}f_t$ (that exists since $f\in\mathcal{AC}_\mathcal{M}$) as the right-hand side in \eqref{E-weak-measure-solution-charact} and it concludes the first part. Conversely, let us assume that $f$ is a solution in the sense of Equation \eqref{E-weak-measure-solution-charact} and consider a test function $\varphi\in C^\infty_c([0,T)\times \mathbb{T}\times \mathbb{R})$. By density, we can assume that it has separate variables, i.e.,
$$\varphi(t,z,\Omega)=\eta(t)\phi(z,\Omega),$$
for $\eta\in C^\infty_c([0,T))$ and $\phi\in C^\infty_c(\mathbb{T}\times \mathbb{R})$. Consider any non-increasing cut-off functions $\xi_\varepsilon\in C^\infty_c([0,\infty))$ like in \eqref{E-scaled-cut-off} for any $\varepsilon>0$ and consider the smooth approximate test functions 
$$\eta_\varepsilon(t):=(1-\xi_\varepsilon(t))\eta(t).$$
It is clear that $\eta_\varepsilon\in C^\infty_c(0,T)$ and, consequently, witting \eqref{E-weak-measure-solution-charact} in weak form amounts to
$$\int_0^T\eta'_\varepsilon(t)\left(\int_{\mathbb{T}\times\mathbb{R}}\phi\,d_{(z,\Omega)}f_t\right)\,dt+\int_0^T\eta_\varepsilon (t)\left(\int_{\mathbb{T}\times \mathbb{R}}\left<\mathcal{V}[f_t],\nabla_{(z,\Omega)}\phi\right>\,d_{(z,\Omega)}f_t\right)\,dt=0,$$
i.e., expanding the derivatives of $\eta_\varepsilon$,
\begin{multline*}
\int_0^T(1-\xi_\varepsilon(t))\eta'(t)\left(\int_{\mathbb{T}\times\mathbb{R}}\phi\,d_{(z,\Omega)}f_t\right)\,dt\\
+\int_0^T(1-\xi_\varepsilon(t))\eta(t)\left(\int_{\mathbb{T}\times \mathbb{R}}\left<\mathcal{V}[f_t],\nabla_{(z,\Omega)}\phi\right>\,d_{(z,\Omega)}f_t\right)\,dt\\
=\int_0^T\xi'_\varepsilon(t)\eta(t)\left(\int_{\mathbb{T}\times \mathbb{R}}\phi\,d_{(z,\Omega)}f_t\right)\,dt.
\end{multline*}
Notice that $1-\xi_\varepsilon\rightarrow 1$ in $C([0,T])$ and $\xi'_\varepsilon\overset{*}{\rightharpoonup}-\delta_0$ in $\mathcal{M}([0,T))$ as $\varepsilon\rightarrow 0$. Then, taking limits as $\varepsilon\rightarrow 0$ in the above identities yields to the weak formulation in Definition \ref{D-weak-measure-solution}.
\end{proof}

Notice that for $f\in \mathcal{C}_\mathcal{M}$ (also for $f\in \widetilde{\mathcal{C}}_\mathcal{M}$) we cannot expect $\mathcal{V}[f]$ to be fully Lipschitz-continuous since $h$ is expected to be barely H\"{o}lder-continuous and this might cause severe problems. Before introducing sharper regularity properties, let us comment on the basic properties that we can infer from the uniform continuity of the kernel $h$.

\begin{lem}\label{L-transport-field-holder}
Consider $\alpha\in (0,\frac{1}{2})$. Then, $h$ is $(1-2\alpha)$-H\"{o}lder continuous in $\mathbb{T}$, namely,
$$\vert h(z_1)-h(z_2)\vert \leq \cosh\pi\,d(z_1,z_2)^{1-2\alpha},$$
for every couple $z_1,z_2\in \mathbb{T}$.
\end{lem}

\begin{proof}
Consider $z_1=e^{i\theta_1},z_2=e^{i\theta_2}\in \mathbb{T}$ for some $\theta_1,\theta_2\in \mathbb{R}$. Taking appropriate representatives for the phases, we can assume that $\theta_1-\theta_2\in (-\pi,\pi]$ without loss of generality. To simplify the proof, we will divide it into two distinguished cases:

$\bullet$ \textit{Case 1}: $\theta_1-\theta_2\in [0,\pi]$. By the Taylor expansion of the sine function, we obtain
$$h(\theta_1)-h(\theta_2)=\frac{\sin\theta_1}{\theta_1^{2\alpha}}-\frac{\sin \theta_2}{\theta_2^{2\alpha}}=\sum_{n=0}^\infty\frac{(-1)^n}{(2n+1)!}(\theta_1^{2n+1-2\alpha}-\theta_2^{2n+1-2\alpha}).$$
Since $\alpha\in (0,\frac{1}{2})$, then $1-2\alpha\in (0,1)$ and, consequently, for $n=0$ we infer the next estimate
$$\vert\theta_1^{1-2\alpha}-\theta_2^{1-2\alpha}\vert\leq \vert \theta_1-\theta_2\vert^{1-2\alpha}.$$
For any other $n\geq 1$ we apply the mean value theorem to obtain
\begin{align*}
\vert \theta_1^{2n+1-2\alpha}-\theta_2^{2n+1-2\alpha}\vert&\leq (2n+1-2\alpha)\pi^{2n-2\alpha}\vert \theta_1-\theta_2\vert\\
&\leq (2n+1-2\alpha)\pi^{2n}\vert \theta_1-\theta_2\vert^{1-2\alpha}\\
&\leq (2n+1)\pi^{2n}\vert \theta_1-\theta_2\vert^{1-2\alpha}.
\end{align*}
Putting everything together, we achieve the desired estimate
$$\vert h(\theta_1)-h(\theta_2)\vert\leq \sum_{n=0}^{\infty}\frac{\pi^{2n}}{(2n)!}\vert \theta_1-\theta_2\vert^{1-2\alpha}=\cosh\pi\, \vert \theta_1-\theta_2\vert^{1-2\alpha}.$$

$\bullet$ \textit{Case 2}: $\theta_1-\theta_2\in (-\pi,0]$. By the antisymmetry of the kernel $h$ with respect to the origin, we can reduce this case to the latter one, thus we omit the proof.

To end the proof notice that $\theta_1-\theta_2\in (-\pi,\pi]$ and, consequently, $\vert \theta_1-\theta_2\vert=\vert \theta_1-\theta_2\vert_o=d(z_1,z_2)$.
\end{proof}

\begin{theo}\label{T-transport-field-holder}
Consider $\alpha\in (0,\frac{1}{2})$, $K>0$ and set $\mu \in\widetilde{\mathcal{C}}_\mathcal{M}$. Then, 
$$\frac{\mathcal{P}[\mu]}{1+\vert \Omega\vert}\in C_b([0,T]\times\mathbb{T}\times \mathbb{R}).$$
In addition, there exists $C>0$, that does not depend on $\mu$, such that
$$\left\vert \mathcal{P}[\mu_t](z_1,\Omega_1)-\mathcal{P}[\mu_t](z_2,\Omega_2)\right\vert \leq \vert \Omega_1-\Omega_2\vert+CK\Vert \mu\Vert_{L^\infty(0,T;\mathcal{M}(\mathbb{T}\times \mathbb{R}))}d(z_1,z_2)^{1-2\alpha},$$
for every $t\in [0,T]$ and $(z_1,\Omega_1),(z_2,\Omega_2)\in\mathbb{T}\times \mathbb{R}$. 
\end{theo}

\begin{proof}
First, let us show the second property. Fix $t\in [0,T]$ and $(z_1,\Omega_1),(z_2,\Omega_2)\in \mathbb{T}\times \mathbb{R}$ and notice that
$$\mathcal{P}[\mu_t](z_1,\Omega_1)-\mathcal{P}[\mu_t](z_2,\Omega_2)=\Omega_1-\Omega_2-K\int_{\mathbb{T}\times \mathbb{R}}(h(z_1\overline{z'})-h(z_2\overline{z'}))d_{(z',\Omega')}\mu_t.$$
Then, the triangle inequality together with Lemma \ref{L-transport-field-holder} imply
$$\vert \mathcal{P}[\mu_t](z_1,\Omega_1)-\mathcal{P}[\mu_t](z_2,\Omega_2)\vert\leq \vert \Omega_1-\Omega_2\vert+K\cosh\pi\int_{\mathbb{T}\times \mathbb{R}}d(z_1\overline{z'},z_2\overline{z'})^{1-2\alpha}d_{(z',\Omega')}\mu_T.$$
Since the Riemannian distance along the torus $\mathbb{T}$ is translation invariant, then $d(z_1\overline{z'},z_2\overline{z'})=d(z_1,z_2)$ for every $z'\in\mathbb{T}$, thus yielding
$$\vert \mathcal{P}[\mu_t](z_1,\Omega_1)-\mathcal{P}[\mu_t](z_2,\Omega_2)\vert\leq \vert \Omega_1-\Omega_2\vert+K\cosh\pi \Vert \mu\Vert_{L^\infty(0,T;\mathcal{M}(\mathbb{T}\times \mathbb{R}))}d(z_1,z_2)^{1-2\alpha}.$$
Second, let us proof the full continuity in all the variables. Consider $t\in [0,T]$ and $(z,\Omega)\in \mathbb{T}\times \mathbb{R}$ ant let $\{t_n\}_{n\in\mathbb{N}}\subseteq [0,T]$ and $\{(z_n,\Omega_n)\}_{n\in \mathbb{N}}\subseteq \mathbb{T}\times\mathbb{R}$ such that $t_n\rightarrow t$ and $(z_n,\Omega_n)\rightarrow (z,\Omega)$. We can split
$$\mathcal{V}[\mu_{t_n}](z_n,\Omega_n)-\mathcal{V}[\mu_{t}](z,\Omega)=A_n+B_n,$$
where each term reads
\begin{align*}
A_n&:=\mathcal{V}[\mu_{t_n}](z_n,\Omega_n)-\mathcal{V}[\mu_{t_n}](z,\Omega),\\
B_n&:=\mathcal{V}[\mu_{t_n}](z,\Omega)-\mathcal{V}[\mu_t](z,\Omega).
\end{align*}
Regarding the first term, the preceding part yields
$$A_n\leq \vert \Omega_n-\Omega\vert+CK\Vert \mu\Vert_{L^\infty(0,T;\mathcal{M}(\mathbb{T}\times \mathbb{R}))}d(z_n,z)^{1-2\alpha},$$
where the convergence $A_n\rightarrow 0$ is clear. On the other hand, the second term reads
$$B_n=-K\int_{\mathbb{T}\times \mathbb{R}}h(z\overline{z'})\,d_{(z',\Omega')}(\mu_{t_n}-\mu_t).$$
Now, the convergence $B_n\rightarrow 0$ follows from the definition of $\widetilde{\mathcal{C}}_\mathcal{M}$ in Definition \ref{D-spaces-measures} and the boundedness and continuity of the map
$$(z',\Omega')\in\mathbb{T}\times \mathbb{R}\longmapsto h(z\overline{z'}).$$
This amounts to the desired continuity property. Regarding the growth estimate notice that
$$\sup_{(z,\Omega)\in\mathbb{T}\times \mathbb{R}}\frac{\left\vert\mathcal{P}[\mu_t](z,\Omega)\right\vert}{1+\vert \Omega\vert}\leq \sup_{(z,\Omega)\in\mathbb{T}\times \mathbb{R}}\frac{\vert \Omega\vert +K\Vert \mu\Vert_{L^\infty(0,T;\mathcal{M}(\mathbb{T}\times \mathbb{R}))}\Vert h\Vert_{C(\mathbb{T})}}{1+\vert \Omega\vert}<\infty.$$
\end{proof}

As a consequence of Lemmas \ref{L-TxR-notation-1} and \ref{L-TxR-notation-2} we can achieve a similar result for $\mathcal{V}[\mu]$. Given that it is a tangent vector field, let us recall the definition of H\"{o}lder-continuity of tangent vector fields along a complete Riemannian manifold.

\begin{defi}
Consider a complete Riemannian manifold $(M,\left<\cdot,\cdot\right>)$, fix $0<\beta\leq 1$ and a tangent vector field $V:M\longrightarrow TM$. $V$ is said to be $\beta$-H\"{o}lder continuous when there exists $C>0$ such that
$$\left\vert \tau[\gamma]_0^1(V_x)-V_y\right\vert \leq C d(x,y)^\beta,$$
for every $x,y\in M$ and any minimizing geodesic $\gamma:[0,1]\longrightarrow M$ joining $x$ to $y$. The smallest such $C$ is called the $\beta$-H\"{o}lder constant of the vector field $V$. Here $d=d(\cdot,\cdot)$ is the Riemannian distance \eqref{E-riemannian-distance} and
$$\tau[\gamma]_{s_1}^{s_2}:T_{\gamma(s_1)}M\longrightarrow T_{\gamma(s_2)}M,$$
stands for the parallel transport from $\gamma(s_1)$ to $\gamma(s_2)$ along the geodesic $\gamma$, see \cite{DC}.
\end{defi}

\begin{cor}\label{C-transport-field-holder}
Consider $\alpha\in (0,\frac{1}{2})$, $K>0$, and set $\mu \in\widetilde{\mathcal{C}}_\mathcal{M}$. Then, 
$$\frac{\mathcal{V}[\mu]}{1+\vert \Omega\vert}\in C([0,T],\mathfrak{X}C_b(\mathbb{T}\times \mathbb{R})).$$
In addition, there exists $C>0$, which does not depend on $\mu$, such that
$$\left\vert \tau_0^1[\widehat{\gamma}](\mathcal{V}[\mu_t](z_1,\Omega_1))-\mathcal{V}[\mu_t](z_2,\Omega_2)\right\vert \leq \vert \Omega_1-\Omega_2\vert+CK\Vert \mu\Vert_{L^\infty(0,T;\mathcal{M}(\mathbb{T}\times \mathbb{R}))}d(z_1,z_2)^{1-\alpha},$$
for every $t\in [0,T]$,  each $(z_1,\Omega_1),(z_2,\Omega_2)\in \mathbb{T}\times \mathbb{R}$ and every minimizing geodesic $\widehat{\gamma}:[0,1]\longrightarrow \mathbb{T}\times \mathbb{R}$ joining $(z_1,\Omega_1)$ to $(z_2,\Omega_2)$.
\end{cor}

\begin{proof}
The first part is clear by the Definition \ref{D-transport-field} along with Lemma \ref{L-TxR-notation-1}. Let us focus on the last part where Lemma \ref{L-TxR-notation-2} will play a role. Set $t\in [0,T]$, $x=(z_1,\Omega_1),\,y=(z_2,\Omega_2)\in \mathbb{T}\times \mathbb{R}$ and a minimizing geodesic $\widehat{\gamma}:[0,1]\longrightarrow \mathbb{T}\times \mathbb{R}$ joining $x$ to $y$. Then,
$$\tau[\widehat{\gamma}]_0^1(\mathcal{V}[\mu_t](z_1,\Omega_1))=(\mathcal{P}[\mu_t](z_1,\Omega_1)\,iz_2,0).$$
Consequently,
\begin{align*}
\vert \tau[\gamma]_0^1(\mathcal{V}[\mu_t](z_1,\Omega_1))-&\mathcal{V}[\mu_t](z_2,\Omega_2)\vert\\
&=\vert ((\mathcal{P}[\mu_t](z_1,\Omega_1)-\mathcal{P}[\mu_t](z_2,\Omega_2))\,iz_2,0) \vert=\vert \mathcal{P}[\mu_t](z_1,\Omega_1)-\mathcal{P}[\mu_t](z_2,\Omega_2)\vert,
\end{align*}
and the result clearly follows from Theorem \ref{T-transport-field-holder}.
\end{proof}

\begin{rem}\label{R-transport-field-holder}
If the tightness condition $\mu\in \widetilde{\mathcal{C}}_\mathcal{M}$ is deprived in Corollary \ref{C-transport-field-holder} and it is replaced by the weaker condition $\mu\in\mathcal{C}_\mathcal{M}$, then time continuity might be lost. However, we still may obtain
$$\frac{\mathcal{V}[f]}{1+\vert \Omega\vert}\in L^\infty(0,T;\mathfrak{X}C_b(\mathbb{T}\times \mathbb{R})).$$
Notice that the time continuity is essential to achieve well-posedness of classical characteristics associated with the transport field $\mathcal{V}[\mu]$. However, in the case $\mu\in \mathcal{C}_\mathcal{M}$, where time continuity is missing (but time boundedness, thus integrability remains), we may still resort on the concept of Caratheodory solution, i.e., absolutely continuous solutions that solve the characteristic system almost everywhere. For simplicity, let us skip it now although it will come into play in the critical regime, see Section \ref{S-critical-regime}.
\end{rem}

\begin{lem}\label{L-existence-characteristic-system}
Consider $\alpha\in (0,\frac{1}{2})$, $K>0$ and fix $\mu\in \widetilde{C}_\mathcal{M}$. For any $x_0=(z_0,\Omega_0)\in \mathbb{T}\times \mathbb{R}$ let us consider the characteristic system issued at $x_0$, i.e.,
\begin{equation}\label{E-characteristic-system-1}
\left\{\begin{array}{l}
\displaystyle\frac{dX}{dt}(t;t_0,x_0)=\mathcal{V}[\mu_t](X(t;t_0,x_0)),\\
\displaystyle \hspace{0.1cm} X(t_0;t_0,x_0)=x_0.
\end{array}\right.
\end{equation}
Then, \eqref{E-characteristic-system-1} enjoys at least one global-in-time $C^1$ solution $X(t;t_0,x_0)=(Z(t;t_0,z_0,\Omega_0),\Omega_0)$. Indeed, if we set $z_0=e^{i\theta_0}$, for some $\theta_0\in \mathbb{R}$, then 
$$Z(t;t_0,z_0,\Omega_0)=e^{i\Theta(t;t_0,\theta_0,\Omega_0)},$$
where $\Theta=\Theta(t;t_0,\theta_0,\Omega_0)$ is a global-in-time $C^1$ solution to
\begin{equation}\label{E-characteristic-system-2}
\left\{
\begin{array}{l}
\displaystyle\frac{d\Theta}{dt}(t;t_0,\theta_0,\Omega_0)=\mathcal{P}[\mu_t](\Theta(t;t_0,\theta_0,\Omega_0),\Omega_0),\\
\displaystyle \hspace{0.1cm}\Theta(t_0;t_0,\theta_0,\Omega_0)=\theta_0.
\end{array}
\right.
\end{equation}
\end{lem}

\begin{proof}
The first part of the result is clear because $\frac{\mathcal{V}[\mu]}{1+\vert \Omega\vert}\in C([0,T],\mathfrak{X}C_b(\mathbb{T}\times\mathbb{R}))$ by virtue of Corollary \ref{C-transport-field-holder}. However, let us comment on the above representation in coordinates. By Theorem \ref{T-transport-field-holder} one also has that $\mathcal{P}[\mu]$, regarded as a function in $\mathbb{R}\times\mathbb{R}$, is a continuous function with sub-linear growth. By the classical Peano theorem, there is at least one global-in-time solution $\Theta=\Theta(t;t_0,\theta_0,\Omega_0)$ to \eqref{E-characteristic-system-2}. Now, let us define $Z(t;t_0,z_0,\Omega_0)$ in terms of $\Theta(t;t_0,z_0,\Omega_0)$ like in the statement of this result. Then,
$$\frac{dZ}{dt}(t;t_0,z_0,\Omega_0)=ie^{i\Theta(t;t_0,\theta_0,\Omega_0)}\frac{d\Theta}{dt}(t;t_0,\theta_0,\Omega_0)=\mathcal{P}[\mu_t](Z(t;t_0,z_0,\Omega_0),\Omega_0)\,iZ(t;t_0,z_0,\Omega_0),$$
and, consequently, so defined $X(t;t_0,x_0)$ is a global-in-time solution to the characteristic system \eqref{E-characteristic-system-1} thanks to the Definition \ref{D-transport-field} of the vector field $\mathcal{V}[\mu]$.
\end{proof}

Notice that in Theorem \ref{T-transport-field-holder}, a $(1-2\alpha)$-H\"{o}lder estimate for $\mathcal{P}[\mu]$ was obtained. Nevertheless, the infinite slope of $h$ at each $\theta\in 2\pi\mathbb{Z}$ prevent us from a full Lipschitz-estimate for $\mathcal{P}[\mu]$. It is well known that H\"{o}lder continuity is not enough for the characteristic system \eqref{E-characteristic-system-2} (or equivalently \eqref{E-characteristic-system-1}) to enjoy a unique global-in-time $C^1$ solution. In the following, we will introduce the key concept that will amount to the forwards uniqueness result, namely, one-sided Lipschitz-continuity.

\begin{defi}\label{D-one-sided-Lipschitz}
Let $(M,\left<\cdot,\cdot\right>)$ be a complete Riemannian manifold and consider a tangent vector field $V$ along $M$. Then, we will say that $V$ is one-sided Lipschitz when there exists a constant $L>0$ such that
\begin{equation}\label{E-one-sided-Lipschitz-1}
\left<V_y,\gamma'(1)\right>-\left<V_x,\gamma'(0)\right>\leq Ld(x,y)^2,
\end{equation}
for every $x,y\in M$ and every minimizing geodesic $\gamma:[0,1]\longrightarrow M$ joining $x$ to $y$.
\end{defi}

\begin{rem}
For such a minimizing geodesic $\gamma$ as in Definition \ref{D-one-sided-Lipschitz}, we can associate the parallel transport from a point $\gamma(s_1)$ to $\gamma(s_2)$ along the geodesic $\gamma$ is a linear isometry between the tangent spaces to $M$ supported at such points
$$\tau[\gamma]_{s_1}^{s_2}:T_{\gamma(s_1)}M\longrightarrow T_{\gamma(s_2)} M.$$
In addition, since $\gamma$ is a geodesic, then the covariant derivative of $\gamma$ vanishes, that is, $\frac{D\gamma'}{ds}=0$ and $\tau[\gamma]_{s_1}^{s_2}(\gamma'(s_1))=\gamma'(s_2)$. Then, the condition \eqref{E-one-sided-Lipschitz-1} can be equivalently restated as follows
\begin{equation}\label{E-one-sided-Lipschitz-2}
-\left<\tau[\gamma]_0^1(V_x)-V_y,\gamma'(1)\right>\leq Ld(x,y)^2.
\end{equation}
\end{rem}

To the best knowledge of the author, such definition has not been clearly proposed previously in the literature as a generalization of the standard one-sided Lipschitz continuity in Euclidean spaces. For the sake of clarity, we list some of the main supporting its definition.

\begin{pro}
Let $(M,g)$ be a complete Riemannian manifold and consider a tangent vector field $V:M\longrightarrow TM$ and a scalar differentiable function $\phi:M\longrightarrow\mathbb{R}$.
\begin{enumerate}
\item If $V$ is Lipschitz-continuous, then $V$ is one-sided Lipschitz.
\item If $M\equiv \mathbb{R}^d$, the $d$-dimensional Euclidean space, then $V$ is one-sided Lipschitz in the sense \eqref{E-one-sided-Lipschitz-1} (equivalently \eqref{E-one-sided-Lipschitz-2}) if, and only if, it is one-sided Lipschitz in the standard sense, i.e.,
$$\left<V_x-V_y,x-y\right>\leq L\vert x-y\vert^2,$$
for every couple of vectors $x,y\in \mathbb{R}^d$.
\item If $\phi$ is $\lambda$-convex, then $-\nabla \phi$ is one-sided Lipschitz with constant $\lambda$.
\end{enumerate}
\end{pro}

\begin{proof}
Consider $x,y\in M$ and any minimizing geodesic $\gamma:[0,1]\longrightarrow M$ joining $x$ to $y$. First, assume that $V$ is Lipschitz-continuous with constant $L$. Then, 
$$-\left<\tau[\gamma]_0^1(V_x)-V_y,\gamma'(1)\right>\leq \vert \tau[\gamma]_0^1(V_x)-V_y\vert\,\vert \gamma'(1)\vert\leq Ld(x,y)^2,$$
where we have used the Cauchy-Schwarz inequality, the Lipschitz continuity of $V$ and that $\vert\gamma'(1)\vert=d(x,y)$ because $\gamma$ is minimizing geodesic between $x$ and $y$. Second, in the Euclidean case it is clear that there is a unique such minimizing geodesic, namely,
$$\gamma(s)=(1-s)x+sy=x+s(y-x),\ s\in [0,1].$$
Notice that $\gamma'(0)=\gamma'(1)=y-x$. Hence, our claim is clear. Finally, assume that $\phi$ is $\lambda$-convex (see, for instance, \cite{A-G-S,F-V,Sa,V}). Then, by definition we obtain
$$\phi(\gamma(s))\leq (1-s)\phi(x)+s\phi(y)+\frac{\lambda}{2}(1-s)sd(x,y)^2,\ s\in [0,1].$$
Then, we find
\begin{align*}
\frac{\phi(\gamma(s))-\phi(x)}{s}&\leq \phi(y)-\phi(x)+\frac{\lambda}{2}(1-s)d(x,y)^2,\\
-\frac{\phi(\gamma(s))-\phi(y)}{s-1}&\leq \phi(x)-\phi(y)+\frac{\lambda}{2}sd(x,y)^2.
\end{align*}
Taking, limits as $s\rightarrow 0^+$ and $s\rightarrow 1^-$ respectively in the first and second expressions yields
\begin{align*}
\left<\nabla \phi(x),\gamma'(0)\right>&\leq \phi(y)-\phi(x)+\frac{\lambda}{2}d(x,y)^2,\\
-\left<\nabla \phi(y),\gamma'(1)\right>&\leq \phi(x)-\phi(y)+\frac{\lambda}{2}d(x,y)^2.
\end{align*}
Taking the sum of both terms, we get that $-\nabla \phi$ is one-sided Lipschitz with constant $\lambda$.
\end{proof}

\begin{figure}
\centering
\begin{subfigure}[b]{0.45\textwidth}
\begin{tikzpicture}[
declare function={
d(\x)= (\x<=-pi) * abs(\x+2*pi) + and(\x>-pi, \x<=pi) * abs(\x) + (\x>pi) * abs(\x-2*pi);
}
]
\begin{axis}[
  axis x line=middle, axis y line=middle,
  xmin=-2*pi, xmax=2*pi, xtick={-6,-4,-2,0,2,4,6}, xlabel=$\theta$,
  ymin=-2, ymax=2, ytick={-1.5,-1,-0.5,0,0.5,1.5},
  legend style={at={(0.8,0.95)}},
]
\addplot [
    domain=-2*pi:2*pi, 
    samples=200, 
    color=black,
]
{-sin(deg(x))/pow(d(x),2*0.25)};
\addlegendentry{$-h(\theta)$}
\end{axis}
\end{tikzpicture}
\caption{$-h(\theta)$}
\label{fig:kuramoto-decomposition-025-1}
\end{subfigure}
\begin{subfigure}[b]{0.45\textwidth}
\begin{tikzpicture}[
  declare function={
 d(\x)= (\x<=-pi) * abs(\x+2*pi) + and(\x>-pi, \x<=pi) * abs(\x) + (\x>pi) * abs(\x-2*pi);
 f(\x)=  (\x<-2*pi+1.165561185207211) * (2*sin(deg(1.165561185207211))/pow(1.165561185207211,2*0.25)-sin(deg(\x))/pow(d(\x),2*0.25)) +
 and(\x>=-2*pi+1.165561185207211, \x<-1.165561185207211) * (sin(deg(1.165561185207211))/pow(1.165561185207211,2*0.25)) +
 and(\x>=-1.165561185207211, \x<=1.165561185207211) * (-sin(deg(\x))/pow(d(\x),2*0.25))+
 and(\x>1.165561185207211, \x<=2*pi-1.165561185207211) * (-sin(deg(1.165561185207211))/pow(1.165561185207211,2*0.25))+
 (\x>2*pi-1.165561185207211) * (-2*sin(deg(1.165561185207211))/pow(1.165561185207211,2*0.25)-sin(deg(\x))/pow(d(\x),2*0.25));
g(\x)= (\x<-2*pi+1.165561185207211) * (-2*sin(deg(1.165561185207211))/pow(1.165561185207211,2*0.25)) +
 and(\x>=-2*pi+1.165561185207211, \x<-1.165561185207211) * (-sin(deg(1.165561185207211))/pow(1.165561185207211,2*0.25)-sin(deg(\x))/pow(d(\x),2*0.25)) +
 and(\x>=-1.165561185207211, \x<=1.165561185207211) * (0)+
 and(\x>1.165561185207211, \x<=2*pi-1.165561185207211) * (sin(deg(1.165561185207211))/pow(1.165561185207211,2*0.25)-sin(deg(\x))/pow(d(\x),2*0.25))+
 (\x>2*pi-1.165561185207211) * (2*sin(deg(1.165561185207211))/pow(1.165561185207211,2*0.25));
  }
]
\begin{axis}[
  axis x line=middle, axis y line=middle,
  xmin=-2*pi, xmax=2*pi, xtick={-6,-4,-2,0,2,4,6}, xlabel=$\theta$,
  ymin=-2, ymax=2, ytick={-1.5,-1,-0.5,0,0.5,1.5},
  legend style={at={(0.76,0.95)}},
]
\addplot[domain=-2*pi:2*pi, samples=200,color=blue]{f(x)};
\addlegendentry{$\Delta(\theta)$}
\addplot[domain=-2*pi:2*pi, samples=200,color=red]{g(x)};
\addlegendentry{$\Lambda(\theta)$}
\end{axis}
\end{tikzpicture} 
\caption{$\Delta(\theta)$ and $\Lambda(\theta)$}
\label{fig:kuramoto-decomposition-025-2}
\end{subfigure}
\caption{Graph of the function $-h(\theta)$ and the functions $\Delta(\theta)$ and $\Lambda(\theta)$ in the decomposition for the value $\alpha=0.25$.}\label{fig:kuramoto-decomposition-025}
\end{figure}
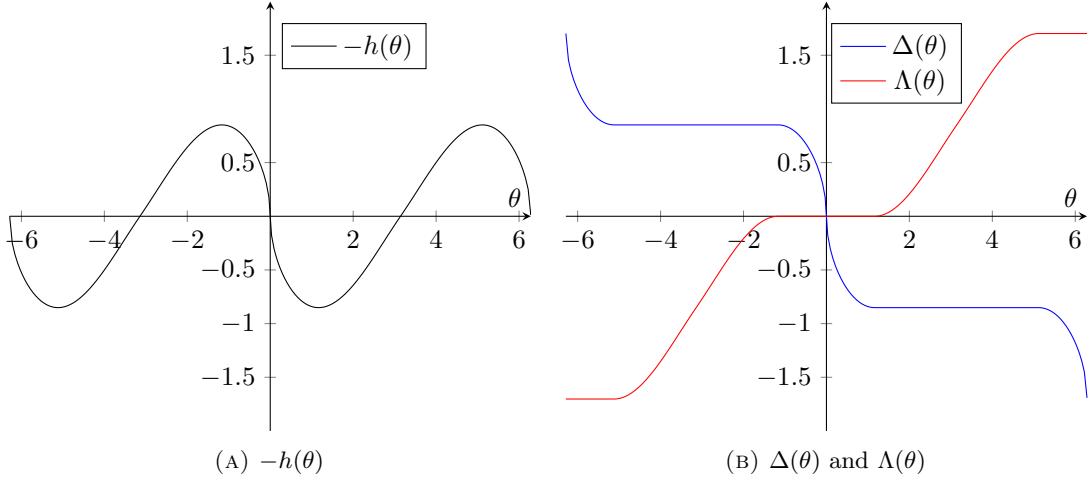

Our next result will show that although $\mathcal{V}[\mu]$ is not fully Lipschitz-continuous with respect to $(z,\Omega)\in \mathbb{T}\times \mathbb{R}$, it is one-sided Lipschitz uniformly in $t\in [0,T]$. The cornerstone in such result is the following split of $-h$ into a decreasing and a Lipschitz-continuous part, see Fig. \ref{fig:kuramoto-decomposition-025}.

\begin{lem}\label{L-split-kernel}
Consider $\alpha\in (0,\frac{1}{2})$ and set $\bar{h}$ and $\tilde{\theta}\in (0,\frac{\pi}{2})$ such that
$$\bar{h}=\max_{0<\theta<\pi}h(\theta)\ \mbox{ and }\ 2\alpha\sin\tilde{\theta}=\tilde{\theta}\cos \tilde{\theta}.$$
Define the couple of functions $\Delta,\Lambda:[-2\pi,2\pi]\longrightarrow\mathbb{R}$ as follows
$$
\Delta(\theta):=\left\{
\begin{array}{ll}
2\bar{h}-h(\theta), & \theta\in [-2\pi,-2\pi+\tilde{\theta}),\\
\bar h, & \theta\in [-2\pi+\tilde{\theta},-\tilde{\theta}),\\
-h(\theta), & \theta\in [-\tilde{\theta},\tilde{\theta}],\\
-\bar{h}, & \theta\in (\tilde{\theta},2\pi-\tilde{\theta}],\\
-h(\theta)-2\bar{h}, & \theta\in (2\pi-\tilde{\theta},2\pi],
\end{array}
\right.
$$
$$
\Lambda(\theta):=
\left\{
\begin{array}{ll}
-2\bar{h}, & \theta\in [-2\pi,-2\pi+\tilde{\theta}),\\
-\bar{h}-h(\theta), & \theta\in [-2\pi+\tilde{\theta},-\tilde{\theta}),\\
0, & \theta\in [-\tilde{\theta},\tilde{\theta}],\\
\bar{h}-h(\theta), & \theta\in (\tilde{\theta},2\pi-\tilde{\theta}],\\
2\bar{h}, & \theta\in (2\pi-\tilde{\theta},2\pi].
\end{array}
\right.
$$
Then, the following properties hold true
\begin{enumerate}
\item $\Delta$ is monotonically decreasing, $\Lambda$ is Lipschitz-continuous and 
$$-h(\theta)=\Delta(\theta)+\Lambda(\theta),\ \forall\,\theta\in [-2\pi,2\pi].$$
\item $-h$ is one-sided Lipschitz in $[-2\pi,2\pi]$, i.e., there exists $L_0>0$ such that
$$\left((-h)(\theta_1)-(-h)(\theta_2)\right)(\theta_1-\theta_2)\leq L_0(\theta_1-\theta_2)^2.$$
\end{enumerate}
\end{lem}

See Figure \ref{fig:kuramoto-decomposition-025} and \cite{P-P-S} for a similar split applied to the agent-based system.

\begin{lem}\label{L-transport-field-sided-Lipschitz}
Consider $\alpha\in (0,\frac{1}{2})$, $K>0$ and set $\mu\in\mathcal{C}_\mathcal{M}$. Then, we have 
$$\left(\mathcal{P}[\mu_t](\theta_1,\Omega_1)-\mathcal{P}[\mu_t](\theta_2,\Omega_2)\right)(\theta_1-\theta_2)\leq (\Omega_1-\Omega_2)(\theta_1-\theta_2)+KL_0\Vert \mu\Vert_{L^\infty(0,T;\mathcal{M}(\mathbb{T}\times \mathbb{R}))}(\theta_1-\theta_2)^2,$$
for every $t\in [0,T]$, each $\theta_1,\theta_2\in \mathbb{R}$ with $\theta_1-\theta_2\in [-\pi,\pi]$ and any $\Omega_1,\Omega_2\in\mathbb{R}$. Here, the constant $L_0$ is the one-sided Lipschitz constant in Lemma \ref{L-split-kernel}.
\end{lem}

\begin{proof}
By the Definition \ref{D-transport-field} we can state
\begin{multline*}
\left(\mathcal{P}[\mu_t](\theta_1,\Omega_1)-\mathcal{P}[\mu_t](\theta_2,\Omega_2)\right)(\theta_1-\theta_2)\\
=(\Omega_1-\Omega_2)(\theta_1-\theta_2)+K\int_{(\theta_1-\pi,\theta_1+\pi]}\int_{\mathbb{R}}((-h)(\theta_1-\theta')-(-h)(\theta_2-\theta'))(\theta_1-\theta_2)\,d_{(\theta',\Omega')}\mu_t.
\end{multline*}
For every $\theta'\in (\theta_1-\pi,\theta_1+\pi]$ we equivalently have $\theta_1-\theta'\in (-\pi,\pi]$. Since $\theta_2-\theta_1\in [-\pi,\pi]$, we also obtain that $\theta_2-\theta'\in (-2\pi,2\pi]$. Hence, we are in the range of applicability of Lemma \ref{L-split-kernel}, that implies
$$\left(\mathcal{P}[\mu_t](\theta_1,\Omega_1)-\mathcal{P}[\mu_t](\theta_2,\Omega_2)\right)(\theta_1-\theta_2)\leq (\Omega_1-\Omega_2)(\theta_1-\theta_2)+KL_0\int_{(\theta_1-\pi,\theta_1+\pi]}\int_{\mathbb{R}}(\theta_1-\theta_2)^2\,d_{(\theta',\Omega')}\mu_t.$$
\end{proof}

\begin{theo}\label{T-transport-field-sided-Lipschitz}
Consider $\alpha\in (0,\frac{1}{2})$, $K>0$ and set $\mu\in\mathcal{C}_\mathcal{M}$. Then, $\mathcal{V}[\mu]$ is one-sided Lipschitz in $\mathbb{T}\times \mathbb{R}$ uniformly in $t\in [0,T]$, i.e., there exists $L=L(\alpha,K,\mu)>0$ such that
$$\left<\mathcal{V}[\mu_t](z_2,\Omega_2),\widehat{\gamma}'(1)\right>-\left<\mathcal{V}[\mu_t](z_1,\Omega_1),\widehat{\gamma}'(0)\right>\leq L\,d((z_1,\Omega_1),(z_2,\Omega_2))^2,$$
for every $t\in [0,T]$, any $(z_1,\Omega_1),(z_2,\Omega_2)\in \mathbb{T}\times\mathbb{R}$ and each minimizing geodesic $\widehat{\gamma}:[0,1]\longrightarrow\mathbb{T}\times \mathbb{R}$ in the manifold $\mathbb{T}\times \mathbb{R}$ joining $(z_1,\Omega_1)$ to $(z_2,\Omega_2)$.
\end{theo}

\begin{proof}
Fix any $t\in [0,T]$ and $(z_1,\Omega_1),(z_2,\Omega_2)\in \mathbb{T}\times \mathbb{R}$. Our first step will be to characterize the minimizing geodesics $\widehat{\gamma}:[0,1]\longrightarrow\mathbb{T}\times \mathbb{R}$ joining $(z_1,\Omega_1)$ to $(z_2,\Omega_2)$. Let us write $z_1=e^{i\theta_1}$ and $z_2=e^{i\theta_2}$ for some $\theta_1,\theta_2\in\mathbb{R}$ and assume that $\theta_2-\theta_1\in (-\pi,\pi]$ without loss of generality.

$\bullet$ \textit{Case 1}: $\theta_2-\theta_1\in (-\pi,\pi)$. In this case, the only minimizing geodesic reads
$$\widehat{\gamma}(s)=(e^{i(\theta_1+s(\theta_2-\theta_1)},\Omega_1+s(\Omega_2-\Omega_1)).$$
Notice that the directions of the geodesic at the endpoints are
$$\widehat{\gamma}'(0)=((\theta_2-\theta_1)\,iz_1,\Omega_2-\Omega_1)\ \mbox{ and }\ \widehat{\gamma}'(1)=((\theta_2-\theta_1)\,iz_2,\Omega_2-\Omega_1).$$
Therefore, we can write
$$\left<\mathcal{V}[\mu_t](z_2,\Omega_2),\widehat{\gamma}'(1)\right>-\left<\mathcal{V}[\mu_t](z_1,\Omega_1),\widehat{\gamma}'(0)\right>=\left(\mathcal{P}[\mu_t](z_2,\Omega_2)-\mathcal{P}[\mu_t](z_1,\Omega_1)\right)(\theta_2-\theta_1).$$

$\bullet$ \textit{Case 2}: $\theta_2-\theta_1=\pi$. Now $z_1$ and $z_2$ are antipodes and there are exactly two different minimizing geodesics, namely,
\begin{align*}
\widehat{\gamma}_+(s)=(e^{i(\theta_1+\pi s)},\Omega_1+s(\Omega_2-\Omega_1)),\\
\widehat{\gamma}_-(s)=(e^{i(\theta_1-\pi s)},\Omega_1+s(\Omega_2-\Omega_1)).
\end{align*}
Its directions at the endpoints read
$$\widehat{\gamma}'_{\pm}(0)=(\pm\pi\,iz_1,\Omega_2-\Omega_1)\ \mbox{ and }\ \widehat{\gamma}'_\pm(1)=(\pm\pi\,iz_2,\Omega_2-\Omega_1).$$
Then, it is clear that
$$\left<\mathcal{V}[\mu_t](z_2,\Omega_2),\widehat{\gamma}_\pm'(1)\right>-\left<\mathcal{V}[\mu_t](z_1,\Omega_1),\widehat{\gamma}'_\pm(0)\right>=\pm\left(\mathcal{P}[\mu_t](z_2,\Omega_2)-\mathcal{P}[\mu_t](z_1,\Omega_1)\right)\pi.$$

No matter the case, we can always use Lemma \ref{L-transport-field-sided-Lipschitz} to arrive at
$$\left<\mathcal{V}[\mu_t](z_2,\Omega_2),\widehat{\gamma}'(1)\right>-\left<\mathcal{V}[\mu_t](z_1,\Omega_1),\widehat{\gamma}'(0)\right>\leq(\Omega_1-\Omega_2)(\theta_1-\theta_2)+KL_0\Vert \mu\Vert_{L^\infty(0,T;\mathcal{M}(\mathbb{T}\times \mathbb{R}))}(\theta_1-\theta_2)^2.$$
Notice that $\theta_1-\theta_2\in [-\pi,\pi]$ and, consequently, $\vert \theta_1-\theta_2\vert=\vert \theta_1-\theta_2\vert_o=d(z_1,z_2)$. Applying Young's inequality we arrive at the desired result for the value
$$L=\frac{1}{2}+KL_0\Vert \mu\Vert_{L^\infty(0,T;\mathcal{M}(\mathbb{T}\times \mathbb{R}))}.$$
\end{proof}

We are now ready to complete the existence part in Lemma \ref{L-existence-characteristic-system} for the characteristic system \eqref{E-characteristic-system-1} with an appropriate notion of uniqueness, namely, one-sided uniqueness. In this way, we obtain the following full well-posedness result.

\begin{theo}\label{T-well-posedness-characteristic-system}
Consider $\alpha\in (0,\frac{1}{2})$, $K>0$ and fix $\mu\in\widetilde{\mathcal{C}}_\mathcal{M}$. The characteristic system \eqref{E-characteristic-system-1} associated with the transport field $\mathcal{V}[\mu]$ enjoys a global-in-time $C^1$ solution that is unique forward-in-time for every given initial data $x_0=(z_0,\Omega_0)=(e^{i\theta_0},\Omega_0)\in\mathbb{T}\times \mathbb{R}$. Indeed, the same representation of the solution holds true, i.e.,
$$X(t;t_0,x_0)=(Z(t;t_0,z_0,\Omega_0),\Omega_0)=(e^{i\Theta(t;t_0,\theta_0,\Omega_0)},\Omega_0),\ t\geq t_0,$$
where $\Theta(t;t_0,\theta_0,\Omega_0)$ is the unique forward-in-time $C^1$ solution to \eqref{E-characteristic-system-2}.
\end{theo}

Although the proof is standard and relies on the one-sided Lipschitz condition in Theorem \ref{T-transport-field-sided-Lipschitz} and the weak differentiability properties of the squared distance in Appendix \ref{appendix-differentiability-distance}, we give a simple proof for the sake of completeness because it involves some delicate points.
\begin{proof}
Let us assume that there are two different solutions $x_1=x_1(t)$ and $x_2=x_2(t)$ to  \eqref{E-characteristic-system-1} with same initial data $x_1(t_0)=x_0=x_2(t_0)$. Define the following function
$$I(t):=\frac{1}{2}d(x_1(t),x_2(t))^2,\ t\geq t_0,$$
where $d(\cdot,\cdot)$ means the Riemannian distance in $\mathbb{T}\times \mathbb{R}$, i.e., \eqref{E-riemannian-distance} (see Appendix \ref{appendix-differentiability-distance}). Recall that the distance function $d(\cdot,\cdot)$ is Lipschitz-continuous, see Proposition \ref{P-derivative-distance-1}. Since $x_1$ and $x_2$ are $C^1$ trajectories, then $I=I(t)$ is absolutely continuous. Hence, its derivative exists for almost every $t$ and it agrees with the one-sided Dini upper derivative. Since Theorem \ref{T-dini-derivative-distance} implies that the one-sided Dini upper directional derivatives ($\frac{d^+ }{dt}$) of the squared distance are finite, then the chain rule yields
\[
\frac{dI}{dt}\equiv\frac{d^+ I}{dt}=d^+\left(\frac{1}{2}d_{x_2(t)}^2\right)_{x_1(t)} \left(\dot{x}_1(t))+d^+\left(\frac{1}{2}d_{x_1(t)}^2\right)_{x_2(t)}(\dot{x}_2(t)\right),
\]
for almost every $t\geq t_0$. Recall that in Theorem \ref{T-dini-derivative-distance} we also got an upper bound for such Dini directional derivatives that reads as follows
\begin{align*}
d^+\left(\frac{1}{2}d_{x_2(t)}^2\right)_{x_1(t)}\Big(\dot{x}_1(t)\Big)&\leq \inf_{\substack{w_1\in \exp_{x_1(t)}^{-1}(x_2(t))\\ \vert w_1\vert =d(x_1(t),x_2(t))}}-\left<\dot{x}_1(t),w_1\right>,\\
d^+\left(\frac{1}{2}d_{x_1(t)}^2\right)_{x_2(t)}\Big(\dot{x}_2(t)\Big)&\leq \inf_{\substack{w_2\in \exp_{x_2(t)}^{-1}(x_1(t))\\ \vert w_2\vert =d(x_1(t),x_2(t))}}-\left<\dot{x}_2(t),w_2\right>.
\end{align*}
Let us fix any minimizing geodesic $\widehat{\gamma}_t:[0,1]\longrightarrow\mathbb{T}\times \mathbb{R}$ joining $x_1(t)$ to $x_2(t)$, for every $t\geq t_0$. Then, we can choose $w_1=\widehat{\gamma}_t'(0)$ and $w_2=-\widehat{\gamma}_t'(1)$ in the above inequalities. Consequently,
$$\frac{dI}{dt}\leq \left<\dot{x}_2(t),\widehat{\gamma}_t'(1)\right>-\left<\dot{x}_1(t),\widehat{\gamma}_t'(0)\right>=\left<\mathcal{V}[\mu_t](x_2(t)),\widehat{\gamma}_t'(1)\right>-\left<\mathcal{V}[\mu_t](x_1(t)),\widehat{\gamma}_t'(0)\right>.$$
Using Theorem \ref{T-transport-field-sided-Lipschitz} we get the estimate
$$\frac{dI}{dt}\leq L d(x_1(t),x_2(t)^2)=2L\,I(t),\ \mbox{ a.e. }t\geq t_0.$$
Since $I(t_0)=\frac{1}{2}d(x_1(0),x_2(0))^2=\frac{1}{2}d(x_0,x_0)^2=0$, Gronwall's Lemma amounts to the desired sided-uniqueness, namely,
$$x_1(t)=x_2(t),\ \forall\,t\geq t_0.$$
\end{proof}

The same ideas as above can be used even if $x_1(0)\neq x_2(0)$ to derive some stability result of the flow. In fact, although $\mathcal{V}[\mu]$ is not fully Lipschitz continuous (but merely one-sided Lipschitz continuous), its characteristic flow is.

\begin{cor}\label{C-Lipschitz-characteristic-flow}
Consider $\alpha\in (0,\frac{1}{2})$, $K>0$ and fix $\mu\in\widetilde{\mathcal{C}}_\mathcal{M}$. Let $X(t;0,z,\Omega)=(Z(t;0,z,\Omega),\Omega)$ be the flow of the characteristic system \eqref{E-characteristic-system-1} associated with the transport field $\mathcal{V}[\mu]$. Then, $X$ is Lipschitz-continuous in $(z,\Omega)$; namely, there exists $L=L(\alpha,K,\mu)$ with
$$d(X(t;0,(z_1,\Omega_1))-X(t;0,(z_2,\Omega_2)))\leq  d((z_1,\Omega_1),(z_2,\Omega_2))e^{Lt},$$
for every $(z_1,\Omega_1),(z_2,\Omega_2)\in \mathbb{T}\times \mathbb{R}$ and each $t\in [0,T]$.
\end{cor}

Again, like in Theorem \ref{T-transport-field-sided-Lipschitz}, the constant $L$ in Corollary \ref{C-Lipschitz-characteristic-flow} is
$$L=\frac{1}{2}+KL_0\Vert \mu\Vert_{L^\infty(0,T;\mathcal{M}(\mathbb{T}\times \mathbb{R}))}.$$

\subsection{Types of measure-valued solutions}

In the above Definition \ref{D-weak-measure-solution} and Proposition \ref{P-weak-measure-solution-charact}, the classical concept of weak measure-valued solutions $f\in \mathcal{AC}_\mathcal{M}$ to \eqref{E-kuramoto-transport-TxR} was revisited. It clearly agrees with the standard concept of measure-valued solution for homogeneous conservative continuity equations (for instance, see \cite{A-C} and references therein). Indeed, such concept requires a very weak regularity of the transport field $\mathcal{V}[f]$ in the continuity equation \eqref{E-kuramoto-transport-TxR}, specifically,
\begin{equation}\label{E-65}
\int_0^T\int_C\vert \mathcal{V}[f_t]\vert\,df_t\,dt<\infty,
\end{equation}
for each compact subset $C\subseteq \mathbb{T}\times \mathbb{R}$. Although such concepts are not usually considered when working on general manifolds like $\mathbb{T}\times \mathbb{R}$, straightforward ideas allow extending them from the Euclidean space $\mathbb{R}^d$ to Riemannian manifolds. Apart from that notion of solution, there are a couple of related concepts that we recall in the following.

\begin{defi}\label{D-measure-solution-flow}
We will say that a time dependent measure $f\in \widetilde{C}_\mathcal{M}$ is a measure-valued solution to \eqref{E-kuramoto-transport-TxR} in the sense of the characteristic flow when
$$f_t=X_f(t;0,\cdot)_{\#}f_0,\ \mbox{ for all }\ t\geq 0,$$
where $X_f$ stands for the flow of the transport field $\mathcal{V}[f]$ as introduced in Theorem \ref{T-well-posedness-characteristic-system}.
\end{defi}

The following result is a straightforward consequence of the definition.

\begin{pro}\label{P-types-solutions-flow-to-weak}
Let $f\in \widetilde{C}_\mathcal{M}$ be a solution to \eqref{E-kuramoto-transport-TxR} in the sense of the characteristic flow. Then
\begin{enumerate} 
\item $f \in \mathcal{AC}_\mathcal{M}.$
\item $f$ is a weak measure-valued solution to \eqref{E-kuramoto-transport-TxR}.
\end{enumerate}
\end{pro}

\begin{defi}\label{D-measure-solution-superposition}
Define the space $\mathcal{C}^1([0,T]):=C^1([0,T],\mathbb{T}\times \mathbb{R})$, consider $f\in \widetilde{C}_\mathcal{M}$ and the set 
$$\mathcal{S}_f([0,T]):=\left\{(x,\gamma)\in (\mathbb{T}\times \mathbb{R})\times \mathcal{C}^1([0,T]):\,\gamma\mbox{ is a characteristic of }\mathcal{V}[f]\mbox{ issued at }x\right\}.$$
Then, $f$ is said to be a superposition solution to \eqref{E-kuramoto-transport-TxR} if $f(t=0)=f_0$ and there exists some probability measure $\eta\in \mathcal{P}((\mathbb{T}\times \mathbb{R})\times \mathcal{C}^1([0,T]))$ such that $\supp\,\eta\subseteq \mathcal{S}_f([0,T])$ and $f_t=f^\eta_t$ for all $t\in [0,T]$, where the measure $f^\eta$ is defined as follows
$$\left<f^\eta_t,\varphi\right>=\int_{(\mathbb{T}\times \mathbb{R})\times \mathcal{C}^1([0,T])}\varphi(\gamma_t)\,d_{(x,\gamma)}\eta,$$
for any test function $\varphi\in C_b(\mathbb{T}\times \mathbb{R})$.
\end{defi}

Mimicking the ideas in \cite[Theorem 4.4]{A-C} we obtain the following result.

\begin{pro}\label{P-type-solutions-weak-to-superposition}
Let $f\in \mathcal{AC}_\mathcal{M}\cap \mathcal{T}_\mathcal{M}$ be any weak measure-valued solution to \eqref{E-kuramoto-transport-TxR}. Then, there exists some probability measure $\eta\in \mathcal{P}((\mathbb{T}\times \mathbb{R})\times \mathcal{C}^1([0,T]))$ supported on $\mathcal{S}_f([0,T])$ such that
$$f_t\equiv f^\eta_t,\ t\geq 0.$$
In other words, $f$ is a superposition solution.
\end{pro}

Recall that in Remark \ref{R-transport-field-holder}, we mentioned that an analogue to Lemma \ref{L-existence-characteristic-system} with $f\in \mathcal{C}_M$ instead of $f\in \widetilde{\mathcal{C}}_{\mathcal{M}}$ can be derived by means of Caratheodory characteristics. Then, Definitions \ref{D-measure-solution-flow} and \ref{D-measure-solution-superposition} along with Propositions \ref{P-types-solutions-flow-to-weak} and \ref{P-type-solutions-weak-to-superposition} have analogues with $f\in \mathcal{C}_\mathcal{M}$ only. Again, we will skip it here.\\

So far, we have not used any special property of the transport field other than
$$\frac{\mathcal{V}[f]}{1+\vert \Omega\vert}\in L^\infty(0,T;\mathfrak{X}C_b(\mathbb{T}\times \mathbb{R})).$$
It guarantees both \eqref{E-65} (giving sense to weak measure-valued solutions) and existence of Caratheodory characteristics (giving sense to solutions in the sense of the flow and superposition solutions). Let us see that under the forward-uniqueness property of the characteristic system, that follows in Theorem \ref{T-well-posedness-characteristic-system} from the one-sided uniqueness property of $\mathcal{V}[f]$, any superposition solution is also a solution in the sense of the characteristic flow.

\begin{pro}\label{P-types-solutions-superposition-to-flow}
Let $f\in \widetilde{C}_\mathcal{M}$ be a superposition solution to \eqref{E-kuramoto-transport-TxR}. Then, $f$ is a solution to \eqref{E-kuramoto-transport-TxR} in the sense of the charactersitic flow.
\end{pro}

\begin{proof}
By definition, $f_t=f^\eta_t$ for $t\in [0,T]$, for some $\eta\in \mathcal{P}((\mathbb{T}\times \mathbb{R})\times \mathcal{C}^1([0,T]))$ with $\supp\eta\subseteq \mathcal{S}_f([0,T])$. Due to Theorem \ref{T-well-posedness-characteristic-system}, the one-sided Lipschitz condition of $\mathcal{V}[f]$ implies one-sided uniqueness of the characteristic system \eqref{E-characteristic-system-1}. Then, $\mathcal{S}_f([0,T])$ agrees with the graph of the flow $X_f$ in the following sense
$$\mathcal{S}_f([0,T])=\left\{(x,X_f(\cdot\,;0,x)):\,x\in \mathbb{T}\times \mathbb{R}\right\}.$$
Consequently, we can write
\begin{equation}\label{E-66}
\left<f_t,\varphi\right>=\left<f^\eta_t,\varphi\right>=\int_{(\mathbb{T}\times \mathbb{R})\times \mathcal{C}^1([0,T])}\varphi(\gamma_t)\,d_{(x,\gamma)}\eta=\int_{(\mathbb{T}\times \mathbb{R})\times \mathcal{C}^1([0,T])}\varphi(X_f(t;0,x))\,d_{(x,\gamma)}\eta,
\end{equation}
for every $t\in [0,T]$ and $\varphi\in C_b(\mathbb{T}\times \mathbb{R})$. Consider the projection
$$\begin{array}{cccc}
\pi_x: & (\mathbb{T}\times \mathbb{R})\times \mathcal{C}^1([0,T])&\longrightarrow &\mathbb{T}\times \mathbb{R}\\
 & (x,\gamma) & \longmapsto & x,
\end{array}$$
along with the marginal $\mu:=(\pi_x)_{\#}\eta$. By the disintegration theorem (see Theorem \ref{T-disintegration} below), let us consider the family $(\eta_x)_{x\in \mathbb{T}\times \mathbb{R}}$ of conditional probabilities or disintegrations of $\eta$. Then,
$$\int_{(\mathbb{T}\times \mathbb{R})\times \mathcal{C}^1([0,T])}\varphi\,d\eta=\int_{\mathbb{T}\times \mathbb{R}}\left(\int_{\mathcal{C}^1([0,T])}\varphi(x,\gamma)\,d_\gamma(\eta_x)\right)\,d\mu,$$
for any $\varphi\in C_b(\mathbb{T}\times \mathbb{R})$. In particular, when applied to \eqref{E-66} we arrive at
$$\left<f_t,\varphi\right>=\int_{\mathbb{T}\times \mathbb{R}}\varphi(X_f(t;0,x))\,d_x\mu,$$
for any $\varphi\in C_b(\mathbb{T}\times \mathbb{R})$. Therefore,
$$f_t=X_f(t;0,\cdot)_{\#}\mu,\ \mbox{ for all }\ t\in [0,T].$$
Taking $t=0$ we get $\mu\equiv f_0$. Then, $f$ is a solution in the sense of the characteristic flow to \eqref{E-kuramoto-transport-TxR}.
\end{proof}

\begin{rem}\label{R-types-solutiones}
The above Propositions \ref{P-types-solutions-flow-to-weak}, \ref{P-type-solutions-weak-to-superposition} and \ref{P-types-solutions-superposition-to-flow} guarantee that all the above three concepts of measure-valued solutions are equivalent in our problem because of the properties:
\begin{enumerate}
\item $\displaystyle\frac{\mathcal{V}[f]}{1+\vert \Omega\vert}\in L^\infty(0,T,\mathfrak{X}C_b(\mathbb{T}\times \mathbb{R}))$.
\item $\mathcal{V}[f]$ is one-sided Lipschitz-continuous uniformly in $t\in [0,T]$.
\end{enumerate}
\end{rem}

\section{Existence of weak measure-valued solutions}\label{S-weak-solutions-existence}

In this part, we shall derive existence of global-in-time measure valued solutions to \eqref{E-kuramoto-transport-TxR} in the subcritical regime $\alpha\in (0,\frac{1}{2})$. The idea does nor rely on any regularization technique of the kernel. Instead, it will rely on a compactness argument as $N\rightarrow\infty$ for any sequences of empirical empirical measures associated with a sequence of solutions to the $N$-particle system \eqref{E-kuramoto-discrete}-\eqref{E-kuramoto-kernel} that initially approximates the given initial datum $f_0$ in Wasserstein distance. Such method will produce weak measure-valued solutions (equivalently, solutions in the sense of the characteristic flow) to \eqref{E-kuramoto-transport-TxR} in the subcritical regime. Similar ideas will be analysed later in Sections \ref{S-critical-regime} and \ref{S-supercritical-regime} for the most singular regimes. Notice that the aforementioned compactness of the empirical measures becomes a first step towards the derivation of the full mean field limit, that will be proved later in Section \ref{S-uniqueness-mean-field} via a Dobrushin-type inequality.

The rest of this section consists of the following three parts. First, we will revisit the concept of empirical measures associated with a solution to the discrete model \eqref{E-kuramoto-discrete}-\eqref{E-kuramoto-kernel} and will show that they automatically are weak measure-valued solutions to \eqref{E-kuramoto-transport-TxR}. In the second part, we will obtain some a priori bounds implying the weak compactness of such sequence of empirical measures. The final step will be to identify the limit as a weak measure-valued solution to the continuous model \eqref{E-kuramoto-transport-TxR}.

\subsection{Empirical measures}
In this part we will recall the definition of empirical measures associated with a discrete solution to \eqref{E-kuramoto-discrete}-\eqref{E-kuramoto-kernel}, see \cite{C-C-H-2,J,L,N,M-P}. We will also inspect whether they provide measure-valued solutions to the macroscopic system \eqref{E-kuramoto-transport-TxR}.

\begin{defi}\label{D-empirical-measures-subcritical}
Fix $N\in\mathbb{N}$ and consider $N$ oscillators with phases and natural frequencies given by the configurations
$$\Theta_0^N=(\theta_{1,0}^N,\ldots,\theta_{N,0}^N)\ \mbox{ and }\ \{\Omega_i^N:\,i=1,\ldots,N\}.$$ 
Let $\Theta^N(t):=(\theta_1^N(t),\ldots,\theta_N^N(t))$ be the unique global-in-time classical solution to the discrete singular Kuramoto model \eqref{E-kuramoto-discrete}-\eqref{E-kuramoto-kernel} according to \cite[Theorem 3.1]{P-P-S}. Then, the associated empirical measures are given by $\mu^N\in \mathcal{P}(\mathbb{T}\times\mathbb{R})$ defined as follows
$$\mu^N_t:=\frac{1}{N}\sum_{i=1}^N\delta_{z_i^N(t)}(z)\otimes \delta_{\Omega_i^N}(\Omega),$$
where $z_i^N(t):=e^{i\theta_i^N(t)}$.
\end{defi}

\begin{theo}\label{T-empirical-measures-subcritical}
Consider $\alpha\in (0,\frac{1}{2})$ and $K>0$. Fix $N\in\mathbb{N}$ and consider $N$ oscillators with initial phases and natural frequencies given by the configurations
$$\Theta_0^N=(\theta_{1,0}^N,\ldots,\theta_{N,0}^N)\ \mbox{ and }\ \{\Omega_i^N:\,i=1,\ldots,N\}.$$ 
Let $\mu^N$ be the associated empirical measure according to Definition \ref{D-empirical-measures-subcritical}. Then, $\mu^N\in \mathcal{AC}_\mathcal{M}\cap \mathcal{T}_\mathcal{M}$ is a weak measure-valued solution to \eqref{E-kuramoto-transport-TxR} and, in addition,
\begin{equation}\label{E-17}
\left\vert\frac{d}{dt}\int_{\mathbb{T}\times \mathbb{R}}\varphi\,d\mu^N_t\right\vert\leq \left(\frac{1}{N}\sum_{i=1}^N\vert \Omega_i^N\vert+K\Vert h\Vert_{C(\mathbb{T})}\right)\Vert \nabla \varphi\Vert_{C_0(\mathbb{T}\times \mathbb{R})},
\end{equation}
for every $t\geq 0$ and every $\varphi\in C^1_0(\mathbb{T}\times\mathbb{R})$.
\end{theo}

\begin{proof}
Let us first prove that $\mu^N\in \mathcal{AC}_\mathcal{M}\cap\mathcal{T}_\mathcal{M}$. Regarding the tightness condition notice that
$$\Vert \vert\Omega\vert\mu^N_t\Vert_{\mathcal{M}(\mathbb{T}\times \mathbb{R})}=\frac{1}{N}\sum_{i=1}^N\vert\Omega_i^N\vert,$$
for every $t\geq 0$. Regarding the absolute continuity in time, set $\varphi\in C^\infty_c(\mathbb{T}\times \mathbb{R})$ and observe that
\begin{equation}\label{E-12}
t\in [0,+\infty)\longmapsto \int_{\mathbb{T}}\int_{\mathbb{R}}\varphi\,d\mu^N_t=\frac{1}{N}\sum_{i=1}^N\varphi(\theta_i^N(t),\Omega_i^N)
\end{equation}
is locally absolutely continuous (it is $C^1$ in fact). Indeed, taking derivatives in \eqref{E-12} yields
\begin{equation}\label{E-13}
\frac{d}{dt}\int_{\mathbb{T}}\varphi\,d\mu^N_t=\frac{1}{N}\sum_{i=1}^N\frac{\partial\varphi}{\partial \theta}(\theta_i^N(t),\Omega_i^N)\dot{\theta}_i^N(t).
\end{equation}
As explained in Remark \ref{rem-complex-derivatives}, we assert that for $(z=e^{i\theta},\Omega)\in\mathbb{T}\times \mathbb{R}$
$$\frac{\partial\varphi}{\partial \theta}(\theta,\Omega)=-ie^{-i\theta}\nabla_z\varphi(\theta,\Omega).$$
Then, the above equation \eqref{E-13} can be restated as follows
\begin{align}\label{E-14}
\begin{split}
\frac{d}{dt}\int_{\mathbb{T}}\varphi\,d\mu^N_t&=\frac{1}{N}\sum_{i=1}^N\Re\left[\nabla_z\varphi(\theta_i^N(t),\Omega_i^N)(-ie^{-i\theta_i^N(t)})\dot{\theta}_i^N(t)\right]\\
&=\frac{1}{N}\sum_{i=1}^N\left<\nabla_z\varphi(\theta_i^N(t),\Omega_i^N),\frac{d}{dt}e^{i\theta_i^N(t)}\right>,
\end{split}
\end{align}
where $\Re$ means real part of the complex number. To describe the associated transport field, notice that
$$
\mathcal{P}[\mu^N_t](\theta,\Omega)=\Omega-K\int_{(-\pi,\pi]}\int_{\mathbb{R}}h(\theta-\theta')\,d_{(\theta',\Omega')}\mu^N_t=\Omega-\frac{K}{N}\sum_{j=1}^N h(\theta-\theta_j^N(t)),
$$
and, consequently,
$$
\mathcal{P}[\mu^N_t](\theta_i^N(t),\Omega_i^N)=\Omega_i-\frac{K}{N}\sum_{j=1}^N h(\theta_i^N(t)-\theta_j^N(t)).
$$
Since $\theta_i^N(t)$ are solutions to the discrete singular Kuramoto model, then we arrive at
\begin{equation}\label{E-15}
\frac{d}{dt}(e^{i\theta_i^N(t)},\Omega_i^N)=(ie^{i\theta_i^N(t)}\dot{\theta}_i^N(t),0)=\mathcal{V}[\mu^N_t](\theta_i^N(t),\Omega_i).
\end{equation}
Putting \eqref{E-15} into \eqref{E-14} implies
$$
\frac{d}{dt}\int_{\mathbb{T}}\int_{\mathbb{R}}\varphi\,d\mu^N_t=\left.\frac{1}{N}\sum_{i=1}^N
\left<\nabla_{(z,\Omega)}\varphi,\mathcal{V}[\mu^N_t]\right>\right\vert_{(z,\Omega)=(e^{i\theta_i^N(t)},\Omega_i^N)}.$$
Then, it becomes apparent that $\mu^N_t$ is a weak solution, namely,
\begin{equation}\label{E-16}
\frac{d}{dt}\int_{\mathbb{T}}\int_{\mathbb{R}}\varphi\,d\mu^N_t=\int_{\mathbb{T}}\int_{\mathbb{R}}\left<\nabla_{(z,\Omega)}\varphi,\mathcal{V}[\mu^N_t]\right>\,d\mu^N_t.
\end{equation}
Notice that all the above computations also makes sense for $\varphi\in C^1_0(\mathbb{T}\times \mathbb{R})$. By putting the sublinear growth of $\mathcal{V}[\mu^N]$ in Theorem \ref{T-transport-field-holder} into \eqref{E-16}, we obtain
\begin{align*}
\left\vert\frac{d}{dt}\int_{\mathbb{T}\times \mathbb{R}}\varphi\,d\mu^N_t\right\vert&\leq \Vert \nabla\varphi\Vert_{C_0(\mathbb{T}\times \mathbb{R})}\int_{\mathbb{T}\times \mathbb{R}}\vert \mathcal{V}[\mu^N_t]\vert\,d\mu^N_t\\
&= \Vert \nabla\varphi\Vert_{C_0(\mathbb{T}\times \mathbb{R})}\int_{\mathbb{T}\times \mathbb{R}}\vert \mathcal{P}[\mu^N_t]\vert\,d\mu^N_t\\
&\leq \Vert\nabla\varphi\Vert_{C_0(\mathbb{T}\times \mathbb{R})}\int_{\mathbb{T}\times \mathbb{R}}\left(\vert \Omega\vert+K\Vert h\Vert_{C(\mathbb{T})}\Vert \mu^N_t\Vert_{\mathcal{M}(\mathbb{T}\times \mathbb{R})}\right)\,d\mu^N_t.
\end{align*}
Then, the estimate \eqref{E-17} becomes true.
\end{proof}

\subsection{A priori estimates and compactness}
Our main goal now is to derive the required compactness allowing us to pass to the limit in the weak formulation in Definition \ref{D-weak-measure-solution}. To such end, we will first derive some necessary estimates in the following result.

\begin{cor}\label{C-mean-field-apriori-subcritical}
Consider $\alpha\in (0,\frac{1}{2})$ and $K>0$. Set $N$ oscillators with phases and natural frequencies given by the configurations
$$\Theta_0^N=(\theta_{1,0}^N,\ldots,\theta_{N,0}^N)\ \mbox{ and }\ \{\Omega_i^N:\,i=1,\ldots,N\},$$
for any $N\in \mathbb{N}$. Assume that there exists a constant $M_1>0$ that does not depend on $N$ such that
\begin{equation}\label{E-Omega-moment-1}
\frac{1}{N}\sum_{i=1}^N\vert\Omega_i^N\vert\leq M_1,
\end{equation}
for all $N\in\mathbb{N}$ and consider the associated empirical measures $\mu^N$ according to Definition \ref{D-empirical-measures-subcritical}. Then,
\begin{align*}
\sup_{t\in [0,T]}\Vert \mu^N_t\Vert_{C^1_0(\mathbb{T}\times \mathbb{R})^*}&\leq 1,\\
\Vert\mu^N_{t_1}-\mu^N_{t_2}\Vert_{C^1_0(\mathbb{T}\times \mathbb{R})^*}&\leq (M_1+K\Vert h\Vert_{C(\mathbb{T})})\vert t_1-t_2\vert,
\end{align*}
for every $N\in\mathbb{N}$ and every $t_1,t_2\geq 0$.
\end{cor}

We skip the proof that is a clear consequence of the assumption \eqref{E-Omega-moment-1} along with the estimate \eqref{E-17} in Theorem \ref{T-empirical-measures-subcritical}. In the sequel, we will need a stronger version of \eqref{E-Omega-moment-1}. The following result introduces the required equi-sumability condition of the natural frequencies along with its relation with condition \eqref{E-Omega-moment-1}.

\begin{pro}\label{P-Omega-moment-equint}
Let us consider a configuration of $N\in\mathbb{N}$ natural frequencies
$$\{\Omega_i^N:\,i=1,\ldots,N\}\subseteq \mathbb{R},$$
for every $N\in \mathbb{N}$.
\begin{enumerate}
\item Assume that the following equi-sumability condition holds true
\begin{equation}\label{E-Omega-moment-1-equiint}
\lim_{R\rightarrow+\infty}\sup_{N\in \mathbb{N}}\frac{1}{N}\sum_{\substack{1\leq i\leq N\\ \vert \Omega_i^N\vert \geq R}}\vert \Omega_i^N\vert=0.
\end{equation}
Then the sumability condition \eqref{E-Omega-moment-1} is also fulfilled.
\item Fix $k>1$ and assume that the following summability condition holds true
$$
\frac{1}{N}\vert \Omega_i^N\vert^k\leq M_k,
$$
for some $N$-independent $M_k>0$ and each $N\in \mathbb{N}$. Then, the following condition takes place
\begin{equation}\label{E-Omega-moment-m-equiint}
\lim_{R\rightarrow 0}\sup_{N\in \mathbb{N}}\frac{1}{N}\sum_{\substack{1\leq i\leq N\\ \vert \Omega_i^N\vert\geq R}}\vert \Omega_i^N\vert^m=0,
\end{equation}
for every $1\leq m<k$. In particular, the equi-sumability condition \eqref{E-Omega-moment-1-equiint} holds.
\end{enumerate}
\end{pro}

\begin{proof}
Regarding the first assertion, fix any arbitrary $R>0$ and notice that
$$
\frac{1}{N}\sum_{i=1}^N\vert \Omega_i\vert^N=\frac{1}{N}\sum_{\substack{1\leq i\leq N\\ \vert \Omega_i^N\vert < R}}\vert \Omega_i^N\vert+\frac{1}{N}\sum_{\substack{1\leq i\leq N\\ \vert \Omega_i^N\vert \geq R}}\vert \Omega_i^N\vert\leq R+\frac{1}{N}\sup_{\substack{1\leq i\leq N\\ \vert \Omega_i^N\vert\geq R}}\vert \Omega_i^N\vert ,
$$
for every $N\in\mathbb{N}$. By virtue of \eqref{E-Omega-moment-1-equiint}, the right hand side is bounded with respect to $N$ and
\begin{equation}\label{E-23}
\sup_{N\in \mathbb{N}}\frac{1}{N}\sum_{i=1}^N\vert \Omega_i^N\vert\leq R+\sup_{N\in \mathbb{N}}\frac{1}{N}\sup_{\substack{1\leq i\leq N\\ \vert \Omega_i^N\vert\geq R}}\vert \Omega_i^N\vert.
\end{equation}
The optimal $M_1$ can be achieved by minimizing the right hand side in \eqref{E-23} with respect to $R$. Regarding the second assertion, fix $1\leq m<k$ and notice that
$$\frac{1}{N}\sum_{\substack{1\leq i\leq N\\ \vert \Omega_i^N\vert \geq R}}\vert \Omega_i^N\vert^m\leq \frac{1}{R^{k-m}}\frac{1}{N}\sum_{i=1}^N\vert \Omega_i^N\vert^k\leq \frac{M_k}{R^{k-m}},$$
for every $R>0$ and $N\in \mathbb{N}$. Taking supreme with respect to $N$ and limits as $R\rightarrow 0$ yields \eqref{E-Omega-moment-m-equiint}.
\end{proof}

\begin{cor}\label{C-mean-field-apriori-2-subcritical}
Let us assume that the hypothesis in Corollary \ref{C-mean-field-apriori-subcritical} and that the equi-summability condition \eqref{E-Omega-moment-1-equiint} are fulfilled. Then, the associated empirical measures $\mu^N$ in Corollary \ref{C-mean-field-apriori-subcritical} are also uniformly equicontinuous in $C([0,+\infty),\mathcal{P}_1(\mathbb{T}\times \mathbb{R})-W_1)$.
\end{cor}

\begin{proof}
Consider any continuous test function $\varphi$ so that
$$\vert\varphi(z,\Omega)\vert \leq C_\varphi(1+\vert \Omega\vert),\ \mbox{ for all }(z,\Omega)\in \mathbb{T}\times \mathbb{R},$$
for some $C_\varphi>0$. Recall the scaled cut-off functions $\xi_R=\xi_R(\Omega)$ in \eqref{E-scaled-cut-off}, for any $R>0$. Fix $\varepsilon>0$ and take $R>0$ large enough so that
\begin{equation}\label{E-31}
C_\varphi\left(\frac{1}{R}+1\right)\frac{1}{N}\sum_{\substack{1\leq i\leq N\\ \vert \Omega_i^N\vert \geq R}}\vert \Omega_i^N\vert\leq \frac{\varepsilon}{6}.
\end{equation}
Notice that $\varphi\,\xi_R\in C_0(\mathbb{T}\times \mathbb{R})$. Then, there exists $\widehat{\varphi}\in C^\infty_c(\mathbb{T}\times \mathbb{R})$ so that
\begin{equation}\label{E-32}
\Vert \varphi\,\xi_R-\widehat{\varphi}\Vert_{C_0(\mathbb{T}\times \mathbb{R})}\leq \frac{\varepsilon}{6}.
\end{equation}
By virtue of Corollary \ref{C-mean-field-apriori-subcritical} there exists $\delta>0$ so that
\begin{equation}\label{E-33}
\left\vert\int_{\mathbb{T}\times \mathbb{R}}\widehat{\varphi}\,d(\mu^N_{t_1}-\mu^N_{t_2})\right\vert\leq \frac{\varepsilon}{3},
\end{equation}
for every $\vert t_1-t_2\vert\leq \delta$ and each $N\in \mathbb{N}$. Now, consider the following split
$$\int_{\mathbb{T}\times \mathbb{R}}\varphi\,d(\mu^N_{t_1}-\mu^N_{t_2})=A_N(t_1,t_2)+B_N(t_1,t_2)+C_N(t_1,t_2),$$
where each term reads
\begin{align*}
A_N(t_1,t_2)&:=\int_{\mathbb{T}\times \mathbb{R}}\varphi\,(1-\xi_R)\,d(\mu^N_{t_1}-\mu^N_{t_2}),\\
B_N(t_1,t_2)&:=\int_{\mathbb{T}\times \mathbb{R}}(\varphi\,\xi_R-\widehat{\varphi})\,d(\mu^N_{t_1}-\mu^N_{t_2}),\\
C_N(t_1,t_2)&:=\int_{\mathbb{T}\times \mathbb{R}}\widehat{\varphi}\,d(\mu^N_{t_1}-\mu^{N}_{t_2}).
\end{align*}
First, \eqref{E-31} along with the equi-sumability condition \eqref{E-Omega-moment-1-equiint} implies
$$\vert A_N(t_1,t_2)\vert\leq \frac{\varepsilon}{3},$$
for all $t_1,t_2\geq 0$ and $N\in \mathbb{N}$. Second, \eqref{E-32} amounts to
$$\vert B_N(t_1,t_2)\vert \leq \frac{\varepsilon}{3},$$
for all $t_1,t_2\geq 0$ and $N\in \mathbb{N}$. Finally, \eqref{E-33} yields
$$\vert C_N(t_1,t_2)\vert\leq \frac{\varepsilon}{3},$$
for all $\vert t_1-t_2\vert\leq \delta$ and $N\in \mathbb{N}$. Putting everything together ends the proof.
\end{proof}

\begin{cor}\label{C-mean-field-compactness-subcritical}
Consider $\alpha\in (0,\frac{1}{2})$ and $K>0$. For any $N\in \mathbb{N}$, set $N$ oscillators with phases and natural frequencies given by the configurations
$$\Theta_0^N=(\theta_{1,0}^N,\ldots,\theta_{N,0}^N)\ \mbox{ and }\ \{\Omega_i^N:\,i=1,\ldots,N\}.$$
Consider the (forward-in-time) unique classical solution $\Theta^N(t)=(\theta_1(t),\ldots,\theta_N(t))$ to \eqref{E-kuramoto-discrete}-\eqref{E-kuramoto-kernel} as given in \cite{P-P-S} and set the corresponding empirical measures $\mu^N$ according to Definition \ref{D-empirical-measures-subcritical}. Assume that the equi-sumability condition \eqref{E-Omega-moment-1-equiint} holds, and take $M_1$ fulfilling \eqref{E-Omega-moment-1} according to Proposition \ref{P-Omega-moment-equint}. Then, for every fixed $T>0$, there exists a subsequence of $\mu^N$, that we denote in the same way for simplicity, and a limiting measure $f\in L^\infty(0,T;\mathcal{M}(\mathbb{T}\times \mathbb{R}))$ such that 
\begin{align}
\sup_{t\in [0,T]}\Vert \vert\Omega\vert\,f_t\Vert_{\mathcal{M}(\mathbb{T}\times \mathbb{R})}&\leq M_1,\label{E-18-1}\\
\lim_{R\rightarrow +\infty}\sup_{t\in [0,T]}\Vert \vert \Omega\vert\,\chi_{\vert \Omega\vert \geq R}\,f_t\Vert_{\mathcal{M}(\mathbb{T}\times \mathbb{R})}&= 0,\label{E-18-2}
\end{align}
and, in addition,
$$f\in W^{1,\infty}([0,T],C^1_0(\mathbb{T}\times \mathbb{R})^*)\cap C([0,T],\mathcal{P}_1(\mathbb{T}\times \mathbb{R})-W_1)\cap \mathcal{AC}_{\mathcal{M}}\cap\mathcal{T}_\mathcal{M},$$
for every fixed $T>0$. Moreover,
\begin{equation}\label{E-20}
\mu^N\rightarrow f\ \mbox{ in }\ C([0,T],\mathcal{P}_1(\mathbb{T}\times \mathbb{R})-W_1),
\end{equation}
where $W_1$ means the Kantorovich-Rubinstein distance.
\end{cor}

\begin{proof}
By virtue of the uniform estimates in the above Corollary \ref{C-mean-field-apriori-subcritical} along with the weak-star version of the Ascoli-Arzel\`a theorem, there exists some subsequence and $f\in L^\infty(0,T;\mathcal{M}(\mathbb{T}\times \mathbb{R}))$ so that
$$\mu^N\rightarrow f\ \mbox{ in }\ C([0,T],C^1_0(\mathbb{T}\times \mathbb{R})^*-\mbox{weak}\,*).$$
Recall that the  embedding
$C^1_0(\mathbb{T}\times \mathbb{R})\hookrightarrow C_0(\mathbb{T}\times \mathbb{R})$
is continuous and dense. Then, we can improve the above convergence into
$$\mu^N\rightarrow f\ \mbox{ in }\ C([0,T],\mathcal{M}(\mathbb{T}\times \mathbb{R})-\mbox{weak}\,*),$$
i.e.,
\begin{equation}\label{E-21}
\lim_{N\rightarrow\infty}\sup_{t\in [0,T]}\left\vert\int_{\mathbb{T}\times \mathbb{R}}\varphi\,d(\mu^N_t-f^N_t)\right\vert=0,
\end{equation}
for all $\varphi\in C_0(\mathbb{T}\times \mathbb{R})$. In order to augment the weak-star convergence into the Wasserstein one, we will first show that the other two properties that can be inherited by the limit. 

$\bullet$ \textit{Step 1}: Corollary \ref{C-mean-field-apriori-subcritical} yields
$$\left\vert \int_{\mathbb{T}\times \mathbb{R}}\varphi\,d(\mu^N_{t_1}-\mu^N_{t_2})\right\vert\leq (M_1+K\Vert h\Vert_{C(\mathbb{T})})\Vert \varphi\Vert_{C_0^1(\mathbb{T}\times \mathbb{R})}\vert t_1-t_2\vert,$$
for every $t_1,t_2\in [0,T]$, each $\varphi\in C^1_0(\mathbb{T}\times \mathbb{R})$, and any $N\in\mathbb{N}$. Taking limits as $N\rightarrow \infty$ and using \eqref{E-21}, we can obtain 
$$
\left\vert \int_{\mathbb{T}\times \mathbb{R}}\varphi\,d(f_{t_1}-f_{t_2})\right\vert\leq (M_1+K\Vert h\Vert_{C(\mathbb{T})})\Vert \varphi\Vert_{C_0^1(\mathbb{T}\times \mathbb{R})}\vert t_1-t_2\vert,
$$
for every $\varphi\in C^1_0(\mathbb{T}\times \mathbb{R})$ and each $t_1,t_2\in [0,T]$. Consequently, $f\in W^{1,\infty}([0,T],C^1_0(\mathbb{T}\times \mathbb{R})^*).$

$\bullet$ \textit{Step 2}: Again, recall the scaled cut-off functions $\xi_R=\xi_R(\Omega)$ in \eqref{E-scaled-cut-off}. Then, notice that the function
$$(z,\Omega)\in \mathbb{T}\times \mathbb{R}\longmapsto \vert\Omega\vert\,\xi_R(\Omega)$$
belongs to $C_0(\mathbb{T}\times \mathbb{R})$. Consequently, \eqref{E-21} implies
\begin{align*}
\Vert \vert\Omega\vert\,f_t\Vert_{\mathcal{M}(\mathbb{T}\times [-R,R])}&\leq \int_{\mathbb{T}\times \mathbb{R}}\vert\Omega\vert\,\xi_R(\Omega)\,d_{(z,\Omega)}f_t=\lim_{N\rightarrow\infty}\int_{\mathbb{T}\times \mathbb{R}}\vert\Omega\vert\,\xi_R(\Omega)\,d_{(z,\Omega)}\mu^N_t\\
&\leq \limsup_{N\rightarrow \infty}\Vert \vert\Omega\vert\mu^N_t\Vert_{\mathcal{M}(\mathbb{T}\times \mathbb{R})}=\limsup_{N\rightarrow \infty}\frac{1}{N}\sum_{i=1}^N\vert\Omega_i^N\vert\leq M_1,
\end{align*}
for every $t\in [0,T]$, and any $R>0$. Taking limit $R\rightarrow +\infty$ entails \eqref{E-18-1}. Similarly, take a couple $0<R<R'$, and consider the test function in $C_0(\mathbb{T}\times \mathbb{R})$ defined as follows
$$(z,\Omega)\in \mathbb{T}\times \mathbb{R}\longmapsto \vert\Omega\vert\,\xi_{R'}(\Omega)(1-\xi_{R/2}(\Omega)).$$
Then, an analogue argument yields
\begin{align*}
\Vert \vert\Omega\vert\,\chi_{R\leq \vert \Omega\vert\leq R'}f_t\Vert_{\mathcal{M}(\mathbb{T}\times \mathbb{R})}&\leq \int_{\mathbb{T}\times \mathbb{R}}\vert \Omega\vert\,\xi_{R'}(\Omega)(1-\xi_{R/2}(\Omega))\,d_{(z,\Omega)}f_t\\
&=\lim_{N\rightarrow\infty}\int_{\mathbb{T}\times \mathbb{R}}\vert \Omega\vert\,\xi_{R'}(\Omega)(1-\xi_{R/2}(\Omega))\,d_{(z,\Omega)}\mu^N_t\\
&\leq \limsup_{N\rightarrow\infty}\Vert \vert \Omega\vert\,\chi_{\vert \Omega\vert \geq R}\mu^N_t\Vert_{\mathcal{M}(\mathbb{T}\times \mathbb{R})}\leq \sup_{N\in \mathbb{N}}\frac{1}{N}\sum_{\substack{1\leq i\leq N\\ \vert \Omega_i^N\vert \geq R}}\vert \Omega_i^N\vert.
\end{align*}
Taking limit when $R'\rightarrow +\infty$ yields
$$\sup_{t\in [0,T]}\Vert \vert\Omega\vert\,\chi_{\vert \Omega\vert \geq R}\,f_t\Vert_{\mathcal{M}(\mathbb{T}\times \mathbb{R})}\leq \sup_{N\in \mathbb{N}}\frac{1}{N}\sum_{\substack{1\leq i\leq N\\ \vert \Omega_i^N\vert \geq R}}\vert \Omega_i^N\vert,$$
for every $R>0$. Finally, we obtain  \eqref{E-18-2}, as $R\rightarrow +\infty$. Recall that $\mu^N$ are uniformly equicontinuous in $C([0,T],\mathcal{P}_1(\mathbb{T}\times \mathbb{R})-W_1)$ thanks to Corollary \ref{C-mean-field-apriori-2-subcritical}. Then, we similarly infer that $f\in C([0,T],\mathcal{P}_1(\mathbb{T}\times \mathbb{R})-W_1)$.

$\bullet$ \textit{Step 3}: Now, take a continuous test function $\varphi\in C(\mathbb{T}\times \mathbb{R})$ with linear growth, that is,
$$\vert \varphi(z,\Omega)\vert \leq C_\varphi(1+\vert \Omega\vert),\ \forall\,(z,\Omega)\in \mathbb{T}\times \mathbb{R},$$
for some $C_\varphi>0$. The integral of interest can be split as follows
$$\left\vert\int_{\mathbb{T}\times \mathbb{R}}\varphi\,d(\mu^N_t-f_t)\right\vert\leq A_N^R(t)+B_N^R(t),$$
where each term reads
\begin{align*}
A_N^R(t)&:=\left\vert\int_{\mathbb{T}\times \mathbb{R}}\varphi\,\xi_R\,d(\mu^N_t-f_t)\right\vert,\\
B_N^R(t)&:=\left\vert\int_{\mathbb{T}\times \mathbb{R}}\varphi(1-\xi_R)\,d(\mu^N_t-f_t)\right\vert.
\end{align*}
Fix $\varepsilon>0$ and consider $R>0$ large enough so that 
$$C_\varphi\left(\frac{1}{R}+1\right)\left(\frac{1}{N}\sum_{\substack{1\leq i\leq N\\ \vert \Omega_i^N\vert \geq R}}\vert \Omega_i^N\vert\right)\leq \frac{\varepsilon}{4}.$$
This can be done by virtue of hypothesis \eqref{E-Omega-moment-1-equiint}. Then, it is clear that
\begin{align*}
B_N^R(t)&\leq C_\varphi\int_{\vert\Omega\vert\geq R}(1+\vert \Omega\vert)\,d_{(z,\Omega)}(\vert \mu^N_t\vert+\vert f_t\vert)\\
&\leq C_\varphi\left(\frac{1}{R}+1\right)\int_{\vert \Omega\vert \geq R}\vert \Omega\vert\,d_{(z,\Omega)}(\vert \mu^N_t\vert+\vert f_t\vert)\leq 2C_\varphi\left(\frac{1}{R}+1\right)\left(\frac{1}{N}\sum_{\substack{1\leq i\leq N\\ \vert \Omega_i^N\vert \geq R}}\vert \Omega_i^N\vert\right)\leq \frac{\varepsilon}{2},
\end{align*}
for every $N\in\mathbb{N}$, and every $t\in [0,T]$. Also, notice that the following function
$$(z,\Omega)\in \mathbb{T}\times \mathbb{R}\longmapsto \varphi(z,\Omega)\,\xi_R(\Omega),$$
belongs to $C_0(\mathbb{T}\times \mathbb{R})$. Then, applying \eqref{E-21} to such function, there exists $N_0\in \mathbb{N}$ so that 
$$A_N^R(t)=\left\vert \int_{\mathbb{T}\times \mathbb{R}}\varphi\,\xi_R\,d(\mu^N_t-f_t)\right\vert\leq \frac{\varepsilon}{2},$$
for every $N\geq N_0$, and every $t\in [0,T]$. Putting everything together implies the uniform-in-time convergence against any continuous function with linear growth, or, equivalently (see \cite[Definition 6.8, Theorem 6.9]{V}), the desired uniform-in-time convergence in the Kantorovich-Rubinstein distance $W_1$.
\end{proof}

\begin{rem}
The above convergence in the Kantorovich-Rubinstein distance $W_1$ can be improved to any other Wasserstein distance $W_p$ with $p>1$ when the equi-sumability condition \eqref{E-Omega-moment-1-equiint} is replaced with the general $p$-equi-sumability condition \eqref{E-Omega-moment-m-equiint}. Indeed, such assumption implies
\begin{align*}
\sup_{t\in [0,T]}\Vert \vert\Omega\vert^p\,f_t\Vert_{\mathcal{M}(\mathbb{T}\times \mathbb{R})}&\leq M_p,\\
\lim_{R\rightarrow 0}\sup_{t\in [0,T]}\Vert \vert \Omega\vert^p\,\chi_{\vert \Omega\vert \geq R}\,f_t\Vert_{\mathcal{M}(\mathbb{T}\times \mathbb{R})}&=0.
\end{align*}
Moreover,
$$
\mu^N\rightarrow f\ \mbox{ in }\ C([0,T],\mathcal{P}_p(\mathbb{T}\times \mathbb{R})-W_p).
$$
\end{rem}

\begin{lem}\label{L-mean-field-uniform-convergence-convolution}
Consider $\alpha\in (0,\frac{1}{2})$ and $K>0$. For every $N\in\mathbb{N}$, consider $N$ oscillators with phases and natural frequencies given by the configurations
$$\Theta_0^N=(\theta_{1,0}^N,\ldots,\theta_{N,0}^N)\ \mbox{ and }\ \{\Omega_i^N:\,i=1,\ldots,N\}.$$
Assume the equi-sumability condition \eqref{E-Omega-moment-1-equiint}, set the associated empirical measures $\mu^N$ according to Definition \ref{D-empirical-measures-subcritical} and any limit $f$ according to Corollary \ref{C-mean-field-compactness-subcritical}. Then,
$$\mathcal{V}[\mu^N]- \mathcal{V}[f]\longrightarrow0\ \mbox{ in }\ C([0,T],\mathfrak{X}C_b(\mathbb{T}\times \mathbb{R})).$$
\end{lem}

\begin{proof}
Since both $\mu^N$ and $f$ belong to $\widetilde{C}_\mathcal{M}$, then Corollary \ref{C-transport-field-holder} guarantees  that
$$\frac{\mathcal{V}[\mu^N]}{1+\vert\Omega\vert},\frac{\mathcal{V}[f]}{1+\vert\Omega\vert}\in C([0,T],\mathfrak{X}C_b(\mathbb{T}\times \mathbb{R})),$$
for every $N\in\mathbb{N}$. Consequently, the continuity of both vector fields is granted. Throughout the rest of the proof, we will show that the convergence is uniform. By the Stone--Weierstrass theorem, there exists $m\in \mathbb{N}$ and $\phi_1,\ldots,\phi_m,\psi_1,\ldots,\psi_m\in C(\mathbb{T})$ depending on $\varepsilon$ and $h$ so that
\begin{equation}\label{E-22}
\left\vert h(\theta-\theta')-\sum_{i=1}^m\phi_i(\theta)\psi_i(\theta')\right\vert\leq \frac{\varepsilon}{2K},\ \forall\,\theta,\theta'\in \mathbb{R},
\end{equation}
for every fixed $\varepsilon>0$. For simplicity, let us define
$$\widehat{h}(\theta,\theta'):=\sum_{i=1}^m\phi_i(\theta)\psi_i(\theta'),\ \forall\,\theta,\theta'\in \mathbb{R}.$$
Then, we have
$$\vert \mathcal{V}[\mu^N_t](\theta,\Omega)-\mathcal{V}[f_t](\theta,\Omega)\vert=\vert \mathcal{P}[\mu^N_t](\theta,\Omega)-\mathcal{P}[f_t](\theta,\Omega)\vert\leq F_N(t,\theta,\Omega)+G_N(t,\theta,\Omega),$$
where each term reads
\begin{align*}
F_N(t,\theta,\Omega)&:=K\left\vert\int_{(-\pi,\pi]\times \mathbb{R}}\widehat{h}(\theta,\theta')\,d_{(\theta',\Omega')}(\mu^N_t-f_t)\right\vert,\\
G_N(t,\theta,\Omega)&:=K\left\vert\int_{(-\pi,\pi]\times \mathbb{R}}(h(\theta-\theta')-\widehat{h}(\theta,\theta'))\,d_{(\theta',\Omega')}(\mu^N_t-f_t)\right\vert.
\end{align*}
Regarding the second term, the bound \eqref{E-22} automatically implies that $$G_N(t,\theta,\Omega)\leq \varepsilon,$$ 
for every $t\in [0,T]$ and $(\theta,\Omega)\in (-\pi,\pi]\times \mathbb{R}$. On the other hand, the first term can be bounded in the following way
$$F_N(t,\theta,\Omega)\leq K\sum_{i=1}^m\Vert \phi_i\Vert_{C(\mathbb{T})}\left\vert\int_{(-\pi,\pi]\times \mathbb{R}}\psi_i(\theta')\,d_{(\theta',\Omega')}(\mu^N_t-f_t)\right\vert.$$
Using \eqref{E-20} in Corollary \ref{C-mean-field-compactness-subcritical}, one obtains that
$$\limsup_{N\rightarrow\infty}\Vert \mathcal{V}[\mu^N]-\mathcal{V}[f]\Vert_{C([0,T],\mathfrak{X}C_b(\mathbb{T}\times \mathbb{R}))}\leq \varepsilon.$$
Since $\varepsilon>0$ is arbitrary, we conclude the proof of this result.
\end{proof}

\subsection{Passing to the limit}
Here we will show that any limit $f$ obtained in Corollary \ref{C-mean-field-compactness-subcritical} yields a weak measure-valued solution to \eqref{E-kuramoto-transport-TxR}, thus solving the initial value problem for \eqref{E-kuramoto-transport-TxR} with any initial data in $\mathcal{P}_1(\mathbb{T}\times \mathbb{R})$. Our first step is to show that any such initial datum can be approximated by an average of Dirac masses supported at a configuration with the equi-sumability condition \eqref{E-Omega-moment-1-equiint}. Although it follows from classical arguments, we will introduce the result in our particular setting for the sake of completeness.

\begin{lem}\label{L-discrete-approximation}
Consider any $\mu \in\mathcal{P}_1(\mathbb{T}\times \mathbb{R})$. Then, there exist $N$ oscillators with phases and natural frequencies given by the configurations
$$\Theta^N=(\theta_{1}^N,\ldots,\theta_{N}^N)\ \mbox{ and }\ \{\Omega_i^N:\,i=1,\ldots,N\},$$
for every $N\in\mathbb{N}$, verifying the equi-sumability condition \eqref{E-Omega-moment-1-equiint} such that the corresponding empirical measures $\mu^N\in \mathcal{P}(\mathbb{T}\times \mathbb{R})$ in Definition \ref{D-empirical-measures-subcritical} verify
$$\lim_{N\rightarrow \infty}W_1(\mu^N,\mu)=0.$$
\end{lem}

\begin{proof}
By a standard application of Kolmogorov's consistency theorem (see \cite{Va2}), there exists a probability space $(E,\mathcal{F},P)$ and a sequence of random variables $\{X_k\}_{k\in \mathbb{N}}$ with values in $\mathbb{T}\times \mathbb{R}$, namely, 
$$\begin{array}{cccc}
X_k=(Z_k,\Omega_k):&E&\longrightarrow& \mathbb{T}\times \mathbb{R},\\
 & \xi & \longmapsto & X_k(\xi),
\end{array}$$
so that $X_k$ are all independent and identically distributed with law $\mu$. Let us define the following random probability measure
$$\mu^N=\frac{1}{N}\sum_{k=1}^N\delta_{X_k}=\frac{1}{N}\sum_{k=1}^N \delta_{Z_k}(z)\otimes \delta_{\Omega_k}(\Omega).$$
A straightforward application of the strong Law of Large Numbers (see \cite{Va1}) shows that
\begin{equation}\label{E-24}
\lim_{N\rightarrow\infty}\int_{\mathbb{T}\times \mathbb{R}}\varphi\,d\mu^N_\xi= \int_{\mathbb{T}\times \mathbb{R}}\varphi\,d\mu,\ P-\mbox{a.s.}
\end{equation}
for every $\varphi\in L^1(\mathbb{T}\times \mathbb{R},d\mu)$. Let us define for each $R>0$ the functions
$$\varphi_R(z,\Omega):=\vert \Omega\vert\,\chi_{\vert \Omega\vert\geq R}.$$
By the assumptions, $\varphi_R\in L^1(\mathbb{T}\times\mathbb{R},d\mu)$ and, consequently, there exits $\{R_n\}_{n\in \mathbb{T}}\subseteq \mathbb{R}^+$ such that
\begin{equation}\label{E-25}
\int_{\mathbb{T}\times \mathbb{R}}\varphi_{R_n}\,d\mu\leq \frac{1}{2n}.
\end{equation}
Let us also set a dense sequence $\{\psi_k\}_{k\in\mathbb{N}}\subseteq C_0(\mathbb{T}\times \mathbb{R})$ in $L^1(\mathbb{T}\times \mathbb{R},d\mu)$. Note that $\psi_k\in L^1(\mathbb{T}\times \mathbb{R},d\mu)$ for every $k\in \mathbb{N}$. Hence, we can apply the strong Law of Large Numbers \eqref{E-24} to the whole family of functions
$$\{\varphi_{R_n}:\,n\in \mathbb{N}\}\cup\{\psi_k:\,k\in \mathbb{N}\},$$
to obtain
\begin{align}
\lim_{N\rightarrow\infty}\int_{\mathbb{T}\times \mathbb{R}}\varphi_{R_n}\,d\mu^N_\xi= \int_{\mathbb{T}\times \mathbb{R}}\varphi_{R_n}\,d\mu,&\ \mbox{ for all }\xi\in E\setminus E_n,\label{E-26}\\
\lim_{N\rightarrow\infty}\int_{\mathbb{T}\times \mathbb{R}}\psi_k\,d\mu^N_\xi=\int_{\mathbb{T}\times \mathbb{R}}\psi_k\,d\mu,&\ \mbox{ for all }\xi\in E\setminus F_k,\label{E-27}
\end{align}
for every $n,k\in \mathbb{N}$, where $E_n,F_k\subseteq E$ are $P$-negligible sets. Let us define the $P$-negligible set 
$$E':=\bigcup_{n,k\in \mathbb{N}}E_n\cup F_k.$$
Then, both \eqref{E-26} and \eqref{E-27} simultaneously hold, for all $n,k\in \mathbb{N}$ and each $\xi\in E\setminus E'$. On the one hand, \eqref{E-26} implies that there exists $N_n\in \mathbb{N}$ such that
\begin{equation}\label{E-28}
\left\vert\int_{\mathbb{T}\times \mathbb{R}}\varphi_{R_n}\,d(\mu^N_\xi-\mu)\right\vert\leq \frac{1}{2n},
\end{equation}
for every $n\in \mathbb{N}$,  every $N\geq N_n$, and each $\xi\in E\setminus E'$. Putting \eqref{E-25} and \eqref{E-28} together yields
\begin{equation}\label{E-29}
\int_{\mathbb{T}\times \mathbb{R}}\varphi_{R_n}\,d\mu^N_\xi\leq \frac{1}{n},
\end{equation}
for every $N\geq N_n$, each $\xi\in E\setminus E'$ and any $n\in \mathbb{N}$. Finally, let us pick a realization $\xi_0\in E\setminus E'$ and set
$$R_n':=\max\left\{R_n,\max_{1\leq k\leq N_n}\vert \Omega_k(\xi_0)\vert\right\},$$
for every $n\in \mathbb{N}$. Note that \eqref{E-29} amounts to
$$\sup_{N\in \mathbb{N}}\int_{\mathbb{T}\times \mathbb{R}}\varphi_R\,d\mu^N_{\xi_0}=\sup_{N\geq N_n}\int_{\mathbb{T}\times \mathbb{R}}\varphi_R\,d\mu^N_{\xi_0}\leq\frac{1}{n},$$ 
for each $R>R_n'$ and any $n\in \mathbb{N}$, thus yielding the equi-sumability condition \eqref{E-Omega-moment-1-equiint} for
$$\{\Omega_k(\xi_0):\,1\leq k\leq N\},\ \ N\in \mathbb{N}.$$
In addition, \eqref{E-27} implies that
\begin{equation}\label{E-30}
\mu^N_{\xi_0}\overset{*}{\rightharpoonup} \mu\ \mbox{ in }\ \mathcal{M}(\mathbb{T}\times \mathbb{R}),
\end{equation}
as a consequence of the density of $\{\psi_k\}_{k\in \mathbb{N}}$ in $C_0(\mathbb{T}\times \mathbb{R})$. Let us improve such weak-star convergence into convergence in the Rubinstein-Kantorovich metric $W_1$. Consider any continuous test function $\varphi$ with
$$\vert \varphi(z,\Omega)\vert\leq C(1+\vert \Omega\vert),\ \mbox{ for all }\,(z,\Omega)\in \mathbb{T}\times \mathbb{R},$$
for some constant $C>0$. Also, recover the cut-off functions $\xi_R=\xi_R(\Omega)$ in \eqref{E-scaled-cut-off} and consider the split
$$\int_{\mathbb{T}\times \mathbb{R}}\varphi\,d(\mu^N_{\xi_0}-\mu)=:A_N+B_N,$$
where each term reads
\begin{align*}
A_N&:=\int_{\mathbb{T}\times \mathbb{R}}\varphi\,\xi_R\,d(\mu^N_{\xi_0}-\mu),\\
B_N&:=\int_{\mathbb{T}\times \mathbb{R}}\varphi\,(1-\xi_R)\,d(\mu^N_{\xi_0}-\mu).
\end{align*}
On the one hand, note that
$$B_N\leq C \int_{\mathbb{T}\times \mathbb{R}}(1+\vert \Omega\vert)\,d\mu^N_{\xi_0}+C\int_{\mathbb{T}\times \mathbb{R}}(1+\vert \Omega\vert)\,d\mu.$$
Fix $\varepsilon>0$. Taking $R>0$ large enough, the assumption $\mu\in \mathcal{P}_1(\mathbb{T}\times \mathbb{R})$ along with the equi-sumability condition \eqref{E-Omega-moment-1-equiint} show that $B_N\leq \frac{\varepsilon}{2}$, for every $N\in \mathbb{N}$. For such $R$, note that $\varphi\,\xi_{R}\in C_0(\mathbb{T}\times \mathbb{R})$. Using the above weak-star convergence \eqref{E-30} of $\mu^N_{\xi_0}$, one obtains $N_0\in \mathbb{N}$ so that $A_N\leq \frac{\varepsilon}{2}$, for every $N\geq N_0$. Putting everything together, we conclude the proof.
\end{proof}

We are now ready to obtain the mean field limit that, in particular, yields the following existence result.

\begin{theo}\label{T-mean-field-existence-subcritical}
Consider $\alpha\in (0,\frac{1}{2})$, $K>0$ and set any initial datum $f_0\in \mathcal{P}_1(\mathbb{T}\times \mathbb{R})$. Then, for every $T>0$ there exists a weak measure-valued solution $f\in \mathcal{AC}_\mathcal{M}\cap \mathcal{T}_\mathcal{M}$ to the initial value problem \eqref{E-kuramoto-transport-TxR}. In addition, \eqref{E-18-1}-\eqref{E-18-2} holds  and
$$f\in C^{0,1}([0,T],C^1_0(\mathbb{T}\times \mathbb{R})^*)\cap C([0,T],\mathcal{P}_1(\mathbb{T}\times \mathbb{R})-W_1).$$
\end{theo}

\begin{proof}
Our first step is to take a discrete approximation like in Lemma \ref{L-discrete-approximation}. Namely, consider $N$ oscillators with phases and natural frequencies given by the configurations
$$\Theta^N_0=(\theta_{1,0}^N,\ldots,\theta_{N,0}^N)\ \mbox{ and }\ \{\Omega_i^N:\,i=1,\ldots,N\},$$
for every $N\in\mathbb{N}$ so that they verify the equi-sumability condition \eqref{E-Omega-moment-1-equiint} and the associated empirical measures $\mu^N_t\in \mathcal{P}(\mathbb{T}\times \mathbb{R})$ in Definition \ref{D-empirical-measures-subcritical} verify
$$\lim_{N\rightarrow \infty}W_1(\mu^N_0,f_0)=0.$$
Using Theorem \ref{T-empirical-measures-subcritical}, we infer that $\mu^N$ are weak-measure valued solutions to \eqref{E-kuramoto-transport-TxR} issued at $\mu^N_0$. Then, they verify the following weak formulation (see Definition \ref{D-weak-measure-solution})
\begin{equation}\label{E-34}
\int_0^T\int_{\mathbb{T\times \mathbb{R}}}\frac{\partial\varphi}{\partial t}\,d_{(z,\Omega)}\mu^N_t\,dt+\int_0^T\int_{\mathbb{T}\times \mathbb{R}}\left<\mathcal{V}[\mu^N_t],\nabla_{(z,\Omega)}\varphi\right>\,d_{(z,\Omega)}\mu^N_t\,dt=-\int_{\mathbb{T}\times \mathbb{R}}\varphi(0,z,\Omega)d_{(z,\Omega)}\mu^N_0,
\end{equation}
for every $\varphi\in C^1_c([0,T)\times \mathbb{T}\times \mathbb{R})$. Using Corollary \ref{C-mean-field-compactness-subcritical}, consider any weak limit $f$ of a subsequence of $\mu^N$, that we still denote in the same way for simplicity. In particular, recall that
$$\mu^N\rightarrow f\ \mbox{ in }\ C([0,T],\mathcal{P}_1(\mathbb{T}\times \mathbb{R})-W_1).$$ 
Now, we can pass to the limit in the weak formulation \eqref{E-34} as $N\rightarrow \infty$. Specifically, regarding the first and third term, the passage to the limit is clear by linearity. Regarding the nonlinear term, let us show that the following term vanishes in the limit $N\rightarrow\infty$
$$I_N:=\int_0^T\int_{\mathbb{T}\times \mathbb{R}}\left<\mathcal{V}[\mu^N_t],\nabla_{(z,\Omega)}\varphi\right>\,d\mu^N_t\,dt-\int_0^T\int_{\mathbb{T}\times \mathbb{R}}\left<\mathcal{V}[f_t],\nabla_{(z,\Omega)}\varphi\right>\,df_t\,dt,$$
for any given $\varphi\in C^1_c([0,T)\times \mathbb{T}\times \mathbb{R})$. Indeed, consider the following split
$$I_N=I_N^1+I_N^2,$$
where each term reads
\begin{align*}
I_N^1&:=\int_0^T\int_{\mathbb{T}\times \mathbb{R}}\left<\mathcal{V}[\mu^N_t]-\mathcal{V}[f_t],\nabla_{(z,\Omega)}\varphi\right>\,d\mu^N_t\,dt,\\
I_N^2&:=\int_0^T\int_{\mathbb{T}\times \mathbb{R}}\left<\mathcal{V}[f_t],\nabla_{(z,\Omega)}\varphi\right>\,d(\mu^N_t-f_t)\,dt.
\end{align*}
On the one hand, note that
$$\vert I_N^1\vert\leq \Vert \varphi\Vert_{C_b([0,T)\times \mathbb{T}\times \mathbb{R})}\Vert \mathcal{V}[\mu^N]-\mathcal{V}[f]\Vert_{C([0,T],\mathfrak{X}C_b(\mathbb{T}\times \mathbb{R}))}\rightarrow 0,$$
by virtue of Corollary \ref{C-mean-field-compactness-subcritical}. On the other hand, notice that $\left<\mathcal{V}[f],\nabla\varphi\right>\in C_c([0,T)\times \mathbb{T}\times \mathbb{R})$ thanks to Corollary \ref{C-transport-field-holder}. Then, we can also pass to the limit in $I_N^2$ to show that $I_n^2\rightarrow 0$, thus ending the proof.
\end{proof}

\section{Uniqueness and rigorous mean field limit}\label{S-uniqueness-mean-field}
The purpose of this part is to derive an upper bound for the growth of some Wasserstein-type distance between any two weak measure-valued solutions of \eqref{E-kuramoto-transport-TxR}. This is a Dobrushin-type inequality that has long been studied to show stability in mean-field equations. It is originally devoted to H. Neunzert when the ambient space is $\mathbb{R}^d$ and the kernel is Lipschitz, see \cite{N}. In that case, the bounded-Lipschitz distance fits with the Lipschitz-continuity property of the kernel. Mimicking the ideas in the above-mentioned paper, an analogue result has been explored for the classical Kuramoto-Sakaguchi equation in \cite{C-C-H-K-K,L}, where the kernel is still Lipschitz-continuous (it agrees with the sine function). However, such approach fails in our case because the singularly-weighted kernel $h$ is no longer Lipszhitz-continuous but barely H\"{o}lder-continuous in the subcritical case. Indeed, it is even discontinuous in the critical and supercritical regimes, that will be studied later in forthcoming sections.

Very recently, the case of non-Lipschitz kernels has been explored for the aggregation equations in $\mathbb{R}^d$. This is a sort of gradient system governed by the negative gradient of a $\lambda$-convex potential, see \cite{C-J-L-V}. In such case, the quadratic Wasserstein distance $W_2$ was considered instead of the bounded-Lipschitz distance for Lipschitz interactions. In our case \eqref{E-kuramoto-transport-TxR}, we will not use any gradient structure. Also, the natural frequencies have been introduced as heterogeneities in the system. That requires turning the ambient Euclidean space $\mathbb{R}^d$ for the aggregation equation into the Riemannian manifold $\mathbb{T}\times \mathbb{R}$, that introduces both the periodicity in $\theta$ and heterogeneity $\Omega$. Despite the above differences with the classical Lipschitz mean-field models and the aggregation equation, we will recover Dobrushin-type estimates. On the one hand, we will prove such estimate for a fibered version of the quadratic Wasserstein distance and general weak measure-valued solutions issued at initial data with the same distribution of natural frequencies. This will be the cornerstone in our uniqueness result. On the other hand, for initial data whose distribution of natural frequencies differ, we will obtain the sub-exponential growth of the classical quadratic Wasserstein distance as long as such distributions of natural frequencies have bounded second order moments. The latter estimate will be used to derive the rigorous mean field limit.

\subsection{Stability of a fibered quadratic Wasserstein distance and uniqueness result}\label{SS-uniqueness}

Before entering into details, we need to revisit some aspects from measure theory and optimal transport that we will need in the sequel. 

\begin{theo}[Disintegration]\label{T-disintegration}
Let $X$ and $Y$ be separable complete metric spaces and define the projection mapping
$$\begin{array}{llll}
\pi_Y: & X\times Y & \longrightarrow & Y,\\
 & (x,y) & \longmapsto & y. 
\end{array}
$$
Consider any Borel probability measure $\mu\in \mathcal{P}(X\times Y,\mathcal{B}(X\times Y))$ and the $Y$-marginal probability measure $\nu:=(\pi_Y)_{\#}\mu$. Then, there exits a Borel measurable map
$$\begin{array}{ccc}
(Y,\mathcal{B}(Y)) & \longrightarrow & \mathcal{P}(X,\mathcal{B}(X)),\\
y & \longmapsto & \mu(\cdot\vert y),
\end{array}
$$
such that the following formula holds true
\begin{equation}\label{E-disintegration}
\int_{X\times Y}\varphi(x,y)d_{(x,y)}\mu=\int_Y\left(\int_X \varphi(x,y)\,d_x\mu(\cdot\vert y)\right)\,d_y\nu,
\end{equation}
for every Borel-measurable map $\varphi:X\times Y\longrightarrow\mathbb{R}$.
\end{theo}

Such a family $\{\mu(\cdot\vert y)\}_{y\in Y}$ is called a disintegration of $\mu$ (or conditional probabilities with respect to $y$). It is uniquely defined $\nu$-a.e. in $Y$, see \cite[Theorem 5.3.1]{A-G-S} and \cite[III.70]{D-M} for more details. 

\begin{lem}\label{L-disintegration-flow-solutions}
Consider $\alpha\in (0,\frac{1}{2})$, $K>0$ and let $f\in \mathcal{AC}_\mathcal{M}\cap \mathcal{T}_\mathcal{M}$ be a weak mesure-valued solution to \eqref{E-kuramoto-transport-TxR} with initial datum $f_0\in \mathcal{P}(\mathbb{T}\times \mathbb{R})$, according to Theorem \ref{T-mean-field-existence-subcritical}. Let $X(t;0,z,\Omega)=(Z(t;0,z,\Omega),\Omega)$ be the flow associated with the transport field $\mathcal{V}[f]$, according to Theorem \ref{T-well-posedness-characteristic-system}. Then,
\begin{enumerate}
\item The solution remains normalized, i.e., 
$$f_t\in \mathcal{P}(\mathbb{T}\times \mathbb{R}),\ \mbox{ for all }t\geq 0.$$
\item The $\Omega$-marginal remains unchanged, i.e., 
$$(\pi_\Omega)_{\#}f_t=(\pi_\Omega)_{\#}f_0\equiv g,\ \mbox{ for all }t\geq 0.$$
\item The disintegrations $\{f_t(\cdot\vert \Omega)\}_{\Omega\in\mathbb{R}}\subseteq \mathcal{P}(\mathbb{T})$ with respect to $\Omega\in \mathbb{R}$ propagate through the flow, i.e.,
$$f_t(\cdot\vert \Omega)=Z(t;0,\cdot,\Omega)_{\#}f_0(\cdot\vert\Omega),\ \mbox{ for all }t\geq 0,\ g\mbox{-a.e. }\Omega\in \mathbb{R}.$$
\end{enumerate}
Here, $\pi_\Omega$ is the projection in $\Omega$, see \eqref{E-projections}.
\end{lem}

\begin{proof}
Recall that, as discussed in Remark \ref{R-types-solutiones}, $f$ is also a solution in the sense of the flow, that is $f_t=X(t;0,\cdot)_{\#}f_0$, for every $t\geq 0$. Recall that it is characterized by the following identity
\begin{equation}\label{E-push-forward}
\int_{\mathbb{T}\times \mathbb{R}}\varphi(z,\Omega)\,d_{(z,\Omega)}f_t=\int_{\mathbb{T}\times \mathbb{R}}\varphi(Z(t;0,z,\Omega),\Omega)\,d_{(z,\Omega)}f_0,
\end{equation}
for any $\varphi\in C_b(\mathbb{T}\times \mathbb{R})$. In particular, taking $\varphi\equiv 1$ the first assertion becomes clear. Regarding the second assertion, consider any text function $\phi\in C_b(\mathbb{R})$ and notice the following chain of identities
$$\int_{\mathbb{R}}\phi(\Omega)\,d_\Omega[(\pi_\Omega)_{\#}f_t]=\int_{\mathbb{T}\times \mathbb{R}}\phi(\Omega)\,d_{(z,\Omega)}f_t=\int_{\mathbb{T}\times \mathbb{R}}\phi(\Omega)\,d_{(z,\Omega)}f_0=\int_{\Omega}\phi(\Omega)d_\Omega[(\pi_\Omega)_{\#}f_0],$$
where we have used the bounded-continuous test function $\varphi(z,\Omega)=\phi(\Omega)$ in Equation \eqref{E-push-forward}. Hence, the second claim is apparent by definition. Finally, consider any $\phi\in C(\mathbb{T})$ and $\psi\in C_b(\mathbb{R})$ and define $\varphi(z,\Omega)=\phi(z)\psi(\Omega)$. Using the disintegration formula \eqref{E-disintegration} in both members of  \eqref{E-push-forward} implies
$$\int_{\mathbb{R}}\psi(\Omega)\left(\int_{\mathbb{T}}\phi(z,\Omega)\,d_zf_t(\cdot \vert \Omega)\right)\,d_\Omega g=\int_{\mathbb{R}}\psi(\Omega)\left(\int_{\mathbb{T}}\phi(Z(t;0,z,\Omega))d_zf_0(\cdot \vert \Omega)\right)\,d_\Omega g.$$
Since $\psi$ is arbitrary we get
$$\int_{\mathbb{T}}\phi(z,\Omega)\,d_zf_t(\cdot \vert \Omega)=\int_{\mathbb{T}}\phi(Z(t;0,z,\Omega))d_zf_0(\cdot \vert \Omega),$$
$g$-a.e. $\Omega\in \mathbb{R}$, for every $t\geq 0$, thus ending the proof of the result.
\end{proof}

In the following result, a new Wasserstein-type distance is introduce in a subspace of $\mathcal{P}(\mathbb{T}\times \mathbb{R})$. This will be the cornerstone in our uniqueness result to avoid the assumption of bounded $\Omega$-moments of the solutions. 

\begin{pro}\label{P-L1(RK)-distance}
Consider any probability measure $g\in \mathcal{P}(\mathbb{R})$ and define
$$\mathcal{P}_g(\mathbb{T}\times \mathbb{R}):=\left\{\mu\in\mathcal{P}(\mathbb{T}\times \mathbb{R}):\,(\pi_\Omega)_{\#}\mu=g\right\},$$
where $\pi_\Omega$ is the projection \eqref{E-projections}. Also, let us define the map $W_{2,g}$ as follows
$$W_{2,g}(\mu^1,\mu^2)=\left(\int_{\mathbb{R}}W_2(\mu^1(\cdot\vert \Omega),\mu^2(\cdot\vert\Omega))^2\,d_\Omega g\right)^{1/2},$$
for every $\mu^1,\mu^2\in\mathcal{P}_g(\mathbb{T}\times \mathbb{R})$, where $\{\mu^1(\cdot\vert\Omega)\}_{\Omega\in\mathbb{R}}$ and $\{\mu^2(\cdot\vert\Omega)\}_{\Omega\in \mathbb{R}}$ stand for their associated families of disintegrations with respect to $\Omega$. Then, $(\mathcal{P}_g(\mathbb{T}\times \mathbb{R}),W_{2,g})$ is a metric space.
\end{pro}

The proof is clear and is a consequence of (disintegration) Theorem  \ref{T-disintegration} along with the fact that $W_2(\cdot,\cdot)$ is a distance in $\mathcal{P}(\mathbb{T})$ and $\Vert\cdot\Vert_{L^2(\mathbb{R},dg)}$ is a norm in $L^2(\mathbb{R},dg)$.

\begin{theo}\label{T-growth-L(RK)}
Consider $\alpha\in (0,\frac{1}{2})$, $K>0$ and let $f^1,f^2\in \mathcal{AC}_\mathcal{M}\cap \mathcal{T}_\mathcal{M}$ be weak measure-valued solutions to  \eqref{E-kuramoto-transport-TxR} with initial data $f^1_0,f^2_0\in \mathcal{P}(\mathbb{T}\times \mathbb{R})$ according to Theorem \ref{T-mean-field-existence-subcritical}. Let us set their distributions of natural frequencies $g^i=\left(\pi_{\Omega}\right)_{\#}f^i_0$. If $g^1\equiv g^2=:g$, then
$$W_{2,g}(f^1_t,f^2_t)\leq W_{2,g}(f^1_0,f^2_0)e^{2KL_0t},$$
for every $t\geq 0$, where $L_0$ is the one-sided Lipschitz constant of $-h$ in Lemma \ref{L-split-kernel}.
\end{theo}

\begin{proof}
Again, $f^1,f^2$ are solutions in the sense of the flow by Remark \ref{R-types-solutiones}. For $g$-a.e. $\Omega\in\mathbb{R}$ fixed, let us consider the corresponding term of the family of disintegrations at the initial time, i.e., $f^1_0(\cdot\vert \Omega)$ and $f^2_0(\cdot\vert \Omega)$. Set an optimal transference plan from the former probability measure in $\mathbb{T}$ to the latter one, i.e.,
$$\mu_{0,\Omega}\in\Pi(f^1_0(\cdot \vert \Omega),f^2_0(\cdot \vert \Omega)):=\left\{\mu\in\mathcal{P}(\mathbb{T}\times \mathbb{T}):\,(\pi_1)_{\#}\mu=f^1_0(\cdot\vert \Omega)\ \mbox{and}\ (\pi_2)_{\#}\mu=f^2_0(\cdot\vert \Omega)\right\},$$
so that the $2$-Wasserstein distance is attained
$$W_2(f^1_0(\cdot\vert\Omega),f^2_0(\cdot\vert\Omega))^2=\int_{\mathbb{T}}\int_{\mathbb{T}}d(z_1,z_2)^2d_{(z_1,z_2)}\mu_{0,\Omega}.$$
Here, we are denoting the projections $\pi_1(z,z')=z$ and $\pi_2(z,z')=z'$. The existence of an optimal transference plan is granted by Brenier's theorem, see \cite{A-G-S,Sa,V}. Then, we can define a competitor transference plan at time $t$ as push-forward of the initial one, that is,
$$\mu_{t,\Omega}:=(Z_{f^1}(t;0,\cdot,\Omega)\otimes Z_{f^2}(t;0,\cdot,\Omega))_{\#}{\mu_{0,\Omega}}\in \mathcal{P}(\mathbb{T}\times \mathbb{T}),$$
where $X_{f^i}(t;0,z,\Omega)=(Z_{f^i}(t;0,z,\Omega),\Omega)$ is the characteristic flow associated with the transport field $\mathcal{V}[f^i]$ according to Theorem \ref{T-well-posedness-characteristic-system}. Notice that by definition
\begin{align*}
(\pi_1)_{\#}\mu_{t,\Omega}&=Z_{f^1}(t;0,\cdot,\Omega)_{\#}((\pi_1)_{\#}\mu_{0,\Omega})=Z_{f^1}(t;0,\cdot,\Omega)_\#f^1_0(\cdot\vert\Omega),\\
(\pi_2)_{\#}\mu_{t,\Omega}&=Z_{f^2}(t;0,\cdot,\Omega)_{\#}((\pi_2)_{\#}\mu_{0,\Omega})=Z_{f^2}(t;0,\cdot,\Omega)_\#f^2_0(\cdot\vert\Omega).
\end{align*}
Using the third statement in Lemma \ref{L-disintegration-flow-solutions} we conclude that
$$(\pi_1)_{\#}\mu_{t,\Omega}=f^1_t(\cdot\vert\Omega)\ \mbox{ and }\ (\pi_2)_{\#}\mu_{t,\Omega}=f^2_t(\cdot\vert\Omega).$$
Thus, $\mu_{t,\Omega}\in\Pi(f^1_t(\cdot\vert\Omega),f^2_t(\cdot\vert\Omega))$ and, consequently,
\begin{align*}
\frac{1}{2}W_2(f^1_t(\cdot\vert\Omega),f^2_t(\cdot\vert\Omega))^2&\leq \int_{\mathbb{T}}\int_{\mathbb{T}}\frac{1}{2}d(z_1,z_2)^2\,d_{(z_1,z_2)}\mu_{t,\Omega}\\
&=\int_{\mathbb{T}}\int_{\mathbb{T}}\frac{1}{2}d(Z_{f^1}(t;0,z_1,\Omega),Z_{f^2}(t;0,z_2,\Omega))^2\,d_{(z_1,z_2)}\mu_{0,\Omega}.
\end{align*}
Integrating the above inequality against $g$ yields
$$\frac{1}{2}W_{2,g}(f^1_t,f^2_t)^2\leq \int_{\mathbb{R}}\int_{\mathbb{T}}\int_{\mathbb{T}}\frac{1}{2}d(Z_{f^1}(t;0,z_1,\Omega),Z_{f^2}(t;0,z_2,\Omega))^2\,d_{(z_1,z_2)}\mu_{0,\Omega}\,d_\Omega g=:I(t).$$
We are interested in proving some Gr\"{o}nwall-type inequality for $I=I(t)$. Notice that two different transport fields $\mathcal{V}[f^1]$ and $\mathcal{V}[f^2]$. Again, for every $z_1=e^{i\theta_1},z_2=e^{i\theta_2}\in \mathbb{T}$ and $g$-a.e. $\Omega\in \mathbb{R}$ fixed, let us define the trajectories
$$
Z_1(t)=Z_{f^1}(t;0,z_1,\Omega)=e^{i\Theta_1(t)}\ \mbox{ and }\ Z_2(t)=Z_{f^2}(t;0,z_2,\Omega)=e^{i\Theta_2(t)},
$$
where $\Theta_1(t)=\Theta_{f^1}(t;0,\theta_1,\Omega)$ and $\Theta_2(t)=\Theta_{f^2}(t;0,\theta_2,\Omega)$ are the unique solutions to \eqref{E-characteristic-system-2} in Lemma \ref{L-existence-characteristic-system} with respective initial datum $\theta_1$ and $\theta_2$. Consider a minimizing geodesic $\gamma_t:[0,1]\longrightarrow\mathbb{T}$ joining $Z_1(t)$ to $Z_2(t)$, for every fixed $t\geq 0$. Since the map
$$t\longmapsto\frac{1}{2}d^2(Z_1(t),Z_2(t)),$$
is clearly absolutely continuous, we can use previous arguments like in the proof of Theorem \ref{T-well-posedness-characteristic-system} to achieve the following estimate
$$\frac{d}{dt}\frac{1}{2}d^2(Z_1(t),Z_2(t))\leq -\left<\mathcal{P}[f^1_t](Z_1(t),\Omega)iZ_1(t),\gamma_t'(0)\right>-\left<\mathcal{P}[f^2_t](Z_2(t),\Omega)iZ_2(t),-\gamma_t'(1)\right>.$$
Recall that it is supported by the one-sided Dini upper directional differentiability of the squared distance in $\mathbb{T}$ (see Appendix \ref{appendix-differentiability-distance}). Now, let us describe such geodesics in $\gamma_t$. The way to go is analogous to that in the proof of Lemma \ref{L-transport-field-sided-Lipschitz}, namely, consider $\theta_{21}(t):=\overline{\Theta_2(t)-\Theta_1(t)}$, the representative of $\Theta_2(t)-\Theta_1(t)$ modulo $2\pi$ that lies in the interval $(-\pi,\pi]$. There are two different cases.

$\bullet$ \textit{Case 1}: $\theta_{21}(t)\in (-\pi,\pi)$. In this case there exists only one such minimizing geodesic and it reads
$$\gamma_t(s)=e^{i(\Theta_1(t)+s\theta_{21}(t))},\ s\in [0,1].$$
Hence, the above inequality can be restated as
$$\frac{d}{dt}\frac{1}{2}d^2(Z_1(t),Z_2(t))\leq \left(\mathcal{P}[f^2_t](\Theta_2(t),\Omega)-\mathcal{P}[f^1_t](\Theta_1(t),\Omega)\right)\theta_{21}(t).$$

$\bullet$ \textit{Case 2}: $\theta_{21}(t)=\pi$. In this second case there are exactly two minimizing geodesics
$$\gamma_{t,\pm}(s)=e^{i(\Theta_1(t)\pm \pi s)},\ s\in [0,1].$$
In such case, the above inequality reads
$$\frac{d}{dt}\frac{1}{2}d^2(Z_1(t),Z_2(t))\leq \left(\mathcal{P}[f^2_t](\Theta_2(t),\Omega)-\mathcal{P}[f^1_t](\Theta_1(t),\Omega)\right)(\pm\pi).$$

Putting everything together, we arrive at the following inequality
\begin{multline*}
\frac{d}{dt}\frac{1}{2}d(Z_{f^1}(t;0,z_1,\Omega),Z_{f^2}(t;0,z_2,\Omega))\\
\leq (\mathcal{P}[f^1_t](\Theta_{f^1}(t;0,\theta_1,\Omega),\Omega)-\mathcal{P}[f^2_t](\Theta_{f^2}(t;0,\theta_2,\Omega),\Omega))\overline{\Theta_{f^1}(t;0,\theta_1,\Omega)-\Theta_{f^2}(t;0,\theta_2,\Omega)},
\end{multline*}
for every $\theta_1,\theta_2$ and almost every $t\geq 0$. By virtue of the dominated convergence theorem, we can take derivatives under the integral sign to show
\begin{multline}\label{E-1}
\frac{dI}{dt}\leq \int_{\mathbb{R}}\int_{(-\pi,\pi]}\int_{(-\pi,\pi]} (\mathcal{P}[f^1_t](\Theta_{f^1}(t;0,\theta_1,\Omega),\Omega)-\mathcal{P}[f^2_t](\Theta_{f^2}(t;0,\theta_2,\Omega),\Omega))\\
\times \overline{\Theta_{f^1}(t;0,\theta_1,\Omega)-\Theta_{f^2}(t;0,\theta_2,\Omega)}\,d_{(\theta_1,\theta_2)}\mu_{0,\Omega}\,d_\Omega g.
\end{multline}
Now, we need to identify $\mathcal{V}[\mu]$ as push-forward of initial data. Indeed, notice that by Theorem \ref{T-disintegration} and Lemma \ref{L-disintegration-flow-solutions}
\begin{align*}
\mathcal{P}[f^i_t](\theta,\Omega)&=\Omega-K\int_{(-\pi,\pi]}\int_{\mathbb{R}}h(\theta-\theta')\,d_{(\theta',\Omega')}f^i_t\\
&=\Omega-K\int_{\mathbb{R}}\left(\int_{(-\pi,\pi]}h(\theta-\theta')\,d_{\theta'}f^i_t(\cdot \vert\Omega')\right)\,d_{\Omega'} g\\
&=\Omega-K\int_{\mathbb{R}}\left(\int_{(-\pi,\pi]}h(\theta-\Theta_{f^i}(t;0,\theta_i',\Omega'))\,d_{\theta'_i}f^i_0(\cdot \vert\Omega')\right)\,d_{\Omega'} g.
\end{align*}
Recall that $(\pi_1)_{\#}\mu_{0,\Omega'}=f^1_0(\cdot \vert \Omega')$ and $(\pi_2)_{\#}\mu_{0,\Omega'}=f^2_0(\cdot \vert \Omega')$. Then, we obtain
\begin{align}
\mathcal{P}[f^1_t](\theta,\Omega)&=\Omega-K\int_{\mathbb{R}}\left(\int_{(-\pi,\pi]}\int_{(-\pi,\pi]}h(\theta-\Theta_{f^1}(t;0,\theta_1',\Omega'))\,d_{(\theta_1',\theta_2')}\mu_{0,\Omega'}\right)\,d_{\Omega'} g,\label{E-2}\\
\mathcal{P}[f^2_t](\theta,\Omega)&=\Omega-K\int_{\mathbb{R}}\left(\int_{(-\pi,\pi]}\int_{(-\pi,\pi]}h(\theta-\Theta_{f^2}(t;0,\theta_2',\Omega'))\,d_{(\theta_1',\theta_2')}\mu_{0,\Omega'}\right)\,d_{\Omega'} g.\label{E-3}
\end{align}
Putting \eqref{E-2}-\eqref{E-3} into \eqref{E-1}, we obtain the following expression
\begin{multline}\label{E-4}
\frac{dI}{dt}\leq -K\int_{((-\pi,\pi]^2\times \mathbb{R})^2} (h(\Theta_{f^1}(t;0,\theta_1,\Omega)-\Theta_{f^1}(t;0,\theta_1',\Omega'))-h(\Theta_{f^2}(t;0,\theta_2,\Omega)-\Theta_{f^2}(t;0,\theta_2',\Omega')))\\
\times\overline{\Theta_{f^1}(t;0,\theta_1,\Omega)-\Theta_{f^2}(t;0,\theta_2,\Omega)} \,d_{(\theta_1,\theta_2)}\mu_{0,\Omega}\,d_\Omega g\,d_{(\theta_1',\theta_2')}\mu_{0,\Omega'}\,d_{\Omega'}g
\end{multline}
Now, let us change variables $(\theta_1,\theta_2,\Omega)$ with $(\theta_1',\theta_2',\Omega')$
\begin{multline}\label{E-5}
\frac{dI}{dt}\leq -K\int_{((-\pi,\pi]^2\times \mathbb{R})^2} -(h(\Theta_{f^1}(t;0,\theta_1,\Omega)-\Theta_{f^1}(t;0,\theta_1',\Omega'))-h(\Theta_{f^2}(t;0,\theta_2,\Omega)-\Theta_{f^2}(t;0,\theta_2',\Omega')))\\
\times\overline{\Theta_{f^1}(t;0,\theta_1',\Omega')-\Theta_{f^2}(t;0,\theta_2',\Omega')} \,d_{(\theta_1,\theta_2)}\mu_{0,\Omega}\,d_\Omega g\,d_{(\theta_1',\theta_2')}\mu_{0,\Omega'}\,d_{\Omega'}g,
\end{multline}
where the antisymmetry of the kernel $h$ around the origin has been used. Taking the mean value of both expressions \eqref{E-4} and \eqref{E-5} yields
\begin{multline}\label{E-6}
\frac{dI}{dt}\leq \frac{K}{2}\int_{((-\pi,\pi]^2\times \mathbb{R})^2} -(h(\Theta_{f^1}(t;0,\theta_1,\Omega)-\Theta_{f^1}(t;0,\theta_1',\Omega'))-h(\Theta_{f^2}(t;0,\theta_2,\Omega)-\Theta_{f^2}(t;0,\theta_2',\Omega')))\\
\times\left(\overline{\Theta_{f^1}(t;0,\theta_1,\Omega)-\Theta_{f^2}(t;0,\theta_2,\Omega)}-\overline{\Theta_{f^1}(t;0,\theta_1',\Omega')-\Theta_{f^2}(t;0,\theta_2',\Omega')} \right)\\
 \times\,d_{(\theta_1,\theta_2)}\mu_{0,\Omega}\,d_\Omega g\,d_{(\theta_1',\theta_2')}\mu_{0,\Omega'}\,d_{\Omega'}g.
\end{multline}
Denote
\begin{align*}
\Theta_1&:=\Theta_{f^1}(t;0,\theta_1,\Omega), & \Theta_1'&:=\Theta_{f^1}(t;0,\theta_1',\Omega'),\\
\Theta_2&:=\Theta_{f^2}(t;0,\theta_2,\Omega), & \Theta_2'&:=\Theta_{f^2}(t;0,\theta_2',\Omega'),
\end{align*}
for almost every $t\geq 0$, each $\theta_1,\theta_2,\theta_1',\theta_2'\in(-\pi,\pi]$ and any $\Omega,\Omega'$ fixed, then the integrand reads
$$\left((-h)(\Theta_1-\Theta_1')-(-h)(\Theta_2-\Theta_2')\right)\left(\overline{\Theta_1-\Theta_2}-\overline{\Theta_1'-\Theta_2'}\right).$$
Let us now make a choice of representatives modulo $2\pi$ for such phases $\widehat{\Theta}_1,\widehat{\Theta}_2,\widehat{\Theta}_1',\widehat{\Theta}_2'\in\mathbb{R}$ with
\begin{align*}
\widehat{\Theta}_1-\widehat{\Theta}_2&\in (-\pi,\pi], & \widehat{\Theta}_1'-\widehat{\Theta}_2'&\in (-\pi,\pi],\\
\widehat{\Theta}_1-\widehat{\Theta}_1'&\in [-2\pi,2\pi], & \widehat{\Theta}_2-\widehat{\Theta}_2'&\in [-2\pi,2\pi].
\end{align*}
Then, the integrand can be rewritten as follows
\begin{align*}
((-h)(\widehat{\Theta}_1-&\widehat{\Theta}_1')-(-h)(\widehat{\Theta}_2-\widehat{\Theta}_2'))\left((\widehat{\Theta}_1-\widehat{\Theta}_2)-(\widehat{\Theta}_1'-\widehat{\Theta}_2')\right)\\
&=((-h)(\widehat{\Theta}_1-\widehat{\Theta}_1')-(-h)(\widehat{\Theta}_2-\widehat{\Theta}_2'))\left((\widehat{\Theta}_1-\widehat{\Theta}_1')-(\widehat{\Theta}_2-\widehat{\Theta}_2')\right).
\end{align*}
Since all the terms lie in $[-2\pi,2\pi]$, then Lemma \ref{L-split-kernel} yields
\begin{multline*}
\left((-h)(\widehat{\Theta}_1-\widehat{\Theta}_1')-(-h)(\widehat{\Theta}_2-\widehat{\Theta}_2')\right)\left((\widehat{\Theta}_1-\widehat{\Theta}_2)-(\widehat{\Theta}_1'-\widehat{\Theta}_2')\right)\\
\leq L_0\,\vert (\widehat{\Theta}_1-\widehat{\Theta}_1')-(\widehat{\Theta}_2-\widehat{\Theta}_2')\vert^2\leq 2L_0\,(\vert \widehat{\Theta}_1-\widehat{\Theta}_2\vert^2+\vert \widehat{\Theta}_1'-\widehat{\Theta}_2'\vert^2)\\
=2L_0\,\left(\vert \Theta_1-\Theta_2\vert_o^2+\vert \Theta_1'-\Theta_2'\vert_o^2\right).
\end{multline*}
Thus, it becomes apparent that
$$\frac{dI}{dt}\leq 4KL_0\,I,\ t\geq 0.$$
By virtue of Gr\"{o}nwall's lemma, we end up with
$$I(t)\leq I(0)e^{4KL_0\,t},\ t\geq 0.$$
Notice that
$$I(0)=\int_{\mathbb{R}}\left(\int_{\mathbb{T}}\int_{\mathbb{T}}\frac{1}{2}d(z_1,z_2)^2\,d_{(z_1,z_2)}\mu_{0,\Omega}\right)\,d_\Omega g=\frac{1}{2}W_{2,g}(f^1_0,f^2_0)^2,$$
which ends the proof.
\end{proof}

As a clear consequence, we obtain the following uniqueness result for general initial data.

\begin{cor}
Consider $\alpha\in (0,\frac{1}{2})$, $K>0$ and let $f^1,f^2\in \mathcal{AC}_\mathcal{M}\cap \mathcal{T}_\mathcal{M}$ be weak measure-valued solutions to \eqref{E-kuramoto-transport-TxR} with initial data $f^1_0,f^2_0\in \mathcal{P}(\mathbb{T}\times \mathbb{R})$. If $f^1_0=f^2_0$, then
$$f^1_t=f^2_t,\ \mbox{ for every }\ t\geq 0.$$
\end{cor}

\subsection{Stability of the quadratic Wasserstein distance and mean field limit}\label{SS-mean-field}

When the distributions of natural frequencies of both solutions do not agree, then the metric space $(\mathcal{P}_g(\mathbb{T}\times \mathbb{R}),W_{2,g})$ cannot be used. In such general case, we will simply resort on the standard quadratic Wasserstein distance $W_2$ in both variables $(z,\Omega)$. Nevertheless, such approach requires the solutions to lie in $\mathcal{P}_2(\mathbb{T}\times \mathbb{R})$. The next result shows that we only need to require that on the initial datum.

\begin{lem}
Consider $\alpha\in (0,\frac{1}{2})$, $K>0$ and let $f\in \mathcal{AC}_\mathcal{M}\cap \mathcal{T}_\mathcal{M}$ be a weak measure-valued solution to \eqref{E-kuramoto-transport-TxR}. Assume that the distribution $g=(\pi_\Omega)_{\#}f$ of natural frequencies has bounded second order moment, i.e., $\Omega^2g\in \mathcal{M}(\mathbb{R})$, then
$$\sup_{t\geq 0}\Vert \Omega^2f_t\Vert_{\mathcal{M}(\mathbb{T}\times \mathbb{R})}\leq \Vert \Omega^2g\Vert_{\mathcal{M}(\mathbb{R})}<\infty.$$
\end{lem}

\begin{proof}
Consider the scaled cut-off functions $\{\xi_R\}_{R>0}$ in \eqref{E-scaled-cut-off} and the compactly supported test functions $\varphi_R(z,\Omega)=\xi_R(\Omega)\Omega^2$. Since $f$ is a solution in the sense of the flow by Remark \ref{R-types-solutiones}, then we claim
$$\int_{\mathbb{T}\times \mathbb{R}}\varphi_R(z,\Omega)\,d_{(z,\Omega)}f_t=\int_{\mathbb{T}\times \mathbb{R}}\varphi_R(Z_{f}(t;0,z,\Omega),\Omega)\,d_{(z,\Omega)}f_0=\int_{\mathbb{T}\times \mathbb{R}}\xi_R(\Omega)\Omega^2\,d_{(R,\Omega)}f_0=\int_\mathbb{R}\xi_R(\Omega)\Omega^2\,d_\Omega g,$$
for every $R>0$. On the one hand, the right hand side has a limit by virtue of the dominated convergence theorem. Then, Fatou's lemma implies
$$\int_{\mathbb{T}\times \mathbb{R}}\Omega^2\,d_{(z,\Omega)}f_t\leq \liminf_{R\rightarrow +\infty}\int_{\mathbb{T}\times \mathbb{R}}\varphi_R(z,\Omega)\,d_{(z,\Omega)}f_t=\lim_{R\rightarrow +\infty}\int_{\mathbb{R}}\xi_R(\Omega)\Omega^2\,d_\Omega g=\int_{\mathbb{R}}\Omega^2\,d_\Omega g,$$
thus,  proving our assertion.
\end{proof}

\begin{theo}\label{T-growth-RK}
Consider $\alpha\in (0,\frac{1}{2})$, $K>0$ and let $f^1,f^2\in \mathcal{AC}_\mathcal{M}\cap \mathcal{T}_\mathcal{M}$ be weak measured-valued solutions to \eqref{E-kuramoto-transport-TxR} with initial data $f^1_0,f^2_0\in \mathcal{P}_2(\mathbb{T}\times \mathbb{R})$. Then, 
$$W_2(f^1_t,f^2_t)\leq e^{\left(\frac{1}{2}+2KL_0\right)t}W_2(f^1_0,f^2_0),$$
for every $t\geq 0$, where $L_0$ is the one-sided Lipschitz constant of $-h$ in Lemma \ref{L-split-kernel}.
\end{theo}

The proof resembles that in Theorem \ref{T-growth-L(RK)}. However, a full quadratic Wasserstein distance in $\mathbb{T}\times \mathbb{R}$ is used instead, what makes the proof comparable to that in \cite{C-J-L-V}. Again, the one-sided Lipschitz property in Lemma \ref{L-split-kernel} will be the key step in the proof. 

\begin{proof}
Since $f^1_0,f^2_0\in\mathcal{P}_2(\mathbb{T}\times \mathbb{R})$, Brenier's theorem ensures the existence of an optimal transference plan $\mu_0$ joining them both, i.e.,
$$\mu_0\in\mathcal{P}(f^1_0,f^2_0):=\left\{\mu\in \mathcal{P}((\mathbb{T}\times \mathbb{R})\times (\mathbb{T}\times \mathbb{R})):\,(\pi_1)_{\#}\mu=f^1_0\ \mbox{ and }\ (\pi_2)_{\#}\mu=f^2_0\right\},$$
such that
$$W_2(f^1_0,f^2_0)^2=\int_{\mathbb{T}\times \mathbb{R}}\int_{\mathbb{T}\times \mathbb{R}}d((z_1,\Omega_1),(z_2,\Omega_2))^2\,d_{((z_1,\Omega_1),(z_2,\Omega_2))}\mu_0.$$
Now the projection are $\pi_1((z,\Omega),(z',\Omega'))=(z,\Omega)$ and $\pi_2((z,\Omega),(z',\Omega'))=(z',\Omega')$. Again, we can construct a competitor at time $t$ via push-forward, namely,
$$\mu_t:=(X_{f^1}(t;0,\cdot)\otimes X_{f^2}(t;0,\cdot))_{\#}\mu_0\in \mathcal{P}((\mathbb{T}\times \mathbb{R})\times (\mathbb{T}\times \mathbb{R})),$$
where $X_{f^i}(t;0,z,\Omega)=(Z_{f^i}(t;0,z,\Omega),\Omega)$ is the characteristic flow associated with the transport field $\mathcal{V}[f^i]$ according to Theorem \ref{T-well-posedness-characteristic-system}. It is clear that $\mu_t\in \Pi(f^1_t,f^2_t)$ and, consequently,
\begin{align*}
\frac{1}{2}W_2(f^1_t,f^2_t)^2&\leq \int_{\mathbb{T}\times\mathbb{R}}\int_{\mathbb{T}\times\mathbb{R}}\frac{1}{2}d((z_1,\Omega_1),(z_2,\Omega_2))^2\,d_{((z_1,\Omega_1),(z_2,\Omega_2))}\mu_t\\
&=\int_{\mathbb{T}\times \mathbb{R}}\int_{\mathbb{T}\times \mathbb{R}}\frac{1}{2}d(X_{f^1}(t;0,z_1,\Omega_1),X_{f^2}(t;0,z_2,\Omega_2))^2\,d_{((z_1,\Omega_1),(z_2,\Omega_2))}\mu_0=:I(t).
\end{align*}
Fix $(z_1=e^{i\theta_1},\Omega_1),(z_2=e^{i\theta_2},\Omega_2)\in \mathbb{T}\times \mathbb{R}$ and define 
$$\Theta_1(t):=\Theta_{f^1}(t;0,\theta_1,\Omega_1)\ \mbox{ and }\ \Theta_2(t):=\Theta_{f^2}(t;0,\theta_2,\Omega_2),$$
the unique forward-in-time solutions to \eqref{E-characteristic-system-2} in Lemma \ref{L-existence-characteristic-system}. Also, consider the following curves in $\mathbb{T}$
\begin{align*}
Z_1(t)&:=Z_{f^1}(t;0,z_1,\Omega_1)=e^{i\Theta_1(t)},\\
Z_2(t)&:=Z_{f^2}(t;0,z_2,\Omega_2)=e^{i\Theta_2(t)},
\end{align*}
and the associated curves in $\mathbb{T}\times \mathbb{R}$,
\begin{align*}
X_1(t)&:=X_{f^1}(t;0,z_1,\Omega_1)=(Z_1(t),\Omega_1),\\
X_2(t)&:=X_{f^2}(t;0,z_2,\Omega_2)=(Z_2(t),\Omega_2).
\end{align*}
Set a minimizing geodesic $\widehat{\gamma}_t:[0,1]\longrightarrow\mathbb{T}\times \mathbb{R}$ joining $X_1(t)$ to $X_2(t)$, for every fixed $t>0$. Again, the following map
$$t\longmapsto \frac{1}{2}d^2(X_1(t),X_2(t)),$$
is absolutely continuous. Taking one-sided upper Dini directional derivatives given in Appendix \ref{appendix-differentiability-distance} entails
$$\frac{d}{dt}\frac{1}{2}d^2(X_1(t),X_2(t))\leq -\left<\mathcal{V}[f^1_t](X_1(t)),\widehat{\gamma}_t'(0)\right>-\left<\mathcal{V}[f^2_t](X_2(t)),-\widehat{\gamma}_t'(1)\right>.$$
Let us now consider $\theta_{21}(t):=\overline{\Theta_2(t)-\Theta_1(t)}$, the representative of $\Theta_2(t)-\Theta_1(t)$ modulo $2\pi$ that lies in $(-\pi,\pi]$. Again, we distinguish two cases:

$\bullet$ \textit{Case 1:} $\theta_{21}(t)\in (-\pi,\pi)$. In this case, the only minimizing geodesic reads
$$\widehat{\gamma}_t(s)=(\gamma_t(s),\Omega_1+s(\Omega_2-\Omega_1))=(e^{i(\Theta_1(t)+s\theta_{21}(t))},\Omega_1+s(\Omega_2-\Omega_1)),\ s\in [0,1].$$
Then, the inequality reads
$$\frac{d}{dt}\frac{1}{2}d^2(X_1(t),X_2(t))\leq (\mathcal{P}[f^2_t](\Theta_2(t),\Omega_2)-\mathcal{P}[f^1_t](\Theta_1(t),\Omega_1))\theta_{21}(t).$$

$\bullet$ \textit{Case 2:} $\theta_{21}(t)=\pi$. In that second case there are exactly two minimizing geodesics
$$\widehat{\gamma}_{t,\pm}(s)=(\gamma_{t,\pm}(s),\Omega_1+s(\Omega_2-\Omega_1))=(e^{i(\Theta_1(t)\pm \pi s)},\Omega_1+s(\Omega_2-\Omega_1)),\ s\in [0,1].$$
Then, we restate the above inequality as follows
$$\frac{d}{dt}\frac{1}{2}d^2(X_1(t),X_2(t))\leq (\mathcal{P}[f^2_t](\Theta_2(t),\Omega_2)-\mathcal{P}[f^1_t](\Theta_1(t),\Omega_1))(\pm \pi).$$

To sum up, we achieve the following estimate
\begin{multline*}
\frac{d}{dt}\frac{1}{2}d^2(X_{f^2}(t;0,z_1,\Omega_1),X_{f^2}(t;0,z_2,\Omega_2))\\
\leq (\mathcal{P}[f^1_t](\Theta_{f^1}(t;0,\theta_1,\Omega_1),\Omega_1)-\mathcal{P}[f^2_t](\Theta_{f^2}(t;0,\theta_2,\Omega_2),\Omega_2))\overline{\Theta_{f^1}(t;0,\theta_1,\Omega_1)-\Theta_{f^2}(t;0,\theta_1,\Omega_2)},
\end{multline*}
for every $\theta_1,\theta_2,\Omega_1,\Omega_2\in \mathbb{R}$ and almost every $t\geq 0$. Using the dominated convergence theorem implies
\begin{multline}\label{E-7}
\frac{dI}{dt}\leq \int_{(-\pi,\pi]\times \mathbb{R}}\int_{(-\pi,\pi]\times \mathbb{R}} (\mathcal{P}[f^1_t](\Theta_{f^1}(t;0,\theta_1,\Omega_1),\Omega_1)-\mathcal{P}[f^2_t](\Theta_{f^2}(t;0,\theta_2,\Omega_2),\Omega_2))\\
\times\overline{\Theta_{f^1}(t;0,\theta_1,\Omega_1)-\Theta_{f^2}(t;0,\theta_1,\Omega_2)}\,d_{((\theta_1,\Omega_1),(\theta_2,\Omega_2)}\mu_0.
\end{multline}
Also, note that
\begin{align*}
\mathcal{P}[f^i_t](\theta,\Omega)&=\Omega-\int_{(-\pi,\pi]}\int_{\mathbb{R}}h(\theta-\theta')\,d_{(\theta',\Omega')}f^i_t\\
&=\Omega-K\int_{(-\pi,\pi]}\int_\mathbb{R}h(\theta-\Theta_{f^i}(t;0,\theta_i',\Omega_i'))\,d_{(\theta_i',\Omega_i')}f^i_0,
\end{align*}
for $i=1,2$. Since $(\pi_1)_{\#}\mu_0=f^1_0$ and $(\pi_2)_{\#}\mu_0=f^2_0$, then
\begin{align}
\mathcal{P}[f^1_t](\Theta,\Omega)&=\Omega-K\int_{(-\pi,\pi]\times \mathbb{R}}\int_{(-\pi,\pi]\times \mathbb{R}}h(\theta-\Theta_{f^1}(t;0,\theta_1',\Omega_1'))\,d_{((\theta_1',\Omega_1'),(\theta_2',\Omega_2'))}\mu_0,\label{E-8}\\
\mathcal{P}[f^2_t](\Theta,\Omega)&=\Omega-K\int_{(-\pi,\pi]\times \mathbb{R}}\int_{(-\pi,\pi]\times \mathbb{R}}h(\theta-\Theta_{f^2}(t;0,\theta_2',\Omega_2'))\,d_{((\theta_1',\Omega_1'),(\theta_2',\Omega_2'))}\mu_0.\label{E-9}
\end{align}
Putting \eqref{E-8}-\eqref{E-9} into \eqref{E-7} amounts to
\begin{align}\label{E-10}
\begin{split}
\frac{d I}{dt}\leq& \int_{((-\pi,\pi]\times \mathbb{R})^4}(\Omega_1-\Omega_2)\,\overline{\Theta_{f^1}(t;0,\theta_1,\Omega_1)-\Theta_{f^2}(t;0,\theta_2,\Omega_2)}\,d_{((\theta_1,\Omega_1),(\theta_2,\Omega_2))}\mu_0\,d_{((\theta_1',\Omega_1'),(\theta_2',\Omega_2'))}\mu_0\\
-K&\int_{((-\pi,\pi]\times \mathbb{R})^4}(h(\Theta_{f^1}(t;0,\theta_1,\Omega_1)-\Theta_{f^1}(t;0,\theta_1',\Omega_1'))-h(\Theta_{f^2}(t;0,\theta_2,\Omega_2)-\Theta_{f^2}(t;0,\theta_2',\Omega_2')))\\
& \hspace{3cm}\times\overline{\Theta_{f^1}(t;0,\theta_1,\Omega_1)-\Theta_{f^2}(t;0,\theta_2,\Omega_2)}\,d_{((\theta_1,\Omega_1),(\theta_2,\Omega_2))}\mu_0\,d_{((\theta_1',\Omega_1'),(\theta_2',\Omega_2'))}\mu_0.
\end{split}
\end{align}
By virtue of the Young inequality in the first term and an analogue symmetrization trick to that in \eqref{E-6} in the second term, we obtain
\begin{multline*}
\frac{dI}{dt}\leq I(t)-\frac{K}{2}\int_{((-\pi,\pi]\times \mathbb{R})^4}(h(\Theta_{f^1}(t;0,\theta_1,\Omega_1)-\Theta_{f^1}(t;0,\theta_1',\Omega_1'))-h(\Theta_{f^2}(t;0,\theta_2,\Omega_2)-\Theta_{f^2}(t;0,\theta_2',\Omega_2')))\\
\times\left(\overline{\Theta_{f^1}(t;0,\theta_1,\Omega_1)-\Theta_{f^2}(t;0,\theta_2,\Omega_2)}-\overline{\Theta_{f^1}(t;0,\theta_1,\Omega_1)-\Theta_{f^2}(t;0,\theta_2,\Omega_2)}\right)\\
\times\,d_{((\theta_1,\Omega_1),(\theta_2,\Omega_2))}\mu_0\,d_{((\theta_1',\Omega_1'),(\theta_2',\Omega_2'))}\mu_0.
\end{multline*}
Mimicking the idea in Theorem \ref{T-growth-L(RK)} that uses the one-sided Lipschitz property of $-h$ in Lemma \ref{L-split-kernel} implies
$$\frac{dI}{dt}\leq (1+4KL_0)I,\ t\geq 0.$$
Using Gr\"{o}nwall's lemma
$$I(t)\leq I(0)e^{(1+4KL_0)t},\ t\geq 0.$$
Finally, notice that
$$I(0)=\int_{\mathbb{T}\times \mathbb{R}}\int_{\mathbb{T}\times \mathbb{R}}\frac{1}{2}d((z_1,\Omega_1),(z_2,\Omega_2))^2\,d_{((z_1,\Omega_1),(z_2,\Omega_2))}\mu_0=\frac{1}{2}W_2(f^1_0,f^2_0)^2,$$
and that ends the proof.
\end{proof}

The above Theorem \ref{T-growth-RK} along with Theorem \ref{T-empirical-measures-subcritical} implies the rigorous mean field limit as depicted in the following result.

\begin{cor}
Consider $\alpha\in (0,\frac{1}{2})$, $K>0$ and let $f\in \mathcal{AC}_\mathcal{M}\cap \mathcal{T}_\mathcal{M}$ be the unique weak measure-valued solution in the sense of the flow to \eqref{E-kuramoto-transport-TxR} with initial datum $f_0\in \mathcal{P}_2(\mathbb{T}\times \mathbb{R})$. Consider $N$ oscillators with phases and natural frequencies given by the configurations
$$\Theta_0^N=(\theta_{1,0}^N,\ldots,\theta_{N,0}^N)\ \mbox{ and }\ \{\Omega_i^N:\,i=1,\ldots,N\},$$ 
for every $N\in\mathbb{N}$. Let $\Theta^N(t):=(\theta_1^N(t),\ldots,\theta_N^N(t))$ be the unique global-in-time classical solution to the discrete singular Kuramoto model according to \cite[Theorem 3.1]{P-P-S} and define the associated empirical measures in $\mathbb{T}\times \mathbb{R}$
$$\mu^N_t:=\frac{1}{N}\sum_{i=1}^N\delta_{z_i^N(t)}(z)\otimes \delta_{\Omega_i^N}(\Omega),$$
where $z_i^N(t):=e^{i\theta_i^N(t)}$. If $\lim_{N\rightarrow\infty}W_2(\mu^N_0,f_0)=0$, then,
$$\lim_{N\rightarrow\infty}\sup_{t\in [0,T]}W_2(\mu^N_t,f_t)=0,\ \mbox{ for all }\ T>0.$$
\end{cor}

\section{Global phase-synchronization of identical oscillators in finite time}\label{S-synchro}

In this section we will focus on the analysis of the dynamics and emergent behavior of solutions to the macroscopic model \eqref{E-kuramoto-transport-TxR} in the subcritical regime of the parameter, i.e., $\alpha\in (0,\frac{1}{2})$. To such end, let us start by recalling the corresponding result for the discrete Kuramoto model with singular weights \eqref{E-kuramoto-discrete}-\eqref{E-kuramoto-kernel} (see \cite[Theorem 5.1]{P-P-S}) 

\begin{theo}\label{T-synchro-discrete-subcritical}
Let $\Theta=(\theta_1, \cdots, \theta_N)$ be the solution to \eqref{E-kuramoto-discrete}-\eqref{E-kuramoto-kernel} with $\alpha\in \left(0,\frac{1}{2}\right)$ for identical oscillators $(\Omega_i = 0)$, for $i=1, \ldots, N$. Assume that the initial configuration $\Theta_0$ is confined in a half circle, i.e., $0<D(\Theta_0) < \pi$. Then, there is complete phase synchronization at a finite time not larger than $T_c$ where
$$
T_c=\frac{D(\Theta_0)^{1-2\alpha}}{2\alpha Kh(D(\Theta_0))}.
$$
\end{theo}

\begin{rem}\label{R-synchro-discrete-subcritical}
Let us consider the average phase of the $N$ oscillator phases when regarded in $\mathbb{R}^N$
$$
\theta_{av}(t):=\frac{1}{N}\sum_{i=1}^N\theta_i(t),\ t\geq 0.
$$
We will sometimes regard it in the unit circle $z_{av}(t)=e^{i\theta_{av}(t)}$. Notice that it is a conserved quantity of the system, i.e., $\theta_{av}(t)=\theta_{av}(0)$ for all $t\geq 0$, when all the natural frequencies $\Omega_i$ vanish. Notice that for some $T_c$ independent on $N$, Theorem \ref{T-synchro-discrete-subcritical} shows that
$$
\theta_i(t)=\theta_{av}(0),\ t\geq T_c.
$$
Hence, one may think that it remains true for measure-valued solutions to the macroscopic system \eqref{E-kuramoto-transport-TxR}.
\end{rem}

Before introducing the main result of this section, we will need to define an analogue concept of average phase at the continuum level for initial configurations whose support is confined to the half circle. For the sake of clarity, capital letters will be used for the corresponding macroscopic quantities.

\begin{defi}\label{D-average-phase}
Consider $f_0\in \mathcal{P}(\mathbb{T}\times \mathbb{R})$. Assume that $D_0=\diam(\supp \rho_0)<\pi$ and let us set
$$C_0=\{e^{i\theta}:\,\theta\in [\theta_*(0),\theta^*(0)]\},$$ 
with $D_0=\theta^*(0)-\theta_*(0)$ to be the geodesic convex hull of $\supp \rho_0$, that is, the smallest geodesically convex set in $\mathbb{T}$ containing $\supp\rho_0$. We will define the average phase of the initial configuration by
$$\Theta_{av}(0):=\int_{[\theta_*(0),\theta^*(0)]}\int_{\mathbb{R}}\theta \,d_{(\theta,\Omega)}f_0=\int_{[\theta_*(0),\theta^*(0)]}\theta\,d_{\theta}\rho_0.$$
Also, we will regard it in the unit torus as follows $Z_{av}(0):=e^{i\Theta_{av}(0)}$.
\end{defi}

Notice that the above $C_0$ exists and is unique by virtue of the assumption $D_0<\pi$. When $D_0\geq\pi$ there might be two (or any) such smallest geodesically convex sets. Moreover, $\Theta_{av}(0)$ is not uniquely defined, but it depends on the choice of representatives that we make for $\theta_*(0)$ and $\theta^*(0)$. However, $Z_{av}(0)$ is uniquely defined since all those representatives agree modulo $2\pi$. For the reader convenience, let us introduce an alternative representation of $Z_{av}(0)$.

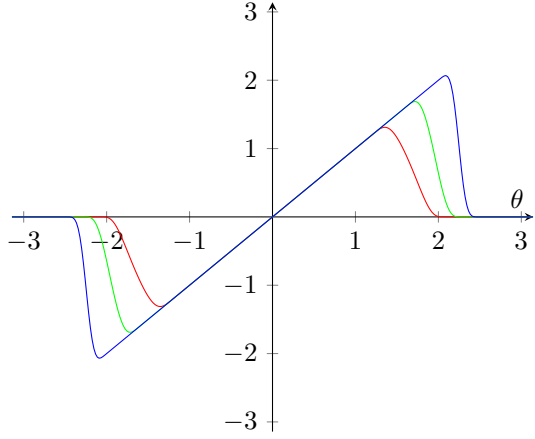
\begin{figure}
\centering
\begin{tikzpicture}[
declare function={
g(\r,\d,\D)= (\r<=\d) * 1 + and(\r>\d, \r<=\D) * 1/(1+exp((2*\r-\d-\D)/((\D-\r)*(\r-\d))))+ (\x>\D) * 0;
fa(\r)=\r*g(abs(\r),pi-2,pi-1);
fb(\r)=\r*g(abs(\r),pi-2*0.8,pi-0.8);
fc(\r)=\r*g(abs(\r),pi-2*0.6,pi-0.6);
}
]
\begin{axis}[
  axis x line=middle, axis y line=middle,
  xmin=-pi, xmax=pi, xtick={-3,-2,-1,0,1,2,3}, xlabel=$\theta$,
  ymin=-pi, ymax=pi, ytick={-3,-2,-1,0,1,2,3},
  legend style={at={(0.79,0.95)}},
]
\addplot [
    domain=-pi:pi, 
    samples=300, 
    color=red,
]
{fa(x)};

\addplot [
    domain=-pi:pi, 
    samples=300, 
    color=green,
]
{fb(x)};

\addplot [
    domain=-pi:pi, 
    samples=300, 
    color=blue,
]
{fc(x)};
\end{axis}
\end{tikzpicture}
\caption{Plot of $\vartheta_{\theta_0,\varepsilon}$ with $\theta_0=-\pi$ and values $\varepsilon=1$ (blue), $\varepsilon=0.8$ (green) and $\varepsilon=0.6$ (red) respectively.}
\label{fig:periodified-phase}
\end{figure}

\begin{lem}[Periodified phase]
Let us consider the following cut-off function
$$\xi_{\delta_1,\delta_2}(r):=\left\{
\begin{array}{ll}
1, & r\in [0,\delta_1),\\
\frac{1}{1+\exp\left(\frac{2r-(\delta_1+\delta_2)}{(\delta_2-r)(r-\delta_1)}\right)}, & r\in [\delta_1,\delta_2),\\
0, & r\in [\delta_2,+\infty),
\end{array}\right.$$
for any $0<\delta_1<\delta_2$. Fix any $\varepsilon>0$ and define (see Figure \ref{fig:periodified-phase})
$$
\vartheta_{\theta_0,\varepsilon}(\theta):=\theta\,\xi_{\pi-2\varepsilon,\pi-\varepsilon}(\vert \theta-\theta_0-\pi\vert),\ \theta\in [\theta_0,\theta_0+2\pi].
$$
Then, $\vartheta_{\theta_0,\varepsilon}(\theta)$ is a periodic function of class $C^\infty$ such that $\vert\vartheta_\varepsilon(\theta)\vert \leq \vert \theta\vert$ and
$$\vartheta_{\theta_0,\varepsilon}(\theta)=\left\{\begin{array}{ll}
\theta, & \vert \theta-\theta_0-\pi\vert \leq \pi-2\varepsilon,\\
0, & \pi-\varepsilon\leq \vert \theta-\theta_0-\pi\vert\leq \pi.
\end{array}
\right.$$
\end{lem}

\begin{lem}\label{L-average-phase-characterization}
Consider $f_0\in \mathcal{P}(\mathbb{T}\times \mathbb{R})$. Assume that $D_0=\diam(\supp \rho_0)<\pi$ and fix any point away from the support $z_0=e^{i\theta_0}\in \mathbb{T}\setminus \supp\rho_0$. Then, there exists $\varepsilon_0>0$ such that
$$\Theta_{av}(0)=\int_{(-\pi,\pi]}\vartheta_{\theta_0,\varepsilon}(\theta)\,d_\theta\rho_0=\int_{\mathbb{T}}\vartheta_{\theta_0,\varepsilon}(z)\,d_z\rho_0,\ \forall\,\varepsilon\in (0,\varepsilon_0).$$
\end{lem}

\begin{theo}\label{T-convergence-equilibrium-continuum-subcritical}
Consider any initial datum $f_0\in \mathcal{P}(\mathbb{T}\times \mathbb{R})$ with identical distribution of natural frequencies, namely, $g=(\pi_\Omega)_{\#}f_0=\delta_0(\Omega)$, where $\pi_\Omega$ is the projection \eqref{E-projections}. Let $f=f_t$ be the unique global-in-time weak measure-valued solution to \eqref{E-kuramoto-transport-TxR} issued at $f_0$ with $\alpha\in (0,\frac{1}{2})$ and assume $D_0:=\diam (\supp\rho_0)<\pi$. Then,
$$f_t=f_\infty\ \mbox{ for all }\ t\geq T_c,$$
where $T_c=\frac{D_0^{1-2\alpha}}{2\alpha Kh(D_0)}$ and the equilibrium $f_\infty$ is given by the monopole $f_\infty:=\delta_{Z_{av}(0)}(z)\otimes \delta_0(\Omega)$.
\end{theo}

\begin{proof}
Consider the smallest geodesically convex subset $C_0$ of $\mathbb{T}$ as given in Definition \ref{D-average-phase}. Let us recall that in the existence Theorem \ref{T-mean-field-existence-subcritical}, the weak measure-valued solution $f$ arose as the limit
\begin{equation}\label{E-62}
\mu^N\rightarrow f\ \mbox{ in }\ \mathcal{P}_1(\mathbb{T}\times \mathbb{R})-W_1,
\end{equation}
of some sequence of empirical measures
$$\mu^N_t=\frac{1}{N}\sum_{i=1}^N\delta_{z_i^N(t)}(z)\delta_{\Omega_i^N}(\Omega),$$
with $z_i^N(t)=e^{i\theta_i^N(t)}$. Recall that $\Theta^N(t)=(\theta_1^N(t),\ldots,\theta_N^N(t))$ are global classical solutions to \eqref{E-kuramoto-discrete}-\eqref{E-kuramoto-kernel}.  Without loss of generality, Lemma \ref{L-discrete-approximation} allows assuming that $\Omega_i^N=0$ for all $i=1,\ldots,N$ and $N\in\mathbb{N}$. Also, the same result allows ensuring that $\mu^N_0$ can be taken so that $\supp\mu^N_0\subseteq C_0^\delta$, where
$$C_0^\delta:=\left\{e^{i\theta}:\,\theta\in \left[\theta_*(0)-\frac{\delta}{2},\theta^*(0)+\frac{\delta}{2}\right]\right\},$$
and $\delta$ is any arbitrary value $\delta\in (0,\pi-D_0)$. By virtue of Theorem \ref{T-synchro-discrete-subcritical} and Remark \ref{R-synchro-discrete-subcritical}, we obtain that
$$\Theta_i^N(t)=z_{av}^N(0),\ t\geq T^\delta_c,$$
for every $i=1,\ldots,N$, every $N\in \mathbb{N}$ and $T^\delta_c$ is given by
$$T_c^\delta=\frac{(D_0+\delta)^{1-2\alpha}}{2\alpha Kh(D_0+\delta)}.$$
Hence, we obtain that
\begin{equation}\label{E-63}
\mu^N_t=\delta_{z_{av}^N(0)}(z)\delta_0(\Omega),\ t\geq T_c^\delta.
\end{equation}
Consider $z_0\in \mathbb{T}\setminus\supp \rho_0$ and $\varepsilon_0>0$ as given in Lemma \ref{L-average-phase-characterization}. Then, it is clear that
\begin{equation}\label{E-64}
\Theta_{av}(0)=\int_{\mathbb{T}\times\mathbb{R}}\vartheta_{\theta_0,\varepsilon}(z)\,d_{(z,\Omega)}f_0=\lim_{N\rightarrow\infty}\int_{\mathbb{T}\times \mathbb{R}}\vartheta_{\theta_0,\varepsilon}(z)\,d_{(z,\Omega)}\mu^N_0=\lim_{N\rightarrow\infty}\theta_{av}^N(0).
\end{equation}
Taking limits in \eqref{E-63} and using \eqref{E-62} and \eqref{E-64} we obtain that
$$f_t=\delta_{Z_{av}(0)}(z)\delta_0(\Omega),\ t\geq T_c^\delta.$$
Since $\delta\in (0,\pi-D_0)$, the result holds true.
\end{proof}

\section{The critical regime}\label{S-critical-regime}
This section is devoted to the critical case, $\alpha=\frac{1}{2}$. Notice that most of the ideas in Sections \ref{S-macroscopic-model}, \ref{S-weak-solutions-existence}, \ref{S-uniqueness-mean-field} and \ref{S-synchro} break down due to the jump discontinuity of the kernel $h$. However, an appropriate concept of measure-valued solution can be considered, thus yielding well-posedness of \eqref{E-kuramoto-transport-TxR} like in Theorem \ref{T-mean-field-existence-subcritical}, analogue Dobrushin-type estimates in Wasserstein distance like in Theorems \ref{T-growth-L(RK)} and \ref{T-growth-RK}, rigorous mean field limit and a similar analysis of the dynamics like in Theorem \ref{T-convergence-equilibrium-continuum-subcritical}.

\subsection{Solutions in the sense of the Filippov flow}
In this section, we will extend the above existence results in Section \ref{S-weak-solutions-existence} to the critical regime. Note that the vector field $\mathcal{V}[\mu]$ in Definition \ref{D-transport-field} does not make sense when $\alpha=\frac{1}{2}$ for general finite Radon measures $\mu\in \mathcal{M}(\mathbb{T}\times \mathbb{R})$. The reason relies on the fact that in such case the kernel $h$ is no longer continuous, but it exhibits a jump discontinuity at $\theta\in 2\pi\mathbb{Z}$. Consequently, Definition \ref{D-transport-field} is wrong due to the possible presence of atoms in the Radon measure $\mu$. Hence, the above regularity theory for $\mathcal{V}[\mu]$ in Section \ref{S-macroscopic-model} fails when $\alpha=\frac{1}{2}$ and it requires to be adapted.

This part will be split in several steps. First, we will adapt Definition \ref{D-transport-field} so that the transport field also makes sense for $\alpha= \frac{1}{2}$. Second, although the classical theory cannot solve the characteristic system due to the discontinuity of the transport field, we can obtain a generalized concept of solutions to the characteristic system. Namely, we will consider Filippov solutions of the characteristic system and will show that they exist, are global and unique forward-in-time. Third, we will show that the unique forward-in-time solutions to the discrete model \eqref{E-kuramoto-discrete} that were obtained in \cite{P-P-S} give rise to solutions in the sense of the Filippov flow to the macroscopic system \eqref{E-kuramoto-transport-TxR} by considering the corresponding empirical measures like in Definition \ref{D-empirical-measures-subcritical}. To end this part, we will derive appropriate a priori bounds on the sequence of empirical measures so that we can pass to the limit and obtain a solutions in the sense of the Filippov flow to the macroscopic system \eqref{E-kuramoto-transport-TxR} issued at any general initial datum.

\begin{defi}\label{D-transport-field-critical}
Consider $\alpha=\frac{1}{2}$ and $K>0$. We will formally define the function $\mathcal{P}[\mu]$ and the tangent vector field $\mathcal{V}[\mu]$ along the manifold $\mathbb{T}\times \mathbb{R}$ given by
\begin{align*}
\mathcal{P}[\mu](\theta,\Omega)&:=\Omega-K\int_{\mathbb{T}\setminus \{e^{i\theta}\}}\int_{\mathbb{R}}h(\theta-\theta')\,d_{(\theta',\Omega')}\mu,\\
\mathcal{V}[\mu](z,\Omega)&:=(\mathcal{P}[\mu](z,\Omega)\,iz,0),
\end{align*}
where $\mu\in \mathcal{M}(\mathbb{T}\times \mathbb{R})$ is any finite Radon measure.
\end{defi}

Notice that in this definition for the critical case, we have removed $z'=e^{i\theta}$. It is consistent with that of the subcritical case in Definition \ref{D-transport-field} because, in such case, $h(0)=0$. Moreover, since the function
$$z'\in \mathbb{T}\setminus\{z\}\longmapsto h(z\bar z'),$$
is continuous, then the above integral also makes sense in the critical case.

\subsubsection{The Filippov set-valued tangent vector field}
In this part we will introduce the concept of \textit{Filippov set-valued tangent vector field} associated with a given single-valued tangent vector field. It will be required in order to define the appropriate concept of generalied characteristic trajectories of \eqref{E-kuramoto-transport-TxR}. The definition is a standard translation of the standard concept for the Euclidean case \cite{Au-Ce,F} to the Riemannian case via local charts. We will first give a general definition that is valid for any complete Riemannian manifold and extends the classical \textit{Filippov convexification} or \textit{Filippov set-valued map} of a given measurable vector field in the Euclidean space \cite{Au-Ce,F}. Our second step will be to identify a clear characterization in our particular Riemannian manifold $\mathbb{T}\times \mathbb{R}$.

\begin{defi}\label{D-Filippov-field}
Let $(M,\left<\cdot,\cdot\right>)$ be a $d$-dimensional Riemannian manifold. Fix $x_0\in M$, and set a local chart $(\mathcal{U},\varphi)$ around $x_0$ and any local frame $U_1,\ldots,U_d\in \mathfrak{X}C^\infty(\mathcal{U})$ on $\mathcal{U}$. For any locally essentially bounded tangent vector field $V\in \mathfrak{X}L^\infty_{loc}(M)$, let us consider its local expression
$$V=f_1U_1+\cdots f_d U_d\ \mbox{ in }\ \mathcal{U},$$
where $f_1,\ldots,f_d\in L^\infty_{loc}(\mathcal{U})$. Define $\mathcal{D}:=\varphi(\mathcal{U})\subseteq \mathbb{R}^d$, $F_i:=f_i\circ \varphi^{-1}$ and the vector field
$$F=(F_1,\ldots,F_d):\,\mathcal{D}\subseteq \mathbb{R}^d\longrightarrow\mathbb{R}^d.$$
Then, we define the Filippov set-valued tangent vector field at $x_0$ by the formula
$$\mathcal{K}[V]_{x_0}=\left\{\sum_{i=1}^d\beta_i U_i(x_0):\,\beta=(\beta_1,\ldots,\beta_d)\in \mathcal{K}[F](\varphi(x_0))\right\},$$
where $\mathcal{K}[F](\varphi(x_0))$ is the Filippov set-valued map associated with $F$ (see \cite{Au-Ce,F,P-P-S}), that is
$$\mathcal{K}[F](\varphi(x_0)):=\bigcap_{\delta>0}\bigcap_{\vert \mathcal{N}\vert=0}\overline{\co}(F(B_\delta(x_0)\setminus \mathcal{N})).$$
\end{defi}

\begin{rem}\label{R-Filippov-tangent-field-1}
Notice that given any local chart $(\mathcal{U},\varphi)$ around $x_0$ with coordinates $\varphi=(x_1,\ldots,x_d)$, we can construct the associated local frame whose basis of tangent vector fields read as follows
$$U_i:=\frac{\partial}{\partial x_i},\ \mbox{ for all }\ i=1,\ldots,d.$$
In such particular local frame we can set the local coordinates of the vector field to be
$$V=f_1\frac{\partial}{\partial x_1}+\cdots f_d\frac{\partial}{\partial x_d}\ \mbox{ in }\ \mathcal{U}.$$
The associated representation of the Filippov set-valued tangent field then reads
$$\mathcal{K}[V]_{x_0}=\left\{\sum_{i=1}^d\beta_i\left.\frac{\partial}{\partial x_i}\right\vert_{x_0}:\,\beta=(\beta_1,\ldots,\beta_d)\in \mathcal{K}[F](\varphi(x_0))\right\}.$$
Our interest in considering general local frames will be clarified later in Lemma \ref{L-one-sided-Lipschitz-multivalued}.
\end{rem}

Let us show that the hypothesis $V\in \mathfrak{X}L^\infty_{loc}(M)$ guarantee the independence of the above definition on the given local charts and frames. We will require the following technical result.

\begin{lem}\label{L-Filippov-chain-rule}
Consider a couple of open set $\mathcal{D}_1,\mathcal{D}_2\subseteq\mathbb{R}^d$, a $C^1$ diffeomorphism $\Phi:\mathcal{D}_1\longrightarrow \mathcal{D}_2$ and a measurable vector field $F:\mathcal{D}_2\longrightarrow \mathbb{R}^d$. Then, the following formula takes place
$$
\mathcal{K}[F\circ \Phi](x)=\mathcal{K}[F](\Phi(x)),
$$
for every $x\in \mathcal{D}_1$.
\end{lem}

\begin{proof}
Since $\Phi$ is a diffeomorphism, then one inclusion follows from the other one when applied to $\Phi^{-1}$. Then, we will just focus on one inclusion; specifically: $\mathcal{K}[F\circ\Phi]\subseteq \mathcal{K}[F]\circ \Phi$. Consider $x\in \mathcal{D}_1$, set any arbitrary $\delta>0$ and any negligible set $\mathcal{N}\subseteq \mathcal{D}_2$. Since $\Phi$ is a diffeomorphism there exists $\delta'>0$ and another negligible set $\mathcal{N}'\subseteq \mathcal{D}_1$ so that
$$\Phi(B_{\delta'}(x)\setminus \mathcal{N'})\subseteq B_\delta(\Phi(x))\setminus \mathcal{N}.$$
Consequently,
$$\mathcal{K}[F\circ \Phi](x)\subseteq \overline{\co}((F\circ \Phi)(B_{\delta'}(x)\setminus \mathcal{N}'))\subseteq \overline{\co}(F(B_\delta(\Phi(x))\setminus \mathcal{N})).$$
Since $\delta$ and $\mathcal{N}$ are arbitrary, then the desired inclusion holds true.
\end{proof}

\begin{lem}\label{L-Filippov-matrix-rule}
Consider an open set $\mathcal{D}\subseteq \mathbb{R}^d$, a continuous map $A:\mathcal{D}\longrightarrow \mathcal{M}_d(\mathbb{R})$ and a locally essentially bounded vector field $F:\mathcal{D}\longrightarrow\mathbb{R}^d$. Then, the following formula takes place
$$
\mathcal{K}[AF](x)=A(x)\mathcal{K}[F][x],
$$
for every $x\in \mathcal{D}$.
\end{lem}

\begin{proof}
Fix $x_0\in \mathcal{D}$. The proof is split into two parts that correspond to  the inclusions $A(x_0)\mathcal{K}[F](x_0)$ $\subseteq  \mathcal{K}[AF](x_0)$ and $\mathcal{K}[AF](x_0)\subseteq A(x_0)\mathcal{K}[F](x_0)$.

$\bullet$ \textit{Step 1}: Consider $z_0\in A(x_0)\mathcal{K}[F](x_0)$ and set $y_0\in \mathcal{K}[F](x_0)$ so that $z_0=A(x_0)y_0$. Fix any arbitrary $\delta>0$ and any negligible set $\mathcal{N}\subseteq \mathcal{D}$. Our goal here is to show that $z_0\in \overline{\co}(AF)(B_\delta(x_0)\setminus \mathcal{N})$. Consider any small enough $\delta'>0$ so that $\delta'<\delta$ and $F\in L^\infty(B_{\delta'}(x_0),\mathbb{R}^d)$. Take any arbitrary $n\in \mathbb{N}$ and note that $y_0\in \overline{\co}(F(B_{\delta'/n}(x_0)\setminus \mathcal{N}))$. Then, there exists $y_n\in \co(F(B_{\delta'/n}(x_0)\setminus \mathcal{N}))$ so that $\vert y_n-y_0\vert\leq \frac{1}{n}$. In addition, Caratheodory's theorem of convex analysis provides
$$q^1_n,\ldots,q^{d+1}_n\in B_{\delta'/n}(x_0)\setminus \mathcal{N}\ \mbox{ and }\ \lambda^{1}_n,\ldots,\lambda^{d+1}_n\in [0,1],$$
so that $\sum_{k=1}^{d+1}\lambda^{k}_n=1$ and $y_n$ can be written as follows
$$y_n=\sum_{k=1}^{d+1}\lambda^{k}_nF(q^k_n),$$
for every $n\in \mathbb{N}$. Let us define the following vectors
$$z_n:=\sum_{k=1}^{d+1}\lambda^k_n A(q^k_n)F(q^k_n),$$
for every $n\in \mathbb{N}$. By definition, $z_n\in \co((AF)(B_{\delta}(x_0)\setminus \mathcal{N}))$. Let us see that $z_n\rightarrow z_0$ and that ends the proof. To that end, let us split as follows
\begin{align*}
z_n=\sum_{k=1}^{d+1}\lambda^k_n (A(q^k_n)-A(x_0))F(q^k_n)+\sum_{k=1}^{d+1}\lambda^k_n A(x_0)F(q^k_n)
=\sum_{k=1}^{d+1}\lambda^k_n (A(q^k_n)-A(x_0))F(q^k_n)+A(x_0)y_n.
\end{align*}
On the one hand, the second term converges towards $A(x_0)y_0=z_0$. On the other hand, let us note that the first one vanishes as $n\rightarrow \infty$. Indeed, taking norms, it can be bounded above by
$$\sum_{k=1}^{d+1}\lambda^k_n \vert A(q^k_n)-A(x_0)\vert \Vert F\Vert_{L^\infty(B_{\delta'}(x_0),\mathbb{R}^d)}.$$
Now the convergence to zero becomes apparent since the coefficients $\lambda^k_n$ sum up to $1$ and $A(q^k_n)\rightarrow A(x_0)$ for every $k=1,\ldots,d+1$ by virtue of the continuity of $A$ along with the convergence $q^k_n\rightarrow x_0$ as $n\rightarrow\infty$.

$\bullet$ \textit{Step 2}: Now, consider $z_0\in \mathcal{K}[AF](x_0)$, fix any arbitrary $\delta>0$ and any negligible set $\mathcal{N}\subseteq \mathcal{D}$. Our goal is to show that there exists $y_0\in \overline{\co}(F(B_\delta(x_0)\setminus \mathcal{N}))$ so that $z_0=A(x_0)y_0$. Consider a small enough $\delta'>0$ so that $\delta'<\delta$ and $F\in L^\infty(B_{\delta'}(x_0),\mathbb{R}^d)$. Take any arbitrary $n\in \mathbb{N}$ and note that $z_0\in\overline{\co}((AF)(B_{\delta'/n}(x_0)\setminus \mathcal{N}))$. Then, there exists $z_n\in \co((AF)(B_{\delta'/n}(x_0)\setminus \mathcal{N}))$ so that $\vert z_n-z_0\vert\leq \frac{1}{n}$. Again, Caratheodory's theorem implies the existence of
$$q^1_n,\ldots,q^{d+1}_n\in B_{\delta'/n}(x_0)\setminus \mathcal{N}\ \mbox{ and }\ \lambda^{1}_n,\ldots,\lambda^{d+1}_n\in [0,1],$$
so that $\sum_{k=1}^{d+1}\lambda^{k}_n=1$ and $y_n$ can be written as follows
$$z_n=\sum_{k=1}^{d+1}\lambda^{k}_nA(a^k_n)F(q^k_n),$$
for every $n\in \mathbb{N}$. Now, consider the following vectors
$$y_n:=\sum_{k=1}^{d+1}\lambda^{k}_nF(q^k_n),$$
for every $n\in \mathbb{N}$. Notice that $y_n\in \co(F(B_\delta(x_0)\setminus \mathcal{N}))$. By the boundedness of $F$ in $B_{\delta'}(x_0)$ it is clear that $\{y_n\}_{n\in \mathbb{N}}$ is a bounded sequence. Hence, Weierstrass theorem, provides a subsequence $y_{\sigma(n)}$ that converges towards some $y_0\in \overline{\co}(F(B_\delta(x_0)\setminus \mathcal{N}))$. In addition, let us note that the following split holds true
\begin{align*}
A(x_0)y_{\sigma(n)}&=\sum_{k=1}^{d+1}\lambda^k_{\sigma(n)}A(x_0)F(q^k_{\sigma(n)})\\
&=\sum_{k=1}^{d+1}\lambda^k_{\sigma(n)}(A(x_0)-A(q^k_{\sigma(n)}))F(q^k_{\sigma(n)})+\sum_{k=1}^{d+1}\lambda^k_{\sigma(n)}A(q^k_{\sigma(n)})F(q^k_{\sigma(n)})\\
&=\sum_{k=1}^{d+1}\lambda^k_{\sigma(n)}(A(x_0)-A(q^k_{\sigma(n)}))F(q^k_{\sigma(n)})+z_{\sigma(n)},
\end{align*}
for every $n\in \mathbb{N}$. Taking limits as $n\rightarrow\infty$ and using the boundedness of $F$ in $B_{\delta'}(x_0)$ and the continuity of $A$ show that $A(x_0)y_0=z_0$ and that ends the proof.
\end{proof}

\begin{rem}
As a simple consequence of the above Lemmas \ref{L-Filippov-chain-rule} and \ref{L-Filippov-matrix-rule} we will show that the above Definition \ref{D-Filippov-field} does not depend on the chosen local frames. To that end, consider a couple of local frames $(\mathcal{U},\varphi,U_1,\ldots,U_d)$ and $(\mathcal{W},\psi,W_1,\ldots,W_d)$. Write $V$ in both basis
$$V=\sum_{i=1}^d f_iU_i=\sum_{j=1}^dg_jW_j.$$
Let us set the matrix $P$ of change from basis $\{W_1,\ldots,W_d\}$ to $\{U_1,\ldots,U_d\}$, i.e., the continuous map $P:\mathcal{U}\cap \mathcal{W}\longrightarrow \mathcal{M}_d(\mathbb{R})$ whose coordinates $P=(p_{ij})_{1\leq i,j\leq d}$ fulfil
$$W_j=\sum_{i=1}^d p_{ij}U_i\ \mbox{ in }\ \mathcal{U}\cap \mathcal{W},$$
for all $j=1,\ldots,d$. Then, the change of basis formula provide the relations
\begin{equation}\label{E-35}
f_i=\sum_{j=1}^d p_{ij}g_j\ \mbox{ in }\ \mathcal{U}\cap \mathcal{W},
\end{equation}
for all $i=1,\ldots,d$. Consider $F_i=f_i\circ \varphi^{-1}$ and $G_i=g_j\circ \psi^{-1}$ and the associated Filippov maps
\begin{align*}
\mathcal{K}[F](\varphi(x_0))&=\bigcap_{\delta>0}\bigcap_{\vert \mathcal{N}\vert=0}\overline{\co}(F(B_\delta(\varphi(x_0))\setminus \mathcal{N})),\\
\mathcal{K}[G](\psi(x_0))&=\bigcap_{\delta>0}\bigcap_{\vert \mathcal{N}\vert=0}\overline{\co}(G(B_\delta(\psi(x_0))\setminus \mathcal{N})).
\end{align*}
Consider the associated Filippov set-valued tangent fields, that we denote upper scripts for distinction
\begin{align*}
\mathcal{K}[V]_{x_0}^\varphi=\left\{\sum_{i=1}^d\beta_i U_i(x_0):\,\beta=(\beta_1,\ldots,\beta_d)\in \mathcal{K}[F](\varphi(x_0))\right\},\\
\mathcal{K}[V]_{x_0}^\psi=\left\{\sum_{j=1}^d\gamma_j V_j(x_0):\,\gamma=(\gamma_1,\ldots,\gamma_d)\in \mathcal{K}[G](\psi(x_0))\right\}.
\end{align*}
Our goal is to show that $\mathcal{K}[V]_{x_0}^\varphi=\mathcal{K}[V]_{x_0}^\psi$. On the one hand, recall that \eqref{E-35} amounts to
\begin{equation}\label{E-38}
F\circ(\varphi\circ\psi^{-1})=(P\circ\psi^{-1})\,G\ \mbox{ in }\ \psi(\mathcal{U}\cap \mathcal{W}),
\end{equation}
Taking $\mathcal{D}_1=\psi(\mathcal{U}\cap \mathcal{V}), \mathcal{D}_2=\varphi(\mathcal{U}\cap \mathcal{V})$ along with $\Phi=\varphi\circ \psi^{-1}$ and $A=P\circ\psi^{-1}$ we can apply Lemmas \ref{L-Filippov-chain-rule} and \ref{L-Filippov-matrix-rule} to \eqref{E-38} and we infer that
$$
\mathcal{K}[F](\varphi(x_0))=P(x_0)\,\mathcal{K}[G](\psi(x_0)).
$$
Then, the claimed independence on the chosen local frame follows.
\end{rem}

We are interested in introducing the Filippov set-valued tangent vector in $\mathbb{T}\times \mathbb{R}$ associated with the transport fields in Definition \ref{D-transport-field-critical} for the critical regime $\alpha=\frac{1}{2}$. Notice that it can be done since $\mathcal{V}[\mu]$ is locally essentially bounded in $\mathbb{T}\times\mathbb{R}$ for any given $\mu \in \mathcal{M}(\mathbb{T}\times \mathbb{R})$.

\begin{pro}\label{P-Filippov-transport-field-critical}
Consider $\alpha=\frac{1}{2}$, $K>0$ and any $\mu\in \mathcal{M}(\mathbb{T}\times \mathbb{R})$. Then, the Filippov set-valued tangent field associated with $\mathcal{V}[\mu]$ reads as follows
$$\mathcal{K}[\mathcal{V}[\mu]](z_0,\Omega_0)=\{(p\,iz_0,0):\,p\in \mathcal{K}[\mathcal{P}[\mu](\cdot,\Omega_0)](\theta_0)\},$$
for every $(z=e^{i\theta_0},\Omega_0)\in \mathbb{T}\times\mathbb{R}$. Here, the set $\mathcal{K}[\mathcal{P}[\mu](\cdot,\Omega_0)](\theta_0)$ is nothing but the standard Filippov set-valued map at $\theta=\theta_0$ associated with the map $\theta\in \mathbb{R}\longmapsto \mathcal{P}[\mu](e^{i\theta},\Omega_0)$ with fixed $\Omega_0$.
\end{pro}

\begin{proof}
Consider any $(z_0=e^{i\theta_0},\Omega_0)$ and take some local chart around it, e.g.,
$$(\theta,\Omega)\in (\theta_0-\pi,\theta_0+\pi)\times \mathbb{R}\longmapsto (e^{i\theta},\Omega).$$
Consider the local frame $\frac{\partial}{\partial\theta}$ and $\frac{\partial}{\partial \Omega}$, see Remark \ref{R-Filippov-tangent-field-1}. Then, we can write the transport field $\mathcal{V}[\mu]$ in local coordinates as follows
$$\mathcal{V}[\mu](z,\Omega)=\mathcal{P}[\mu](z,\Omega)(iz,0)+0\,(0,1),$$
for every $z\in \mathbb{T}\setminus\{\bar z_0\}$ and $\Omega\in \mathbb{R}$. Set the vector field
$$\begin{array}{cccc}
F: & (\theta_0-\pi,\theta_0+\pi)\times \mathbb{R} & \longrightarrow & \mathbb{R}^2,\\
 & (\theta,\Omega) & \longmapsto & (\mathcal{P}[\mu](e^{i\theta},\Omega),0).
\end{array}$$
Then, Definition \ref{D-Filippov-field} shows that
$$\mathcal{K}[\mathcal{V}[\mu]](z_0,\Omega_0)=\{(p\,iz_0,q):\,(p,q)\in \mathcal{K}[F](\theta_0,\Omega_0)\}.$$
Notice that the $q$-component in any vector in $(p,q)\in\mathcal{K}[F](\theta_0,\Omega_0)$ vanishes by virtue of the definition of $F$ and that $F$ is continuous in the second variable. A straightforward computation then shows that
$$\mathcal{K}[F_1](\theta_0,\Omega_0)=\mathcal{K}[F_1(\cdot,\Omega_0)](\theta_0)\equiv \mathcal{K}[\mathcal{P}[\mu](\cdot,\Omega_0)](\theta_0),$$
and that ends the proof.
\end{proof}

\subsubsection{The Filippov flow of the transport field}
In this part we will recover some regularity properties of the Filippov set-valued vector field that extend those in Subsection \ref{SS-transport-field-properties} to the critical case and will be useful throughout the next parts. Our final goal will be to show the well-posedness of a Filippov flow associated with the transport field $\mathcal{V}[\mu]$ for $\alpha=\frac{1}{2}$. 

\begin{lem}\label{L-transport-field-bounded-critical}
Consider $\alpha=\frac{1}{2}$, $K>0$ and $\mu \in L^\infty(0,T;\mathcal{M}(\mathbb{T}\times \mathbb{R}))$. Then,
\begin{align*}
\frac{\mathcal{P}[\mu]}{1+\vert \Omega\vert}&\in L^\infty([0,T]\times \mathbb{T}\times \mathbb{R}),\\
\frac{\mathcal{V}[\mu]}{1+\vert \Omega\vert}&\in L^\infty(0,T;\mathfrak{X}L^\infty(\mathbb{T}\times \mathbb{R}).
\end{align*}
\end{lem}

\begin{proof}
By analogue arguments to those in Corollary \ref{C-transport-field-holder} we show that
$$\esssup\limits_{(t,z,\Omega)\in [0,T]\times\mathbb{T}\times \mathbb{R}}\frac{\left\vert\mathcal{P}[\mu_t](z,\Omega)\right\vert}{1+\vert \Omega\vert}\leq \sup_{(z,\Omega)\in\mathbb{T}\times \mathbb{R}}\frac{\vert \Omega\vert +K\Vert \mu\Vert_{L^\infty(0,T;\mathcal{M}(\mathbb{T}\times \mathbb{R}))}\Vert h\Vert_{L^\infty(\mathbb{T})}}{1+\vert \Omega\vert}<\infty,$$
what proves the claimed results.
\end{proof}

Recall that in Subsection \ref{SS-transport-field-properties} the existence of classical solutions in the subcritical case when $\mu\in \widetilde{\mathcal{C}}_\mathcal{M}$ relied on the following continuity property (see Corollary \ref{C-transport-field-holder})
$$\frac{\mathcal{V}[\mu]}{1+\vert \Omega\vert}\in C([0,T],\mathfrak{X}C_b(\mathbb{T}\times \mathbb{R})).$$
There, time continuity came from the continuity of $\mu$ in $C([0,T],\mathcal{M}(\mathbb{T}\times \mathbb{R})-\mbox{narrow})$ and was lost when $\mu\in \mathcal{C}_\mathcal{M}$ only (see Remark \ref{R-transport-field-holder}). In our setting for $\alpha=\frac{1}{2}$, we expect to lose continuity both in $t$ in $z$ due to the jump discontinuity of $h$. Indeed, note that since the microscopic system shows global phase synchronization in finite time under certain regimes of the natural frequencies (see \cite{P-P-S}), it is possible that Dirac masses emerge and gain mass as times goes on. It justifies that time continuity might also be lost at certain times and the velocity field becomes discontinuous at such phase value. Our next result clarifies to what extend  time continuity can be preserved (see \cite[Lemma A.1]{C-J-L-V} for similar results).

\begin{theo}\label{T-transport-field-convergence-critical}
Consider $\alpha=\frac{1}{2}$, $K>0$ and let $\{\mu_n\}_{n\in \mathbb{N}}$ and $\mu$ be in $\mathcal{P}(\mathbb{T}\times \mathbb{R})$ so that
$$\mu_n\rightarrow \mu\ \mbox{ in }\ \mathcal{P}(\mathbb{T}\times \mathbb{R})-\mbox{narrow}.$$
Then, the following convergence takes place
$$\lim_{n\rightarrow \infty}\sup_{\Omega\in \mathbb{R}}\vert \mathcal{V}[\mu_n](z,\Omega)-\mathcal{V}[\mu](z,\Omega)\vert=0,$$
for each continuity point $z\in \mathbb{T}$ of the marginal measure $(\pi_z)_{\#}\mu$, where $\pi_z$ is the projection \eqref{E-projections}. In particular, it happens a.e. in $\mathbb{T}$.
\end{theo}

\begin{proof}
Let us define $\rho_n:=(\pi_z)_{\#}\mu_n$ and $\rho:=(\pi_z)_{\#}\mu$. Then, notice that 
$$
\left\vert\mathcal{V}[\mu_n](z,\Omega)-\mathcal{V}[\mu](z,\Omega)\right\vert=\left\vert\mathcal{P}[\mu_n](z,\Omega)-\mathcal{P}[\mu](z,\Omega)\right\vert=\left\vert K\int_{(-\pi,\pi]\setminus \{\theta\}}h(\theta-\theta')\,d_{\theta'}(\rho_n-\rho)\right\vert,
$$
for $z=e^{i\theta}\in \mathbb{T}$ such that $\theta\in (-\pi,\pi]$ and $\Omega\in \mathbb{R}$. Define
$$F_n(\theta):=K\int_{(-\pi,\pi]\setminus \{\theta\}}h(\theta-\theta')\,d_{\theta'}(\rho_n-\rho)=F_n^1(\theta)+F_n^2(\theta)+F_n^3(\theta),$$
for every $\theta\in (-\pi,\pi]$ and $n\in \mathbb{N}$, where each term reads
\begin{align*}
F_n^1(\theta)&:=K\int_{(-\pi,\pi]\setminus\{\theta\}}h_\varepsilon(\theta-\theta')\,d_{\theta'}(\rho_n-\rho),\\
F_n^2(\theta)&:=K\int_{(-\pi,\pi]\setminus \{\theta\}}\chi_{\vert \theta'-\theta\vert_o\geq \varepsilon^{1/2}}(h(\theta-\theta')-h_\varepsilon(\theta-\theta'))\,d_{\theta'}(\rho_n-\rho),\\
F_n^3(\theta)&:=K\int_{(-\pi,\pi]\setminus \{\theta\}}\chi_{\vert \theta'-\theta\vert_o< \varepsilon^{1/2}}(h(\theta-\theta')-h_\varepsilon(\theta-\theta'))\,d_{\theta'}(\rho_n-\rho).
\end{align*}
Here $\varepsilon>0$ is any fix but arbitrary parameter. First, notice that $\rho_n\rightarrow\rho$ narrow in $\mathcal{P}(\mathbb{T})$. Since $h_\varepsilon$ is continuous, then,
$$\lim_{n\rightarrow \infty}F_n^1(\theta)=0,$$
for every $\theta\in (-\pi,\pi]$. Second, by a clear application of the mean value theorem, we obtain
\begin{equation}\label{E-52}
\vert h(\theta)-h_\varepsilon(\theta)\vert\leq \frac{1}{2}\frac{\varepsilon}{\vert \theta\vert_o},
\end{equation}
for all $\theta \in \mathbb{R}\setminus 2\pi\mathbb{Z}$. As an application of \eqref{E-52} we obtain the upper bound
\begin{equation}\label{E-53}
\limsup_{n\rightarrow\infty}\vert F_n^2(\theta)\vert\leq \varepsilon^{1/2},
\end{equation}
for all $\theta\in (-\pi,\pi]$. Finally, consider any cut-off function $\xi\in C^\infty_c([0,+\infty))$ like in \eqref{E-scaled-cut-off}. Notice that
\begin{align*}
\vert F_n^3(\theta)\vert&\leq 2\int_{(-\pi,\pi]\setminus \{\theta\}}\xi\left(\frac{\vert \theta'-\theta\vert_o}{\varepsilon^{1/2}}\right)\,d_{\theta'}(\rho_n+\rho)\\
&\leq 2\int_{(-\pi,\pi]}\xi\left(\frac{\vert \theta'-\theta\vert_o}{\varepsilon^{1/2}}\right)\,d_{\theta'}(\rho_n-\rho)+4\int_{(-\pi,\pi]}\chi_{\vert \theta'-\theta\vert_o\leq \varepsilon^{1/2}}\,d_{\theta'}\rho.
\end{align*}
Consequently,
\begin{equation}\label{E-54}
\limsup_{n\rightarrow\infty}\vert F_n^3(\theta)\vert\leq 4\int_{(-\pi,\pi]}\chi_{\vert \theta'-\theta\vert_o\leq \varepsilon^{1/2}}\,d_{\theta'}\rho.
\end{equation}
Putting \eqref{E-52}, \eqref{E-53} and \eqref{E-54} together yields
$$\limsup_{n\rightarrow\infty}\vert F_n(\theta)\vert\leq \varepsilon^{1/2}+4\int_{(-\pi,\pi]}\chi_{\vert \theta'-\theta\vert_o\leq \varepsilon^{1/2}}\,d_{\theta'}\rho,$$
for every $\varepsilon>0$. The fist term vanishes as $\varepsilon\rightarrow 0$. Since $\theta$ is a continuity point of $\rho$, i.e., $\rho(\{\theta\})=0$, so does the second term. Thus,
$$\limsup_{n\rightarrow\infty}\vert F_n(\theta)\vert=0,$$
and that ends the first part. The second part is clear and follows, for instance, from the Lebesgue decomposition theorem. 
\end{proof}

\begin{cor}\label{C-transport-field-time-convergence-critical-1}
Consider $\alpha=\frac{1}{2}$, $K>0$ and let $\mu$ be in $C([0,T],\mathcal{P}(\mathbb{T}\times \mathbb{R})-\mbox{narrow})$. Then,
\begin{equation}\label{E-55}
\lim_{\tau\rightarrow t}\sup_{\Omega\in\mathbb{R}}\vert \mathcal{V}[\mu_\tau](z,\Omega)-\mathcal{V}[\mu_t](z,\Omega)\vert=0,
\end{equation}
for every almost every $(t,z)\in [0,T]\times \mathbb{T}$.
\end{cor}

\begin{proof}
Theorem \ref{T-transport-field-convergence-critical} provides a negligible set $\mathcal{N}_t\subseteq \mathbb{T}$ so that \eqref{E-55} takes place, for every $z\in \mathbb{T}\setminus\mathcal{N}_t$ and each $t\in [0,T]$. Consider the set
$$\mathcal{N}:=\{(t,z):\,t\in [0,T],\,z\in \mathcal{N}_t\}.$$
It is a negligible set of $\mathbb{T}\times \mathbb{R}$ and the claimed convergence takes place in $([0,T]\times \mathbb{T})\setminus \mathcal{N}$.
\end{proof}

Similarly we obtain.

\begin{cor}\label{C-transport-field-time-convergence-critical-2}
Consider $\alpha=\frac{1}{2}$, $K>0$ and let $\{\mu^n\}_{n\in \mathbb{N}}$ and $\mu$ be in $C([0,T],\mathcal{P}(\mathbb{T}\times \mathbb{R})-\mbox{narrow})$ with
$$\mu^n\rightarrow\mu\ \mbox{ in }\ C([0,T],\mathcal{P}(\mathbb{T}\times \mathbb{R})-\mbox{narrow}).$$
Then, the following convergence takes place
$$
\lim_{n\rightarrow\infty}\sup_{\Omega\in\mathbb{R}}\vert \mathcal{V}[\mu^n_t](z,\Omega)-\mathcal{V}[\mu_t](z,\Omega)\vert=0,
$$
for almost every $(t,z)\in [0,T]\times \mathbb{T}$.
\end{cor}

In the sequel, we will resort on the concept of solutions in the sense of Filippov, to be understood in the sense of absolutely continuous solutions that solve the differential inclusion into the Filippov set-valued tangent field almost everywhere. The following Lemmas summarize the main results regarding the existence of Caratheodory solution to differential inclusions, see \cite{Au-Ce,F,Lo,P}.

\begin{lem}\label{L-Fillipov-map-nonautonomous-1}
Let $F:\mathbb{R}_0^+\times \mathbb{R}^d\longrightarrow\mathbb{R}^d$ be any measurable map. Assume that for every compact subset $K\subseteq \mathbb{R}_0^+\times \mathbb{R}^d$ there exists a nonnegative function $m_K\in L^1_{loc}(\mathbb{R}_0^+)$ such that
\begin{equation}\label{E-42}
\vert F(t,x)\vert \leq m_K(t),\ \forall\,(t,x)\in K.
\end{equation}
Consider its associated Filippov set-valued map with respect to the variable $x$, namely,
\begin{equation}\label{E-Filippov-map-x}
\mathcal{K}[F(t,\cdot)](x)=\bigcap_{\delta>0}\bigcap_{\vert\mathcal{N}\vert =0}\overline{\co}(F(t,B_\delta(x)\setminus \mathcal{N})),
\end{equation}
for all $(t,x)\in \mathbb{R}_0^+\times \mathbb{R}^N$. Then, it satisfies the following properties:
\begin{enumerate}
\item $\mathcal{K}[F(t,\cdot)](x)$ is non-empty, convex and compact for a.e. $t\geq 0$ and every $x\in \mathbb{R}^d$.
\item The set-valued map $t\in\mathbb{R}_0^+ \longmapsto \mathcal{K}[F(t,\cdot)](x)$ is Effros-measurable for every $x\in \mathbb{R}^d$.
\item $\mathcal{K}[F(t,\cdot)]$ is upper semicontinuous for a.e. $t\geq 0$.
\item For every compact subset $K\subseteq \mathbb{R}_0^+\times \mathbb{R}^d$ and the above $m_K\in L^1_{loc}(\mathbb{R}_0^+)$ one has
\begin{equation}\label{E-43}
\vert m(\mathcal{K}[F(t,\cdot)](x)\vert\leq m_K(t),\ \forall\,(t,x)\in K.
\end{equation}
\end{enumerate}
Here, $m(C)=\pi_C(0)$ and $\pi_C$ is the orthogonal projection operator onto the convex set $C$. Then, the above condition \eqref{E-43} equivalently reads
$$\dist(0,\mathcal{K}[F(t,\cdot)](x)))\leq m_K(t),\ \forall\,(t,x)\in K.$$
\end{lem}

\begin{proof}
All the assertion are clear by definition and the measurability of $F$ except, at most, the third one whose proof is given in \cite[Proposition 2.1.1]{Au-Ce} when applied to the autonomous field $F(t,\cdot)$ for fixed $t$.
\end{proof}

\begin{lem}\label{L-Fillipov-map-nonautonomous-2}
Let $\mathcal{F}:\mathbb{R}_0^+\times \mathbb{R}^d\longrightarrow\mathbb{R}^d\longrightarrow\mathcal{P}(\mathbb{R}^d)$ be any set-valued map. Assume that it fulfills the following assumptions:
\begin{enumerate}
\item $\mathcal{F}$ takes non-empty, convex and compact values.
\item $\mathcal{F}(\cdot,x)$ is Effros-measurable for every $x\in \mathbb{R}^d$.
\item $\mathcal{F}(t,\cdot)$ is upper semicontinuous for a.e. $t\geq 0$.
\item For every compact subset $K\subseteq \mathbb{R}_0^+\times \mathbb{R}^d$ there exists a nonnegative $m_K\in L^1_{loc}(\mathbb{R}_0^+)$ such that
\begin{equation}\label{E-44}
\vert m(\mathcal{F}(t,x))\vert\leq m_K(t),\ \forall\,(t,x)\in K.
\end{equation}
\end{enumerate}
Consider the following initial value problem (IVP) issued at any initial datum $x_0\in \mathbb{R}^d$
$$
\left\{\begin{array}{l}
\dot{x}\in \mathcal{F}(t,x),\\
x(0)=x_0.
\end{array}\right.
$$
Then, it has a local-in-time absolutely continuous solution. If the locally integably boundedness condition \eqref{E-44} holds true globally, i.e. with $K$ replaced with the whole $\mathbb{R}_0^+\times \mathbb{R}^d$, then a global-in-time absolutely continuous solution exists.
\end{lem}

A similar result can be found in \cite[Theorems 2.1.3, 2.1.4]{Au-Ce} when $\mathcal{F}$ satisfies a stronger assumption, namely, $\mathcal{F}$ is upper semicontinuous in the joint variables $(t,x)$ (also see \cite[Lemmas 3.2, 3.3, 3.4]{P-P-S} where the autonomous case is reviewed). When time upper semicontinuity is missing such result does not longer apply. Fortunately, Lema \ref{L-Fillipov-map-nonautonomous-2} provides a solution (see \cite[Theorem 2.7.5]{F} and \cite{Lo,P} for the detailed proofs). Note that such result becomes a literal translation to the multivalued case of the classical Caratheodory's existence theorem for single-valued dynamical system. Finally, putting Lemmas \ref{L-Fillipov-map-nonautonomous-1} and \ref{L-Fillipov-map-nonautonomous-2} together, we arrive at the next result (see \cite[Theorem 2.7.8]{F}).

\begin{lem}\label{L-Fillipov-map-nonautonomous-3}
Let $F:\mathbb{R}_0^+\times \mathbb{R}^d\longrightarrow\mathbb{R}^d$ be any measurable map and assume that the local integrably boundedness condition \eqref{E-42} is satisfied. Let $\mathcal{K}[F(t,\cdot)]$ be the Filippov set-valued map with respect to $x$ according to \eqref{E-Filippov-map-x}. Consider the following initial value problem (IVP) issued at any initial datum $x_0\in \mathbb{R}^d$
$$
\left\{\begin{array}{l}
\dot{x}\in \mathcal{K}[F(t,\cdot)](x),\\
x(0)=x_0.
\end{array}\right.
$$
Then, it has a local-in-time absolutely continuous solution. If the locally integably boundedness condition \eqref{E-42} holds true globally, i.e. with $K$ replaced with the whole $\mathbb{R}_0^+\times \mathbb{R}^d$, then a global-in-time absolutely continuous solution exists.
\end{lem}

Solutions to such differential inclusion are often called \textit{Filippov solutions} to the discontinuous non-autonomous dynamical system
$$\left\{\begin{array}{l}
\dot{x}=F(t,x),\\
x(0)=x_0.
\end{array}\right.$$

\begin{lem}\label{L-existence-characteristic-system-critical}
Consider $\alpha=\frac{1}{2}$, $K>0$ and fix $\mu\in L^\infty(0,T;\mathcal{M}(\mathbb{T}\times \mathbb{R}))$. For any $x_0=(z_0,\Omega_0)\in \mathbb{T}\times \mathbb{R}$ let us consider the characteristic system issued at $x_0$, i.e.,
\begin{equation}\label{E-characteristic-system-critical-1}
\left\{\begin{array}{l}
\displaystyle\frac{dX}{dt}(t;t_0,x_0)\in\mathcal{K}[\mathcal{V}[\mu_t]](X(t;t_0,x_0)),\\
\displaystyle \hspace{0.1cm} X(t_0;t_0,x_0)=x_0.
\end{array}\right.
\end{equation}
Then, \eqref{E-characteristic-system-critical-1} has at least one global-in-time Filippov solution $X(t;t_0,x_0)=(Z(t;t_0,z_0,\Omega_0),\Omega_0)$. Indeed, if we set $z_0=e^{i\theta_0}$ for some $\theta_0\in \mathbb{R}$, then 
$$Z(t;t_0,z_0,\Omega_0)=e^{i\Theta(t;t_0,\theta_0,\Omega_0)},$$
where $\Theta=\Theta(t;t_0,\theta_0,\Omega_0)$ is a global-in-time Filippov solution to
\begin{equation}\label{E-characteristic-system-critical-2}
\left\{
\begin{array}{l}
\displaystyle\frac{d\Theta}{dt}(t;t_0,\theta_0,\Omega_0)\in \mathcal{K}[\mathcal{P}[\mu_t](\cdot,\Omega_0)](\Theta(t;t_0,\theta_0,\Omega_0)),\\
\displaystyle \hspace{0.1cm}\Theta(t_0;t_0,\theta_0,\Omega_0)=\theta_0.
\end{array}
\right.
\end{equation}
\end{lem}

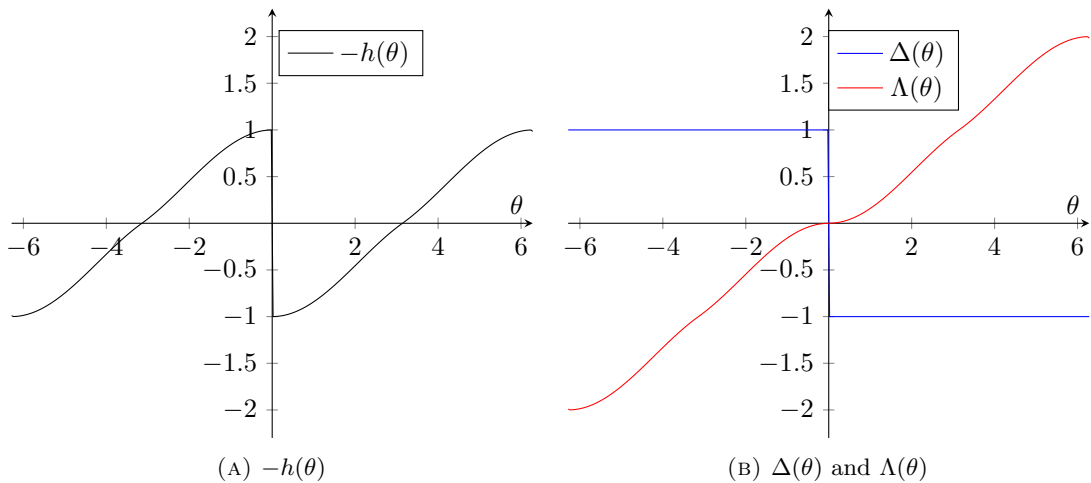
\begin{figure}
\centering
\begin{subfigure}[b]{0.45\textwidth}
\begin{tikzpicture}[
declare function={
d(\x)= (\x<=-pi) * abs(\x+2*pi) + and(\x>-pi, \x<=pi) * abs(\x) + (\x>pi) * abs(\x-2*pi);
}
]
\begin{axis}[
  axis x line=middle, axis y line=middle,
  xmin=-2*pi, xmax=2*pi, xtick={-6,-4,-2,0,2,4,6}, xlabel=$\theta$,
  ymin=-2.3, ymax=2.3, ytick={-2,-1.5,-1,-0.5,0,0.5,1,1.5,2},
  legend style={at={(0.79,0.95)}},
]
\addplot [
    domain=-2*pi+0.001:2*pi-0.001, 
    samples=300, 
    color=black,
]
{-sin(deg(x))/pow(d(x),2*0.5)};
\addlegendentry{$-h(\theta)$}
\end{axis}
\end{tikzpicture}
\caption{$-h(\theta)$}
\label{fig:kuramoto-decomposition-05-1}
\end{subfigure}
\begin{subfigure}[b]{0.45\textwidth}
\begin{tikzpicture}[
  declare function={
 d(\x)= (\x<=-pi) * abs(\x+2*pi) + and(\x>-pi, \x<=pi) * abs(\x) + (\x>pi) * abs(\x-2*pi);
 f(\x)= (\x<0) *  (1) + (\x>0) * (-1); 
 g(\x)= (\x<0) *  (-1-sin(deg(\x))/pow(d(\x),2*0.5))) + (\x>0) * (1-sin(deg(\x))/pow(d(\x),2*0.5))); 
  }
]
\begin{axis}[
  axis x line=middle, axis y line=middle,
  xmin=-2*pi, xmax=2*pi, xtick={-6,-4,-2,0,2,4,6}, xlabel=$\theta$,
  ymin=-2.3, ymax=2.3, ytick={-2,-1.5,-1,-0.5,0,0.5,1,1.5,2},
  legend style={at={(0.75,0.95)}},
]
\addplot[domain=-2*pi+0.001:2*pi-0.001, samples=300,color=blue]{f(x)};
\addlegendentry{$\Delta(\theta)$}
\addplot[domain=-2*pi+0.001:2*pi-0.001, samples=300,color=red]{g(x)};
\addlegendentry{$\Lambda(\theta)$}
\end{axis}
\end{tikzpicture} 
\caption{$\Delta(\theta)$ and $\Lambda(\theta)$}
\label{fig:kuramoto-decomposition-05-2}
\end{subfigure}
\caption{Graph of the function $-h(\theta)$ and the functions $\Delta(\theta)$ and $\Lambda(\theta)$ in the decomposition for the value $\alpha=0.5$.}\label{fig:kuramoto-decomposition-05}
\end{figure}

Our next step is to show the uniqueness of the above Filippov trajectories. Our approach mimics that in Subsection \ref{SS-transport-field-properties}, specifically, we will show an analogue decomposition to that in Lemma \ref{L-split-kernel} (see also Figure \ref{fig:kuramoto-decomposition-05}) that implies an one-sided Lipschitz condition.

\begin{lem}\label{L-split-kernel-critical}
Consider $\alpha=\frac{1}{2}$ and define the couple of functions $\Delta,\Lambda:[-2\pi,2\pi]\longrightarrow\mathbb{R}$ as follows
$$
\Delta(\theta):=\left\{
\begin{array}{ll}
-1-h(\theta), & \theta\in [-2\pi,0),\\
1-h(\theta), & \theta\in (0,2\pi].
\end{array}
\right.
$$
$$
\hspace{-1cm} \Lambda(\theta):=
\left\{
\begin{array}{ll}
1, & \theta\in [-2\pi,0),\\
-1, & \theta\in (0,2\pi].
\end{array}
\right.
$$
Then, the following properties hold true
\begin{enumerate}
\item $\Delta$ is monotonically decreasing, $\Lambda$ is Lipschitz-continuous and 
$$-h(\theta)=\Delta(\theta)+\Lambda(\theta),\ \forall\theta\in [-2\pi,2\pi].$$
\item $-h$ is one-sided Lipschitz in $[-2\pi,2\pi]$, i.e., there exists $L_0>0$ such that
$$\left((-h)(\theta_1)-(-h)(\theta_2)\right)(\theta_1-\theta_2)\leq L_0(\theta_1-\theta_2)^2.$$
\end{enumerate}
\end{lem}

Hence, we are ready to introduce the following results.

\begin{lem}\label{L-transport-field-sided-Lipschitz-critical}
Consider $\alpha\in \frac{1}{2}$, $K>0$ and set $\mu\in L^\infty(0,T;\mathcal{M}(\mathbb{T}\times \mathbb{R}))$. Then, for $L_0$ given in Lemma \ref{L-split-kernel-critical} we obtain that
$$\left(\mathcal{P}[\mu_t](\theta_1,\Omega_1)-\mathcal{P}[\mu_t](\theta_2,\Omega_2)\right)(\theta_1-\theta_2)\leq (\Omega_1-\Omega_2)(\theta_1-\theta_2)+KL_0\Vert \mu\Vert_{L^\infty(0,T;\mathcal{M}(\mathbb{T}\times \mathbb{R}))}(\theta_1-\theta_2)^2,$$
for every $t\in [0,T]$, each $\theta_1,\theta_2\in \mathbb{R}$ with $\theta_1-\theta_2\in [-\pi,\pi]$ and any $\Omega_1,\Omega_2\in\mathbb{R}$.
\end{lem}

\begin{theo}\label{T-transport-field-sided-Lipschitz-critical}
Consider $\alpha=\frac{1}{2}$, $K>0$ and set $\mu\in L^\infty(0,T;\mathcal{M}(\mathbb{T}\times \mathbb{R}))$. Then, $\mathcal{V}[\mu]$ is one-sided Lipschitz in $\mathbb{T}\times \mathbb{R}$ uniformly in $t\in [0,T]$, i.e., there exists $L=L(\alpha,K,\mu)>0$ such that
$$\left<\mathcal{V}[\mu_t](z_2,\Omega_2),\widehat{\gamma}'(1)\right>-\left<\mathcal{V}[\mu_t](z_1,\Omega_1),\widehat{\gamma}'(0)\right>\leq L\,d((z_1,\Omega_1),(z_2,\Omega_2))^2,$$
for every $t\in [0,T]$, any $(z_1,\Omega_1),(z_2,\Omega_2)\in \mathbb{T}\times\mathbb{R}$ and each minimizing geodesic $\widehat{\gamma}:[0,1]\longrightarrow\mathbb{T}\times \mathbb{R}$ in the manifold $\mathbb{T}\times \mathbb{R}$ joining $(z_1,\Omega_1)$ to $(z_2,\Omega_2)$.
\end{theo}

The proofs follow similar arguments to those in Lemma \ref{L-transport-field-sided-Lipschitz} and Theorem \ref{T-transport-field-sided-Lipschitz}, then we omit them. Let us finally transform the point-wise one-sided Lipschitz condition into an analogue one-sided condition in the multivalued sense for the associated Filippov set-valued map. 

\begin{defi}\label{D-one-sided-Lipschitz-multivalued}
Let $(M,\left<\cdot,\cdot\right>)$ be a complete Riemannian manifold and $V:M\longrightarrow 2^{TM}$ a set-valued tangent vector field along $M$. Then, we will say that $V$ is one-sided Lipschitz (in multivalued sense) when there exists a constant $L>0$ such that
$$
\left<v_y,\gamma'(1)\right>-\left<v_x,\gamma'(0)\right>\leq Ld(x,y)^2,
$$
for every $x,y\in M$, each $v_x\in X_x$, $v_y\in X_y$ and any minimizing geodesic $\gamma:[0,1]\longrightarrow M$ joining $x$ to $y$.
\end{defi}

\begin{lem}\label{L-one-sided-Lipschitz-multivalued}
Let $(M,\left<\cdot,\cdot\right>)$ be a complete Riemannian manifold and consider an essentially locally bounded tangent field $V\in \mathfrak{X}L^\infty_{loc}(M)$. If $V$ is one-sided Lipschitz a.e. with constant $L>0$, then so is its Filippov set-valued tangent field $\mathcal{K}[V]$ given by Definition \ref{D-Filippov-field} with same constant $L$.
\end{lem}

The starting point of the proof is that the analogue result with the general $(M,\left<\cdot,\right>)$ replaced by the Euclidean space with the flat metric is clearly true and, in particular, was shown in \cite[Lemma 3.5]{P-P-S}. Then, we will use the local description in coordinates appearing in Definition \ref{D-Filippov-field} to augment our local result into a global one. 

\begin{proof}
Consider $x,y\in M$, set any minimizing geodesic $\gamma:[0,1]\longrightarrow M$ joining $x$ to $y$ and fix any $v_x\in \mathcal{K}[V]_{x}$ and $v_y\in \mathcal{K}[V]_y$. Our goal is to show that
\begin{equation}\label{E-45}
\left<v_y,\gamma'(1)\right>-\left<v_y,\gamma'(0)\right>\leq L d(x,y)^2.
\end{equation}

$\bullet$ \textit{Step 1: Local normal orthonormal frame around $x$}. We can construct a local frame $(\mathcal{U},\varphi,E_1,\ldots,E_d)$ centered at $x$ where $\{E_1,\ldots,E_d\}$ is an orthonormal basis at each point and $(\mathcal{U},\varphi)$ is a normal neighborhood. The existence of normal neighborhoods around any points is classical in Riemannian Geometry (see \cite{DC,K,Pe}) and follows from the inverse function theorem. Indeed, let us consider the injectivity radius of $M$ at $x$
$$\iota(x):=\dist(x,\cut(x))=\sup\{\delta>0:\,\exp_x:\mathbb{B}_\delta(0)\longrightarrow \exp_x(\mathbb{B}_\delta(0))\ \mbox{ is a diffeomorphism}\},$$
that is a positive number. Then, any $\delta<\iota(x)$ provides a geodesic ball $\exp_x(\mathbb{B}_\delta(0))$. Take $\mathcal{U}=\exp_x(\mathbb{B}_\delta(0))$ as normal neighborhood of $x$. 

Regarding the local chart, let us consider $\varphi:=(\exp_x\circ \mathcal{L})^{-1}$ for any linear isometry $\mathcal{L}:\mathbb{R}^d\longrightarrow T_x M$, that maps the standard basis $\{e_1,\ldots,e_d\}$ of $\mathbb{R}^d$ into some orthonormal basis $\{u_1,\ldots,u_d\}$ of $T_x M$, i.e., $\mathcal{L}(e_i)=u_i$ for every $i=1,\ldots,d$. Obviously, the associated local frame $\frac{\partial}{\partial x_1},\ldots,\frac{\partial}{\partial x_d}$ is not necessarily orthonormal in the full $\mathcal{U}$. In fact, only locally flat Riemannian manifolds (like $M=\mathbb{T}\times \mathbb{R}$) can enjoy such property. Nevertheless, we can augment $\{u_1,\ldots,u_d\}$ into a complete local frame on $\mathcal{U}$ through the next procedure. Consider the unique minimizing geodesic $\gamma_{x,z}:[0,1]\longrightarrow M$ joining $x$ to $z$, i.e.,
$$\gamma_{x,z}(s):=\exp_x(s\exp_x^{-1}(z)),\ s\in [0,1],$$
for any $z\in \mathcal{U}$. Then, we can define the local frame at $z$ via paralleling transporting $\{u_1,\ldots,u_d\}$ along $\gamma_{x,z}$, that is,
$$U_i(z):=\tau[\gamma_{x,z}]_0^1(u_i),$$
for all $i=1,\ldots,d$. Since the parallel transport is a linear isometry between tangent spaces we recover the local orthonormal character of the frame $U_1,\ldots,U_d$. In addition, the following property holds true by definition
$$
\nabla_{\gamma'_{x,z}(s)} U_i(\gamma_{x,z}(s))=0,
$$
for all $s\in [0,1]$ and every $z\in \mathcal{U}$. In other words, the tangent fields $U_i$ are all parallel along any radial geodesic issued at $x$ within the normal neighborhood $\mathcal{U}$. 

$\bullet$ \textit{Step 2: One-sided Lipschitz property of $F$}. Using such local frame $(\mathcal{U},\varphi,E_1,\ldots,E_d)$, we can locally write the tangent field $V$ as follows
$$V=\sum_{i=1}^d f_iE_i\ \mbox{ in }\ \mathcal{U}.$$
Again, define $F_i=f_i\circ\varphi^{-1}$ and the vector field $F=(F_1,\ldots,F_d)$. Our goal is to see that $F$ is one-sided-Lipschitz at $\varphi(x)=0$ with constant $L$. To that end, set any $z\in \mathcal{U}$ and consider the local coordinates
\begin{align*}
V_x&=\sum_{i=1}^d F_i(\varphi(x)) E_i(x)=\sum_{i=1}^d F_i(\varphi(x))e_i,\\
V_z&=\sum_{i=1}^d F_i(\varphi(z)) E_i(z)=\sum_{i=1}^d F_i(\varphi(z))\tau[\gamma_{x,z}]_0^1(e_i)=\tau[\gamma_{x,z}]\left(\sum_{i=1}^d F_i(\varphi(z))e_i\right).
\end{align*}
Since $V$ is one-sided Lipschitz at $x$ with constant $L$, then
$$\left<V_z,\gamma_{x,z}'(1)\right>-\left<V_x,\gamma_{x,z}'(0)\right>\leq L\,d(x,z)^2.$$
Notice that $\gamma_{x,z}'(0)=\exp_x^{-1}(z)=\mathcal{L}(\varphi(z))$ and, consequently,
\begin{align*}
(F(\varphi(z))-F(\varphi(x)))\cdot (\varphi(z)-\varphi(0))&=\left<\tau[\gamma_{x,z}]_1^0(V_z)-V_x,\gamma'_{x,z}(0)\right>\\
&=\left<V_z,\gamma_{x,z}'(1)\right>-\left<V_x,\gamma_{x,z}'(0)\right>\\
&\leq L\,d(x,z)^2=L\,\vert \varphi(z)-\varphi(x)\vert^2,
\end{align*}
for every $z\in \mathcal{U}$. As advanced in the above comment, \cite[Lemma 3.5]{P-P-S} shows that $\mathcal{K}[F]$ is one-sided Lipschitz continuous at $\varphi(x)$ with constant $L$, i.e.,
\begin{equation}\label{E-47}
(\beta^z-\beta^x)\cdot (\varphi(z)-\varphi(x))\leq L\vert \varphi(z)-\varphi(x)\vert^2,
\end{equation}
for every $z\in \mathcal{U}$ and every $\beta^x\in \mathcal{K}[F](\varphi(x))$ and $\beta^z\in \mathcal{K}[F](\varphi(z))$.

$\bullet$ \textit{Step 3: Local result}. In this step we will assume that the whole minimizing geodesic $\gamma$ joining $x$ to $y$ lies in $\mathcal{U}$, i.e.,  $\gamma(s)\in \mathcal{U}$ for all $s\in [0,1]$. By the uniqueness, we obtain $\gamma=\gamma_{x,y}$. In that case, we write
\begin{align*}
v_x&=\sum_{i=1}^d \beta_i U_i(x)=\sum_{i=1}^d\beta_i u_i,\\
v_y&=\sum_{i=1}^d\gamma_i U_i(y)=\sum_{i=1}^d\gamma_i\tau[\gamma_{x,y}]_0^1(u_i)=\tau[\gamma_{x,y}]_0^1\left(\sum_{i=1}^d \gamma_i u_i\right),
\end{align*}
for some coefficients $\beta^x=(\beta_1^x,\ldots,\beta_d^x)\in \mathcal{K}[F](\varphi(x))$ and $\beta^y=(\beta_1^y,\ldots,\beta_d^y)\in \mathcal{K}[F](\varphi(y))$. Note that $\gamma'(0)=\exp_x^{-1}(y)=\mathcal{L}(\varphi(y))$ and $\gamma$ is a geodesic. Then, we infer
$$\left<v_y,\gamma'(1)\right>-\left<v_x,\gamma'(0)\right>=\left<\sum_{i=1}^d (\beta_i^y-\beta_i^x)u_i,\gamma'(0)\right>=(\beta^y-\beta^x)\cdot (\varphi(y)-\varphi(0)).$$
Hence, \eqref{E-45} follows from \eqref{E-47}.
 
$\bullet$ \textit{Step 4: Global result}. Now, let us assume that the minimizing geodesic $\gamma$ does not necessarily lies in the normal neighborhood $\mathcal{U}$. Notice that for every $s\in [0,1]$ we can repeat \textit{Step 1} to find a normal orthonormal frame $(\mathcal{U}^s,\varphi^s,E_1^s,\ldots,E_d^s)$ around $\gamma(s)$. Indeed, recall that $\mathcal{U}^s=\exp_{\gamma(s)}(\mathbb{B}_{\delta_s}(0))$, for any $0<\delta<\iota(\gamma(s))$ and each $s\in [0,1]$. Notice that we can choose all the $\delta_s$ to be the same $\delta$ independently on $s$. To such end, recall that the injectivity radius $\iota:M\longrightarrow\mathbb{R}^+$ is a continuous function (see \cite[Proposition 2.1.10]{K}). Since $\gamma([0,1])$ is compact, then
$$\delta_0:=\min_{s\in [0,1]}\iota(\gamma(s))>0,$$
and we choose such any $0<\delta<\delta_0$ as the radii of each geodesic ball about $\gamma(s)$. Fix $n\in \mathbb{N}$ the smallest integer with $d(x,y)\leq k\delta$ and define the numbers
$$s_k=\left\{\begin{array}{ll} k\delta, & k=0,\ldots,n-1,\\ 1, & k=n.\end{array}\right.$$
Since $\gamma$ is a minimizing geodesic we infer
$$\gamma([s_k,s_{k+1}])\subseteq \mathcal{U}^{s_k},$$
for all $k=0,\ldots,n-1$. Hence, for all $k=0,\ldots,n-1$ the piece of geodesic $\gamma([s_k,s_{k+1}])$ satisfies the same properties as in \textit{Step 1} for the normal local orthonormal frame $(\mathcal{U}^{s_k},\varphi^{s_k},E_1^{s_k},\ldots,E_d^{s_k})$. Let us consider any $v_{\gamma(s_k)}\in \mathcal{K}[V]_{\gamma(s_k)}$, for every $k=1,\ldots,n-1$. Thus, we have
$$\left<v_{\gamma(s_{k+1})},\gamma'(s_{k+1})\right>-\left<v_{\gamma(s_{k})},\gamma'(s_{k})\right>\leq L\,d(\gamma(s_k),\gamma(s_{k+1}))^2,$$
for every $k=0,\ldots,n-1$. Then,
\begin{align*}
\left<v_y,\gamma'(1)\right>-\left<v_x,\gamma'(0)\right>&=\sum_{k=0}^{n-1}\left<v_{\gamma(s_{k+1})},\gamma'(s_{k+1})\right>-\left<v_{\gamma(s_{k})},\gamma'(s_{k})\right>\\
&\leq L\sum_{k=0}^{n-1}d(\gamma(s_k),\gamma(s_{k+1}))^2\leq L\left(\sum_{k=0}^{n-1}d(\gamma(s_k),\gamma(s_{k+1}))\right)^2=L\,d(x,y)^2,
\end{align*}
where we have used that $\gamma$ is minimizing in the last step.
\end{proof}

As a direct consequence of Theorem \ref{T-transport-field-sided-Lipschitz-critical} and Lemma \ref{L-one-sided-Lipschitz-multivalued} we obtain the following result.

\begin{cor}\label{T-set-valued-Filippov-field-sided_lipschitz-critical}
Consider $\alpha=\frac{1}{2}$, $K>0$ and set $\mu\in L^\infty(0,T;\mathcal{M}(\mathbb{T}\times \mathbb{R}))$. Then, the Filippov set-valued vector field $\mathcal{K}[\mathcal{V}[\mu]]$ associated with the transport field $\mathcal{V}[\mu]$ is one-sided Lipschitz in $\mathbb{T}\times \mathbb{R}$ in multivalued sense uniformly in $t\in [0,T]$, i.e., there exists $L=L(\alpha,K,\mu)>0$ such that
$$\left<v_{(z_2,\Omega_2)},\widehat{\gamma}'(1)\right>-\left<v_{(z_1,\Omega_1)},\widehat{\gamma}'(0)\right>\leq Ld((z_1,\Omega_1),(z_2,\Omega_2))^2,$$
for every $t\in [0,T]$, each $(z_1,\Omega_1),\,(z_2,\Omega_2)\in \mathbb{T}\times \mathbb{R}$, any couple $v_{(z_1,\Omega_1)}\in \mathcal{K}[\mathcal{V}[\mu_t]](z_1,\Omega_1)$ and $v_{(z_2,\Omega_2)}\in \mathcal{K}[\mathcal{V}[\mu_t]](z_2,\Omega_2)$, and each minimizing geodesic $\widehat{\gamma}:[0,1]\longrightarrow\mathbb{T}\times \mathbb{R}$ joining $(z_1,\Omega_1)$ to $(z_2,\Omega_2)$.
\end{cor}

We are now ready to recover a one-sided uniqueness property for the Filippov characteristics in \eqref{E-characteristic-system-critical-1}.

\begin{theo}\label{T-well-posedness-characteristic-system-critical}
Consider $\alpha= \frac{1}{2}$, $K>0$ and fix $\mu\in L^\infty(0,T;\mathcal{M}(\mathbb{T}\times \mathbb{R}))$. The characteristic system \eqref{E-characteristic-system-critical-1} associated with the transport field $\mathcal{V}[\mu]$ enjoys a global-in-time absolutely continuous Filippov solution that is unique forward-in-time for every given initial data $x_0=(z_0,\Omega_0)=(e^{i\theta_0},\Omega_0)\in\mathbb{T}\times \mathbb{R}$. Indeed, the same representation of the solution holds true, i.e.,
$$X(t;t_0,x_0)=(Z(t;t_0,z_0,\Omega_0),\Omega_0)=(e^{i\Theta(t;t_0,\theta_0,\Omega_0)},\Omega_0),\ t\geq t_0,$$
where $\Theta(t;t_0,\theta_0,\Omega_0)$ is the unique forward-in-time absolutely continuous Filippov solution to \eqref{E-characteristic-system-critical-2}.
\end{theo}

We will say that $X(\cdot\,;t_0,x_0)$ is the \textit{Filippov characteristic} issued at $x_0$ at time $t_0$ and $X(t;t_0,\cdot)$
is the \textit{Filippov flow} from $t_0$ to $t$. The proof is a simple adaptation to that of Theorem \ref{T-well-posedness-characteristic-system}, when classical solutions are replaced with Filippov solutions. Again, the proof relies on the one-sided Lipschitz condition in multivalued sense for the Fillippov set-valued tangent field associated with the transport field (recall Theorem  \ref{T-transport-field-sided-Lipschitz-critical}). Also, the weak differentiability properties of the squared distance in Appendix \ref{appendix-differentiability-distance} will be used. For clarity, we just provide a sketch.

\begin{proof}
Consider two different solutions in the sense of Filippov $x_1=x_1(t)$ and $x_2=x_2(t)$ to the characteristic system \eqref{E-characteristic-system-critical-1} issued at $x_1(t_0)=x_0=x_2(t_0)$ and define
$$I(t):=\frac{1}{2}d(x_1(t),x_2(t))^2,\ t\geq t_0.$$
Since the Filippov trajectories are at least locally absolutely continuous, then so is $I=I(t)$. Again, we compute it in terms of the one-sided Dini upper derivative
$$\frac{dI}{dt}\equiv\frac{d^+ I}{dt}=d^+\left(\frac{1}{2}d_{x_2(t)}^2\right)_{x_1(t)}(\dot{x}_1(t))+d^+\left(\frac{1}{2}d_{x_1(t)}^2\right)_{x_2(t)}(\dot{x}_2(t)),$$
for almost every $t\geq t_0$. By virtue of Theorem \ref{T-dini-derivative-distance} 
\begin{align*}
d^+\left(\frac{1}{2}d_{x_2(t)}^2\right)_{x_1(t)}(\dot{x}_1(t))&\leq \inf_{\substack{w_1\in \exp_{x_1(t)}^{-1}(x_2(t))\\ \vert w_1\vert =d(x_1(t),x_2(t))}}-\left<\dot{x}_1(t),w_1\right>,\\
d^+\left(\frac{1}{2}d_{x_1(t)}^2\right)_{x_2(t)}(\dot{x}_2(t))&\leq \inf_{\substack{w_2\in \exp_{x_2(t)}^{-1}(x_1(t))\\ \vert w_2\vert =d(x_1(t),x_2(t))}}-\left<\dot{x}_2(t),w_2\right>.
\end{align*}
Fix a minimizing geodesic $\widehat{\gamma}_t:[0,1]\longrightarrow\mathbb{T}\times \mathbb{R}$ joining $x_1(t)$ to $x_2(t)$, for almost every $t\geq t_0$. Then, we can choose $w_1=\widehat{\gamma}_t'(0)$ and $w_2=-\widehat{\gamma}_t'(1)$ in the above inequalities. Consequently,
$$\frac{dI}{dt}\leq \left<\dot{x}_2(t),\widehat{\gamma}_t'(1)\right>-\left<\dot{x}_1(t),\widehat{\gamma}_t'(0)\right>,$$
where, recall that $\dot{x}_1(t)\in \mathcal{K}[\mathcal{V}[\mu_t]](x_1(t))$ and $\dot{x}_2(t)\in \mathcal{K}[\mathcal{V}[\mu_t]](x_2(t))$. Then, Theorem \ref{T-transport-field-sided-Lipschitz} implies
$$\frac{dI}{dt}\leq L d(x_1(t),x_2(t)^2)=2L\,I(t),\ \mbox{ a.e. }t\geq t_0,$$
and Gr\"{o}nwall's lemma concludes the proof.
\end{proof}

\subsubsection{Empirical measures as solutions in the sense of the Filippov flow}\label{SSS-empirical-measures-critical} 
Before entering into details, we need to recall some notation that was introduced in \cite[Subsection 2.4]{P-P-S} for the agent-based system \eqref{E-kuramoto-discrete}. Specifically, we will need appropriate notation for cluster of oscillators that arise after collisions in the microscopic dynamics. Take any phase configuration of the $N$ oscillators, i.e.,
\[
\Theta=(\theta_1,\ldots,\theta_N)\in\mathbb{R}^N.
\]
We will say that the $i$-th oscillator \textit{collides} with $j$-th oscillator when $\bar\theta_i=\bar\theta_j$. That induces the relation
\[
i\overset{\Theta}{\sim}j\ \mbox{ when }\ \bar\theta_i=\bar\theta_j.
\]
Since it is an equivalence relation, we can denote its equivalence classes by
$$
\mathcal C_i(\Theta):=\{j\in\{1,\ldots,N\}:\,i\overset{\Theta}{\sim}j\}.
$$
As it is apparent from the definition, $\mathcal C_i(\Theta)$ is the set of \textit{indices of collision} with the $i$-th oscillator. Each of the above equivalence classes can be regarded as a \textit{cluster}. Let us denote by $\mathcal{E}(\Theta)$ the family of all the different clusters. Then, $\mathcal{E}(\Theta)$ establish a partition of $\{1,\ldots,N\}$ that we will call the \textit{collisional type} of $\Theta$. For simplicity of notation, we will enumerate the different equivalence classes
\[
\mathcal{E}(\Theta)=\{E_1(\Theta),\ldots,E_{\kappa(\Theta)}(\Theta)\},\]
in such a way that the minimal representatives in each of them, i.e., $\iota_k(\Theta):=\min E_k(\Theta)$, are increasingly ordered. $\kappa(\Theta):=\#\mathcal{E}(\Theta)$ will denote the total amount of clusters in such a phase configuration $\Theta$ and we will denote the size of the $k$-th cluster, that is the number of particles which form the $k$-th cluster, by $n_k(\Theta):=\# E_k(\Theta)$, for each $k=1,\ldots,\kappa(\Theta)$. If we are not only given a static phase configuration but a whole absolutely continuous trajectory $t\mapsto\Theta(t)=(\theta_1(t),\ldots\theta_N(t))\in\mathbb{R}^N$ governing the dynamics of the $N$ oscillators. Then, as long as it is clear from the context, we will simplify the notation and will denote
\begin{equation}\label{E-cluster-notation}
\mathcal C_i(t):=\mathcal C_i(\Theta(t)),\ \mathcal{E}(t):=\mathcal{E}(\Theta(t)),\ \kappa(t):=\kappa(\Theta(t)),\ n_k(t):=n_k(\Theta(t)).
\end{equation}
Similarly, time may be omitted in our notation for simplicity. 

Our next goal is to recover an analogue to Theorem \ref{T-empirical-measures-subcritical} for the critical regime $\alpha=\frac{1}{2}$, where the classical solutions to \eqref{E-kuramoto-discrete}-\eqref{E-kuramoto-kernel} are replace with the Filippov solutions obtained in \cite{P-P-S}. We will no longer obtain weak-measure valued solutions like in the preceding theorem but measure valued solutions in the sense of the Filippov flow, to be introduced in the sequel.

\begin{theo}\label{T-empirical-measures-critical}
Consider $\alpha=\frac{1}{2}$ and $K>0$. Fix $N\in\mathbb{N}$ and consider $N$ oscillators with initial phases and natural frequencies given by the configurations
$$\Theta_0^N=(\theta_{1,0}^N,\ldots,\theta_{N,0}^N)\ \mbox{ and }\ \{\Omega_i^N:\,i=1,\ldots,N\}.$$ 
Consider the (forward-in-time) unique solution $\Theta^N(t)=(\theta_1^N(t),\ldots,\theta^N_t)$ to \eqref{E-kuramoto-discrete}-\eqref{E-kuramoto-kernel} in the sense of Filippov associated with such initial configuration as given in \cite[Theorem 3.3]{P-P-S} and set the corresponding empirical measures $\mu^N$ according to Definition \ref{D-empirical-measures-subcritical}. Then, $\mu^N\in \mathcal{AC}_\mathcal{M}\cap \mathcal{T}_\mathcal{M}$ and
\begin{equation}\label{E-solutions-Filippov-flow-empirical}
X^N(t;0,\cdot)_{\#}\mu^N_0=\mu^N_t,
\end{equation}
for every $t\geq 0$, where $X^N(t;0,\cdot)$ is the Filippov flow associated with the transport field $\mathcal{V}[\mu^N]$. In addition, an analogue of \eqref{E-17} holds true, namely,
\begin{equation}\label{E-17-critical}
\left\vert\frac{d}{dt}\int_{\mathbb{T}\times \mathbb{R}}\varphi\,d\mu^N_t\right\vert\leq \left(\frac{1}{N}\sum_{i=1}^N\vert \Omega_i^N\vert+K\right)\Vert \nabla \varphi\Vert_{C_0(\mathbb{T}\times \mathbb{R})},
\end{equation}
for almost every $t\geq 0$ and every $\varphi\in C^1_0(\mathbb{T}\times\mathbb{R})$.
\end{theo}

\begin{proof}
The proof that $\mu^N\in \mathcal{AC}_\mathcal{M}\cap \mathcal{T}_\mathcal{M}$ is parallel to that of Theorem \ref{T-empirical-measures-subcritical}. However, the time regularity is much tighter now. Specifically, fix any $\varphi\in C^\infty_c(\mathbb{T}\times \mathbb{R})$ and consider the map
\begin{equation}\label{E-56}
t\in [0,+\infty)\longmapsto\int_{\mathbb{T}\times \mathbb{R}}\varphi\,d\mu^N_t=\frac{1}{N}\sum_{i=1}^N\varphi(\theta_i^N(t),\Omega_i^N).
\end{equation}
Since $\varphi$ is Lipschitz-continuous and $\Theta^N=\Theta^N(t)$ is a Filippov solution (then, it is merely absolutely continuous), such map is not necessarily $C^1$ like in the subcritical case, but it is locally absolutely continuous at least. The tightness condition is again clear since $\mu^N$ has a bounded first order $\Omega$-moment uniformly-in-time like in the subcritical regime. Our next goal is to prove that \eqref{E-solutions-Filippov-flow-empirical} holds true. Take any $\varphi\in C_b(\mathbb{T}\times \mathbb{R})$ and note that, by definition,
$$\int_{\mathbb{T}\times \mathbb{R}}\varphi\,d(X^N(t;0,\cdot)_{\#}\mu^N_0)=\frac{1}{N}\sum_{i=1}^N\varphi(X^N(t;0,z_i^N(0),\Omega_i^N)).$$
Consequently, we just must check that
$$X^N(t;0,z_i^N(0),\Omega_i^N)=(z_i^N(t),\Omega_i^N),$$
for every $t\geq 0$. By the one-sided uniqueness of Filippov trajectories (recall Theorem \ref{T-well-posedness-characteristic-system-critical}), we only need to show that
$$\frac{d}{dt}(z_i^N,\Omega_i^N)\in \mathcal{K}[\mathcal{V}[\mu^N_t]](z_i^N,\Omega_i^N).$$
Using the representation of the Filippov set-valued map in Proposition \ref{P-Filippov-transport-field-critical}, notice that we equivalently need to show that
\begin{equation}\label{E-48}
\dot{\theta}_i^N\in \mathcal{K}[\mathcal{P}[\mu^N_t](\cdot,\Omega_i^N)](\theta_i^N),
\end{equation}
for almost every $t\geq 0$. Note that, bearing in mind the above notation \eqref{E-cluster-notation} implies
\begin{align*}
\mathcal{P}[\mu^N_t](\theta,\Omega_i^N)=\Omega_i^N-\frac{K}{N}\sum_{\substack{j=1\\ \bar\theta_j^N(t)\neq \bar\theta}}^N h(\theta-\theta_j^N(t))
=\Omega_i^N-\frac{K}{N}\sum_{\substack{k=1\\ \bar\theta_{\iota_k}^N(t)\neq \bar \theta}}^{\kappa^N(t)}n_k(t)h(\theta-\theta_{\iota_k}^N(t)),
\end{align*}
where $\bar\theta$ is the representative of $\theta$ in $(-\pi,\pi]$ modulo $2\pi$. Then, one clearly obtains that
\begin{equation}\label{E-49}
\mathcal{K}[\mathcal{P}[\mu^N_t](\cdot,\Omega_i^N)](\theta)=\left\{\Omega_i^N-\frac{K}{N}\sum_{\substack{j=1\\ \bar\theta_j^N(t)\neq \bar \theta}}^{N}h(\theta-\theta_j^N(t))-\frac{K}{N}\sum_{\substack{j=1\\ \bar\theta_j^N(t)= \bar \theta}}^{N} \widehat{y}_i:\,\widehat{y}_i\in [-1,1]\right\}.
\end{equation}
In particular, if we evaluate \eqref{E-49} at $\theta=\theta_i^N(t)$ we obtain the set
\begin{eqnarray}\label{E-50}
\mathcal{K}[&&\hspace{-0.7cm}\mathcal{P}[\mu^N_t](\cdot,\Omega_i^N)](\theta_i^N(t)) \nonumber\\
&&=\left\{\Omega_i^N-\frac{K}{N}\sum_{j\notin \mathcal{C}_i(\Theta^N(t))}h(\theta_i^N(t)-\theta_j^N(t))-\frac{K}{N}\#\mathcal{C}_i(\Theta^N(t))\widehat{y}_i:\,\widehat{y}_i\in [-1,1]\right\}.
\end{eqnarray}
On the other hand, since $\Theta^N=\Theta^N(t)$ is a Filippov solution to \eqref{E-kuramoto-discrete}-\eqref{E-kuramoto-kernel}, the characterization of the Filippov set-valued map associated with its right hand side (see \cite[Proposition 3.2]{P-P-S}) shows that there exists $Y^N=(y^N_{ij})_{1\leq i,j\leq N}$ with $y_{ij}^N\in L^\infty(0,+\infty)$ and $Y^N(t)\in \Skew_N([-1,1])$ for a.e. $t\geq 0$ so that
\begin{equation}\label{E-51}
\dot{\theta}_i^N(t)=\Omega_i^N-\frac{K}{N}\sum_{j\notin \mathcal{C}_i(\Theta^N(t))}h(\theta_i^N(t)-\theta_j^N(t))-\frac{K}{N}\sum_{j\in \mathcal{C}_i(\Theta^N(t))}y_{ij}^N(t),
\end{equation}
for almost every $t\geq 0$ and every $i=1,\ldots,N$. Let us define the following mean value
$$\displaystyle\widehat{y}^N_i(t):=\frac{1}{\# \mathcal{C}_i(\Theta^N(t))}\sum_{j\in \mathcal{C}_i(\Theta^N(t))}y_{ij}^N(t)\in [-1,1].$$
Then, it now becomes apparent that from such choice and the explicit expressions \eqref{E-50} and \eqref{E-51} we have that \eqref{E-48} is verified. Finally, let us conclude \eqref{E-17-critical}. Fix any $\varphi\in C^1_0(\mathbb{T}\times \mathbb{R})$ and take derivatives almost everywhere in \eqref{E-56}
$$\frac{d}{dt}\int_{\mathbb{T}\times \mathbb{R}}\varphi\,d\mu^N_t=\frac{d}{dt}\left(\frac{1}{N}\sum_{i=1}^N\varphi(\theta_i^N(t),\Omega_i^N)\right)=\frac{1}{N}\sum_{i=1}\frac{\partial\varphi}{\partial \theta}(\theta_i^N(t),\Omega_i^N)\dot{\theta}_i^N(t).$$
Since $\theta_i^N(t)$ are solutions in the sense of Filippov to the discrete singular system \eqref{E-kuramoto-discrete}-\eqref{E-kuramoto-kernel}, then 
$$\vert \dot{\theta}_i^N(t)\vert\leq \vert \Omega_i^N\vert+K,\ \mbox{ for a.e. }t\geq 0,$$
and that ends the proof.
\end{proof}

Notice that we can repeat all the ideas of the compactness result in Corollary \ref{C-mean-field-compactness-subcritical} for $\alpha=\frac{1}{2}$.

\begin{cor}\label{C-mean-field-compactness-critical}
Consider $\alpha=\frac{1}{2}$ and $K>0$. For every $N\in\mathbb{N}$, et $N$ oscillators with initial phases and natural frequencies given by the configurations
$$\Theta_0^N=(\theta_{1,0}^N,\ldots,\theta_{N,0}^N)\ \mbox{ and }\ \{\Omega_i^N:\,i=1,\ldots,N\}.$$
Consider the (forward-in-time) unique solution $\Theta^N(t)=(\theta_1(t),\ldots,\theta_N(t))$ to \eqref{E-kuramoto-discrete}-\eqref{E-kuramoto-kernel} in the sense of Filippov associated with such initial configurations as given in \cite{P-P-S} and set the corresponding empirical measures $\mu^N$ according to Definition \ref{D-empirical-measures-subcritical}. Assume that the equi-sumability condition \eqref{E-Omega-moment-1-equiint} holds true and take $M_1$ fulfilling \eqref{E-Omega-moment-1} according to Proposition \ref{P-Omega-moment-equint}. Then, for every fixed $T>0$, there exists a subsequence of $\mu^N$, that we denote in the same way for simplicity, and a limiting measure $f\in L^\infty(0,T;\mathcal{M}(\mathbb{T}\times \mathbb{R}))$ such that \eqref{E-18-1} and \eqref{E-18-2} hold true. Moreover,
$$f\in W^{1,\infty}([0,T],C^1_0(\mathbb{T}\times \mathbb{R})^*)\cap C([0,T],\mathcal{P}_1(\mathbb{T}\times \mathbb{R})-W_1)\cap \mathcal{AC}_{\mathcal{M}}\cap\mathcal{T}_\mathcal{M},$$
and the following convergence takes place
$$
\mu^N\rightarrow f\ \mbox{ in }\ C([0,T],\mathcal{P}_1(\mathbb{T}\times \mathbb{R})-W_1),
$$
where $W_1$ means the Kantorovich-Rubinstein distance.
\end{cor}

Notice that the strong uniform convergence in Lemma \ref{L-mean-field-uniform-convergence-convolution} cannot hold in the critical regime. Nevertheless, the sequence $\vert\mathcal{V}[\mu^N]-\mathcal{V}[f]\vert$ is essentially uniformly bounded. Consequently it enjoys a subsequence, that we denote in the same way, so that it converges weakly-star in $L^\infty$. Using a standard application of Banach--Saks' theorem, we claim that the weak-star limit agrees with the a.e. limit in Corollary \ref{C-transport-field-time-convergence-critical-2}. 

\begin{cor}\label{C-mean-field-weak-star-convergence-fields-critical}
Under the assumptions in Corollary \ref{C-mean-field-compactness-subcritical} the following convergence takes place 
$$\vert\mathcal{V}[\mu^N]- \mathcal{V}[f]\vert\overset{*}{\rightharpoonup}0\ \mbox{ in }\ L^\infty((0,+\infty)\times \mathbb{T}\times \mathbb{R}).$$
\end{cor}

Before showing our final result, that allows concluding the mean field limit and the existence of measure-valued solutions in the sense of the Filippov flow, we need to guarantee that the above weak-star convergence implies uniform convergence of the flows. In the sequel, we adapt the result in \cite[Theorem 4.3]{B-G} to vector fields in Riemannian manifolds.

\begin{lem}\label{L-convergence-Filippov-flow}
Let $(M,\left<\cdot,\cdot\right>)$ be a complete Riemannian manifold and consider $\{V^n\}_{n\in \mathbb{N}}$ and $V$ in $L^\infty((0,+\infty), \mathfrak{X}L^\infty_{loc}(M))$ that are one-sided Lipschitz with same constant $L>0$. Assume that
$$\vert V^n-V\vert \rightharpoonup 0\ \mbox{ in }\ L^1_{loc}((0,+\infty)\times \mathbb{T}\times \mathbb{R}).$$
Then, the associated Filippov flows $X^n=X^n(t;0,x)$ and $X=X(t;0,x)$ verify
$$X^n\rightarrow X\ \mbox{ in }\ C_{loc}((0,+\infty)\times \mathbb{T}\times \mathbb{R}).$$
\end{lem}

Since the proof is clear, we omit it. It simply relies on Definition \ref{D-Filippov-field}, where the Filippov set-valued tangent field is introduced using local coordinates, and the analogue result in \cite[Theorem 4.3]{B-G} for quasi-monotone operators in the Euclidean spaces. 

\begin{theo}\label{T-mean-field-existence-critical}
Consider $\alpha=\frac{1}{2}$, $K>0$ and set any initial datum $f_0\in \mathcal{P}_1(\mathbb{T}\times \mathbb{R})$. Then, for every $T>0$ there exists a measure-valued solution $f\in \mathcal{AC}_\mathcal{M}\cap \mathcal{T}_\mathcal{M}$ in the sense of the Filippov flow to the initial value problem \eqref{E-kuramoto-transport-TxR}, i.e.,
$$X(t;0,\cdot)_{\#}f_t=f_0,\ \mbox{ for all }t\in [0,T],$$
where $X=X(t;0,z,\Omega)$ is the Filippov flow associated with $\mathcal{V}[f]$. In addition, \eqref{E-18-1}-\eqref{E-18-2} holds true and
$$f\in W^{1,\infty}([0,T],C^1_0(\mathbb{T}\times \mathbb{R})^*)\cap C([0,T],\mathcal{P}_1(\mathbb{T}\times \mathbb{R})-W_1).$$
\end{theo}

\begin{proof}
Consider a discrete approximation like in Lemma \ref{L-discrete-approximation}. Namely, consider $N$ oscillators with phases and natural frequencies given by the configurations
$$\Theta_0^N=(\theta_{1,0}^N,\ldots,\theta_{N,0}^N)\ \mbox{ and }\ \{\Omega_i^N:\,i=1,\ldots,N\},$$
for every $N\in\mathbb{N}$ so that they verify the equi-sumability condition \eqref{E-Omega-moment-1-equiint} and the associated empirical measures $\mu^N_t\in \mathcal{P}(\mathbb{T}\times \mathbb{R})$ in Definition \ref{D-empirical-measures-subcritical} verify
\begin{equation}\label{E-60}
\lim_{N\rightarrow \infty}W_1(\mu^N_0,f_0).
\end{equation}
By virtue of Theorem \ref{T-empirical-measures-critical}, $\mu^N$ are measure-valued solutions to \eqref{E-kuramoto-transport-TxR} in the sense of the Filippov flow issued at $\mu^N_0$, i.e.,
\begin{equation}\label{E-57}
X^N(t;0,\cdot)_{\#}\mu^N_0=\mu^N_t,\ \mbox{ for all }t\geq 0.
\end{equation}
By Corollary \ref{C-mean-field-compactness-critical}, there exists a limiting measure $f$ so that
$$\mu^N\rightarrow f \ \mbox{ in }\ C([0,T],\mathcal{P}_1(\mathbb{T}\times \mathbb{R})-W_1).$$
Using Corollary \ref{C-mean-field-weak-star-convergence-fields-critical} we claim that
$$\vert \mathcal{V}[\mu^N]-\mathcal{V}[\mu]\vert \overset{*}{\rightharpoonup} 0\ \mbox{ in }\ L^\infty((0,+\infty)\times \mathbb{T}\times \mathbb{R}).$$
In particular, the convergence takes place weakly in $L^1_{loc}((0,+\infty)\times \mathbb{T}\times \mathbb{R})$. Using Lemma \ref{L-convergence-Filippov-flow} we obtain
\begin{equation}\label{E-59}
X^N\rightarrow X\ \mbox{ in }\ C_{loc}((0,+\infty)\times \mathbb{T}\times\mathbb{R}).
\end{equation}
Let us finally show that we can pass to the limit in \eqref{E-57}. Writing \eqref{E-57} in weak form against any test function $\varphi\in C_c(\mathbb{T}\times \mathbb{R})$ we can write
$$
\int_{\mathbb{T}\times \mathbb{R}}\varphi(X^N(t;0,z,\Omega))\,d_{(z,\Omega)}\mu^N_0=\int_{\mathbb{T}\times \mathbb{R}}\varphi(z,\Omega)\,d_{(z,\Omega)}\mu^N_t,
$$
for every $t\geq 0$. First, it is clear that by the above convergence of the empirical measures we can pass to the limit in the right hand side. Regarding the left hand side, we need to prove that the following sequence
$$I^N(t):=\int_{\mathbb{T}\times \mathbb{R}}\varphi(X^N(t;0,z,\Omega))\,d_{(z,\Omega)}\mu^N_0-\int_{\mathbb{T}\times \mathbb{R}}\varphi(X(t;0,z,\Omega))\,d_{(z,\Omega)}f_0,$$
vanishes as $N\rightarrow 0$ for every $t\geq 0$. Consider the following split
$$I^N(t)=A^N_R(t)+B^N_R(t)+C^N(t),$$
for any $R>0$, where each of the terms reads
\begin{align*}
A^N_R(t)&:=\int_{\mathbb{T}\times \mathbb{R}}\xi\left(\frac{\vert \Omega\vert}{R}\right)(\varphi(X^N(t;0,x,\Omega))-\varphi(X(t;0,z,\Omega)))\,d_{(z,\Omega)}\mu^N_0,\\
B^N_R(t)&:=\int_{\mathbb{T}\times \mathbb{R}}\left(1-\xi\left(\frac{\vert \Omega\vert}{R}\right)\right)(\varphi(X^N(t;0,x,\Omega))-\varphi(X(t;0,z,\Omega)))\,d_{(z,\Omega)}\mu^N_0,\\
C^N(t)&:=\int_{\mathbb{T}\times \mathbb{R}}\varphi(X(t;0,z,\Omega))\,d_{(z,\Omega)}(\mu^N_0-f_0).
\end{align*}
Here $\xi\in C_c([0,+\infty))$ is any cut-off function such that $0\leq \xi\leq 1$, $\left.\xi\right\vert_{[0,1]}\equiv 1$ and $\left.\xi\right\vert_{[2,+\infty)}=0$. First, \eqref{E-59} guarantees that
$$\lim_{N\rightarrow \infty}A^N_R(t)=0,$$
for each $R>0$ and every $t\geq 0$. Second, notice that
$$B^N_R(t)\leq 2\Vert \varphi\Vert_{C_0(\mathbb{T}\times \mathbb{R})}\Vert \vert \Omega\vert\chi_{\vert \Omega\vert\geq R}\mu^N_0\Vert_{\mathcal{M}(\mathbb{T}\times \mathbb{R})}\leq 2\Vert \varphi\Vert_{C_0(\mathbb{T}\times \mathbb{R})}\frac{1}{N}\sum_{\substack{1\leq i\leq N\\ \vert \Omega_i^N\vert\geq R}}\vert \Omega_i^N\vert,$$
for every $N\in \mathbb{N}$ and each $t\geq 0$. Third  \eqref{E-60} shows that
$$\lim_{N\rightarrow\infty}C^N(t)=0,$$
for all $t\geq 0$. Putting everything together yields,
$$\limsup_{N\rightarrow \infty}I^N(t)\leq 2\Vert \varphi\Vert_{C_0(\mathbb{T}\times \mathbb{R})}\frac{1}{N}\sum_{\substack{1\leq i\leq N\\ \vert \Omega_i^N\vert\geq R}}\vert \Omega_i^N\vert,$$
for every $R>0$. Therefore, the equi-sumability condition \eqref{E-Omega-moment-1-equiint} ends the proof.
\end{proof}

\subsection{Uniqueness of solutions in the sense of the Filippov flow}\label{SS-critical-regime-solutions-filippov-flow-uniqueness}
In this part, we shall show that the ideas in Theorem \ref{T-growth-L(RK)} of Subsection \ref{SS-uniqueness} for the bound on the fibered quadratic Wasserstein distance can be extended to the critical regime. As a byproduct, we will recover an uniqueness result for solutions in the sense of the Filippov flow to the non-linear transport equation \eqref{E-kuramoto-transport-TxR}. Our proof relies on an approximation of the discontinuous kernel $h$ through the regularized kernels $h_\varepsilon$. In the sequel, we provide some technical lemma that will be used along the proof.

\begin{defi}\label{D-transport-field-regularized-critical}
Consider $\alpha=\frac{1}{2}$, $K>0$ and $\varepsilon>0$. We will formally define the function $\mathcal{P}_\varepsilon[\mu]$ and the tangent vector field $\mathcal{V}_\varepsilon[\mu]$ along the manifold $\mathbb{T}\times \mathbb{R}$ given by
\begin{align*}
\mathcal{P}_\varepsilon[\mu](\theta,\Omega)&:=\Omega-K\int_{\mathbb{T}}\int_{\mathbb{R}}h_\varepsilon(\theta-\theta')\,d_{(\theta',\Omega')}\mu,\\
\mathcal{V}_\varepsilon[\mu](z,\Omega)&:=(\mathcal{P}_\varepsilon[\mu](z,\Omega)\,iz,0),
\end{align*}
where $\mu\in \mathcal{M}(\mathbb{T}\times \mathbb{R})$ is any finite Radon measure.
\end{defi}

In a similar way to the decompositions in Lemmas \ref{L-split-kernel} and \ref{L-split-kernel-regularized-critical} for the kernels of the subcritical and critical regimes, we need an appropriate split of the regularized kernels of the critical case in a consistent way so that we obtain a common one-sided-Lipschitz constant both for $h$ and $h_\varepsilon$. Notice that the standard Lipschitz constant of $h_\varepsilon$, obtained via a uniform bound of the first derivative $h_\varepsilon'$, should be avoided as it blows up when $\varepsilon\rightarrow 0$.

\begin{lem}\label{L-split-kernel-regularized-critical}
Consider $\alpha=\frac{1}{2}$, $\varepsilon>0$ and set $\bar{h}_\varepsilon$ and $\tilde{\theta}_\varepsilon\in (0,\frac{\pi}{2})$ such that
$$\bar{h}_\varepsilon:=\max_{0\leq \theta\leq \pi}h_\varepsilon(\theta)\ \mbox{ and }\ \tilde{\theta}_\varepsilon\tan(\tilde{\theta}_\varepsilon)-\tilde{\theta}_\varepsilon^2=\varepsilon^2.$$
Define the couple of functions $\Delta_\varepsilon,\Lambda_\varepsilon:[-2\pi,2\pi]\longrightarrow\mathbb{R}$ as follows
$$
\Delta_\varepsilon(\theta):=\left\{
\begin{array}{ll}
2\bar{h}_\varepsilon-h_\varepsilon(\theta), & \theta\in [-2\pi,-2\pi+\tilde{\theta}_\varepsilon),\\
\bar h_\varepsilon, & \theta\in [-2\pi+\tilde{\theta}_\varepsilon,-\tilde{\theta}_\varepsilon),\\
-h_\varepsilon(\theta), & \theta\in [-\tilde{\theta}_\varepsilon,\tilde{\theta}_\varepsilon],\\
-\bar{h}_\varepsilon, & \theta\in (\tilde{\theta}_\varepsilon,2\pi-\tilde{\theta}_\varepsilon],\\
-h_\varepsilon(\theta)-2\bar{h}_\varepsilon, & \theta\in (2\pi-\tilde{\theta}_\varepsilon,2\pi],
\end{array}
\right.
$$
$$
\Lambda_\varepsilon(\theta):=
\left\{
\begin{array}{ll}
-2\bar{h}_\varepsilon, & \theta\in [-2\pi,-2\pi+\tilde{\theta}_\varepsilon),\\
-\bar{h}_\varepsilon-h_\varepsilon(\theta), & \theta\in [-2\pi+\tilde{\theta}_\varepsilon,-\tilde{\theta}_\varepsilon),\\
0, & \theta\in [-\tilde{\theta}_\varepsilon,\tilde{\theta}_\varepsilon],\\
\bar{h}_\varepsilon-h_\varepsilon(\theta), & \theta\in (\tilde{\theta}_\varepsilon,2\pi-\tilde{\theta}_\varepsilon],\\
2\bar{h}_\varepsilon, & \theta\in (2\pi-\tilde{\theta}_\varepsilon,2\pi].
\end{array}
\right.
$$
Then, there exists $0<L_\varepsilon\leq -\inf_{\theta\in (0,\pi)}h'(\theta)$ such that following properties hold true
\begin{enumerate}
\item $\Delta_\varepsilon$ is monotonically decreasing, $\Lambda_\varepsilon$ is Lipschitz-continuous with constant $L_\varepsilon$ and
$$-h_\varepsilon(\theta)=\Delta_\varepsilon(\theta)+\Lambda_\varepsilon(\theta),\ \forall\theta\in [-2\pi,2\pi].$$
\item $-h_\varepsilon$ is one-sided Lipschitz in $[-2\pi,2\pi]$ with same constant $L_\varepsilon$, namely,
$$\left((-h_\varepsilon)(\theta_1)-(-h_\varepsilon)(\theta_2)\right)(\theta_1-\theta_2)\leq L_\varepsilon(\theta_1-\theta_2)^2.$$
\end{enumerate}
\end{lem}

\begin{proof}
Everything is clear except, at most, the estimate of the Lipschitz constant of $\Lambda_\varepsilon$. Since such function is piecewise smooth, it is enough to compute the Lipschitz constant on any of the pieces of its domain. We will only focus on the interval $(\tilde{\theta}_\varepsilon,2\pi-\tilde{\theta}_\varepsilon)$ where the functions is not constant. In the other non-constant piece the estimate follows from similar arguments. Let us compute
$$\Lambda_\varepsilon'(\theta)=\frac{d}{d\theta}\left(\bar h_\varepsilon-h_\varepsilon(\theta)\right)=-h_\varepsilon'(\theta)=\frac{1}{(\varepsilon^2+\vert \theta\vert_o^2)^{1/2}}\left[\frac{\vert \theta\vert_o\sin\vert \theta\vert_o}{\varepsilon^2+\vert \theta\vert_o^2}-\cos\theta\right],$$
for every $\theta\in (\tilde{\theta}_\varepsilon,2\pi-\tilde{\theta}_\varepsilon)$. Since $\Lambda_\varepsilon$ is increasing in the whole interval $(\tilde{\theta}_\varepsilon,2\pi-\tilde{\theta}_\varepsilon)$ as a consequence of the definition of $\tilde{\theta}_\varepsilon$, then $g_\varepsilon'$ is non-negative along it and we conclude that
\begin{multline*}
\sup_{\theta\in (\tilde{\theta}_\varepsilon,2\pi-\tilde{\theta}_\varepsilon)}\vert \Lambda_\varepsilon'(\theta)\vert=\sup_{\theta\in (\tilde{\theta}_\varepsilon,2\pi-\tilde{\theta}_\varepsilon)}\Lambda_\varepsilon'(\theta)=\sup_{\theta\in (\tilde{\theta}_\varepsilon,2\pi-\tilde{\theta}_\varepsilon)}\frac{1}{(\varepsilon^2+\vert \theta\vert_o^2)^{1/2}}\left[\frac{\vert \theta\vert_o\sin\vert \theta\vert_o}{\varepsilon^2+\vert \theta\vert_o^2}-\cos\theta\right]\\
\leq \sup_{\theta\in (\tilde{\theta}_\varepsilon,2\pi-\tilde{\theta}_\varepsilon)}\frac{1}{\vert \theta\vert_o}\left[\frac{\vert \theta\vert_o\sin\vert \theta\vert_o}{\vert \theta\vert_o^2}-\cos\theta\right]=\sup_{\theta\in (\tilde{\theta}_\varepsilon,2\pi-\tilde{\theta}_\varepsilon)}(-h'(\theta))=-\inf_{\theta\in (\tilde{\theta}_\varepsilon,2\pi-\tilde{\theta}_\varepsilon)}h'(\theta).
\end{multline*}
The proof then follows from the mean value theorem.
\end{proof}

\begin{theo}\label{T-transport-field-convergence-regularized-critical}
Consider $\alpha=\frac{1}{2}$, $K>0$ and let $\{\mu^\varepsilon\}_{\varepsilon>0}$ and $\mu$ be in $\mathcal{P}(\mathbb{T}\times \mathbb{R})$ so that
$$\mu^\varepsilon\rightarrow \mu\ \mbox{ in }\ \mathcal{P}(\mathbb{T}\times \mathbb{R})-\mbox{narrow}.$$
Then, the following convergence takes place
$$\lim_{\varepsilon\rightarrow 0}\sup_{\Omega\in \mathbb{R}}\vert \mathcal{V}_\varepsilon[\mu^\varepsilon](z,\Omega)-\mathcal{V}[\mu](z,\Omega)\vert=0,$$
for each continuity point $z\in \mathbb{T}$ of the marginal measure $(\pi_z)_{\#}\mu$. In particular, it happens a.e. in $\mathbb{T}$.
\end{theo}

\begin{proof}
Let us consider 
$$F_\varepsilon(\theta):=\vert \mathcal{V}_\varepsilon[\mu^\varepsilon](\theta,\Omega)-\mathcal{V}[\mu](\theta,\Omega)\vert=\left\vert K\int_{(-\pi,\pi]\setminus\{\theta\}}h_\varepsilon(\theta-\theta')\,d_{\theta'}\rho_\varepsilon-K\int_{(-\pi,\pi]\setminus\{\theta\}}h(\theta-\theta')\,d_{\theta'}\rho\right\vert,$$
for every $\theta\in (-\pi,\pi]$ and $\varepsilon>0$, where $\rho_\varepsilon=(\pi_z)_{\#}\mu_\varepsilon$ and $\rho=(\pi_z)_{\#}\mu$. Let us consider the following spit 
$$F_\varepsilon(\theta)\leq F^1_\varepsilon(\theta)+F_{\varepsilon,\delta}^2(\theta)+F_{\varepsilon,\delta}^3(\theta),$$
for every $\delta>\varepsilon^{1/2}$, where each term looks like
\begin{align*}
F_{\varepsilon}^1(\theta)&:=\left\vert K\int_{(-\pi,\pi]\setminus \{\theta\}}h(\theta-\theta')\,d_{\theta'}(\rho_\varepsilon-\rho)\right\vert,\\
F_{\varepsilon,\delta}^2(\theta)&:=\left\vert K\int_{(-\pi,\pi]\setminus\{\theta\}}\chi_{\vert \theta-\theta'\vert_o\geq \delta}(h_\varepsilon(\theta-\theta')-h(\theta-\theta'))\,d\rho_\varepsilon\right\vert,\\
F_{\varepsilon,\delta}^3(\theta)&:=\left\vert K\int_{(-\pi,\pi]\setminus\{\theta\}}\chi_{\vert \theta-\theta'\vert<\delta} (h_\varepsilon(\theta-\theta')-h(\theta-\theta'))\,d\rho_\varepsilon \right\vert.
\end{align*}
Let us fix $\theta\in (-\pi,\pi]$ any continuity point of $\rho$. On the one hand, Theorem \ref{T-transport-field-convergence-critical} implies that
$$\lim_{\varepsilon\rightarrow 0}F_\varepsilon^1(\theta)=0.$$
Second, the estimate \eqref{E-52} in the proof of such result entails
$$F_{\varepsilon,\delta}^2(\theta)\leq \frac{K}{2}\frac{\varepsilon}{\delta}\leq \frac{K\varepsilon^{1/2}}{2},$$
for every $\delta>\varepsilon^{1/2}$. Then, taking limits $\varepsilon\rightarrow 0$ we obtain
$$\lim_{\varepsilon\rightarrow 0}F_{\varepsilon,\delta}^2(\theta)=0.$$
To deal with the third term, let us consider any cut-off function $\xi\in C_c([0,+\infty))$ like in \eqref{E-scaled-cut-off}. Then,
\begin{align*}
F_{\varepsilon,\delta}^3(\theta)&\leq K\int_{(-\pi,\pi]\setminus\{\theta\}}\xi\left(\frac{\vert \theta-\theta'\vert_o}{\delta}\right)\,d\rho_\varepsilon\\
&=K\int_{(-\pi,\pi]\setminus\{\theta\}}\xi\left(\frac{\vert \theta-\theta'\vert_o}{\delta}\right)\,d(\rho_\varepsilon-\rho)+K\int_{(-\pi,\pi]\setminus\{\theta\}}\xi\left(\frac{\vert \theta-\theta'\vert_o}{\delta}\right)\,d\rho\\
&\leq K\int_{(-\pi,\pi]\setminus\{\theta\}}\xi\left(\frac{\vert \theta-\theta'\vert_o}{\delta}\right)\,d(\rho_\varepsilon-\rho)+K\int_{(-\pi,\pi]\setminus\{\theta\}}\chi_{\vert \theta-\theta'\vert_o\geq 2\delta}\,d\rho.
\end{align*}
Taking limits yields
$$\limsup_{\varepsilon\rightarrow 0}F_{\varepsilon,\delta}^3(\theta)\leq K\int_{(-\pi,\pi]\setminus\{\theta\}}\chi_{\vert \theta-\theta'\vert_o\geq 2\delta}\,d\rho.$$
Putting everything together, we obtain
$$\limsup_{\varepsilon\rightarrow 0}F_\varepsilon(\theta)\leq K\int_{(-\pi,\pi]\setminus\{\theta\}}\chi_{\vert \theta-\theta'\vert_o\geq 2\delta}\,d\rho,$$
for any arbitrary $\delta>0$. Since $\theta$ is a continuity point of $\rho$ we conclude the desired result by taking limits $\delta\rightarrow 0$.
\end{proof}

\begin{theo}\label{T-growth-L(RK)-critical}
Consider $\alpha=\frac{1}{2}$, $K>0$ and let $f^1,f^2\in \mathcal{AC}_\mathcal{M}\cap \mathcal{T}_\mathcal{M}$ be solutions in the sense of the Filippov flow to \eqref{E-kuramoto-transport-TxR} with initial data $f^1_0,f^2_0\in \mathcal{P}(\mathbb{T}\times \mathbb{R})$ according to Theorem \ref{T-mean-field-existence-critical}. Let us set their distributions of natural frequencies $g^i=\left(\pi_{\Omega}\right)_{\#}f^i_0$. If  $g:=g^1\equiv g^2$, then
$$W_{2,g}(f^1_t,f^2_t)\leq W_{2,g}(f^1_0,f^2_0)e^{2KL_0t},$$
for every $t\geq 0$, where $L_0$ is the one-sided Lipschitz constant of $-h$ in Lemma \ref{L-split-kernel-critical} and $W_{2,g}$ is the fibered quadratic Wasserstein distance in Proposition \ref{P-L1(RK)-distance}.
\end{theo}

\begin{proof}
Using similar arguments to those in Section \ref{S-weak-solutions-existence} for the Lipschitz-continuous regularized kernel $h_\varepsilon$, we can construct global classical solutions $f^{1,\varepsilon},f^{2,\varepsilon}\in\mathcal{AC}_\mathcal{M}\cap \mathcal{T}_\mathcal{M}$ to the regularized systems
$$\left\{\begin{array}{l}
\displaystyle\frac{\partial f^{1,\varepsilon}}{\partial t}+\divop(\mathcal{V}_\varepsilon[f^{1,\varepsilon}]f^{1,\varepsilon})=0,\\
f^{1,\varepsilon}_0=f^1_0,
\end{array}\right. \ \mbox{ and }\  \left\{\begin{array}{l}
\displaystyle\frac{\partial f^{2,\varepsilon}}{\partial t}+\divop(\mathcal{V}_\varepsilon[f^{2,\varepsilon}]f^{2,\varepsilon})=0,\\
f^{2,\varepsilon}_0=f^2_0.
\end{array}\right.$$
For $g$-a.e. $\Omega\in\mathbb{R}$ fixed, let us consider the corresponding term of the family of disintegrations at the initial time, i.e., $f^1_0(\cdot\vert \Omega)$ and $f^2_0(\cdot\vert \Omega)$. Set an optimal transference plan from the former probability measure in $\mathbb{T}$ to the latter one, i.e.,
$$\mu_{0,\Omega}\in\Pi(f^1_0(\cdot \vert \Omega),f^2_0(\cdot \vert \Omega)):=\left\{\mu\in\mathcal{P}(\mathbb{T}\times \mathbb{T}):\,(\pi_1)_{\#}\mu=f^1_0(\cdot\vert \Omega)\ \mbox{and}\ (\pi_2)_{\#}\mu=f^2_0(\cdot\vert \Omega)\right\},$$
so that the $2$-Wasserstein distance is attained
$$W_2(f^1_0(\cdot\vert\Omega),f^2_0(\cdot\vert\Omega))^2=\int_{\mathbb{T}}\int_{\mathbb{T}}d(z_1,z_2)^2d_{(z_1,z_2)}\mu_{0,\Omega}.$$
Again, we are denoting the projections $\pi_1(z,z')=z$ and $\pi_2(z,z')=z'$. We then can define the following competitor transference plans at time $t$
\begin{align*}
\mu_{t,\Omega}&:=(Z_{f^1}(t;0,\cdot,\Omega)\otimes Z_{f^2}(t;0,\cdot,\Omega))_{\#}{\mu_{0,\Omega}}\in \mathcal{P}(\mathbb{T}\times \mathbb{T}),\\
\mu_{t,\Omega}^\varepsilon&:=(Z_{f^{1,\varepsilon}}(t;0,\cdot,\Omega)\otimes Z_{f^{2,\varepsilon}}(t;0,\cdot,\Omega))_{\#}{\mu_{0,\Omega}}\in \mathcal{P}(\mathbb{T}\times \mathbb{T}),
\end{align*}
where $X_{f^i}(t;0,z,\Omega)=(Z_{f^i}(t;0,z,\Omega),\Omega)$ is the Filippov flow associated with the transport field $\mathcal{V}[f^i]$ according to Theorem \ref{T-transport-field-sided-Lipschitz-critical} and $X_{f^{i,\varepsilon}}(t;0,z,\Omega)=(Z_{f^{i,\varepsilon}}(t;0,z,\Omega),\Omega)$ is the classical flow of $\mathcal{V}_\varepsilon[f^{i,\varepsilon}]$ according to Theorem \ref{T-transport-field-sided-Lipschitz}. Notice that $\mu_{t,\Omega},\mu_{t,\Omega}^\varepsilon\in\Pi(f^1_t(\cdot\vert\Omega),f^2_t(\cdot\vert\Omega))$. Hence,
\begin{align*}
\frac{1}{2}W_2(f^1_t(\cdot\vert\Omega),f^2_t(\cdot\vert\Omega))^2&\leq \int_{\mathbb{T}}\int_{\mathbb{T}}\frac{1}{2}d(z_1,z_2)^2\,d_{(z_1,z_2)}\mu_{t,\Omega}\\
&=\int_{\mathbb{T}}\int_{\mathbb{T}}\frac{1}{2}d(Z_{f^1}(t;0,z_1,\Omega),Z_{f^2}(t;0,z_2,\Omega))^2\,d_{(z_1,z_2)}\mu_{0,\Omega}.
\end{align*}
Integrating the above inequality against $g$ yields
$$\frac{1}{2}W_{2,g}(f^1_t,f^2_t)^2\leq \int_{\mathbb{R}}\int_{\mathbb{T}}\int_{\mathbb{T}}\frac{1}{2}d(Z_{f^1}(t;0,z_1,\Omega),Z_{f^2}(t;0,z_2,\Omega))^2\,d_{(z_1,z_2)}\mu_{0,\Omega}\,d_\Omega g=:I(t).$$
Let us define
$$I_\varepsilon(t):=\int_{\mathbb{R}}\int_{\mathbb{T}}\int_{\mathbb{T}}\frac{1}{2}d(Z_{f^{1,\varepsilon}}(t;0,z_1,\Omega),Z_{f^{2,\varepsilon}}(t;0,z_2,\Omega))^2\,d_{(z_1,z_2)}\mu_{0,\Omega}\,d_\Omega g.$$
Since the kernel $h_\varepsilon$ is globally Lipschitz-continuous, we can mimic the proof of Theorem \ref{T-growth-L(RK)} and obtain
$$I_\varepsilon(t)\leq I_\varepsilon(0)e^{4KL_\varepsilon t},$$
for every $t\geq 0$, where $L_\varepsilon$ is the one sided-Lipschitz constant of $-h_\varepsilon$ in Lemma \ref{L-split-kernel-regularized-critical}. On the one hand, $I_\varepsilon(0)=I(0)=\frac{1}{2}W_{2,g}(f^1_0,f^2_0)$. On the other hand, recall that
$$L_\varepsilon\leq -\inf_{\theta\in (0,\pi)}h'(\theta)=L_0,$$
for every $\varepsilon>0$. Consequently,
$$
I_\varepsilon(t)\leq \frac{1}{2}W_{2,g}(f^1_0,f^2_0)^2e^{4KL_0 t},
$$
for every $t\geq 0$ and $\varepsilon>0$. Our last goal is to show that $\lim_{\varepsilon\rightarrow 0}I_\varepsilon(t)=I(t)$. Consider the scaled cut-off functions $\xi_=\xi_R(\Omega)$ in \eqref{E-scaled-cut-off} and define the decomposition
$$I_\varepsilon(t)-I(t)=A_{\varepsilon}^R(t)+B_\varepsilon^R(t),$$
where both terms are given by the formulas
\begin{align*}
A_\varepsilon^R(t)&:=\int_{\mathbb{R}}\int_{\mathbb{T}}\int_{\mathbb{T}}\xi_R(\Omega)\frac{1}{2}\left(d(Z_{f^{1,\varepsilon}}(t;0,z_1,\Omega),Z_{f^{2,\varepsilon}}(t;0,z_2,\Omega))^2\right.\\
&\hspace{4.5cm}\left.-d(Z_{f^{1}}(t;0,z_1,\Omega),Z_{f^{2}}(t;0,z_2,\Omega))^2\right)\,d_{(z_1,z_2)}\mu_{0,\Omega}\,d_\Omega g,\\
B_\varepsilon^R(t)&:=\int_{\mathbb{R}}\int_{\mathbb{T}}\int_{\mathbb{T}}(1-\xi_R(\Omega))\frac{1}{2}\left(d(Z_{f^{1,\varepsilon}}(t;0,z_1,\Omega),Z_{f^{2,\varepsilon}}(t;0,z_2,\Omega))^2\right.\\
&\hspace{4.5cm}\left.-d(Z_{f^{1}}(t;0,z_1,\Omega),Z_{f^{2}}(t;0,z_2,\Omega))^2\right)\,d_{(z_1,z_2)}\mu_{0,\Omega}\,d_\Omega g.
\end{align*}
Since the vector fields $\mathcal{V}_\varepsilon[f^{i,\varepsilon}]$ and $\mathcal{V}[f^i]$ are all essentially uniformly bounded, then Theorem \ref{T-transport-field-convergence-regularized-critical}, the Alaoglu--Bourbaki theorem and a standard application of Banach--Saks' theorem show that
$$\vert \mathcal{V}_\varepsilon[f^{i,\varepsilon}]-\mathcal{V}[f^i]\vert \overset{*}{\rightharpoonup}0\ \mbox{ in }\ L^\infty((0,+\infty)\times \mathbb{T}\times \mathbb{R}),$$
for every $i=1,2$. Then, Lemma \ref{L-convergence-Filippov-flow} implies that
$$X_{f^{i,\varepsilon}}\rightarrow X_{f^i}\ \mbox{ in }\ C_{\textnormal{loc}}([0,+\infty)\times \mathbb{T}\times \mathbb{R}),$$
for every $i=1,2$. Since the squared distance is at least uniformly continuous, then we claim that
$$\lim_{\varepsilon\rightarrow 0}A_\varepsilon^R(t)=0,$$
for every $R>0$ and each $t\geq 0$. On the other hand,
$$\vert B_\varepsilon^R(t)\vert\leq \pi^2\Vert \chi_{\vert \Omega\vert \geq 2R}g\Vert_{\mathcal{M}(\mathbb{R})},$$
for every $R>0$ and each $t\geq 0$. Putting everything together yields
$$\limsup_{\varepsilon\rightarrow 0}\vert I_\varepsilon(t)-I(t)\vert\leq\pi^2\Vert \chi_{\vert \Omega\vert \geq 2R}g\Vert_{\mathcal{M}(\mathbb{R})}.$$
Taking limits $R\rightarrow+\infty$ concludes the proof by tightness of $g$.
\end{proof}

\begin{cor}\label{C-uniqueness-flow-solutions-critical}
Consider $\alpha=\frac{1}{2}$, $K>0$ and let $f^1,f^2\in \mathcal{AC}_\mathcal{M}\cap \mathcal{T}_\mathcal{M}$ be solutions in the sense of the Filippov flow to the non-linear transport equation \eqref{E-kuramoto-transport-TxR} with initial data $f^1_0,f^2_0\in \mathcal{P}(\mathbb{T}\times \mathbb{R})$. If $f^1_0=f^2_0$, then
$$f^1_t=f^2_t,\ \mbox{ for every }\ t\geq 0.$$
\end{cor}

\subsection{Mean field limit towards solutions in the sense of the Filippov flow}\label{SS-critical-regime-mean-field}
This part is devoted to adapt the bound of the quadratic Wasserstein distance in \ref{T-growth-RK} of Subsection \ref{SS-mean-field} to the critical case. As a byproduct, we will recover a quantitative version of the mean field limit towards solutions in the sense of the Filippov flow.
 
\begin{theo}\label{T-growth-RK-critical}
Consider $\alpha=\frac{1}{2}$, $K>0$ and let $f^1,f^2\in \mathcal{AC}_\mathcal{M}\cap \mathcal{T}_\mathcal{M}$ be solutions in the sense of the Filippov flow to the non-linear transport equation \eqref{E-kuramoto-transport-TxR} with initial data $f^1_0,f^2_0\in \mathcal{P}_2(\mathbb{T}\times \mathbb{R})$. Then, there exists $C=C(\alpha,K,f^1_0,f^2_0)>0$ such that
$$W_2(f^1_t,f^2_t)\leq e^{\left(\frac{1}{2}+2KL_0\right)t}W_2(f^1_0,f^2_0),$$
for every $t\geq 0$, where $L_0$ is the one-sided Lipschitz constant of $-h$ in Lemma \ref{L-split-kernel-critical}.
\end{theo}

The proof follows a similar approximation argument like in Theorem \ref{T-growth-L(RK)-critical}. Then, we omit it.

\begin{cor}\label{C-growth-RK-mean-field-critical}
Consider $\alpha=\frac{1}{2}$, $K>0$ and let $f\in \mathcal{AC}_\mathcal{M}\cap \mathcal{T}_\mathcal{M}$ be the unique solution in the sense of the Filippov flow to \eqref{E-kuramoto-transport-TxR} with initial datum $f_0\in \mathcal{P}_2(\mathbb{T}\times \mathbb{R})$. Consider $N$ oscillators with initial phases and natural frequencies given by the configurations
$$\Theta_0^N=(\theta_{1,0}^N,\ldots,\theta_{N,0}^N)\ \mbox{ and }\ \{\Omega_i^N:\,i=1,\ldots,N\},$$ 
for every $N\in\mathbb{N}$. Let $\Theta^N(t):=(\theta_1^N(t),\ldots,\theta_N^N(t))$ be the unique global-in-time Filippov solution to the discrete singular Kuramoto model according to \cite[Theorem 3.1]{P-P-S} and define the associated empirical measures in $\mathbb{T}\times \mathbb{R}$
$$\mu^N_t:=\frac{1}{N}\sum_{i=1}^N\delta_{z_i^N(t)}(z)\otimes \delta_{\Omega_i^N}(\Omega),$$
where $z_i^N(t):=e^{i\theta_i^N(t)}$. If $\lim_{N\rightarrow\infty}W_2(\mu^N_0,f_0)=0$, then
$$\lim_{N\rightarrow\infty}\sup_{t\in [0,T]}W_2(\mu^N_t,f_t)=0,\ \mbox{ for all }\ T>0.$$
\end{cor}

\subsection{Global phase-synchronization of identical oscillators in finite time}\label{SS-critical-regime-synchro}
Recall that in the agent-based system \eqref{E-kuramoto-discrete} enjoys finite-time global phase synchronization of identical oscillators. Such result was derived in \cite[Theorem 5.5]{P-P-S}.

\begin{theo}\label{T-synchro-discrete-critical}
Let $\Theta=(\theta_1, \cdots, \theta_N)$ be the Filippov solution to \eqref{E-kuramoto-discrete}-\eqref{E-kuramoto-kernel} with $\alpha =\frac{1}{2}$ for identical oscillators $(\Omega_i = 0)$, for $i=1, \ldots, N$. Assume that the initial configuration $\Theta_0$ is confined in a half circle, i.e., $0<D(\Theta_0) < \pi$. Then, there is complete phase synchronization at a finite time not larger than $T_c$, where
$$
T_c=\frac{D(\Theta_0)}{Kh(D(\Theta_0))}.
$$
\end{theo}

Mimicking the ideas in Theorem \ref{T-convergence-equilibrium-continuum-subcritical}, we obtain the following analogue.

\begin{theo}\label{T-convergence-equilibrium-continuum-critical}
Consider any initial datum $f_0\in \mathcal{P}(\mathbb{T}\times \mathbb{R})$ with identical distribution of natural frequencies, namely, $g=(\pi_\Omega)_{\#}f_0=\delta_0(\Omega)$, where $\pi_\Omega$ is the projection \eqref{E-projections}. Let $f=f_t$ be the unique global-in-time measure-valued solution in the sense of the Filippov flow to \eqref{E-kuramoto-transport-TxR} issued at $f_0$ with $\alpha =\frac{1}{2}$ and assume $D_0:=\diam (\supp\rho_0)<\pi$. Then,
$$f_t=f_\infty\ \mbox{ for all }\ t\geq T_c,$$
where $T_c=\frac{D_0}{Kh(D_0)}$ and the equilibrium $f_\infty$ is given by the monopole $f_\infty:=\delta_{Z_{av}(0)}(z)\otimes \delta_0(\Omega)$.
\end{theo}

\section{The supercritical regime}\label{S-supercritical-regime}
This part is devoted to the derivation of weak measure-valued solutions via a different technique. Namely, we will explore a singular hyperbolic limit of vanishing inertia type on a kinetic second order regularized system. Once the regularized model with inertia is presented, we will introduce an appropriate scaling where the inertia term is neglected and singularization of the weights emerges as the scaling parameter $\varepsilon$ tends to zero. As a consequence of the rigorous hydrodynamic limit $\varepsilon\rightarrow 0$, we obtain weak measure-valued solutions of the macroscopic singular system in the supercritical regime $\alpha\in (\frac{1}{2},1)$ for identical oscillators, i.e., $g=\delta_0$. Of course, the idea also works in the most regular regime $\alpha\in (0,\frac{1}{2})$ even for non-identical oscillators, that is $g\neq \delta_0$, thus recovering the above weak measure-valued solutions in Theorems \ref{T-mean-field-existence-subcritical}.

\subsection{Second order regularized system}\label{SS-second-order-regularized}
This part is structured as follows. First, we will introduce the agent-based second order Kuramoto model with inertia endowed with regular weighted coupling, frequency damping and noise. Second, we will recall the derivation of the Vlasov--McKean kinetic equation associated with the second order stochastic regular system.

\subsubsection{The second order agent-based system}\label{SSS-second-order-agent-based}
Let us consider the following scaled second order stochastic system for the dynamics of the $N$ oscillators under the effect of inertia, frequency damping and noise:
\begin{equation}\label{E-second-order-kuramoto-discrete}
\left\{
\begin{array}{l}
\displaystyle d\theta_i=\omega_i dt,\\
\displaystyle \varepsilon d\omega_i=\Omega_idt+\frac{\nu}{N}\sum_{j=1}^N h_\varepsilon(\theta_j-\theta_i)dt-\omega_idt+\sqrt{2\varepsilon}\,dW^i_t,\\
\theta_i(0)=\theta_{i,0},\ \omega_i(0)=\omega_{i,0}.
\end{array}
\right.
\end{equation}
for $i=1,\ldots,N$. Here, $W^i=W^i_t$ are independent Brownian motions. Again, $\theta_i=\theta_i(t)$ are phase values of the signals whilst $\omega_i=\omega_i(t)$ yields the evolution of their frequencies, to be distinguished from the static natural frequencies $\Omega_i$. We have introduced an inertia term modulated by the ``inertia parameter'' $\varepsilon>0$ that makes the transient to the original first-order dynamics faster as $\varepsilon\rightarrow 0$. In turns, the noise disappears and singularity in the coupling functions $h_\varepsilon$ emerges as $\varepsilon\rightarrow 0$. Hence, we formally recover the singular Kuramoto model as the reduced first order dynamics when $\varepsilon=0$. Introducing inertia is not an artificial method and one can indeed find a large literature concerning the original second order Kuramoto model with inertia as a suitable model of synchronization of coupled oscillators, see \cite{C-H-L-X-Y,C-H-M,C-H-Y-1,C-H-Y-2} and references therein.

\subsubsection{The second order Vlasov-McKean equation}\label{SSS-second-order-scaling}
As mentioned in Subsection \ref{SS-formal-mean-field-limit}, classical mean field and propagation of chaos methods \cite{C-C-H-2,H-J,H-M,J,J-W,J-W-2,Ka,M-M,M-M-W,S} allow deriving the kinetic equation as $N\rightarrow\infty$ associated with the stochastic agent based model \eqref{E-second-order-kuramoto-discrete}:
\begin{equation}\label{E-second-order-kuramoto-kinetic-pre}
\frac{\partial P^\varepsilon}{\partial t}+\omega\frac{\partial P^\varepsilon}{\partial \theta}=\frac{1}{\varepsilon}\frac{\partial}{\partial \omega}\left((\omega-\Omega)P^\varepsilon+\frac{\partial P^\varepsilon}{\partial \omega}+K\int_{(-\pi,\pi]\times \mathbb{R}\times\mathbb{R}}h_\varepsilon(\theta-\theta')P^\varepsilon(\theta',\omega',\Omega') P^\varepsilon\,d\theta'\,d\omega'\,d\Omega'\right),
\end{equation}
for every $t\geq 0$, $\theta\in (-\pi,\pi]$, $\omega\in\RR$ and $\Omega\in\RR$. Here, $P^\varepsilon(t,\theta,\omega,\Omega)$ describes the probability distribution of finding an oscillator at time $t$ with phase $\theta$, frequency $\omega$ and natural frequency $\Omega$, respectively. Again, we endow \eqref{E-second-order-kuramoto-kinetic-pre} with periodic boundary conditions
\begin{equation}\label{E-second-order-kuramoto-kinetic-periodicity}
P^\varepsilon(t,-\pi,\omega,\Omega)=P^\varepsilon(t,\pi,\omega,\Omega).
\end{equation}
Using Appendix \ref{appendix-relation-measure-spaces} to identify $\theta\in (-\pi,\pi]$ with $z\in\mathbb{T}$ via the formula $z=e^{i\theta}$. An important fact to be remarked is that the interaction term in \eqref{E-second-order-kuramoto-kinetic-pre} can be simplified. Specifically, consider the associated macroscopic quantities
\begin{align}\label{E-hierarchy-f-rho-g}
\begin{split}
f^\varepsilon(t,\theta,\Omega)&:=\int_{\RR} P^\varepsilon\,d\omega,\\
\rho^\varepsilon(t,\theta)&:=\int_{\RR^2}P^\varepsilon\,d\omega\,d\Omega=\int_{\RR}f^\varepsilon\,d\Omega,\\
g(\Omega)&:=\int_{\mathbb{T}\times \mathbb{R}} P^\varepsilon\,d\theta\,d\omega=\int_{\mathbb{T}} f^\varepsilon\,d\theta,
\end{split}
\end{align}
and note that such term is the following convolution
$$\int_{(-\pi,\pi]\times \mathbb{R}\times \mathbb{R}}h_\varepsilon(\theta-\theta')P^\varepsilon(t,\theta',\omega',\Omega')\,d\theta'\,d\omega'\,d\Omega'=(h_\varepsilon*\rho^\varepsilon)(t,\theta,\omega,\Omega),$$
where the convolution is considered as periodic functions, that is, as functions in $\mathbb{T}$. Then, \eqref{E-second-order-kuramoto-kinetic-pre}-\eqref{E-second-order-kuramoto-kinetic-periodicity} can be restated as follows
\begin{equation}\label{E-second-order-kuramoto-kinetic}
\frac{\partial P^\varepsilon}{\partial t}+\omega\frac{\partial P^\varepsilon}{\partial \theta}=\frac{1}{\varepsilon}\frac{\partial}{\partial \omega}\left((\omega-\Omega)P^\varepsilon+\frac{\partial P^\varepsilon}{\partial \omega}+K(h_\varepsilon*\rho^\varepsilon) P^\varepsilon\right).
\end{equation}

Let us now set the desired hypothesis on the sequence of initial data $f_\varepsilon^0$ and the distribution $g_\varepsilon$ of natural frequencies. Depending on the degree of integrability that we wish to recover on the limiting distribution, i.e.,  $g\in L^p(\mathbb{R})$ or just for $g\in\mathcal{M}(\mathbb{R})$, we respectively assume that either
\begin{equation}\label{kura-kin-2nd-order-hypothesis-lp}\tag{$H_p$}
\left\{\begin{array}{l}
\displaystyle f_\varepsilon^0=f_\varepsilon^0(\theta,\omega,\Omega)\geq 0\ \mbox{ and }\ f_\varepsilon^0\in C^\infty_c(\mathbb{T}\times\RR\times \RR),\\ 
\displaystyle \Vert f_\varepsilon^0\Vert_{L^1(\mathbb{T}\times \RR\times \RR)}=1\ \mbox{ and }\ f_\varepsilon^0\overset{*}{\rightharpoonup}f^0\,\mbox{ in }\mathcal{M}(\mathbb{T}\times \RR),\\
\displaystyle \Vert f_\varepsilon^0\Vert_{L^p_\Omega(\mathbb{R};L^1_{(\theta,\omega)}(\mathbb{T}\times \mathbb{R}))}\leq C_0\ \mbox{ and }\ f_\varepsilon^0\overset{*}{\rightharpoonup}f^0\,\mbox{ in }\,L^p(\RR,\mathcal{M}(\mathbb{T})),\\
\Vert\Omega^2 g_\varepsilon\Vert_{L^1(\RR)}\leq V_0\ \mbox{ and }\ \frac{1}{2}\Vert \omega^2 f_\varepsilon^0\Vert_{L^1(\mathbb{T}\times \RR)}\leq E_0,
\end{array}\right.
\end{equation}
fulfil for some $1<p<\infty$, or the limiting assumptions with $p\rightarrow 1$ hold true, i.e.,
\begin{equation}\label{kura-kin-2nd-order-hypothesis}\tag{$H_1$}
\left\{\begin{array}{l}
\displaystyle f_\varepsilon^0=f_\varepsilon^0(\theta,\omega,\Omega)\geq 0\ \mbox{ and }\ f_\varepsilon^0\in C^\infty_c(\mathbb{T}\times\RR\times \RR),\\ 
\displaystyle \Vert f_\varepsilon^0\Vert_{L^1(\mathbb{T}\times \RR\times \RR)}=1\ \mbox{ and }\ f_\varepsilon^0\overset{*}{\rightharpoonup}f^0\,\mbox{ in }\,\mathcal{M}(\mathbb{T}\times \RR),\\
\Vert\Omega^2 g_\varepsilon\Vert_{L^1(\RR)}\leq V_0\ \mbox{ and }\ \frac{1}{2}\Vert \omega^2 f_\varepsilon^0\Vert_{L^1(\mathbb{T}\times \RR)}\leq E_0.
\end{array}\right.
\end{equation}

\begin{rem}\label{R-kura-kin-2nd-order-hypothesis-g}
Note that both \eqref{kura-kin-2nd-order-hypothesis} and \eqref{kura-kin-2nd-order-hypothesis-lp} impose some properties on $g_\varepsilon$. Specifically, recall that $g_\varepsilon$ emerges as a marginal distribution of $f_\varepsilon^0$, namely,
$$g_\varepsilon(\Omega):=\int_{\mathbb{T}}f_\varepsilon^0(\theta,\Omega)\,d\theta.$$
Similarly, one can define the $\Omega$-marginal measure of $f^0\in \mathcal{P}(\mathbb{T}\times \RR)$, i.e., $g=(\pi_\Omega)_{\#}f^0$ defined by
$$\left<g,\varphi\right>:=\int_{\mathbb{T}}\int_{\RR}\varphi(\Omega)(df^0)(\theta,\Omega),$$
for every test function $\varphi\in C_0(\RR)$.

\begin{enumerate}
\item If \eqref{kura-kin-2nd-order-hypothesis} fulfils, then the following properties of $g_\varepsilon$ hold true
\begin{equation}\label{kura-kin-2nd-order-hypothesis-g}
\left\{
\begin{array}{l}
\displaystyle g_\varepsilon=g_\varepsilon(\Omega)\geq 0,\\
\displaystyle g_\varepsilon\in C^\infty_c(\Omega),\\
\displaystyle \Vert g_\varepsilon\Vert_{L^1(\RR)}=1,\\
\displaystyle g_\varepsilon\overset{*}{\rightharpoonup} g\ \mbox{ in }\ \mathcal{M}(\RR).
\end{array}
\right.
\end{equation}
\item Similarly, if in addition \eqref{kura-kin-2nd-order-hypothesis-lp} are verified, then not only do we recover \eqref{kura-kin-2nd-order-hypothesis-g} but also
\begin{equation}\label{kura-kin-2nd-order-hypothesis-g-lp}
\left\{
\begin{array}{l}
\Vert g_\varepsilon\Vert_{L^p(\RR)}\leq C_0,\\
\displaystyle g_\varepsilon\rightharpoonup g\ \mbox{ in }\ L^p(\RR).
\end{array}
\right.
\end{equation}
\end{enumerate}
\end{rem}

Notice that under the assumptions \eqref{kura-kin-2nd-order-hypothesis} or \eqref{kura-kin-2nd-order-hypothesis-lp}, for the above compactly supported initial data $f_\varepsilon^0$,  classical techniques assure the existence of a global-in-time classical solution $f_\varepsilon=f_\varepsilon(t,\theta,\omega,\Omega)$ to such system \eqref{E-second-order-kuramoto-kinetic}. Our main goal is to show that via the asymptotic limit $\varepsilon\searrow 0$ we can rigorously derive weak measure-valued solutions to the the macroscopic singularly weighted Kuramoto model of interest in the subcritical regime $\alpha\in (0,\frac{1}{2})$. The remaining parts are structured as follows. First, we introduce the hierarchy of frequency moments that for positive $\varepsilon$ is not a closed system, as usual. Second, we introduce some a priori bounds that allows passing to the limit in a weak sense, closing such hierarchy of frequency moments and obtaining weak measure-valued solutions to the singular Kuramoto model in the subcritical regime $\alpha\in (0,\frac{1}{2})$.

\subsection{A priori estimates}\label{SS-second-order-regularized-apriori}
Apart from \eqref{E-hierarchy-f-rho-g}, we will be concerned with the following set of frequency moments of the distribution function $P^\varepsilon=P^\varepsilon(t,\theta,\omega,\Omega)$
\begin{align*}
j^\varepsilon(t,\theta,\Omega)&:=\int_{\RR} \omega P^\varepsilon(t,\theta,\omega,\Omega)\,d\omega,\\
\mathcal{S}^\varepsilon(t,\theta,\Omega)&:=\int_{\RR}\omega^2 P^\varepsilon(t,\theta,\omega,\Omega)\,d\omega,\\
\mathcal{T}^\varepsilon(t,\theta,\Omega)&:=\int_{\RR}\omega^3 P^\varepsilon(t,\theta,\omega,\Omega)\,d\omega.
\end{align*}

The corresponding hierarchy of frequency moments can be easily derived from \eqref{E-second-order-kuramoto-kinetic} if one multiplies it by $1$, $\omega$ and $\omega^2$ and integrates with respect to $\omega$. The regularity of the global-in-time solution $P^\varepsilon$ along with the periodicity conditions with respect to $\theta$ yields the equations
\begin{align}
\frac{\partial f^\varepsilon}{\partial t}+\frac{\partial j^\varepsilon}{\partial \theta}&=0,\label{E-second-order-kuramoto-kinetic-moments-0}\\
\varepsilon\frac{\partial j^\varepsilon}{\partial t}+\varepsilon\frac{\partial \mathcal{S}^\varepsilon}{\partial \theta}+j^\varepsilon-\Omega f^\varepsilon+(h_\varepsilon*\rho^\varepsilon)f^\varepsilon&=0,\label{E-second-order-kuramoto-kinetic-moments-1}\\
\varepsilon\frac{\partial \mathcal{S}^\varepsilon}{\partial t}+\varepsilon\frac{\partial \mathcal{T}_\varepsilon}{\partial \theta}+2(\mathcal{S}^\varepsilon-\Omega j^\varepsilon)-2f^\varepsilon+2(h_\varepsilon*\rho^\varepsilon)j^\varepsilon&=0.\label{E-second-order-kuramoto-kinetic-moments-2}
\end{align}

Let us now focus on deriving some a priori estimates of the system \eqref{E-second-order-kuramoto-kinetic}. To such end, we first introduce the primitive function of the kernel, that will give some insight about the inter-particle interactions of our system.

\begin{defi}\label{D-h-scaled-primitive}
Let us define
$$
W_\varepsilon(\theta):=\int_0^{\theta} h_\varepsilon(\theta')\,d\theta'=\int_0^{\bar \theta} h_\varepsilon(\theta')\,d\theta',
$$
for every $\alpha\in (0,1)$ and $\varepsilon>0$. As usual, $\bar\theta$ denotes the representative modulo $2\pi$ of $\theta$ in $(-\pi,\pi]$.
\end{defi}

By definition, $W_\varepsilon$ enjoy nice regularity properties in the whole range of the parameter $\alpha\in (0,1)$ by virtue of the mild singularity of $h_\varepsilon$, although we will focus on $\alpha\in (\frac{1}{2},1)$ in this sections.

\begin{pro}\label{P-h-scaled-primitive}
The following properties hold true
\begin{enumerate}
\item $W_\varepsilon$ is $2\pi$-periodic.
\item If $\alpha\in \left(0,\frac{1}{2}\right)$ and $1\leq p<\frac{1}{2\alpha}$, there exists a positive constant $M_{\alpha,p}$ such that
$$\Vert W_\varepsilon\Vert_{W^{2,p}(\mathbb{T})}\leq M_{\alpha,p},\ \forall\varepsilon>0.$$
\item If $\alpha=\frac{1}{2}$, there exists a positive constant $M$ such that 
$$\Vert W_\varepsilon\Vert_{W^{1,\infty}(\mathbb{T})}\leq M,\ \forall\varepsilon>0.$$
\item If $\alpha\in \left(\frac{1}{2},1\right)$ and $1\leq p<\frac{1}{2\alpha-1}$, there exists a positive constant $M_{\alpha,p}$ such that
$$\Vert W_\varepsilon\Vert_{W^{1,p}(\mathbb{T})}\leq M_{\alpha,p},\ \forall\varepsilon>0.$$
\item $W_\varepsilon$ is a primitive function of $h_\varepsilon$.
\item $W_\varepsilon\geq 0$ and the identity only holds at $\theta\in 2\pi\mathbb{Z}$.
\end{enumerate}
\end{pro}

\begin{rem}\label{R-h-scaled-primitive}
Recall that the Sobolev embedding theorem entails the compact inclusion
$$W^{1,p}(\mathbb{T})\subset\subset C(\mathbb{T}),\ \forall\,p>1.$$
In particular, there is some constant $M_\alpha>0$ that does not depend on $\varepsilon>0$ such that
$$
\Vert W_\varepsilon\Vert_{C(\mathbb{T})}\leq M_\alpha,\ \forall\,\varepsilon>0.
$$
\end{rem}

\begin{lem}\label{L-estimates-kura-2nd-order}
Consider the strong solution $P^\varepsilon$ to \eqref{E-second-order-kuramoto-kinetic} whose initial data $P^\varepsilon_0$ fulfil the assumptions \eqref{kura-kin-2nd-order-hypothesis}. Then, the next formula holds true for every $\varepsilon>0$
$$\frac{d}{dt}\left(\varepsilon\int_{\mathbb{T}}\int_{\mathbb{R}}\mathcal{S}^\varepsilon\,d\theta\,d\Omega+K\int_{\mathbb{T}}(W_\varepsilon*\rho^\varepsilon)\rho^\varepsilon\,d\theta\right)+2\int_{\mathbb{T}}\int_{\mathbb{R}}\mathcal{S}^\varepsilon\,d\theta\,d\Omega=2\int_{\mathbb{T}}\int_{\mathbb{R}}\Omega j^\varepsilon\,d\theta\,d\Omega+2.$$
\end{lem}

\begin{proof}
Let us integrate \eqref{E-second-order-kuramoto-kinetic-moments-2} with respect to $\theta$ and $\Omega$ to obtain
$$\frac{d}{dt}\left(\varepsilon\int_{\mathbb{T}}\int_{\mathbb{R}}\mathcal{S}^\varepsilon\,d\theta\,d\Omega\right)+2\int_{\mathbb{T}}\int_{\mathbb{R}}(h_\varepsilon*\rho^\varepsilon)j^\varepsilon\,d\theta\,d\Omega+2\int_{\mathbb{T}}\int_{\mathbb{R}}\mathcal{S}^\varepsilon\,d\theta\,d\Omega=2\int_{\mathbb{T}}\int_{\mathbb{R}}\Omega j^\varepsilon\,d\theta\,d\Omega+2.$$
The only effort to be done is to identify the second term. To such end, we notice that $\frac{\partial W_\varepsilon}{\partial \theta}=h_\varepsilon$ and consequently,
\begin{align*}
\int_{\mathbb{T}}\int_{\mathbb{R}}(h_\varepsilon*\rho^\varepsilon)j^\varepsilon\,d\theta\,d\Omega&=\int_{\mathbb{T}}\int_{\mathbb{R}}\left(\frac{\partial W_\varepsilon}{\partial \theta}*\rho^\varepsilon\right)j^\varepsilon\,d\theta\,d\Omega=\int_{\mathbb{T}}\int_{\mathbb{R}}\frac{\partial}{\partial\theta}(W_\varepsilon*\rho^\varepsilon)j^\varepsilon\,d\theta\,d\Omega\\
&=-\int_{\mathbb{T}}\int_{\mathbb{R}}(W_\varepsilon*\rho^\varepsilon)\frac{\partial j^\varepsilon}{\partial \theta}\,d\theta\,d\Omega=\int_{\mathbb{T}}\int_{\mathbb{R}}\left(W_\varepsilon*\rho^\varepsilon\right)\frac{\partial f^\varepsilon}{\partial t}\,d\theta\,d\Omega,
\end{align*}
where an integration by parts and the continuity equation \eqref{E-second-order-kuramoto-kinetic-moments-0} have been used. Now, let us restate the last term
\begin{multline*}
\int_{\mathbb{T}}\int_{\mathbb{R}}\left(W_\varepsilon*\rho^\varepsilon\right)\frac{\partial f^\varepsilon}{\partial t}\,d\theta\,d\Omega=\int_{\mathbb{T}^2}\int_{\mathbb{R}^2}W_{\varepsilon}(\theta-\theta')\frac{\partial f^\varepsilon}{\partial t}(t,\theta,\Omega)f^\varepsilon(t,\theta',\Omega')\,d\theta\,d\theta'\,d\Omega\,d\Omega'\\
=\frac{d}{dt}\int_{\mathbb{T}}(W_\varepsilon*\rho^\varepsilon)\rho^\varepsilon\,d\theta-\int_{\mathbb{T}^2}\int_{\mathbb{R}^2}W_\varepsilon(\theta-\theta')f^\varepsilon(t,\theta,\Omega)\frac{\partial f^\varepsilon}{\partial t}(t,\theta',\Omega')\,d\theta\,d\theta'\,d\Omega\,d\Omega'.
\end{multline*}
Since $W_\varepsilon$ is an even function, then the left hand side agrees with the second term of the right hand side and consequently,
$$2\int_{\mathbb{T}}(h_\varepsilon*\rho^\varepsilon)j^\varepsilon\,d\theta\,d\Omega=\frac{d}{dt}\int_{\mathbb{T}}(W_\varepsilon*\rho^\varepsilon)\rho^\varepsilon\,d\theta.$$
This ends the proof of this lemma.
\end{proof}

\begin{theo}\label{T-estimates-kura-kin-2nd-order-1}
Consider the strong solution $P^\varepsilon$ to \eqref{E-second-order-kuramoto-kinetic} whose initial data $P^\varepsilon_0$ fulfill the assumptions \eqref{kura-kin-2nd-order-hypothesis-lp}, for some $p\in [1,+\infty)$. Then, the following estimates 
\begin{align*}
\Vert f^\varepsilon\Vert_{L^\infty(0,T;L^p(\mathbb{R},L^1(\mathbb{T})))}&\leq  C_0,\\
\Vert j^\varepsilon\Vert_{L^2(0,T;L^q(\mathbb{R},L^1(\mathbb{T})))}&\leq C_0^{1/2}\left(2\varepsilon E_0+KM_\alpha+T(V_0+2)\right)^{1/2},\\
\Vert \mathcal{S}^\varepsilon\Vert_{L^1(0,T;L^1(\mathbb{T}\times\mathbb{R}))}&\leq 2\varepsilon E_0+KM_\alpha+T(V_0+2),
\end{align*}
hold  for every $\varepsilon>0$, where  $q:=\frac{2p}{1+p}$.
\end{theo}

\begin{proof}

$\bullet$ \textit{Step 1:} By integration with respect to $\theta$ in the continuity equation \eqref{E-second-order-kuramoto-kinetic} we achieve
$$\frac{d}{dt}\int_{\mathbb{T}}f^\varepsilon\,d\theta=0\ \Longrightarrow\ \int_{\mathbb{T}}f^\varepsilon\,d\theta=\int_\mathbb{T}f^\varepsilon_0\,d\theta=g^\varepsilon.$$
Then, the first a priori estimate for $f^\varepsilon$ follows from Remark \ref{R-kura-kin-2nd-order-hypothesis-g} (with $C_0=1$ if $p=1$).\\

$\bullet$ \textit{Step 2:} Using the Cauchy--Schwarz and Young inequalities on the first term in the right hand side of the formula obtained in Lemma \ref{L-estimates-kura-2nd-order} we arrive at
$$
\frac{d}{dt}\left(\varepsilon\int_{\mathbb{T}}\int_{\mathbb{R}}\mathcal{S}^\varepsilon\,d\theta\,d\Omega+K\int_{\mathbb{T}}(W_\varepsilon*\rho^\varepsilon)\rho^\varepsilon\,d\theta\right)+2\int_{\mathbb{T}}\int_{\mathbb{R}}\mathcal{S}^\varepsilon\,d\theta\,d\Omega\leq \Vert \Omega^2g^\varepsilon\Vert_{L^1(\mathbb{R})}+\int_{\mathbb{T}}\int_{\mathbb{R}}\mathcal{S}^\varepsilon\,d\theta\,d\Omega+2.
$$
Now, let us integrate with respect to time in $[0,T]$. Using the fundamental theorem of calculus and neglecting the terms corresponding to time $t=T$ (notice that $W_\varepsilon\geq 0$ by Remark \ref{R-h-scaled-primitive}), we obtain
$$\int_0^T\int_{\mathbb{T}}\int_{\mathbb{R}}\mathcal{S}^\varepsilon\,dt\,d\theta\,d\Omega\leq \varepsilon_0\int_{\mathbb{T}}\int_{\mathbb{R}}\mathcal{S}^\varepsilon_0\,d\theta\,d\Omega+K\int_{\mathbb{T}}(W_\varepsilon*\rho^\varepsilon_0)f^\varepsilon_0\,d\theta+T\left(\Vert \Omega^2g^\varepsilon\Vert_{L^p(\mathbb{R})}+2\right).$$
Using the assumptions \eqref{kura-kin-2nd-order-hypothesis-lp} and the uniform-in-$\varepsilon$ bound of $W_\varepsilon$ in Remark \ref{R-h-scaled-primitive}, we obtain the third a priori estimate.\\

$\bullet$ \textit{Step 3:} The second a priori estimate is a consequence of the first one that we can obtain by interpolation in $L^p$ spaces. Indeed,
$$\int_{\mathbb{T}}\vert j^\varepsilon\vert \,d\theta\leq \int_{\mathbb{T}}\int_{\mathbb{R}}\vert \omega\vert f^\varepsilon\,d\omega\leq\left(\int_{\mathbb{T}} f^\varepsilon\,d\theta\right)^{1/2}\left(\int_{\mathbb{T}}\mathcal{S}^\varepsilon\,d\theta\right)^{1/2}=(g^\varepsilon)^{1/2}\left(\int_{\mathbb{T}}\mathcal{S}^\varepsilon\,d\theta\right)^{1/2},$$
where the Cauchy-Schwarz inequality has been used. Define the exponent $q=\frac{2p}{1+p}$ and notice that
$$\frac{1}{2p}+\frac{1}{2}=\frac{1}{q}.$$
Then, we can take $L^q$ norms in the above expression and use the generalized H\"{o}lder inequality with exponents $2p$ and $2$ to obtain
$$\Vert j^\varepsilon(t)\Vert_{L^q(\mathbb{R},L^1(\mathbb{T}))}\leq \Vert g^\varepsilon\Vert_{L^p(\mathbb{R})}^{1/2}\Vert \mathcal{S}^\varepsilon(t)\Vert_{L^1(\mathbb{T}\times\mathbb{R})}^{1/2},$$
for every $t\in [0,T]$. Finally, let us take $L^2$ norms with respect to time to achieve the second a priori estimate for $j^\varepsilon$ in $L^2(0,T;L^q(\mathbb{R},L^1(\mathbb{T})))$.
\end{proof}

\subsection{Compactness of the regularized system}\label{SS-second-order-regularized-compactness}

In this part, we will derive the corresponding weak-star compactness as a consequence of Theorem \ref{T-estimates-kura-kin-2nd-order-1}. We will do it for both type of assumptions \eqref{kura-kin-2nd-order-hypothesis} and \eqref{kura-kin-2nd-order-hypothesis-lp} that we can assume on the initial data. Then, we obtain the following two Corollaries for each of the two cases:
\begin{cor}\label{C-compactness-kura-kin-2nd-order}
Consider the strong solution $P^\varepsilon$ to \eqref{E-second-order-kuramoto-kinetic} whose initial data $P^\varepsilon_0$ fulfill the assumptions \eqref{kura-kin-2nd-order-hypothesis}. Then, there exists $f\in L^\infty(0,T;\mathcal{M}(\mathbb{T}\times \mathbb{R}))$ and $j\in L^2(0,T;L^1(\mathbb{T}\times\mathbb{R}))$ such that
\begin{align*}
f^\varepsilon&\overset{*}{\rightharpoonup}f\ \mbox{ in }\ L^\infty(0,T;\mathcal{M}(\mathbb{T}\times\mathbb{R})),\\
j^\varepsilon&\overset{*}{\rightharpoonup} j\ \mbox{ in }\ L^2(0,T;\mathcal{M}(\mathbb{T}\times \mathbb{R})),
\end{align*}
up to a subsequence that we denote the same for simplicity.
\end{cor}

\begin{cor}\label{C-compactness-kura-kin-2nd-order-lp}
Consider the strong solution $P^\varepsilon$ to \eqref{E-second-order-kuramoto-kinetic} whose initial data $P^\varepsilon_0$ fulfill the assumptions \eqref{kura-kin-2nd-order-hypothesis-lp}, for some $1<p<\infty$. Then, the above weak-star limits $f\in L^\infty(0,T;\mathcal{M}(\mathbb{T}\times \mathbb{R}))$ and $j\in L^2(0,T;L^1(\mathbb{T}\times\mathbb{R}))$ in Corollary \ref{C-compactness-kura-kin-2nd-order} also satisfy
\begin{align*}
f^\varepsilon&\overset{*}{\rightharpoonup}f\ \mbox{ in }\ L^\infty(0,T;L^p(\mathbb{R},\mathcal{M}(\mathbb{T}))),\\
j^\varepsilon&\overset{*}{\rightharpoonup} j\ \mbox{ in }\ L^2(0,T;L^{\frac{2p}{1+p}}(\mathbb{R},\mathcal{M}(\mathbb{T}))).
\end{align*}
\end{cor}

We skip the proof of both results since it follows from standard arguments supported by the Alaoglu theorem along with the Riesz representation theorem of Banach-valued Lebesgue spaces.

Apart from the above weak convergence in time, we can recover a stronger convergence result of the density $f^\varepsilon$. This is the content of the next result.

\begin{theo}\label{T-compactness-kura-kin-2nd-order-stronger}
Consider the strong solution $P^\varepsilon$ to \eqref{E-second-order-kuramoto-kinetic} with initial data $P^\varepsilon_0$.
\begin{enumerate}
\item If $P^\varepsilon_0$ fulfils the assumptions \eqref{kura-kin-2nd-order-hypothesis}, then
$$
f^\varepsilon\rightarrow f\ \mbox{ in }\ C([0,T],\mathcal{M}(\mathbb{T}\times\mathbb{R})-\mbox{narrow}).
$$
\item If $P^\varepsilon_0$ fulfils the assumptions \eqref{kura-kin-2nd-order-hypothesis-lp}, for some $1<p<\infty$, then
$$
f^\varepsilon\rightarrow f\ \mbox{ in }\ C([0,T],L^p(\mathbb{R},\mathcal{M}(\mathbb{T}))-\mbox{weak}\,*).
$$
\end{enumerate}
\end{theo}

\begin{proof}
Since both proofs are similar, we just focus on the most involved one, that is, the first one. The second result follows a parallel train of thoughts. Let use the continuity equation \eqref{E-second-order-kuramoto-kinetic-moments-0} that we write in weak form against a test function $\varphi(t,\theta,\Omega)=\eta(t)\phi(\theta,\Omega)$, where $\eta$ and $\phi$ are smooth and compactly supported
\begin{multline*}
\int_0^T\frac{\partial\eta}{\partial t}\int_{\mathbb{T}}\int_{\mathbb{R}}f^\varepsilon(t,\theta,\Omega)\phi(\theta,\Omega)\,dt\,d\theta\,d\Omega\\
=-\int_0^T\int_{\mathbb{T}}\int_{\mathbb{R}}j^\varepsilon(t,\theta,\Omega)\eta(t)\frac{\partial\phi}{\partial\theta}(\theta,\Omega)\,dt\,d\theta\,d\Omega\\
\leq\Vert j^\varepsilon\Vert_{L^2(0,T;L^1(\mathbb{T}\times \mathbb{R}))}\Vert \phi\Vert_{W^{1,\infty}_0(\mathbb{T}\times \mathbb{R})}.
\end{multline*}
Then, the standard characterization of Sobolev spaces yields
$$
\left\Vert\int_{\mathbb{T}}\int_{\mathbb{R}}f^\varepsilon(\cdot,\theta,\Omega)\phi(\theta,\Omega)\,d\theta\,d\Omega\right\Vert_{H^1(0,T)}\\
\leq (T^{1/2}\Vert f^\varepsilon\Vert_{L^\infty(0,T;L^1(\mathbb{T}\times \mathbb{R}))}+\Vert j^\varepsilon\Vert_{L^2(0,T;L^1(\mathbb{T}\times\mathbb{R}))}) \Vert \phi\Vert_{W^{1,\infty}_0(\mathbb{T}\times\mathbb{R})}.
$$
By the Sobolev embedding theorem and the first and second estimates in Theorem \ref{T-estimates-kura-kin-2nd-order-1} we can also claim that there exists $C>0$ that does not depend neither on $\varepsilon$ nor in the chosen test functions so that
$$
\left\vert\int_{\mathbb{T}}\int_{\mathbb{R}}(f^\varepsilon(t_1,\theta,\Omega)-f^\varepsilon(t_2,\theta,\Omega))\phi(\theta,\Omega)\,d\theta\,d\Omega\right\vert\leq C\vert t_1-t_2\vert^{1/2}\Vert\phi\Vert_{W^{1,\infty}_0(\mathbb{T}\times\mathbb{R})},
$$
for every $t_1,t_2\in [0,T]$ and each $\varepsilon>0$. Since the chosen test functions are arbitrary, then we  obtain
$$
\Vert f^\varepsilon(t_1,\cdot,\cdot)-f^\varepsilon(t_2,\cdot,\cdot)\Vert_{W^{-1,1}(\mathbb{T}\times \mathbb{R})}\leq C\vert t_1-t_2\vert^{1/2},
$$
for every $t_1,t_2\in [0,T]$ and each $\varepsilon>0$. Also notice that
$$
\Vert f^\varepsilon\Vert_{C([0,T],W^{-1,1}(\mathbb{T}\times \mathbb{R}))}\leq C,
$$
for every $\varepsilon>0$ for another constant $C$ that, without loss of generality, we denote in the same way. Such assertion is nothing but a consequence of the first a priori estimate in Theorem \ref{T-estimates-kura-kin-2nd-order-1} and the chain of continuous embeddings
$$L^1(\mathbb{T}\times\mathbb{R})\hookrightarrow\mathcal{M}(\mathbb{T}\times\mathbb{R})\hookrightarrow W^{-1,1}(\mathbb{T}\times\mathbb{R}),$$
where the last embedding is a consequence of the dense embedding of $W^{1,\infty}_0(\mathbb{T}\times \mathbb{R})$ into $C_0(\mathbb{T}\times \mathbb{R})$. Then, we can use the Ascoli-Arzel\`a theorem in weak form to the space $C([0,T],(W^{1,\infty}_0(\mathbb{T}\times \mathbb{R}))^*)\cong C([0,T],W^{-1,1}(\mathbb{T}\times \mathbb{R}))$ and obtain a subsequence of $f_\varepsilon$ that we denote in the same way so that
$$f^\varepsilon\rightarrow f\ \mbox{ in }\ C([0,T],W^{-1,1}(\mathbb{T}\times \mathbb{R})-\mbox{weak}*).$$
By the firs estimate in Theorem \ref{T-estimates-kura-kin-2nd-order-1} along with the above-mentioned density of $W^{1,\infty}_0(\mathbb{T}\times \mathbb{R})$ into $C_0(\mathbb{T}\times \mathbb{R})$ we can actually improve the above convergence into
$$f^\varepsilon\rightarrow f\ \mbox{ in }\ C([0,T],\mathcal{M}(\mathbb{T}\times \mathbb{R})-\mbox{weak}*).$$
Let us finally show that the above weak-star convergence can be improved into narrow convergence in the spaces of finite Radon measures $\mathcal{M}(\mathbb{T}\times \mathbb{R})$. To such end, note that
$$\int_{\mathbb{T}}\int_{\mathbb{R}}\vert\Omega\vert f^\varepsilon\,d\theta\,d\Omega=\int_{\mathbb{T}}\int_{\mathbb{R}}\vert \Omega\vert f^\varepsilon_0\,d\theta\,d\Omega\leq C_0,$$
Then, the strengthened version of the Prokhorov's compactness theorem in Lemma \ref{L-Prokhorov} (see Appendix \ref{appendix-Prokhorov}) yields the desired strong-narrow convergence:
$$f^\varepsilon\rightarrow f\ \mbox{ in }C([0,T],\mathcal{M}(\mathbb{T}\times\mathbb{R})-\mbox{narrow}).$$
\end{proof}

As a simple consequence, we can prove that the tensor product $f^\varepsilon\otimes f^\varepsilon$ also enjoys good compactness properties.

\begin{cor}\label{C-compactness-kura-kin-2nd-order-stronger}
Consider the strong solution $P^\varepsilon$ to \eqref{E-second-order-kuramoto-kinetic} with initial data $P^\varepsilon_0$.
\begin{enumerate}
\item If $P^\varepsilon_0$ fulfils the assumptions \eqref{kura-kin-2nd-order-hypothesis}, then
$$
f^\varepsilon\otimes f^\varepsilon\overset{*}{\rightharpoonup}f\otimes f\ \mbox{ in }\ L^\infty(0,T;\mathcal{M}(\mathbb{T}^2\times\mathbb{R}^2)-\mbox{narrow}).
$$
\item If $f_\varepsilon^0$ fulfils the assumptions \eqref{kura-kin-2nd-order-hypothesis-lp}, for some $1<p<\infty$, then
$$
f^\varepsilon\otimes f^\varepsilon\overset{*}{\rightharpoonup} f\otimes f\ \mbox{ in }\ L^\infty(0,T;L^p(\mathbb{R}^2,\mathcal{M}(\mathbb{T}^2))-\mbox{weak}\,*).
$$
\end{enumerate}
\end{cor}

Our next step is to show that we can pass to the limit in the balance laws for the phase density and phase current \eqref{E-second-order-kuramoto-kinetic-moments-0}-\eqref{E-second-order-kuramoto-kinetic-moments-1}. To such end, let us write them in weak form by multiplication against any test function $\varphi\in C^1_0([0,T)\times \mathbb{T}\times\RR)$ and integrate by parts:

\begin{align}
\int_0^T\int_{\mathbb{T}}\int_{\RR}f^\varepsilon\frac{\partial \varphi}{\partial t}\,dt\,d\theta\,d\Omega&+\int_0^T\int_{\mathbb{T}}\int_{\RR} j^\varepsilon\frac{\partial \varphi}{\partial \theta}\,dt\,d\theta\,d\Omega=-\int_{\mathbb{T}}\int_{\RR}f^\varepsilon_0\,\varphi(0,\cdot,\cdot)\,d\theta\,d\Omega,\label{E-second-order-kuramoto-kinetic-moments-0-weak}\\
\varepsilon\int_0^T\int_{\mathbb{T}}\int_{\RR}j^\varepsilon \frac{\partial \varphi}{\partial t}\,dt\,d\theta\,d\Omega&+\varepsilon\int_0^T\int_{\mathbb{T}}\int_{\RR}\mathcal{S}^\varepsilon\frac{\partial\varphi}{\partial \theta}\,dt\,d\theta\,d\Omega=-\varepsilon\int_{\mathbb{T}}\int_{\RR}j^\varepsilon_0\,\varphi(0,\cdot,\cdot)\,d\theta\,d\Omega\nonumber\\
&+\int_0^T\int_{\mathbb{T}}\int_{\RR}(j^\varepsilon-\Omega f^\varepsilon)\,\varphi\,dt\,d\theta\,d\Omega+\int_0^T\int_{\mathbb{T}}\int_{\RR}(h_\varepsilon*\rho^\varepsilon)f^\varepsilon\varphi\,dt\,d\theta\,d\Omega.\label{E-second-order-kuramoto-kinetic-moments-1-weak}
\end{align}
Note that under the weak assumptions \eqref{kura-kin-2nd-order-hypothesis}, Corollary \ref{C-compactness-kura-kin-2nd-order} and the a priori estimate of $\mathcal{S}_\varepsilon$ in Theorem \ref{T-estimates-kura-kin-2nd-order-1} allow passing to the limit all the terms except at most two of them; namely, the nonlinear term and the term involving $\Omega f^\varepsilon$. Let us first address the latter one first and we will discuss about the most difficult convergence result for the nonlinear term later.

\begin{pro}\label{P-estimates-kura-kin-2nd-order-1}
Consider the strong solution $P^\varepsilon$ to \eqref{E-second-order-kuramoto-kinetic} with initial data $P^\varepsilon_0$.
\begin{enumerate}
\item If $P^\varepsilon_0$ fulfil the assumptions \eqref{kura-kin-2nd-order-hypothesis}, then
\begin{align*}
\Vert \Omega f^\varepsilon\Vert_{L^\infty(0,T;L^1(\mathbb{T}\times \RR))}&\leq V_0^{1/2},\\
\Vert \Omega f\Vert_{L^\infty(0,T;\mathcal{M}(\mathbb{T}\times\RR))}&\leq V_0^{1/2}.
\end{align*}
Moreover, the following weak-star convergence takes place
$$\Omega f^\varepsilon\overset{*}{\rightharpoonup}\Omega f\ \mbox{ in }\ L^\infty(0,T;\mathcal{M}(\mathbb{T}\times\RR)).$$
\item If $P^\varepsilon_0$ fulfil the assumptions \eqref{kura-kin-2nd-order-hypothesis-lp}, for some $1<p<\infty$, then
\begin{align*}
\Vert \Omega f^\varepsilon\Vert_{L^\infty(0,T;L^\frac{2p}{1+p}(\mathbb{R},L^1(\mathbb{T})))}&\leq V_0^{1/2}C_0^{1/2},\\
\Vert \Omega f\Vert_{L^\infty(0,T;L^\frac{2p}{1+p}(\mathbb{R},\mathcal{M}(\mathbb{T})))}&\leq V_0^{1/2}C_0^{1/2}.
\end{align*}
Moreover, the following weak-star convergence takes place
$$\Omega f^\varepsilon\overset{*}{\rightharpoonup}\Omega f\ \mbox{ in }\ L^\infty(0,T;L^\frac{2p}{1+p}(\mathbb{R},\mathcal{M}(\mathbb{T}))).$$
\end{enumerate}
\end{pro}

\begin{proof}
Since both proofs are similar, we just focus on the second one. Notice that by the continuity equation we obtain that
$$\frac{d}{dt}\int_{\mathbb{T}}\vert \Omega\vert f^\varepsilon\,d\theta=0\ \Longrightarrow\ \int_{\mathbb{T}}\vert \Omega\vert f^\varepsilon\,d\theta=\int_{\mathbb{T}}\vert \Omega\vert f^\varepsilon_0\,d\theta=\vert \Omega\vert g^\varepsilon.$$
Again, we can take $L^\frac{2p}{1+p}$-norms and use the generalized H\"{o}lder inequality to arrive at
$$\Vert f^\varepsilon(t)\Vert_{L^{\frac{2p}{1+p}}(\mathbb{R},L^1(\mathbb{T}))}\leq \Vert \Omega^2 g^\varepsilon\Vert_{L^1(\mathbb{R})}^{1/2}\Vert g^\varepsilon\Vert_{L^p(\mathbb{R})}^{1/2},$$
for every $t\in [0,T]$. Taking supreme yields the desired estimate of $\Omega f^\varepsilon$ by virtue of the assumptions \eqref{kura-kin-2nd-order-hypothesis-lp}. Regarding the limiting estimate, let us first note that 
$$(t,\Omega)\in (0,T)\times \mathbb{R}\longmapsto \Vert \Omega f_t(\cdot,\Omega)\Vert_{\mathcal{M}(\mathbb{T})}=\vert \Omega\vert \Vert f_t(\cdot,\Omega)\Vert_{\mathcal{M}(\mathbb{T})}$$ 
belongs to $L^1_{\textnormal{loc}}((0,T)\times\mathbb{R})$. In order to get the $L^\infty(0,T;L^q(\mathbb{R}))$-estimate with $q=\frac{2p}{1+p}$, we consider any $\varphi\in C_c((0,T)\times\mathbb{R})$ and compute
\begin{align*}
\int_0^T\int_{\mathbb{R}}&\varphi(t,\Omega)\Vert \Omega f_t(\cdot,\Omega)\Vert_{\mathcal{M}(\mathbb{T})}\,dt\,d\Omega=\int_0^T\int_{\mathbb{R}}\vert \Omega\vert\varphi(t,\Omega)\Vert f_t(\cdot,\Omega)\Vert_{\mathcal{M}(\mathbb{T})}\,dt\,d\Omega\\
&=\int_0^T\int_{\mathbb{R}}\int_{\mathbb{T}}\left<f_t(\cdot,\Omega),\vert \Omega\vert\varphi(t,\Omega)\right>_{\mathcal{M}(\mathbb{T})}\,dt\,d\Omega=\lim_{\varepsilon\rightarrow 0}\int_0^T\int_{\mathbb{T}}\int_{\mathbb{R}}\vert\Omega\vert \varphi(t,\Omega)f^\varepsilon(t,\theta,\Omega)\,dt\,d\theta\,d\Omega\\
&\leq \limsup_{\varepsilon\rightarrow 0}\Vert\Omega f^\varepsilon\Vert_{L^\infty(0,T;L^q(\mathbb{R},L^1(\mathbb{T})))}\Vert\varphi\Vert_{L^1(0,T;L^{q'}(\mathbb{R}))}\leq V_0^{1/2}C_0^{1/2}\Vert\varphi\Vert_{L^1(0,T;L^{q'}(\mathbb{R}))},
\end{align*}
where we have used the above estimate for $\Omega f^\varepsilon$ and the weak-star convergence in Corollary \ref{C-compactness-kura-kin-2nd-order-lp}, notice that $\vert \Omega\vert \varphi(t,\Omega)$ belongs to space of test functions $L^1(0,T;L^{q'}(\mathbb{R};C(\mathbb{T})))$. Since $\varphi$ is any compactly supported test function we can conclude the desired estimate on $\Omega f$ by the Riesz representation theorem.

Finally, let us show the convergence result. To this end, we fix a test function $\varphi\in L^1(0,T;L^{q'}(\mathbb{R},C(\mathbb{T})))$. By density, we can assume that $\varphi\in C_c((0,T)\times \mathbb{T}\times\mathbb{R})$. Then,
\begin{align*}
&\int_0^T\int_{\mathbb{T}}\int_{\RR}\varphi(t,\theta,\Omega)(\Omega f^\varepsilon(t,\theta,\Omega)-\Omega f_t(\theta,\Omega))\,dt\,d\theta\,d\Omega\\
 & \qquad =\int_0^T\int_{\mathbb{T}}\int_{\RR}\Omega\varphi(t,\theta,\Omega)(f^\varepsilon(t,\theta,\Omega)-f_t(\theta,\Omega))\,dt\,d\theta\,d\Omega.
\end{align*}
Since $\Omega\varphi(t,\theta,\Omega)$ belongs to $L^1(0,T;L^{q'}(\mathbb{R},C(\mathbb{T}))$, then we can apply Corollary \eqref{C-compactness-kura-kin-2nd-order-lp} and obtain that the last integral converges towards zero when $\varepsilon\rightarrow 0$, thus concluding the proof.
\end{proof}

Notice that the above result allows passing to the limit in the term $\Omega f^\varepsilon$ in \eqref{E-second-order-kuramoto-kinetic-moments-1-weak}. Then, the only term that remains to be studied is the nonlinear one.

\subsection{Convergence of the nonlinear term and hydrodynamic limit}\label{SS-second-order-regularized-convergence}
In this part, we will discuss about the nonlinear term \eqref{E-second-order-kuramoto-kinetic-moments-1-weak}. On the one hand, we will show that we cannot pass to the limit for general $g$. The main reason is that the proposed cancellation property would fail for $\alpha\in (\frac{1}{2},1)$. On the other hand, we will show that it proves useful in the identical case $g=\delta_0$. Although we will not comment on the more regular regime $\alpha\in (0,\frac{1}{2})$ here, it is clear that $h_\varepsilon*\rho^\varepsilon$ converges strongly for general $g$ although $f^\varepsilon$ is just narrowly convergent. In such way, we can recover the weak measure-valued solutions to \eqref{E-kuramoto-transport-TxR} like in Section \ref{S-weak-solutions-existence}, but we will skip it here for simplicity.

The main idea is to write the nonlinear term in \eqref{E-second-order-kuramoto-kinetic-moments-1-weak} appropriately using a well known symmetrization idea. Specifically, note that
\begin{multline*}
\int_0^T\int_\mathbb{T}\int_\mathbb{R}(h_\varepsilon*\rho^\varepsilon)f^\varepsilon\varphi\,dt\,d\theta\,d\Omega=\int_0^T\int_{\mathbb{T}^2\times\mathbb{R}^2}h_\varepsilon(\theta-\theta')\varphi(t,\theta,\Omega)f^\varepsilon(t,\theta,\Omega)f^\varepsilon(t,\theta',\Omega')\,dt\,d\theta\,d\theta'\,d\Omega\,d\Omega'\\
=-\int_0^T\int_{\mathbb{T}^2\times \mathbb{R}^2}h_\varepsilon(\theta-\theta')\varphi(t,\theta',\Omega')f^\varepsilon(t,\theta,\Omega)f^\varepsilon(t,\theta',\Omega')\,dt\,d\theta\,d\theta'\,d\Omega\,d\Omega',
\end{multline*}
where we have changed variables $(\theta,\Omega)$ with $(\theta',\Omega')$ and we have used the antisymmetry of the kernel $h_\varepsilon$ in the last line. Taking the mean value of both expressions we obtain
\begin{equation}\label{E-nonlinear-term-symmetrized}
\int_0^T\int_\mathbb{T}\int_\mathbb{R}(h_\varepsilon*\rho^\varepsilon)f^\varepsilon\varphi\,dt\,d\theta\,d\Omega=\int_0^T\int_{\mathbb{T}^2\times \mathbb{R}^2}H^\varphi_\varepsilon(t,\theta,\theta',\Omega,\Omega')f^\varepsilon(t,\theta,\Omega)f^\varepsilon(t,\theta',\Omega')\,dt\,d\theta\,d\theta'\,d\Omega\,d\Omega',
\end{equation}
where the function $H_\varepsilon^\varphi$ reads
$$
H_\varepsilon^\varphi(t,\theta,\theta',\Omega,\Omega'):=\frac{1}{2}h_\varepsilon(\theta-\theta')(\varphi(t,\theta,\Omega)-\varphi(t,\theta',\Omega')).
$$
Notice that when $\varepsilon\rightarrow 0$, the above function $H_0^\varphi$ is not continuous at $\theta=\theta'$ unless $\Omega=\Omega'$. Here, the Lipschitz continuity of the test function $\varphi$ in the variable $\theta$ plays a role to cancel the full singularity of $h$ in the whole range $\alpha\in (\frac{1}{2},1)$. Then, it is not clear at all that \eqref{E-nonlinear-term-symmetrized} makes sense in the limit $\varepsilon\rightarrow 0$, as the limiting measure $f$ is expected to have atoms. However, we can solve the problem at least for identical oscillators, i.e., $g=\delta_0$. This is the content of the main result in this part.

\begin{theo}\label{T-hydrodynamic-field-existence-supercritical}
Fix $\alpha\in (\frac{1}{2},1)$ and consider the strong solution $P^\varepsilon$ to \eqref{E-second-order-kuramoto-kinetic} whose initial data $P^\varepsilon_0$ fulfil the assumptions \eqref{kura-kin-2nd-order-hypothesis}. Assume that the limiting oscillators are all identical, that is, $g=\delta_0$. Then, there is a subsequence of $f^\varepsilon$, denoted in the same way, and $f\in \mathcal{AC}_\mathcal{M}\cap \mathcal{T}_\mathcal{M}$ so that 
$$f^\varepsilon\rightarrow f\ \mbox{ in }\ C([0,T],\mathcal{M}(\mathbb{T}\times \mathbb{R})-\mbox{narrow}),$$
and its associated macroscopic density $\rho$ is a weak measure-valued solution (in the symmetrized sense) to \eqref{E-kuramoto-transport-TxR} in the identical case
$$
\left\{
\begin{array}{l}
\displaystyle\frac{\partial \rho}{\partial t}-\frac{\partial}{\partial \theta}\left((h*\rho)\rho\right)=0,\\
\displaystyle \rho(0,\cdot,\cdot)=\rho_0.
\end{array}
\right.
$$
Specifically, in weak form
$$\int_0^T\int_\mathbb{T}\frac{\partial\phi}{\partial t}\,d_\theta\rho-\int_0^T\int_{\mathbb{T}^2}\frac{1}{2}h(\theta-\theta')\left(\frac{\partial \phi}{\partial\theta}(t,\theta)-\frac{\partial \phi}{\partial\theta}(t,\theta')\right)\,d_{(\theta,\theta')}(\rho\otimes \rho)=-\int_\mathbb{T}\phi(0,\theta)d_\theta\rho_0,$$
for any test function $\phi\in C^1_c([0,T)\times \mathbb{T})$. Here, $\rho_0\in\mathcal{P}(\mathbb{T})$ is any prescribed initial distribution of oscillators as given in the assumptions \eqref{kura-kin-2nd-order-hypothesis}. If in addition \eqref{kura-kin-2nd-order-hypothesis-lp} holds true, for some $1<p<\infty$, then
$$f^\varepsilon\rightarrow f\ \mbox{ in }\ C([0,T],L^p(\mathbb{R},\mathcal{M}(\mathbb{T}))-\mbox{weak}\,*).$$
\end{theo}

\begin{proof}
Let us recover the weak formulation \eqref{E-second-order-kuramoto-kinetic-moments-1-weak} and use the symmetrized version \eqref{E-nonlinear-term-symmetrized} of the nonlinear term. Then, we obtain
\begin{multline}\label{E-67}
\varepsilon\int_0^T\int_{\mathbb{T}}\int_{\RR}j^\varepsilon \frac{\partial \varphi}{\partial t}\,dt\,d\theta\,d\Omega+\varepsilon\int_0^T\int_{\mathbb{T}}\int_{\RR}\mathcal{S}^\varepsilon\frac{\partial\varphi}{\partial \theta}\,dt\,d\theta\,d\Omega\\
=-\varepsilon\int_{\mathbb{T}}\int_{\RR}j^\varepsilon_0\,\varphi(0,\cdot,\cdot)\,d\theta\,d\Omega+\int_0^T\int_{\mathbb{T}}\int_{\RR}(j^\varepsilon-\Omega f^\varepsilon)\,\varphi\,dt\,d\theta\,d\Omega\\
+\int_0^T\int_{\mathbb{T}^2}\int_{\mathbb{R}^2}H^\varphi_\varepsilon(t,\theta,\theta',\Omega,\Omega')f^\varepsilon(t,\theta,\Omega)f^\varepsilon(t,\theta',\Omega')\,dt\,d\theta\,d\theta'\,d\Omega\,d\Omega',
\end{multline}
for any test function $\varphi\in C^1_c([0,T)\times\mathbb{T}\times \mathbb{R})$. The main idea is to get rid of the $\Omega$ variable in \eqref{E-67}, as it becomes irrelevant in the limit $\varepsilon\rightarrow 0$ due to the assumed hypothesis $g=\delta_0$. To simplify notation, denote
\begin{align*}
\widehat{j}^\varepsilon(t,\theta)&:=\int_{\mathbb{R}}j^\varepsilon(t,\theta,\Omega)\,d\Omega=(\pi_\theta)_{\#}j^\varepsilon, & \widehat{\mathcal{S}}^\varepsilon(t,\theta)&:=\int_{\mathbb{R}}\mathcal{S}^\varepsilon(t,\theta,\omega)\,d\Omega=(\pi_\theta)_{\#}\mathcal{S}^\varepsilon,\\
\widehat{j}_t(\theta)&:=(\pi_\theta)_{\#}j_t. & &
\end{align*}
Recall that Theorem \ref{T-compactness-kura-kin-2nd-order-stronger}, Corollaries \ref{C-compactness-kura-kin-2nd-order-stronger}, \ref{C-compactness-kura-kin-2nd-order} along with Theorem \ref{T-estimates-kura-kin-2nd-order-1} respectively imply
\begin{equation}\label{E-68}
\begin{array}{l}
\displaystyle \rho^\varepsilon\rightarrow \rho \ \mbox{ in }\ C([0,T],\mathcal{M}(\mathbb{T})-\mbox{narrow}),\\
\displaystyle \widehat{j}^\varepsilon\overset{*}{\rightharpoonup} \widehat{j}\ \mbox{ in }\ L^2(0,T;\mathcal{M}(\mathbb{T})-\mbox{narrow}),\\
\displaystyle\Vert \widehat{\mathcal{S}}^\varepsilon\Vert_{L^1(0,T;L^1(\mathbb{T}))}\leq 2\varepsilon E_0+KM_\alpha+T(V_0+2),\\
\displaystyle\rho^\varepsilon\otimes \rho^\varepsilon\overset{*}{\rightharpoonup}\rho\otimes \rho\ \mbox{ in }\ L^\infty(0,T;\mathcal{M}(\mathbb{T}^2)).
\end{array}
\end{equation}
Now, we can take $\varphi(t,\theta,\Omega)=\phi(t,\theta)$ in \eqref{E-67} with $\phi\in C^1_c([0,T)\times \mathbb{T})$. Notice that it can be done without loss of generality by virtue of the tightness a priori estimates in Proposition \ref{P-estimates-kura-kin-2nd-order-1}, thus yielding
\begin{multline}\label{E-69}
\varepsilon\int_0^T\int_{\mathbb{T}}\widehat{j}^\varepsilon \frac{\partial \phi}{\partial t}\,dt\,d\theta+\varepsilon\int_0^T\int_{\mathbb{T}}\widehat{\mathcal{S}}^\varepsilon\frac{\partial\phi}{\partial \theta}\,dt\,d\theta\\
=-\varepsilon\int_{\mathbb{T}}\widehat{j}^\varepsilon_0\,\phi(0,\cdot)\,d\theta+\int_0^T\int_{\mathbb{T}}\widehat{j}^\varepsilon\,\phi\,dt\,d\theta-\int_0^T\int_{\mathbb{T}}\int_{\mathbb{R}}\Omega f^\varepsilon\,\phi\,dt\,d\theta\,d\Omega\\
+\int_0^T\int_{\mathbb{T}^2}\widehat{H}^\phi_\varepsilon(t,\theta,\theta')\rho^\varepsilon(t,\theta)\rho^\varepsilon(t,\theta')\,dt\,d\theta\,d\theta',
\end{multline}
for the bounded and continuous function
$$\widehat{H}^\phi_\varepsilon(t,\theta,\theta')=\frac{1}{2}h_\varepsilon(\theta-\theta')(\phi(t,\theta)-\phi(t,\theta')).$$
Then, \eqref{E-68} clearly allows passing to the limit in all the terms of \eqref{E-69} (including the nonlinear term)
\begin{equation}\label{E-70}
\int_0^T\int_\mathbb{T}\phi\,d_\theta \widehat{j}\,dt=\int_0^T\int_{\mathbb{T}}\int_{\mathbb{R}}\Omega \,d_{(\theta,\Omega)}f\,dt+\int_0^T\int_{\mathbb{T}^2}H_0^\phi(t,\theta,\theta')d_{(\theta,\theta')}\rho\otimes \rho.
\end{equation}
To finish, let us just identify the first term in the right hand side of \eqref{E-70}, which has not been closed in terms of the macroscopic quantities $\rho$ and $\widehat{j}$ yet. To such end, notice that $\rho^\varepsilon\rightarrow \rho$ in $C([0,T],\mathcal{M}(\mathbb{T})-\mbox{narrow})$ and also $g^\varepsilon\overset{*}{\rightharpoonup} g\equiv \delta_0$ narrowly. Similarly, notice that $f^\varepsilon\rightarrow f$ in $C([0,T],\mathcal{M}(\mathbb{T}\times \mathbb{R})-\mbox{narrow})$. Hence, uniqueness implies
$$f_t(\theta,\Omega)=\rho_t(\theta)\otimes \delta_0(\Omega),\ \mbox{ for all }\ t\in [0,T].$$
Consequently,
$$\int_0^T\int_{\mathbb{T}}\int_{\mathbb{R}}\Omega\, d_{(\theta,\Omega)}f\,dt=\int_0^T\int_{\mathbb{T}}\int_{\mathbb{R}}\Omega\, d_{(\theta,\Omega)}(\rho_t\otimes \delta_0)=0,$$
thus ending the proof.
\end{proof}

\appendix

\section{Periodic finite Radon measures}\label{appendix-relation-measure-spaces}
In this appendix we recall the relations between some spaces of finite Radon measures, that are used throughout the whole paper. Our aim is to appropriately set the space of measures with periodic properties to be used. First, note that the spaces of densities $L^1(\mathbb{T}),\ L^1(-\pi,\pi),\ L^1((-\pi,\pi])$ and $L^1([-\pi,\pi])$ can be mixed up since such functions are defined on the whole torus $\mathbb{T}$ except at most at $-1+0\,i\equiv (-1,0)$. However, it is no longer valid for the spaces of finite Radon measures:
$$\mathcal{M}(\mathbb{T}),\ \mathcal{M}(-\pi,\pi),\ \mathcal{M}((-\pi,\pi])\ \mbox{ and }\ \mathcal{M}([-\pi,\pi]).$$
The main reason is that all such spaces contain the Dirac mass $\delta_{\pi}$ except the second one. Also, the last measure space might duplicate the Dirac masses at $\delta_{-\pi}$ and $\delta_{\pi}$, that can be identified in $\mathcal{M}(\mathbb{T})$ though. Naturally, one has to rule out such doubling of point masses by appropriately setting the good spaces. The main idea is to note that if one unfolds the torus $\mathbb{T}$ into the interval $(-\pi,\pi]$, then each measure in $\mathcal{M}(\mathbb{T})$ can be identified with a measure in $\mathcal{M}((-\pi,\pi])$ and conversely. This is the content of the next straightforward result:

\begin{theo}\label{theo-relation-measure-spaces}
The next Banach spaces of finite Radon measures are topologically isomorphic when endowed with the total variation norm:
$$
\mathcal{M}(\mathbb{T})\cong \mathcal{M}_p([-\pi,\pi])\cong \RR\,\oplus_1\,\mathcal{M}(-\pi,\pi)\cong \mathcal{M}((-\pi,\pi]).
$$
\end{theo}

Although it is straightforward, we will sketch the proof of such result for the readers convenience. 

\subsection{Periodic functions}
Before we sketch the proof, we will first introduce some natural identification between Banach spaces of regular functions with analogue periodicity properties.

\begin{defi}\label{def-complex-derivatives}
For any derivable function $f:\mathbb{T}\longrightarrow\RR$ along $\mathbb{T}$, we define the associated derivatives
\begin{align*}
\frac{\partial f}{\partial z}(e^{i\theta})&:=-ie^{-i\theta}\frac{d}{d\theta}f(e^{i\theta}),\\
\frac{\partial f}{\partial \bar z}(e^{i\theta})&:=ie^{i\theta}\frac{d}{d\theta}f(e^{i\theta}),
\end{align*}
\end{defi}

\begin{rem}\label{rem-complex-derivatives}
The above derivatives have long been used in complex analysis. Consider $z=e^{i\theta}$, $\bar z=e^{-i\theta}$ and define $g(z)=f(\bar z)$. Then, the motivation underlying the above definition is simply a formal chain rule, namely,
\begin{align*}
\frac{d}{d\theta}f(e^{i\theta})&=\frac{\partial f}{\partial z}\frac{dz}{d\theta}=i e^{i\theta}\frac{\partial f}{\partial z} ,\\
\frac{d}{d\theta} f(e^{i\theta})&=\frac{d}{d\theta}g(e^{-i\theta})=\frac{\partial g}{\partial \bar z}\frac{d\bar z}{d\theta}=\frac{\partial f}{\partial \bar{z}}\frac{\partial \bar z}{\partial \bar z}\frac{d\bar z}{d\theta}=-ie^{-i\theta}\frac{\partial f}{\partial \bar{z}}.
\end{align*}
Notice that one can go from one to the other by taking complex conjugation, i.e., $\frac{\partial f}{\partial \bar z}=\overline{\frac{\partial f}{\partial z}}$. Using only one of the above derivatives we can recover the other one (along with the full differential map), thus avoiding redundancy of information. Indeed, $\mathbb{T}$ is a Riemannian manifold with the standard metric and one can easily check that
$$\left<\nabla f(z),iz\right>=df_{z}(iz)=\frac{\partial f}{\partial z} iz\ \Longrightarrow\ \nabla f(z)=\left<\nabla f(z),iz\right>iz=-z^2\frac{\partial f}{\partial z}=\frac{\partial f}{\partial \bar z}.$$
\end{rem}

\begin{defi}\label{def-space-periodic-functions}
We will set the next Banach spaces of test functions:
\begin{enumerate}
\item $C^1(\mathbb{T})$ will denote the Banach space of continuously differentiable functions $f:\mathbb{T}\longrightarrow\RR$. It is endowed with the complete norm
$$\Vert f\Vert_{C^1(\mathbb{T})}:=\left\Vert f\right\Vert_{C(\mathbb{T})}+\left\Vert \frac{\partial f}{\partial z}\right\Vert_{C(\mathbb{T})}.$$
\item $C^1_p([-\pi,\pi])$ will denote the Banach space of continuously differentiable functions $g:[-\pi,\pi]\longrightarrow\RR$ such that $g$ and its derivatives have same values at the endpoints of the interval. It is endowed with the complete norm
$$\Vert g\Vert_{C^1_p([-\pi,\pi])}:=\left\Vert g\right\Vert_{C([-\pi,\pi])}+\left\Vert\frac{d g}{d\theta}\right\Vert_{C([-\pi,\pi])}.$$
\item As it is usual,  $C^1_0((-\pi,\pi])$ (respectively, $C^1_0(-\pi,\pi)$) denote the Banach space of continuously differentiable functions such that $g$ and its derivatives vanish at $-\pi$ (respectively at $-\pi$ and $\pi$). It is endowed with the preceding complete norm.
\end{enumerate}
The spaces $C(\mathbb{T})$, $C_p([-\pi,\pi])$, $C_0((-\pi,\pi])$ and $C_0(-\pi,\pi)$ are also Banach spaces when endowed with the uniform norm. Similarly, spaces with higher order derivatives can also be considered.
\end{defi}

By definition, the next results are clear:

\begin{pro}\label{pro-space-functions-periodic-1}
For every $f\in C^1(\mathbb{T})$, let us define $\Phi[f](\theta):=f(e^{i\theta})$. Then, the following map is an isometric isomorphism
$$\begin{array}{cccc}
\Phi: & C(\mathbb{T}) & \longrightarrow & C_p([-\pi,\pi]),\\
 & f & \longmapsto & \Phi[f].
\end{array}$$
\end{pro}

\begin{pro}\label{pro-space-functions-periodic-2}
For every $g\in C_p([-\pi,\pi])$, let us define 
\begin{align*}
\Psi_1[g]&:=g(-\pi)=g(\pi)\in\RR,\\
\Psi_2[g]&:=g-\Psi_1[g]\in C_0(-\pi,\pi).
\end{align*}
Then, the next map is a topological isomorphism
$$\begin{array}{cccc}
\Psi: & C_p([-\pi,\pi]) & \longrightarrow & \RR\,\oplus_\infty\,C_0(-\pi,\pi),\\
 & g & \longmapsto & \Psi[g]:=(\Psi_1[g],\Psi_2[g]).
\end{array}$$
\end{pro}

\begin{pro}\label{pro-space-functions-periodic-3}
Let us set any cut-off function $\eta\in C([-\pi,\pi])$ such that $\eta(-\pi)=0,\,\eta(\pi)=1$ and $0\leq \eta\leq 1$. For every $g\in C_0((-\pi,\pi])$, let us define
\begin{align*}
\Lambda_1^\eta[g]&:=g(\pi)\in\RR,\\
\Lambda_2^\eta[g]&:=g-\Lambda_1^\eta[g] \eta\in C_0(-\pi,\pi).
\end{align*}
Then, the next map is a topological isomorphism
$$\begin{array}{cccc}
\Lambda^\eta: & C_0((-\pi,\pi]) & \longrightarrow & \RR\,\oplus_\infty\,C_0(-\pi,\pi),\\
 & g & \longmapsto & (\Lambda_1^\eta[g],\Lambda_2^\eta).
\end{array}$$
\end{pro}

\subsection{Periodic measures}

\begin{defi}\label{def-space-periodic-measures}
As it is usual, taking duals of the Banach spaces in Definition \ref{def-space-periodic-functions} we arrive at the next Banach spaces of finite Radon measures endowed with the (dual) total variation norm:
\begin{align*}
\mathcal{M}(\mathbb{T})&:=C(\mathbb{T})^*,\\
\mathcal{M}_p([-\pi,\pi])&:=C_p([-\pi,\pi])^*,\\
\mathcal{M}(-\pi,\pi)&:=C_c(-\pi,\pi)^*,\\
\mathcal{M}((-\pi,\pi])&:=C_c((-\pi,\pi])^*.
\end{align*}
\end{defi}

\begin{proof}[{Proof of Theorem \ref{theo-relation-measure-spaces}}]
The proof is just a simple consequence of Propositions \ref{pro-space-functions-periodic-1}, \ref{pro-space-functions-periodic-2} and \ref{pro-space-functions-periodic-3} that follows from taking dual operators to $\Phi$, $\Psi$ and $\Lambda^\eta$
$$\begin{array}{cccc}
\Phi^*: & \mathcal{M}_p([-\pi,\pi]) & \longrightarrow & \mathcal{M}(\mathbb{T}),\\
 & \mu & \longmapsto & \Phi^*(\mu),\\
 \Psi^*: & \RR\,\oplus_1\,\mathcal{M}(-\pi,\pi) & \longrightarrow & \mathcal{M}_p([-\pi,\pi]),\\
 & (b,\nu) & \longrightarrow & \Psi^*(b,\nu),\\
 (\Lambda^\eta)^*: & \RR\,\oplus_1\,\mathcal{M}(-\pi,\pi) & \longmapsto & \mathcal{M}((-\pi,\pi]),\\
 & (b,\nu) & \longmapsto & (\Lambda^\eta)^*(b,\nu).
\end{array}$$
Indeed, given $\mu\in \mathcal{M}_p([-\pi,\pi])$, $b\in \RR$ and $\nu\in \mathcal{M}(-\pi,\pi)$ and setting $f\in C(\mathbb{T})$ and $g\in C_p([-\pi,\pi])$, the duality read
\begin{align*}
\left<\Phi^*(\mu),f\right>&=\left<\mu,\Phi[f]\right>,\\
\left<\Psi^*(b,\nu),g\right>&=\Psi_1[g] b+\left<\nu,\Psi_2[g]\right>,\\
\left<(\Lambda^\eta)^*(b,\nu),g\right>&=\Lambda^\eta_1[g]b+\left<\nu,\Lambda^\eta_2[g]\right>.
\end{align*}
 \end{proof}
 
\begin{rem}\label{isomorphism_measure_spaces_basis}
In particular, the next identification takes place under the above topological isomorphisms:
$$\begin{array}{ccccc}
\mathcal{M}(\mathbb{T}) &\cong &\RR\,\oplus\,\mathcal{M}(-\pi,\pi)&\cong& \mathcal{M}((-\pi,\pi]),\\
\delta_{(-1,0)} & \equiv & (1,0\,d\theta) & \equiv & \delta_\pi.
\end{array}$$
\end{rem}

Notice that we have only provided topological isomorphisms between the above spaces of measures in Theorem \ref{theo-relation-measure-spaces}. In particular, we can identify the spaces $\mathcal{M}(\mathbb{T})$ with $\mathcal{M}((-\pi,\pi])$ via the composition 
$$\mathfrak{I}^\eta:=(\Lambda^\eta)^*\circ(\Psi^*)^{-1}\circ (\Phi^*)^{-1}=((\Psi\circ\Phi)^{-1}\circ \Lambda^\eta)^*.$$
Specifically, it means that, given $\mu \in\mathcal{M}(\mathbb{T})$, the measure $\mathfrak{J}^\eta[\mu]\in \mathcal{M}((-\pi,\pi])$ acts as follows
$$\left<\mathfrak{J}^\eta[\mu],\phi\right>=\left<\mu,f_\phi^\eta\right>,$$
for any $\phi\in C_0((-\pi,\pi])$, where the continuous function $f_\phi^\eta:=((\Psi\circ\Phi)^{-1}\circ \Lambda^\eta)[\phi]\in C(\mathbb{T})$ reads
$$f_\phi^\eta(e^{i\theta})=\phi(\theta)+(1-\eta(\theta))\phi(\pi),\ \theta\in (-\pi,\pi].$$
Although $\mathfrak{J}^\eta$ is not an isometry, we still can do better and introduce a simpler isometry between the particular spaces $\mathcal{M}(\mathbb{T})$ and $\mathcal{M}((-\pi,\pi])$. This will be the content of our last result.

\begin{theo}\label{theo-relation-measure-spaces-isometry}
Let us consider the bijective and continuous mapping
$$
\begin{array}{cccc}
\iota: & (-\pi,\pi] & \longrightarrow &\mathbb{T}\\
 & \theta & \longmapsto & e^{i\theta}. 
\end{array}
$$
Then, the associated push-forward mapping is a surjective isometry
$$\begin{array}{cccc}
\iota_{\#}: & \mathcal{M}((-\pi,\pi]) & \longrightarrow & \mathcal{M}(\mathbb{T})\\
 & \mu & \longmapsto & \iota_{\#}\mu.
\end{array}$$
\end{theo}

\begin{proof}
First, it is clear that $\iota_{\#}$ is bijective because so is $\iota$. Let us now prove that it is a linear isometry. On the one hand, consider any $\mu \in \mathcal{M}((-\pi,\pi])$. Then,
\begin{align*}
\Vert \iota_{\#}\mu\Vert_{\mathcal{M}(\mathbb{T})}&=\sup_{\Vert \varphi\Vert_{C(\mathbb{T})}\leq 1}\int_{\mathbb{T}}\varphi\,d(\iota_{\#}\mu)=\sup_{\Vert \varphi\Vert_{C(\mathbb{T})}\leq 1}\int_{(-\pi,\pi]}(\varphi\circ\iota) \,d\mu\\
&\leq \sup_{\Vert \phi\Vert_{C_0((-\pi,\pi])}\leq 1}\int_{(-\pi,\pi]}\phi\,d\mu=\Vert\mu\Vert_{\mathcal{M}((-\pi,\pi])}.
\end{align*}
On the other hand, let us show the reverse inequality. Take any $\phi\in C_0((-\pi,\pi])$ with $\Vert \phi\Vert_{C_0((-\pi,\pi])}\leq 1$. By density, we can assume that $\phi\in C_c((-\pi,\pi])$. Then, there exists some $M_0\in (-\pi,\pi)$ such that
$$\phi(\theta)=0,\ \mbox{ for all }\ \theta\in (-\pi,M_0].$$
For every $\varepsilon\in (0,\pi+M_0)$ let us consider a cut-off function $\eta_\varepsilon\in C([-\pi,\pi])$  with $0\leq \eta_\varepsilon\leq 1$ such that
$$\eta_\varepsilon(-\pi)=0\ \mbox{ and }\ \left.\eta_\varepsilon\right\vert_{[-\pi+\varepsilon,\pi]}\equiv 1.$$
For any $\varepsilon>0,$ let us define
$$f_\varepsilon(e^{i\theta}):=f_\phi^{\eta_\varepsilon}(e^{i\theta})=\phi(\theta)+(1-\eta_\varepsilon)\phi(\pi),\ \theta\in (-\pi,\pi],$$
as in the proof of Theorem \ref{theo-relation-measure-spaces}. Then, it is clear that
$$\int_{(-\pi,\pi]}\phi\,d\mu=\int_{(-\pi,\pi]}(f_\varepsilon\circ\iota)\,d\mu-\phi(\pi)\int_{(-\pi,\pi]}(1-\eta_\varepsilon)\,d\mu=\int_{\mathbb{T}}f_\varepsilon\,d(\iota_{\#}\mu)-\phi(\pi)\int_{(-\pi,\pi]}(1-\eta_\varepsilon)\,d\mu.$$
Notice that due to the boundedness of $1-\eta_\varepsilon$ we have that it belongs to $L^1(\mu)$ and the above terms make sense. In fact
$$\int_{(-\pi,\pi]}\phi\,d\mu\leq \int_{\mathbb{T}}\vert f_\varepsilon\vert\,d\vert\iota_{\#}\mu\vert+\Vert\phi\Vert_{C_0((-\pi,\pi])}\left\vert\int_{(-\pi,\pi]}(1-\eta_\varepsilon)\,d\mu\right\vert,$$
for every $\varepsilon>0$. Regarding the first term, it is clear that $\Vert f_\varepsilon\Vert_{C(\mathbb{T})}=\Vert \varphi\Vert_{C(\mathbb{T})}\leq 1$ for every $\varepsilon\in (0,\pi+M_0)$. On the other hand, notice that the second term vanishes as $\varepsilon\rightarrow 0$ due to the dominated convergence theorem. Then, putting both facts together we can conclude
$$\int_{(-\pi,\pi]}\phi\,d\mu\leq \Vert \iota_{\#}\mu\Vert_{\mathcal{M}(\mathbb{T})},$$
and it ends the proof.
\end{proof}

\section{Ascoli-Arzel\`a-Prokhorov theorem}\label{appendix-Prokhorov}
The standard Ascoli-Arzel\`a theorem states that any uniformly bounded and equicontinuous sequence of functions $\{f_n\}_{n\in\mathbb{N}}\subseteq C(K)$, where $K$ is a compact metric space, admits an uniformly convergent subsequence $\{f_{\sigma(n)}\}_{n\in\mathbb{N}}$. The main idea relies on a Cantor diagonal argument along with the Bolzano-Weierstrass theorem in finite dimensional spaces. Regarding spaces like $C(K,B)$ with $B$ a Banach space we can mimic the proof when we substitute the Bolzano-Weierstrass theorem by some infinite-dimensional compactness theorem. For instance:
\begin{enumerate}
\item When $B$ is a dual space, boundedness implies weak-star compactness by Alaoglu's theorem. Then, any uniformly bounded and equicontinuous sequence $\{\mu_n\}_{n\in\mathbb{N}}\subseteq C(K,B)$ has a convergent subsequence in $C(K,B-\mbox{weak}\,*)$.
\item When $B$ is a reflexive space, boundedness implies weak compactness. Then, any uniformly bounded and equicontinuous sequence $\{\mu_n\}_{n\in\mathbb{N}}\subseteq C(K,B)$ has a convergent subsequence in $C(K,B-\mbox{weak})$.
\end{enumerate}
In particular, if $\{\mu_n\}_{n\in\mathbb{N}}\subseteq C([0,T],\mathcal{M}(X))$ is an uniformly bounded and equicontinuous sequence, where $\mathcal{M}(X)$ is the space of finite Radon measures over a locally compact and Hausdorff topological space $X$, then there exists some subsequence $\{\mu_{\sigma(n)}\}_{n\in\mathbb{N}}$ that converges in $C([0,T],\mathcal{M}(\mathbb{T}\times \mathbb{R})-\mbox{weak}\,*)$. The space of finite Radon measures $\mathcal{M}(X)$ can also be endowed with the narrow topology. Recall that $\{\nu_n\}_{n\in\mathbb{N}}\subseteq \mathcal{M}(X)$ is said to converge narrowly to $\nu\in\mathcal{M}(X)$ when
$$\int_{X}\varphi\,d\nu_n\rightarrow\int_{X}\varphi\,d\nu,$$
for every $\varphi\in C_b(X)$, where $C_b(X)$ stands for the space of continuous and bounded functions on $X$. A nice characterization of the narrow compactness follows from the well known Prokhorov's compactness theorem. Namely, a the sequence $\{\nu_n\}_{n\in\mathbb{N}}\subseteq \mathcal{M}(X)$ that is uniformly bounded (in total variation norm) and uniformly tight enjoys a narrowly convergent subsequence. Combining such compactness result with the Cantor diagonal argument we obtain the following compactness result.
\begin{lem}\label{L-Prokhorov}
Let $K$ be a compact metric space and $X$ be a locally compact and Hausdorff topological space. Consider $\{\mu_n\}_{n\in\mathbb{N}}\subseteq C(K,\mathcal{M}(X))$ such that
it is uniformly bounded and equicontinuous in $C(K,\mathcal{M}(X))$. In addition, assume that it is uniformly tight in $n$ and $K$, i.e., for every $\varepsilon>0$ there exists a compact subset $K_\varepsilon\subseteq X$ so that
$$\vert\mu_n(k)\vert(X\setminus K_\varepsilon)<\varepsilon,\ \forall k\in K,\,\forall n\in\mathbb{N}.$$
Then, there exists some $\mu\in C(K,\mathcal{M}(X))$ and a subsequence $\{\mu_{\sigma(n)}\}_{n\in\mathbb{N}}$ such that
$$\{\mu_{\sigma(n)}\}_{n\in\mathbb{N}}\rightarrow \mu\ \mbox{ in }\ C(K,\mathcal{M}(X)-\mbox{narrow}).$$
\end{lem}

\section{Differentiability properties of the squared distance}\label{appendix-differentiability-distance}

In this Appendix we will revisit some well known results about (non)-differentiability of the squared distance in a Riemannian manifold. Most of them are folklore in Riemannian Geometry and require no comment. Nevertheless, we will comment on the appropriate concept of differentiability that we are interested in, namely, the one-sided Dini upper derivative.

\begin{pro}\label{P-derivative-distance-1}
Let $(M,\left<\cdot,\cdot\right>)$ be a $d$-dimensional complete Riemannian manifold and consider the Riemannian distance as recalled in \eqref{E-riemannian-distance}. Fix any $y\in M$ and define the distance function $d_y:M\rightarrow \mathbb{R}$ towards some point $y\in M$ by the rule
$$d_y(x)=d(x,y),\ \forall\,x\in M.$$
Then, following properties hold true:
\begin{enumerate}
\item $d_y$ is Lipschitz in $M$.
\item $d_y$ is derivable at almost every $y\in M$.
\item $d_y$ is derivable in $M\setminus\{\cut(y)\cup\{y\}\}$ and 
$$(\nabla d_y)(x)=-\frac{\exp_x^{-1}(y)}{\vert \exp_x^{-1}(y)\vert}.$$
\item $\frac{1}{2}d_y^2$ fails to be everywhere directionally derivable unless $M$ is diffeomorphic to the flat space $\mathbb{R}^d$. However, it is derivable at any $x\in M\setminus(\cut(x)\cup\{x\})$ and
$$\nabla\left(\frac{1}{2}d_y^2\right)(x)=-\exp_x^{-1}(y).$$
\end{enumerate}
Here, $\cut(x)$ denotes the cut locus of the point $x$ in $M$ and $\exp_x:T_xM\rightarrow M$ is nothing but the Riemannian exponential map at such $x$.
\end{pro}

Since the proofs are standard, we do not provide proofs here. They can be found in any textbook of Riemannian Geometry. Instead, we just focus on the later assertion that is the less apparent one. According to such result not only $\frac{1}{2}d_y^2$ fails to be derivable at some points $x\in \cut(y)$, but also the lateral directional derivatives might not agree at points $x\in M$ that can be joined with $y$ through several minimizing geodesic. The proof is apparently hidden in the litterature; however, one can find a short proof following simple arguments in \cite{W}. 

\begin{rem}\label{R-derivative-distance}
Just to illustrate a meaningful example, consider $M=\mathbb{T}$ with the (standard) induced metric. Recall that for $z_1=e^{i\theta_1}$ and $z_2=w^{i\theta_2}$ where $\theta_1,\theta_2\in\mathbb{R}$ one has the clear identity
$$d(z_1,z_2)=\vert \theta_1-\theta_2\vert_o.$$
Also, consider some $z=e^{i\theta}$, where $\theta\in\mathbb{R}$, and its antipode $\bar z=-z=e^{i(\theta+\pi)}$. For every $\omega>0$, one can define a geodesic with speed $\omega$ as follows
$$\gamma_{z,\omega}(s):=e^{i(\theta+\omega s)},\ s\in \mathbb{R},$$
Indeed, $\gamma_{z,\omega}(0)=z$ and it is minimizing in any interval whose length is not larger than $\frac{\pi}{\omega}$. Consequently,
$$\frac{1}{2}d_{\bar z}^2(\gamma_\omega(s))=\left\{\begin{array}{ll}
\frac{1}{2}(\pi+\omega s)^2, & s\in (-\frac{\pi}{\omega},0],\\
\frac{1}{2}(\pi-\omega s)^2, & s\in [0,\frac{\pi}{\omega}).
\end{array}\right.$$
Hence, both one-sided derivatives exist but they differ from each other, namely,
$$\left.\frac{d}{ds}\right\vert_{s=0^\mp}\frac{1}{2}d_{\bar{z}}^2(\gamma_{z,\omega}(s))=\pm\omega\pi.$$
\end{rem}

Apart from the above differentiability properties, others have been explored in the literature with applications to optimal mas transportation and Wasserstein distances. We address some of them in the following result (see \cite[Proposition 2.9]{F-V}, \cite[Proposition 6]{M}, \cite[Third Appendix of Chapter 10]{V} for more details).

\begin{pro}\label{P-derivative-distance-2}
Let $(M,\left<\cdot,\cdot\right>)$ be a complete Riemannian manifold, fix $y\in M$ and the squared distance function $\frac{1}{2}d_y^2$ towards the point $y$. Then,
\begin{enumerate}
\item $\frac{1}{2}d_y^2$ is superdifferentiable at every $x\in M$ and for every $w\in \exp_x^{-1}(y)$ with $\vert w\vert =d(x,y)$ one has that $-w$ is an upper gradient, i.e.,
$$\frac{1}{2}d_y^2(\exp_x(v))\leq \frac{1}{2}d_y^2(x)-\left<w,v\right>+o(\vert v\vert), \mbox{ as }\ \vert v\vert\rightarrow 0.$$
\item If in addition $M$ has non-negative sectional curvatures then $\frac{1}{2}d_y^2$ is $1$-semiconcave, i.e.,
$$\frac{1}{2}d_y^2(\gamma(s))\geq (1-s)\frac{1}{2}d_y^2(x_1)+s\frac{1}{2}d_y^2(x_2)+s(1-s)\frac{1}{2}d^2(x_1,x_2),\ s\in [0,1],$$
for any couple $x_1,x_2\in M$ and any geodesic $\gamma:[0,1]\longrightarrow M$ joining $x_1$ to $x_2$.
\end{enumerate}

\end{pro}

Both superdifferentiability and semiconcavity are locally equivalent (see \cite{V} for more details). On the one hand, the non-negativity condition on the sectional curvatures is required in order to obtain uniform estimates in term of the quadratic modulus of semiconcavity and it is certainly a very imposing hypothesis that we will not assume. On the other hand, the superdifferentiability is not enough for our purpose since the $o(\vert v\vert)$ term in the right hand-side is not necessarily uniform. For the purposes in this paper we will resort on an slightly different tool coming form non-smooth analysis, namely, the \textit{one-sided upper Dini directional derivative} of a function.

\begin{defi}\label{D-Dini}
Let $(M,\left<\cdot,\cdot\right>)$ be a finite-dimensional complete Riemannian manifold and consider some function $f:M\rightarrow \mathbb{R}$, any $x\in M$ and any direction $v\in T_xM$.  Then, the one-sided upper and lower Dini derivatives of $f$ at $x$ in the direction $v$ stand for
\begin{align*}
(d^+f)_x(v)&:=\left.\frac{d^+}{ds}\right\vert_{s=0} f(\exp_x(sv))=\limsup_{s\rightarrow 0^+}\frac{f(\exp_x(sv))-f(x)}{s},\\
(d_+f)_x(v)&:=\left.\frac{d_+}{ds}\right\vert_{s=0} f(\exp_x(sv))=\liminf_{s\rightarrow 0^+}\frac{f(\exp_x(sv))-f(x)}{s}.
\end{align*}
By definition both derivatives are ordered
$$-\infty\leq (d_+f)_x(v)\leq (d^+f)_x(v)\leq +\infty,$$
and we will say that $f$ is one-sided upper (respectively lower) Dini derivable at $x$ in the direction $v$ if the corresponding one-sided upper (respectively lower) Dini derivative is finite. If in addition both derivatives agree, then $f$ is also one-sided directionally derivable at $x$ in the direction $v$ in the standard sense and all the derivatives agree.
\end{defi}

An important fact to be remarked is that the geodesic $s\mapsto \exp_x(sv)$ has been chosen as a representative of a curve with direction $v$ at $x$. However, one might have taken any other $C^1$ curve and apparently it would have provided a ``different'' definition of directional derivative. Since it will be of interest for our purposes, let us show that any such curve representing the direction $v$ at $x$ could have been chosen, yielding the same definition.

\begin{lem}\label{L-Dini-characterization}
Let $(M,\left<\cdot,\cdot\right>)$ be a Riemannian manifold and fix any $x\in M$ and any couple $v,w\in T_xM$. Consider a couple of $C^1$ curves $\gamma_1:(-s_0,s_0)\rightarrow M$ and $\gamma_2:(-s_0,s_0)\rightarrow M$ such that $\gamma_1(0)=x=\gamma_2(0)$, $\gamma_1'(0)=v$ and $\gamma_2'(0)=w$. Then,
$$\limsup_{s\rightarrow 0}\frac{d(\gamma_1(s),\gamma_2(s))}{s}\leq\vert v-w\vert.$$
\end{lem}

\begin{proof}
Although the proof is standard, we provide a simple proof for the sake of completeness. Consider some $R>0$ smaller enough than the radius of injectivity at $x$ (e.g. half of it) and consider the ball $B_R(0)\subseteq T_xM$ along with the associated geodesic ball $\mathbb{B}_R(x):=\exp_x(B_R(0))$. By definition one has that $\exp_x:B_R(0)\longrightarrow\mathbb{B}_R(x)$ is a diffeomorphism. Without loss of generality, we will assume that $s_0$ is small enough so that $\gamma_i(s)\in \mathbb{B}_R(x)$ for all $s\in(-s_0,s_0)$ and $i=1,2$. Hence, we can define the curves in $B_R(0)$
$$\bar\gamma_i(s):=\exp_x^{-1}(\gamma_i(s)),\ s\in (-s_0,s_0),\ i=1,2.$$
Equivalently, $\gamma_i(s)=\exp_x(\bar \gamma_i(s))$ and taking derivatives one has
$$\gamma_i'(0)=(d\exp_x)_0(\bar\gamma_i'(0))\Longrightarrow\ \bar \gamma_1'(0)=v\ \mbox{ and }\ \bar \gamma_2'(0)=w,$$
where we have used that $(d\exp_x)_0$ is nothing but the identity map in $T_x M$ and $\gamma_1'(0)=v$, $\gamma_2'(0)=w$.

Also, consider the interpolating curves between $\bar\gamma_1(s)$ and $\bar\gamma_2(s)$
$$\bar\gamma_s(\varepsilon)=(1-\varepsilon)\bar\gamma_1(s)+\varepsilon\gamma_2(s),\ \varepsilon\in [0,1],$$
for every $s\in (-s_0,s_0)$. They have associated interpolating curves between $\gamma_1(s)$ and $\gamma_2(s)$
$$\gamma_s(\varepsilon)=\exp_x(\bar\gamma_s(\varepsilon))=\exp_x((1-\varepsilon)\gamma_1(s)+\varepsilon\gamma_2(\varepsilon)),\ \varepsilon\in [0,1].$$
Then, we can estimate
$$\frac{d(\gamma_1(s),\gamma_2(s))}{s}=\frac{d(\exp_x(\bar{\gamma}_1(s)),\exp_x(\bar{\gamma}_2(s)))}{s}\leq [\exp_x]_{C^{0,1}(\mathbb{B}_R(x))}\frac{\vert \bar\gamma_1(s)- \bar\gamma_2(s)\vert}{s},$$
where $[\,\cdot\,]_{C^{0,1}(B_R(x))}$ stands for the Lipschitz constant in $B_R(x)$. Taking $\limsup$ and using that $\bar\gamma_1$ and $\bar\gamma_2$ are both derivable at $s=0$ with derivatives $v$ and $w$ respectively we obtain
$$\limsup_{s\rightarrow 0}\frac{d(\gamma_1(s),\gamma_2(s))}{s}\leq [\exp_x]_{C^{0,1}(\mathbb{B}_R(x))}\vert v-w\vert.$$
Now, notice that $R>0$ can be chosen arbitrarily small. Taking $\liminf$ with respect to such radius we obtain
$$\limsup_{s\rightarrow 0}\frac{d(\gamma_1(s),\gamma_2(s))}{s}\leq \liminf_{R\rightarrow 0}\,[\exp_x]_{C^{0,1}(\mathbb{B}_R(x))}\vert v-w\vert.$$
Also, the mean value theorem implies
$$\liminf_{R\rightarrow 0}\,[\exp_x]_{C^{0,1}(\mathbb{B}_R(x))}\leq \liminf_{R\rightarrow 0}\sup_{\xi\in B_R(0)}\vert (d\exp_x)_\xi\vert_{T^*_xM}.$$
Since $(d\exp_x)_0$ agrees with the identity map in $T_x M$, that has operator norm equals to $1$, we infer
$$\liminf_{R\rightarrow 0}\,[\exp_x]_{C^{0,1}(\mathbb{B}_R(x))}\leq 1,$$
and that ends the proof of our result.
\end{proof}

As a consequence we obtain the following characterizations of the one-sided Dini derivatives.

\begin{theo}\label{T-Dini-characterization}
Let $(M,\left<\cdot,\cdot\right>)$ be a complete Riemannian manifold, consider some $x\in M$ and $v\in T_x M$ and choose any $C^1$ curve $\gamma:(-s_0,s_0)\longrightarrow M$ such that $\gamma(0)=x$ and $\gamma'(0)=v$. If $f:M\longrightarrow\mathbb{R}$ is locally Lipschitz map around $x$, then
\begin{align*}
(d^+ f)_x (v)&=\limsup_{s\rightarrow 0^*}\frac{f(\gamma(s))-f(x)}{s},\\
(d_+ f)_x (v)&=\liminf_{s\rightarrow 0^*}\frac{f(\gamma(s))-f(x)}{s}.
\end{align*} 
\end{theo}

\begin{proof}
Let us define the auxiliary $C^1$ curve $\widetilde{\gamma}(s)=\exp_x(sv)$ for $s\in (-s_0,s_0)$. According to the above Definition \ref{D-Dini} we only need to prove that
$$\lim_{s\rightarrow 0^+}\frac{f(\gamma(s))-f(\widetilde{\gamma}(s))}{s}=0.$$
Consider $R$ smaller enough that the radius of injectivity at $x$ and set $B_R(0)\subseteq T_x M$ along with the geodesic ball $\mathbb{B}_R(x):=\exp_x(B_R(0))$ that is a relatively compact set. Consider $L_R$ the Lipschitz constant of $f$ in $\overline{\mathbb{B}}_R(x)$ and assume that $s_0$ is small enough so that $\gamma(s),\widetilde{\gamma}(s)\in \overline{\mathbb{B}}_R(x)$ for every $s\in (-s_0,s_0)$. Thus,
$$\left\vert\frac{f(\gamma(s))-f(\widetilde{\gamma}(s))}{s}\right\vert\leq L_R\frac{d(\gamma(s),\widetilde{\gamma}(s))}{s},\ s\in (-s_0,s_0).$$
Since $\gamma'(0)=v=\widetilde{\gamma}'(0)$ we conclude the proof of this result by virtue of Lemma \ref{L-Dini-characterization}.
\end{proof}

We are now ready to give a simple proof of the one-sided uper Dini directional differentiability of the distance.

\begin{theo}\label{T-dini-derivative-distance}
Let $(M,\left<\cdot,\cdot\right>)$ be a complete Riemannian manifold, fix $y\in M$ and the squared distance $\frac{1}{2}d_y^2$ towards the point $y$. Then, $\frac{1}{2}d_y^2$ is one-sided upper Dini directionally derivable in all $M$ and
$$d^+\left(\frac{1}{2}d_y^2\right)_x(v)\leq \inf_{\substack{w\in \exp_x^{-1}(y)\\ \vert w\vert=d(x,y)}} -\left<v,w\right>,$$
for any $x\in M$ and any direction $v\in T_x M$.
\end{theo}

\begin{proof}
Consider $x\in M$ and $v\in T_x M$ and set $w\in \exp_x^{-1}(y)$ with $\vert w\vert=d(x,y)$. Also, consider a minimizing geodesic $\gamma_0:[0,1]\longrightarrow M$ such that $\gamma_0(0)=x$, $\gamma_0(1)=y$ and $\gamma_0'(0)=w$. Then, there exists some smooth variation $\{\gamma_s\}_{s\in (-s_0,s_0)}$ of $\gamma_0$ such that
$$\gamma_s(0)=\exp_x(sv)\ \mbox{ and }\ \gamma_s(1)=y,$$
for every $s\in (-s_0,s_0)$. Such a variation enjoys an associated variational field
$$V(\xi):=\frac{\partial \gamma}{\partial s}(s=0,\xi)=\left.\frac{d}{ds}\right\vert_{s=0}\gamma_s(\xi),\ \xi\in [0,1].$$
Also, we can define the associated energy functional
$$\mathcal{E}[\gamma_s]:=\frac{1}{2}\int_0^1\vert \gamma_s'(\xi)\vert^2\,d\xi.$$
Then, the first variation formula of energy (see \cite[Proposition 2.4]{DC}) around the geodesic $\gamma_0$ amounts to
$$
\left.\frac{d}{ds}\right\vert_{s=0}\mathcal{E}[\gamma_s]=\left<\gamma_0'(1),V(1)\right>-\left<\gamma_0'(0),V(0)\right>=\left<\gamma_0'(1),0\right>-\left<w,v\right>=-\left<w,v\right>.
$$
Since $\gamma_0$ is minimizing then we can equivalently restate
\begin{equation}\label{E-App3-1}
\lim_{s\rightarrow 0}\frac{\mathcal{E}[\gamma_s]-\frac{1}{2}d_y^2(x)}{s}=-\left<w,v\right>,
\end{equation}
and it is clear that
\begin{align*}
\frac{\frac{1}{2}d_y^2(\exp_x(sv))-\frac{1}{2}d_y^2(x)}{s}\leq \frac{\mathcal{E}[\gamma_s]-\frac{1}{2}d_y^2(x)}{s},\ s\in (-s_0,s_0),
\end{align*}
Taking $\limsup$ as $s\rightarrow 0^+$ and using the above formula \eqref{E-App3-1} implies
$$d^+\left(\frac{1}{2}d_y^2\right)_x(v)\leq -\left<v,w\right>,$$
where $w\in \exp_x^{-1}(y)$ such that $\vert w\vert=d(x,y)$ is arbitrary. This ends the proof. 
\end{proof}


\section*{Acknowledgment}
This work has been partially supported by the the MECD (Spain) research grant FPU14/06304, the MINECO-Feder (Spain) research grant number MTM2014-53406-R and the Junta de Andaluc\'ia (Spain) Project FQM 954.

The author thanks Prof. Juan Soler and Dr. Jinyeong Park for the fruitful discussions.


\end{document}